\documentclass[12pt]{report}

\usepackage{amssymb}
\usepackage{amsmath}
\usepackage{amsthm}
\usepackage{fullpage}
\usepackage{makeidx}

\usepackage{xypic}
\input xy
\xyoption{all}
\xyoption{2cell}
\xyoption{arrow}
\UseTwocells
\UseHalfTwocells
\UseCompositeMaps
\UseAllTwocells

\newtheorem{theorem}{Theorem}[section]

\newtheorem{prop}[theorem]{Proposition}

\newtheorem{lemma}[theorem]{Lemma}
\newtheorem{coro}[theorem]{Corollary}

\theoremstyle{definition}
\newtheorem{definition}[theorem]{Definition}
\newtheorem{example}[theorem]{Example}
\newtheorem{remark}[theorem]{Remark}
\newtheorem{note}[theorem]{Note}

\newcommand{\Int}{\mathop{\rm Int}\nolimits}

\newcommand{\uvw}[1]{u_{i{#1}}^{\alpha{#1}}\times
v_{j{#1}}^{\beta{#1}}\times w_{k{#1}}^{\varepsilon{#1}}}

\newcommand{\uvws}[1]{u[i{#1},\alpha{#1}]\times
v[j{#1},\beta{#1}]\times w[k{#1},\varepsilon{#1}]}
\newcommand{\huvws}[1]{u[\hat{i}{#1},\hat{\alpha}{#1}]\times
v[\hat{j}{#1},\hat{\beta}{#1}]\times
w[\hat{k}{#1},\hat{\varepsilon}{#1}]}
\newcommand{\lmns}[1]{u[l{#1},\sigma{#1}]\times
v[m{#1},\tau{#1}]\times w[n{#1},\omega{#1}]}
\newcommand{\ui}[1]{u[i{#1},\alpha{#1}]}
\newcommand{\ul}[1]{u[l{#1},\sigma{#1}]}
\newcommand{\vj}[1]{v[j{#1},\beta{#1}]}
\newcommand{\vm}[1]{v[m{#1},\tau{#1}]}
\newcommand{\wk}[1]{w[k{#1},\varepsilon{#1}]}
\newcommand{\wn}[1]{w[n{#1},\omega{#1}]}
\newcommand{\uvwt}{u\times v\times w}
\newcommand{\dpg}{d_{p}^{\gamma}}
\newcommand{\dna}{d_{n}^{\alpha}}

\newcommand{\vIj}[1]{v^{I}[j{#1},\beta{#1}]}
\newcommand{\wIk}[1]{w^{I}[k{#1},\varepsilon{#1}]}
\newcommand{\wJk}[1]{w^{J}[k{#1},\varepsilon{#1}]}
\newcommand{\A}{\mathcal{A}}
\newcommand{\M}{\mathcal{M}}
\newcommand{\PS}{\mathcal{P}}
\newcommand{\C}{\mathcal{C}}
\newcommand{\fdim}{\text{\rm fr }\dim}
\newcommand{\Cl}{\text{\rm Cl}}
\newcommand{\pr}{\text{\rm pr}}

\newcommand{\uaaa}[3]{u_{1}[#1_{1}#3,#2_{1}#3]\times
u_{2}[#1_{2}#3,#2_{2}#3]\times u_{3}[#1_{3}#3,#2_{3}#3]\times
u_{4}[#1_{4}#3,#2_{4}#3]}
\newcommand{\ubbb}{u_{1}\times u_{2}\times u_{3}\times u_{4}}
\newcommand{\ua}[1]{u_{#1}[i_{#1},\alpha_{#1}]}
\newcommand{\uaa}[4]{u_{#4}[#1_{#4}#3,#2_{#4}#3]}
\index{
\includeonly{
 }}

\title{ \huge {\bf The $\omega$-Categories
Associated With Products of Infinite-Dimensional Globes}}
\author{{\small by}\vspace{0.5cm}\\
Hongbin Cui\vspace{0.3cm}\\ {\small Department of Mathematics,
University of Glasgow}\vspace{1.5cm}\\ A thesis submitted to\\ the
Faculty of Science\\ at the University of Glasgow\\ for the degree
of\\ Doctor of Philosophy\vspace{0.5cm}\\
 {\small\copyright Hongbin Cui}
}
\date{November, 2000}

\makeindex
\begin{document}
\maketitle \sloppy \setlength{\evensidemargin}{0.2cm}
\setlength{\oddsidemargin}{0.2cm} \setlength{\textwidth}{17cm}
\setlength{\textheight}{24cm}

\pagenumbering{roman}
\chapter*{Summary}
\addcontentsline{toc}{chapter}{Summary}

The results in this thesis are organised in four chapters.

Chapter 1 is  preliminary. We state the necessary definitions and
results in $\omega$-complexes, atomic complexes and products of
$\omega$-complexes. Some definitions are restated to meet the
requirement for the following chapters. There is a new  proof for
the existence of {\em `natural homomorphism'} (Theorem
\ref{natural_homo}) and a new result for the decomposition of
molecules in loop-free $\omega$-complexes (Theorem
\ref{decom_mole}).

In Chapter 2, we study the product of three infinite dimensional
globes. The main result in this chapter is that a subcomplex in
the product of three infinite dimensional globes is a molecule if
and only if it is pairwise molecular (Theorem \ref{3glb_main}).
The definition for pairwise molecular subcomplexes is given in
section 1. One direction of the main theorem,  molecules are
necessarily pairwise molecular, is proved in section 2. Some
properties of pairwise molecular subcomplexes are studied in
section 3. These properties are the preparation for a more
explicit description of pairwise molecular subcomplexes, which is
given in section 4. The properties for the sources and targets of
pairwise molecular subcomplexes are studied in section 5, where we
prove that the class of pairwise molecular subcomplexes is closed
under source and target operation; there are also algorithms to
calculate the sources and targets of a pairwise molecular
subcomplex. Section 6 deals with the composition of pairwise
molecular subcomplexes. The proof of the main theorem is completed
in section 7, where an algorithm for decomposing molecules into
atoms is implied in the proof.

The construction of molecules in the product of three infinite
dimensional globes is studied in Chapter 3. The main result is
that any molecule can be constructed inductively by a systematic
approach. Section 1 gives another description for molecules in the
product of three infinite dimensional globes which is the
theoretical basis for the construction. Section 2 states the
inductive process of constructing molecules. The justification for
the construction is given in section 3.

The main result in Chapter 4 is that a subcomplex in the product
of four infinite dimensional globes is a molecule if and only if
it is pairwise molecular (Theorem \ref{4glb_main}).  In the first
four sections, some basic concepts and properties have to be
reestablished to suit more general case. The organisation for the
last three sections is parallel to that in Chapter 2. The
corresponding results for sources, targets, composition and
decomposition of pairwise molecular subcomplexes are also
achieved.


  \chapter*{Acknowledgements}
  I would like to express my deep appreciation to my supervisor,
   Doctor R. J. Steiner for his guidance,
   help and encouragement throughout the course of this research,
   without which this thesis would not have been possible; and to
   all the staff in the department for their
  help during my study at the Department of Mathematics,
   University of Glasgow.

  I should also like to thank the University
  of Glasgow and the Committee of Vice-Chancellors and Principals of the
  Universities of the United Kingdom for financing this research through a
University Scholarship and an ORS Award, respectively, from
October 1997 to September 2000.


\chapter*{Statement}
\addcontentsline{toc}{chapter}{Statement}

  This thesis is submitted according to the regulations for the degree
   of Doctor of Philosophy in the University of Glasgow. It presents part
    of research results carried out by the author during the academic
    years 1997-2000.

  All the results of this thesis are the original work of the author
  except for the instances indicated within the text. Some results
  are obtained jointly with  Doctor R. J. Steiner.

\tableofcontents


\chapter*{Introduction}
\addcontentsline{toc}{chapter}{Introduction}

In this work, we study the $\omega$-complexes of  products of
infinite-dimensional globes.

An $n$-category is an algebraic structure consisting of objects,
morphisms between objects, $2$-morphisms between morphisms, and so
on up to $n$-morphisms, subject to various composition rules.

The study of $n$-categories started from $2$-categories which
generalise the idea that natural transformations can certainly be
thought as morphisms between morphisms. The theory of bicategories
(a generalisation of $2$-categories) has successfully been
established by the wonderful work of, for example, Eilenberg and
Kelly \cite{e-k}, Kelly \cite{kelly}, Kelly and Street
\cite{kelly-street}, and Mac Lane and Par\'{e} \cite{m-p}.

The concept of $\omega$-category or $\infty$-category (\cite{b-h,
street2}) is a generalisation  of  $n$-category with no
restriction of `up to n'. It was originated by Brown and Higgins
in \cite{b-h}, in connection with homotopy theory. It was not long
after the concept was introduced that the researchers realised
that a sort of pasting diagrams representing compositions in
multiple categories should be introduced. There are several
approaches in the study of such pasting diagrams with different
names such as parity complexes, pasting schemes, directed
complexes or $\omega$-complexes. See Al-Agl and Steiner
\cite{al_rjs1}, Johnson \cite{johnson1}, Kapranov and Voevodsky
\cite{k-v1}, Power \cite{power1}, Steiner \cite{rjs2, rjs} and
Street \cite{street1, street3}. We follow the approach in Steiner
\cite{rjs} because the concept of $\omega$-complex is certainly
the most general one.

There is a  concept of products of $\omega$-complexes defined in
Steiner \cite{rjs}. It is natural because the products of
$\omega$-complexes give the tensor product of the underlying
$\omega$-categories. (For the study of tensor products of multiple
categories, see the work of Gray \cite{gray1}, Al-Agl and Steiner
\cite{al_rjs1}, Crans \cite{crans1}, Joyal and Street
\cite{joyal-street1}, and Brown and Higgins \cite{b-h2}). It is
shown in paper \cite{rjs} that the products of $\omega$-complexes
are still $\omega$-complexes. Since the definition for the product
is given by generators and relations, it is natural to seek
explicit descriptions for the products of $\omega$-complexes. This
problem is difficult in general, since the molecules, which are
representatives of elements in the underlying $\omega$-categories,
in the products are difficult to recognise. We thus consider the
solution for the products of the simplest $\omega$-complexes,
globes.

An $n$-dimensional globe $u$ is the an $\omega$-complex
representing the $n$-category with exactly one $n$-morphism and
two $m$-morphisms $u_{m}^{-}$ and $u_{m}^{+}$ for every
non-negative integer $m<n$ such that the $l$-source
$d_{l}^{-}u_{m}^{\gamma}$ and $l$-target $d_{l}^{+}u_{m}^{\gamma}$
of $u_{m}^{\gamma}$ are $u_{l}^{-}$ and $u_{l}^{+}$ respectively
for $l<m\leq n$ and $\gamma=\pm$. The infinite dimensional globe
is the obvious generalisation of $n$-globes. The globes are basic
$\omega$-complexes because they serves as the generators in the
standard representation of $\omega$-categories. (See Crans
\cite{crans2}.) The product of, for example, three infinite
dimensional globes $u\times v\times w$ is generated by elements of
the form $u_{i}^{\alpha}\times v_{j}^{\beta}\times
w_{k}^{\varepsilon}$ (called atoms) with relations resembling
those in homological algebra. Thus an element (called molecule) in
the product of three globes is a union of atoms (called
subcomplex). One of the main result in this thesis characterises
molecules in the product of three infinite dimensional globes, in
terms of such subcomplexes.

The study for the product of infinite dimensional globes is
important not only because infinite dimensional globe is a basic
$\omega$-complex, but also because it may help to understand the
products of general $\omega$-complexes. According to the
approaches used in paper \cite{rjs}, it seems that the product of
infinite dimensional globes has a sort of universal property which
may be used to study product of general $\omega$-complexes,
although we have not yet been able to  describe this universal
property precisely. Moreover, the explicit description of products
of infinite dimensional globes may also help in better
understanding some work in weak $n$-categories. (See Baez and
Neuchl \cite{baez-neuchl}, and Kapranov and Voevodsky
\cite{k-v2}.)

For the product of two infinite dimensional globes, there are
descriptions in paper \cite{street1} and \cite{rjs}. The
description in \cite{rjs} is more explicit in the way that the
molecules are easily recognised and constructed, and there are
explicit algorithms to calculate sources and targets of a molecule
and the composites of molecules, there is also an algorithm to
decompose a molecule into atoms. The main work in this thesis is
to extend these results to products of three and four infinite
dimensional globes.

As stated above, the construction for the product of two infinite
dimensional globes is fairly clear. So it is natural to reduce the
problem for the product of three infinite dimensional globes to
that of two infinite dimensional globes. This consideration leads
to the idea of describing a molecule in product of three infinite
dimensional globes by projecting it to the (twisted) products of
two infinite dimensional ones. This results in the definition of
{\em pairwise molecular subcomplexes} in the product of three
infinite dimensional globes. It is proved that molecules are
exactly pairwise molecular subcomplexes.

A more explicit description for molecules in the product of three
infinite dimensional globes is influenced by \cite{street1} and
\cite{rjs}. Some conditions in this description come from the
requirement that a molecule should be {\em well-formed}, while
some come from the comparison with the description of molecules in
the product of two infinite dimensional globes. A crucial concept
is the {\em adjacency} of maximal atoms in a subcomplex. This
description has some new features distinguished from that for the
product of two infinite dimensional globes. Some restrictions must
be given because of the middle factor. For example, if there is a
pair of distinct maximal atoms $\uvw{_{1}}$ and $\uvw{_{2}}$ in a
pairwise molecular subcomplex such that $i_{1}>i_{2}$,
$\min\{j_{1},j_{2}\}>0$ and $k_{1}<k_{2}$ , it is required that
there is a maximal atom $\uvw{}$ such that $i>i_{2}$, $j\geq
\min\{j_{1},j_{2}\}-1$ and $k>k_{1}$.

After the descriptions of molecules in the product of three
infinite dimensional globes are proposed, we have to prove that
pairwise molecular subcomplexes are closed under source and target
operations, and they are also closed under composition operations.
The algorithms for calculating the sources, targets and composites
of a pairwise molecular subcomplex are also studied.

To prove that pairwise molecular subcomplexes are molecules, we
have to show that they can be decomposed into atoms. To do this, a
total order, called {\em natural order}, on the set of atoms in
the product of three infinite dimensional globes is introduced.
The natural order is designed so that the maximal atoms of
dimensions greater than the {\em frame dimension $p$} (see paper
\cite{rjs2}) in a pairwise molecular subcomplex can be listed as
$\lambda_{1}$, $\lambda_{2}$, $\dots$, $\lambda_{S}$ with
$\lambda_{s}\cap\lambda_{t}\subset d_{p}^{+}\lambda_{s}\cap
d_{p}^{-}\lambda_{t}$ for all $1\leq s<t\leq S$. This means that
the decomposition approach in paper \cite{rjs2} applies. In the
proof that pairwise molecular subcomplexes are molecules, there is
also an explicit algorithm to decompose a molecule into atoms.

At this stage, we have satisfactory descriptions for the product
of three infinite dimensional globes. However, these are still
descriptive. From these descriptions, it is fairly easy to check
whether a given subcomplex is a molecule. But we still cannot
construct all the molecules. Our next goal is to find a systematic
way to construct all molecules. The approach is based on the
middle factor.  According our results, we can construct any
molecule, inductively, by listing its maximal atoms as
$\lambda_{1}$, $\lambda_{2}$, $\dots$, $\lambda_{R}$ with
$\lambda_{r}=\uvw{_{r}}$ such that $j_{1}\geq\cdots\geq j_{R}$ and
such that $i_{r}>i_{r+1}$ when $1\leq r<R$ and $j_{r}=j_{r+1}$. In
more detail, let maximal atoms $\lambda_{1}$, $\dots$,
$\lambda_{r}$ be an initial segment of the list. We can easily
determine whether $\lambda_{1}\cup\cdots\cup\lambda_{r}$ is
already a molecule and determine the set of possible next maximal
atoms $\lambda_{r+1}$, so that all the molecules can be
constructed inductively.

Up to now, we have a completely satisfactory theory for the
product of three infinite dimensional globes.

Our discussion for the product of four infinite dimensional globes
is roughly parallel to that for the product of three infinite
dimensional globes. Since the construction of the product of three
infinite dimensional globes, by our results, is thought to be
clear, we propose that the molecules in the product of four
infinite dimensional globes should be the subcomplexes such that
they are projected to the molecules in the (twisted) products of
three infinite dimensional globes. This results in the basic
definition for pairwise molecular subcomplexes in the product of
four infinite dimensional globes.

To work out the more explicit description (the one without using
projection), some basic concepts, for example, the definition of
{\em adjacency} and {\em projection maximal} must be reestablished
because of another middle factor. Compared with the description
for the molecules in the product of three infinite dimensional
globes, this description is less explicit. However, it is good
enough to check whether a given subcomplex in the product of four
infinite dimensional globes is a molecule. Best of all, both
descriptions for the molecules in the product of four infinite
dimensional globes can easily be stated for those in the product
of more infinite dimensional globes. This may lead us to further
study the product of more infinite dimensional globes.

After the basic concepts and tools are properly established, the
rest of the work very much resembles  that for dealing with the
product of three infinite dimensional globes: the closedness of
molecular subcomplexes under the source, target and composition
operations are proved, algorithms for the calculations of
sources, targets and composites are given, and in the proof that
molecular subcomplexes are exactly molecules, an algorithm for
decomposing molecules into atoms is also established.

Unfortunately, we have not been able to work out the construction
of molecules in the product of four infinite dimensional globes.
The difficulty remains that there are two `middle' factors. Thus
our theory for the product of four infinite dimensional globes is
not as satisfactory as that for the product of three infinite
dimensional globes.

We end the introduction by raising some questions following this
work.

1. What are the explicit descriptions for the product of $n$
infinite dimensional globes with $n>4$?

We have proposed some fairly reasonable explicit descriptions for
the molecules in the product of $n$ infinite dimensional globes
which resembles very much that for the product of four infinite
dimensional globes. Some proofs in the study of the product of
four infinite dimensional globes are already quite complicated,
and the problem is how to generalise them. We feel pretty
confident about working this out.

2. How can one construct the molecules in the product of four
infinite dimensional globes?

As stated above, the construction of the product of three infinite
dimensional globes is satisfactory because there is systematic way
to construct any molecule in the product of three infinite
dimensional globes. However, we have not yet worked out the
analogue for the product of four infinite dimensional globes. The
difficulty remains how to handle the two `middle' factors. We
still have no idea of what the construction should look like.

3. What about the explicit descriptions for the product of general
$\omega$-complexes.

As stated at the beginning of the introduction, the study of the
product of infinite dimensional globes may help to understand
products of general $\omega$-complexes. Following this idea, for
example, the construction of such $\omega$-complexes as
$(u_{1}\#_{p}u_{2})\times v\times w$, where $u_{1}$, $u_{2}$, $v$
and $w$ are infinite dimensional globes, must firstly be studied
before one can carry on the study for the general problem.

The following are  two questions which we have not had time to
think of deeply.

4. What about the construction of the joins of infinite
dimensional globes (simplexes)?

5. What about  the product of globes in the weak $n$-categories or
weak $\omega$-categories? (for the definition of weak
$n$-categories and weak $\omega$-categories, see \cite{baez},
\cite{baez-dolan}, \cite{bat} \cite{street2} and \cite{tam}).

\pagenumbering{arabic}
\chapter{Preliminaries}

In this chapter, we give some basic definitions and discuss some
properties of $\omega$-complexes and products of
$\omega$-complexes which are used throughout the thesis. All the
results are based on papers \cite{rjs} and \cite{rjs2}, although
some treatments are different from those in these two papers. In
the last section, we give a new decomposition theorem which will
be used later in the thesis.

Throughout the thesis, non-negative integers are denoted by $i$,
$j$, $k$, $l$, $m$, $n$, $p$, $q$, etc. We also use $\alpha$,
$\beta$, $\gamma$, $\sigma$, $\tau$, $\varepsilon$, $\omega$, etc,
to denote signs $\pm$.

\section{$\omega$-complexes}

In this section, we define $\omega$-complexes and give some local
descriptions of  $\omega$-complexes.

It is well known that a small category can be described purely by
its morphism set by regarding  objects as  identities.

Informally, an $\omega$-category is a set $X$ which forms the
morphism set of a small category $C_{n}$ for every non-negative
integer $n$ such that every element $x$ in $X$ is an identity in
some $C_{n}$ and $ob(C_{0})\subset ob(C_{1})\subset\dots$, where
$ob(C_{n})$ denote the set of objects (identities) of $C_{n}$. We
also require that the categorical structures commute for every
pair of non-negative integers. The formal definition is as
follows.

\begin{definition}
A partial $\omega$-category is a set $X$ together with unary
operations $d_{0}^{-}$, $d_{0}^{+}$, $d_{1}^{-}$, $d_{1}^{+}$,
$\dots$ and not everywhere defined binary operations $\#_{0}$,
$\#_{1}$, $\dots$ on $X$ such that the following conditions hold
for all elements $x$, $x'$, $y$, $y'$ and $z$ in $X$, all
non-negative integers $m$ and~$n$, and all signs $\alpha$ and
$\beta$:
\begin{enumerate}
\item
if $x\#_{n}y$ is defined, then $d_{n}^{+}x=d_{n}^{-}y$;
\item
$$ d_{m}^{\beta}d_{n}^{\alpha}x=
  \begin{cases}
    d_{m}^{\beta}x & \text{if $m<n$}, \\
    d_{n}^{\alpha}x & \text{if $m\geq n$};
  \end{cases}
$$
\item
$d_{n}^{-}x\#_{n}x=x\#_{n}d_{n}^{+}x=x$;
\item
if $x\#_{n}y$ is defined, then $$
\begin{array}{rl}
&d_{m}^{\alpha}(x\#_{n}y)=d_{m}^{\alpha}x=d_{m}^{\alpha}y \text{
for $m<n$},\\ &d_{n}^{-}(x\#_{n}y)=d_{n}^{-}x,\,\,\,
d_{n}^{+}(x\#_{n}y)=d_{n}^{+}y,\\
&d_{m}^{\alpha}(x\#_{n}y)=d_{m}^{\alpha}x\#_{n}d_{m}^{\alpha}y
\text{ for $m>n$};
\end{array}
$$
\item
$(x\#_{n}y)\#_{n}z=x\#_{n}(y\#_{n}z)$ if either side is defined;
\item
$(x\#_{n}y)\#_{m}(x'\#_{n}y')=(x\#_{m}x')\#_{n}(y\#_{m}y')$ if
$m<n$ and the left side is defined;
\item \label{dim}
for every $x\in X$ there is a non-negative integer $p$ such that
$d_{n}^{\alpha}x=x$ if and only if $n\geq p$.

\end{enumerate}

\end{definition}

The unique non-negative integer $p$ in condition \ref{dim} is
called the {\em dimension} of $x$, denoted by $\dim x$.

\begin{example}\label{partial_ex}
There is a partial $\omega$-category $X=\{a, b, x, y\}$ such that
$\dim a =\dim b=0$, $\dim x=\dim y=1$, $d_{0}^{-}x=d_{0}^{+}y=a$
and $d_{0}^{+}x=d_{0}^{-}y=b$.
\end{example}

\begin{definition}
Let $X$ be a partial $\omega$-category. If
$d_{n}^{+}x=d_{n}^{-}y$ implies that $x\#_{n}y$ is defined for all
elements $x$,~$y$ in $X$ and for all non-negative integers $n$,
then $X$ is an {\em $\omega$-category}.
\end{definition}

From Example \ref{partial_ex}, a partial $\omega$-category is not
necessarily an $\omega$-category.

It is natural to consider representing a partial $\omega$-category
$X$ by a suitable `pasting diagram'. The `pasting diagram' is a
sort of cell complex such that the indecomposable elements of $X$
are represented by atoms, the operations $d_{n}^{\alpha}$ are
represented by parts of boundaries, composites are represented by
well behaved unions, elements in the $\omega$-category are
represented by  subcomplexes which are well-behaved unions of
atoms.

\begin{example}\label{first_ex}
There is an $\omega$-category $X$ with the following presentation:
there are generators $a$, $x$,~$y$ and relations $$\dim a=1,\ \dim
x=\dim y=2,\ d^+_1 x=d^-_1 y,\ d^+_0 a=d^-_0 x=d^-_0 y.$$ Then $X$
has $16$~elements which can be represented by subcomplexes of the
diagram in the following figure: $$
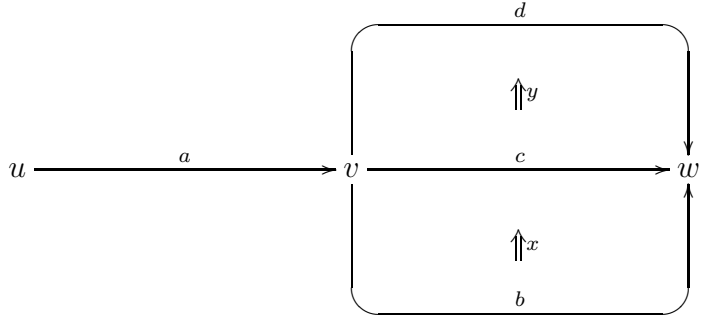
\xymatrix@R=4mm{ &&&&&&&&&
\\
&&&&&&
\\
&&&&&& \ar @{=>} [u] _y
\\
u \ar [rrrr] ^a &&&& v \ar `u [uuurrrr] `[rrrr] ^d [rrrr] \ar
[rrrr] ^c \ar `d [dddrrrr] `[rrrr] ^b [rrrr] &&&& w
\\
&&&&&&
\\
&&&&&& \ar @{=>} [u] _x
\\
&&&&&&&& }$$ There are three cells $a$, $x$,~$y$ representing the
generators; three additional $0$-cells $u$, $v$,~$w$ representing
$d^-_0 a$, $d^+_0 a=d^-_0 x=d^-_0 y$ and $d^+_0 x=d^+_0 y$; three
additional $1$-cells $b$, $c$,~$d$ representing $d^-_1 x$, $d^+_1
x=d^-_1 y$ and $d^+_1 y$; and the seven subcomplexes $$x\cup y,\
a\cup b,\ a\cup c,\ a\cup d,\ a\cup x,\ a\cup y,\ a\cup x\cup y$$
representing $$x\#_1 y,\ a\#_0 b,\ a\#_0 c,\ a\#_0 d,\ a\#_0 x,\
a\#_0 y,\ a\#_0(x\#_1 y).$$ In this figure, $d^-_0$, $d^+_0$,
$d^-_1$, $d^+_1$ are represented by left end, right end, bottom
and top respectively; for example, $d^+_1 a=a$ because $\dim a=1$,
and $$d^+_1[a\#_0(x\#_1 y)]=d^+_1 a\#_0 d^+_1(x\#_1 y)=a\#_0 d^+_1
y =a\#_0 d.$$

\end{example}

Suppose that $x\#_{n}y$ is a composite in a partial
$\omega$-category, and suppose that $x$ and $y$ are represented by
complexes in a pasting diagram. We then have $d^+_n x=d^-_n y=z$,
say, and $z$ must be represented by a subcomplex of the
intersection $x\cap y$. In fact our intuition requires $z$ to be
the whole of $x\cap y$. For we want $z$ to be at one extreme
of~$x$ and at the opposite extreme of~$y$, so $x\setminus z$ and
$y\setminus z$ should be on opposite sides of~$z$, and therefore
disjoint.

For an example of what can go wrong if this requirement is not
satisfied, let $X$ be the partial $\omega$-category in Example
\ref{partial_ex}. This partial $\omega$-category can be
represented by the following diagram $$
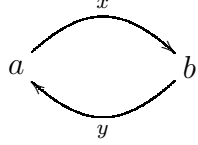
\xymatrix{ a \ar @(ur,ul)
[rr] ^x && b \ar @(dl,dr) [ll] ^y }$$ where the composites
$x\#_{0}y$ and $y\#_{0}x$ do not exist. We argue that the
composites like these, if exist, would lead to an unsatisfactory
behaviour in such pasting diagrams. Suppose otherwise that the
composites $x\#_{0}y$ and $y\#_{0}x$ both exist. They are distinct
because $d^-_0(x\#_{0}y)\neq d^-_0(y\#_{0}x)$, so it is not
satisfactory to have them both represented by the union $x\cup y$.
This unsatisfactory behaviour arises because $x\cap y$ strictly
contains~$a$ and strictly contains~$b$.

These considerations lead to the following definition.

\begin{definition}
An {\em $\omega$-complex} is a set $K$ together with a family of
subsets called atoms and a family of subsets called molecules such
that the following conditions hold.
\begin{enumerate}
\item
The molecules form a partial $\omega$-category.
\item
Let $x$ and $y$ be molecules. Then $x\#_{n}y$ is defined if and
only if $x\cap y=d_{n}^{+}x=d_{n}^{-}y$; if $x\#_{n}y$ is defined,
then $x\#_{n}y=x\cup y$.
\item
Every atom is an molecule; every molecule is generated from some
atoms by applying composition operations $\#_{0}$, $\#_{1}$,
$\dots$.
\item
The set $K$ is the union of its atoms.
\item
For an atom $a$ and a sign $\alpha$, let $\partial^{\alpha} a$ be
given by $$\partial^{\alpha} a=
  \begin{cases}
    d_{p-1}^{\alpha}a & \text{if $\dim a=p>0$}, \\
    \emptyset & \text{if $\dim a=0$};
  \end{cases}
$$ let the {\em interior} of $a$ be the subset $\Int a$ given by
$$\Int a= a\setminus (\partial^{-}a\cup\partial^{+}a).$$ Then
interiors of atoms are non-empty and disjoint.

\end{enumerate}

\end{definition}

\begin{example}\label{globe23}
There is an $\omega$-complex $u_{2}$ called $2$-dimensional globe.
It is a closed disk. The boundary of the disk consists of two
semicircles $u_{1}^{-}$ and $u_{1}^{+}$ intersecting at two
distinct points $u_{0}^{-}$ and $u_{0}^{+}$. The atoms are $u_{2}$
itself, the two semicircles $u_{1}^{-}$ and $u_{1}^{+}$, and the
two distinct points $u_{0}^{-}$ and $u_{0}^{+}$. The operators
$d_{m}^{\beta}$ are such that $d_{m}^{\beta}u_{2}=u_{m}^{\beta}$
for $m<2$ and $d_{0}^{\beta}u_{1}^{\alpha}=u_{0}^{\beta}$. It is
easy to see that all the molecules in $u_{2}$ are atoms, and they
form an $\omega$-category. The $\omega$-complex $u_{2}$ can be
represented by the following diagram. $$ \xymatrixrowsep{3pc}
\xymatrixcolsep{3pc} 
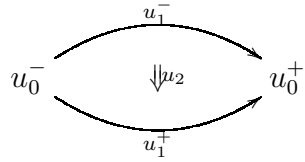
\diagram u_{0}^{-}
\rrtwocell<8>^{u_{1}^{-}}_{u_{1}^{+}}{u_{2}} &&u_{0}^{+}
\enddiagram
$$

Similarly, there is an $\omega$-complex $u_{3}$ called
$3$-dimensional globe. It is a closed $3$-dimensional ball. The
boundary sphere consists of two hemispheres $u_{2}^{-}$ and
$u_{2}^{+}$ intersecting in a circle, and the circle consists of
two semicircles $u_{1}^{-}$ and $u_{1}^{+}$ intersecting in two
distinct points $u_{0}^{-}$ and $u_{0}^{+}$. The atoms are the
ball $u_{3}$ itself, the two hemispheres $u_{2}^{-}$ and
$u_{2}^{+}$, the two semicircles $u_{1}^{-}$ and $u_{1}^{+}$, and
the two distinct points $u_{0}^{-}$ and $u_{0}^{+}$. The operators
$d_{m}^{\beta}$ are such that $d_{m}^{\beta}u_{3}=u_{m}^{\beta}$
for $m<3$ and $d_{m}^{\beta}u_{n}^{\alpha}=u_{m}^{\beta}$ for
$m<n<3$. It is easy to see that all the molecules in $u_{3}$ are
atoms, and they form an $\omega$-category.

As another example of $\omega$-complex, let $K$ be a $7$ element
set $\{e_{3}, e_{2}^{-}, e_{2}^{+}, e_{1}^{-}, e_{1}^{+},
e_{0}^{-}, e_{0}^{+}\}$. The atoms are $\bar{e}_{3}=\{e_{3},
e_{2}^{-}, e_{2}^{+},e_{1}^{-}, e_{1}^{+}, e_{0}^{-},e_{0}^{+}\}$,
$\bar{e}_{2}^{-}=\{e_{2}^{-},e_{1}^{-},e_{1}^{+},e_{0}^{-},e_{0}^{+}\}$,
$\bar{e}_{2}^{+}=\{e_{2}^{+},e_{1}^{-},e_{1}^{+},e_{0}^{-},e_{0}^{+}\}$,
$\bar{e}_{1}^{-}=\{e_{1}^{-},e_{0}^{-},e_{0}^{+}\}$,
$\bar{e}_{1}^{+}=\{e_{1}^{+},e_{0}^{-},e_{0}^{+}\}$,
$\bar{e}_{0}^{-}=\{e_{0}^{-}\}$, and
$\bar{e}_{0}^{+}=\{e_{0}^{+}\}$. The operators $d_{m}^{\beta}$ are
such that $d_{m}^{\beta}\bar{e}_{3}=\bar{e}_{m}^{\beta}$ for $m<3$
and $d_{m}^{\beta}\bar{e}_{n}^{\alpha}=\bar{e}_{m}^{\beta}$ for
$m<n<3$. It turns out that all the molecules in $K$ are atoms, and
they indeed form an $\omega$-category.
\end{example}

\begin{example}\label{globe_def}

There is an $\omega$-complex $u$ called $p$-dimensional globe such
that the atoms in $u$ can be listed as $u_{p}$, $u_{p-1}^{-}$,
$u_{p-1}^{+}$, $\dots$, $u_{0}^{-}$, $u_{0}^{+}$ such that
$d_{m}^{\beta}u_{p}=u_{m}^{\beta}$ for $m<p$ and
$d_{m}^{\beta}u_{n}^{\alpha}=u_{m}^{\beta}$ for $m<n<p$. It is
easy to check that all the molecules are atoms in $p$-dimensional
globes. We also denote the $p$-dimensional globe by $u_{p}$.

For instance, both of the subcomplexes $u_{3}$ and $K$ described
in Example \ref{globe23}  represent the $3$-dimensional globe. We
are going to see that they are equivalent.

Similarly, there is an $\omega$-complex $u$ called {\em infinite
dimensional globe} with exactly two $n$-dimensional atoms
$u_{n}^{-}$ and $u_{n}^{+}$ for every non-negative integer $n$,
such that $d_{m}^{\beta}u_{n}^{\alpha}=u_{m}^{\beta}$ for $m<n$.
It is easy to see that all the molecules in a globe are atoms, and
they form an $\omega$-category.

In the thesis, an atom $u_{n}^{\alpha}$ in an infinite dimensional
globe $u$ is also denoted by $u[n,\alpha]$.

\end{example}

We now state some results about local description of
$\omega$-complexes in \cite{rjs}.

\begin{prop}
\begin{enumerate}
\item
Let $x$ be a molecule in an $\omega$-complex. Then
$d_{n}^{\alpha}x\subset x$ for every sign $\alpha$ and every
non-negative integer $n$.

\item
Let $a$ be an atom in an $\omega$-complex. If
$\partial^{\alpha}a\neq\emptyset$, then  $\partial^{\alpha}a$ is a
molecule and $\dim\partial^{\alpha}a<\dim a$ for every sign
$\alpha$.

\end{enumerate}
\end{prop}

\begin{prop}\label{axi}
Let $\xi$ be an element in an $\omega$-complex. If $a$ is an atom
of minimal dimension such that $\xi\in a$, then $\xi\in \Int a$.
\end{prop}

\begin{prop}\label{axi2}
Let $x$ be a molecule and $a$ be an atom in an $\omega$-complex.
Then $a\subset x$ if and only if $\Int a\cap x\neq\emptyset$.
\end{prop}

\begin{prop}\label{dpgx}
Let $x$ be a molecule in an $\omega$-complex. Then
$$d_{n}^{\alpha}x =\bigcup\{a:\,\, a\subset x \text{ and } \dim
a\leq
n\}\setminus\bigcup\{b\setminus\partial^{\alpha}b:\,\,b\subset x
\text{ and } \dim b=n+1\},$$ where $a$ and $b$ are atoms.

\end{prop}

According to this proposition, we can see that an element $\xi\in
d_{n}^{\alpha}x$ if and only if (1) $\xi\in a$ for some atom
$a\subset x$ with $\dim a\leq n$; and (2) for every atom $b\subset
x$ with $\dim b=n+1$, if $\xi\in b$, then $\xi\in
d_{n}^{\alpha}b$.

As an example, we use Proposition \ref{dpgx} to verify that, in
Example, \ref{first_ex} $d_{1}^{-}(a\#_{0}(x\#_{1}y))=a\#_{0}b$.
By the above theorem, $d_{1}^{-}(a\#_{0}(x\#_{1}y))$ is the
difference of the union $a\cup b\cup c\cup d$ and $\Int d\cup \Int
c$. Thus $d_{1}^{-}(a\#_{0}(x\#_{1}y))=a\cup b=a\#_{0}b$.

\begin{coro}\label{by_atoms}
An $\omega$-complex is determined by its atoms , their dimensions
and the functions $\partial^{-}$ and $\partial^{+}$.
\end{coro}

\begin{definition}
Let $X$ and $Y$ be partial $\omega$-categories. A {\em
homomorphism} $f:X\to Y$  is a map such that
\begin{enumerate}
\item
$f(d_{n}^{\gamma}x)=d_{n}^{\gamma}f(x)$ for all $x\in X$, all
non-negative integers n and all signs $\gamma$;
\item
$f(x\#_{n}y)=f(x)\#_{n}f(y)$ whenever $x\#_{n}y$ is defined.
\end{enumerate}
\end{definition}

\begin{example}
Let $u$  be a infinite dimensional globes. Let $u_{p}$ be a  $p$
dimensional globe. It is evident that there is a homomorphism
$f_{p}^{u}:\M(u)\to \M(u_{p})$ of $\omega$-categories such that
for all atom $u_{i}^{\alpha}\in\M(u)$ $$f_{p}^{u}(u_{i}^{\alpha})=
  \begin{cases}
    u_{i}^{\alpha} & \text{when $i<p$}, \\
    u_{p} & \text{when $i\geq p$}.
  \end{cases}
$$
\end{example}

We end this section by introducing a definition of equivalence of
$\omega$-complexes.

Let $K$ be an $\omega$-complex. A {\em subcomplex} is a subset of
$K$ which can be written as a union of atoms. The set of all
subcomplexes of $K$ is denoted by $\C(K)$; the set of all atoms of
$K$ is denoted by $\A(K)$; The set of all molecules of $K$ is
denoted by $\M(K)$.

\begin{definition}\label{equivalent}
Let $K$ and $L$ be $\omega$-complexes. We say $K$ and $L$ are {\em
equivalent} if there exists a map  $f: \C(K)\to \C(L)$ called an
{\em equivalence of $\omega$-complexes} such that the following
conditions hold:

\begin{enumerate}
\item
If $a\in\A(K)$, then $f(a)\in \A(L)$. Moreover, $f|_{\A(K)}:
\A(K)\to\A(L)$ is a bijection.

\item
If $A$ is a set of atoms, then $f(\bigcup A)=\bigcup
\{f(a):\,\,a\in A\}$.

\item
If $a\in \A(K)$, then $\dim f(a)=\dim a$.

\item
If $a\in\A(K)$, then $f(\partial^{\alpha} a)=\partial^{\alpha}
f(a)$ for $\alpha=\pm$.

\end{enumerate}

\end{definition}

It is easy to check that the geometric description and
combinatorial descriptions for $3$-dimensional globes in Example
\ref{globe_def} are equivalent. From this, we may use the
geometric model to understand the combinatorial model and vice
versa.

We are going to prove that an equivalence of $\omega$-complexes
preserves molecules. We need several technical lemmas.

\begin{lemma}\label{intersection_subcomplex}
Let $K$ be an $\omega$-complex. If $c_{1},c_{2}\in\C(K)$, then
$c_{1}\cap c_{2}\in \C(K)$.
\end{lemma}

\begin{proof}
It suffices to prove that $a\cap b$ is a subcomplex of $K$  for
every pair $a$ and $b$ of atoms in $K$.

Let $\xi\in a\cap b$. Let $a_{\xi}$ be the atom of the minimal
dimension with $\xi\in a_{\xi}$. According to Propositions
\ref{axi} and \ref{axi2}, we have $a_{\xi}\subset a\cap b$. It
follows that $a\cap b=\bigcup\{a_{\xi}:\,\,\xi\subset a\cap b\}$.
Thus $a\cap b$ is a subcomplex of $K$, as required.
\end{proof}

\begin{lemma}
Let $f:\C(K)\to \C(L)$ be an equivalence of $\omega$-complexes. If
$c_{1}, c_{2}\in\C(K)$, then $f(c_{1}\cap c_{2})=f(c_{1})\cap
f(c_{2})$.
\end{lemma}
\begin{proof}
Since $f$ preserves unions of subcomplexes, we have
$f(c_{1})\subset f(c_{2})$ for a pair of subcomplexes $c_{1}$ and
$c_{2}$ in $K$ with $c_{1}\subset c_{2}$. Note that
$f:\C(K)\to\C(L)$ is a bijection, it follows easily from Lemma
\ref{intersection_subcomplex} that $f(c_{1}\cap
c_{2})=f(c_{1})\cap f(c_{2})$, as required.
\end{proof}

\begin{lemma}
Let $f:\C(K)\to \C(L)$ be an equivalence of $\omega$-complexes. If
$x\in\M(K)$ and $f(x)\in M(L)$, then
$f(d_{p}^{\gamma}x)=d_{p}^{\gamma}f(x)$.
\end{lemma}

\begin{proof}
Suppose that $b\in\A(L)$ and $\Int b\subset f(d_{p}^{\gamma}x)$.
Then there exists $a\in\A(K)$ with $b=f(a)$ such that $a\subset
d_{p}^{\gamma}x$. Thus $a\subset x$ and $\dim a\leq p$. It follows
that $b=f(a)\subset f(x)$ and $\dim b=\dim f(a)=\dim a\leq p$. Now
suppose that $b'\in\A(L)$ with $\dim b'=p+1$ such that $\Int
b\subset b'$. Then there exists $a'\in\A(K)$ such that $f(a')=b'$.
It is evident that $\dim a'=p+1$ and $a\subset a'$. So we have
$a\subset
\partial^{\gamma} a'$. This implies that
$b=f(a)\subset f(\partial^{\gamma} a')=\partial^{\gamma}
f(a')=\partial^{\gamma}b'$. According to Proposition \ref{dpgx},
we have $\Int b\subset d_{p}^{\gamma}f(x)$. It follows that
$f(d_{p}^{\gamma}x)\subset d_{p}^{\gamma}f(x)$.

By a similar argument, we can prove the reverse inclusion.

This completes the proof.
\end{proof}

\begin{prop}
Let $f:\C(K)\to \C(L)$ be an equivalence of $\omega$-complexes. If
$x\in\M(K)$, then $f(x)\in\M(L)$. Moreover, $f|_{\M(K)}:
\M(K)\to\M(L)$ is a homomorphism of partial $\omega$-categories.
\end{prop}
\begin{proof}
We give the proof by induction.

Firstly, if $a\in\M(K)$, then $f(a)\in \M(L)$ by the definition of
equivalence.

Suppose that $m>1$ and $f(x)\in\M(K)$ if $x$ can be written as a
composite of less then $m$ atoms. Let $x$ be an atom which can be
written as a composite of $m$ atoms. We must prove that
$f(x)\in\M(L)$.

Indeed, it is evident that $x$ has a proper decomposition
$x=y\#_{p}z$ into molecules such that $y$ and $z$ are molecules
which can be written as a composite of less than $p$ atoms. By the
inductive hypothesis, we have $f(y)\in \M(L)$ and $f(z)\in\M(L)$.
To prove $f(x)\in\M(L)$, it suffices to show that the composite
$f(y)\#_{p}f(z)$ exists and that $f(y\#_{p}z)=f(y)\#_{p}f(z)$.

Since $x=y\#_{p}z$, we have $d_{p}^{+}y=d_{p}^{-}z=y\cap z$. By
the previous lemmas, we get
$d_{p}^{+}f(y)=f(d_{p}^{+}y)=f(d_{p}^{-}z)=d_{p}^{-}f(z)=f(y)\cap
f(z)$. Therefore $f(y)\#_{p}f(z)$ is defined and
$f(y\#_{p}z)=f(y\cup z)=f(y)\cup f(z)=f(y)\#_{p}f(z)$, as
required.

By a similar argument, we can show that $f$ preserves composition
operation. Thus $f:\M(K)\to\M(L)$ is a homomorphism of partial
$\omega$-complexes.

This completes the proof.
\end{proof}

\section{Atomic Complexes}

Corollary \ref{by_atoms} shows that it is possible to describe an
$\omega$-complex by its atoms and the boundary operators
$\partial^{-}$ and $\partial^{+}$. This leads us to a concept
consisting of atoms and boundary operators which we call an atomic
complex.

In this section, we first define atomic complexes and state some
properties of atomic complexes. Then we state a necessary and
sufficient condition for an atomic complex to be an
$\omega$-complex. From this theorem, we will see that the results
in paper \cite{rjs2} for loop-free directed complexes can be
generalised to loop-free $\omega$-complexes. We shall discuss this
in section \ref{section1.4}.

\begin{definition}
An {\em atomic complex} is a set $K$ together with a family of
subsets $\A(K)$ called atoms and functions $\dim$, $\partial^{-}$
and $\partial^{+}$ defined on $\A(K)$ such that the following
conditions hold.

\begin{enumerate}
\item
For every atom $a$, $\dim a$ is an non-negative integer called the
{\em dimension} of $a$.

\item
If $a$ is an atom and $\alpha$ is a sign, then
$\partial^{\alpha}a$ is a subset of $a$ such that
$\partial^{\alpha}a$ is a union of $atoms$ of dimensions less than
$\dim a$.

\item
$K=\bigcup \A(K)$.

\item
For an atom $a$, let $\Int a=a\setminus
(\partial^{-}a\cup\partial^{+}a)$. Then the interiors of atoms are
non-empty and disjoint.

\end{enumerate}

\end{definition}

\begin{prop}
$\omega$-complexes are atomic complexes.
\end{prop}

To give the necessary and sufficient conditions for an atomic
complex to be an $\omega$-complex, we need to define operators
$d_{n}^{\alpha}$ on an arbitrary subset of $K$. This can be given
by generalising Proposition \ref{dpgx}.

\begin{definition}
Let $K$ be an atomic complex.
\begin{itemize}
\item \label{atomic_dpg}
If $x\subset K$ and $\alpha=\pm$, then $$d_{n}^{\alpha}x
=\bigcup\{a:\,\, a\subset x \text{ and } \dim a\leq
n\}\setminus\bigcup\{b\setminus\partial^{\alpha}b:\,\,b\subset x
\text{ and } \dim b=n+1\},$$ where $a$ and $b$ are atoms.

\item
If $x\subset K$ and $y\subset K$, then the composite $x\#_{n}y$ is
defined if and only if $x\cap y=d_{n}^{+}x=d_{n}^{-}y$; If
$x\#_{n}y$ is defined, then $x\#_{n}y=x\cup y$.

\item
A molecule is a subset generated from atoms by finitely applying
the composition operations $\#_{n}$ ($n=0,1,\dots$).

\end{itemize}
\end{definition}

With the definition of the operators $d_{n}^{\alpha}$ on an
arbitrary subsets of  $\omega$-complexes, we can define {\em
finite dimensional subcomplexes}.

\begin{definition}
Let $K$ be an atomic complex.
\begin{itemize}
\item
If $x$ is a union of atoms in $K$, then $x$ is a {\em subcomplex}
of $K$.
\item
Let $x$ be a subcomplex of $K$. If there exists an integer $n$
such that $x=d_{n}^{-}x=d_{n}^{+}x$, then $x$ is {\em finite
dimensional}.
\end{itemize}
\end{definition}

\begin{prop}
Let $a$ and $c$ be distinct atoms in an atomic complex. If $\Int
a\cap c\neq\emptyset$, then $\dim a\leq \dim c$ and $a\subset
\partial^{\alpha}c$ for some sign $\alpha$.
\end{prop}

\begin{prop}\label{dpg_sub}
Let $x$ and $y$ be subcomplexes of an atomic complex. If $y\subset
x$, then $d_{n}^{\alpha}x\cap y\subset d_{n}^{\alpha}y$.
\end{prop}

\begin{prop}\label{dpgunion}
Let $x$ and $y$ be subcomplexes of an atomic complex. Then
$d_{n}^{\alpha}(y\cup z)=(d_{n}^{\alpha}y\cap d_{n}^{\alpha}z)\cup
(d_{n}^{\alpha}y\setminus z)\cup (d_{n}^{\alpha}z\setminus y)$.
\end{prop}

Now we can state the necessary and sufficient conditions for an
atomic complex to be an $\omega$-complex.

\begin{theorem}\label{atomic_omega}
Let $K$ be an atomic complex. Then $K$ is an $\omega$-complex if
and only if the following conditions hold.
\begin{enumerate}
\item
If $a$ is an atom and $\dim a>0$, then $\partial^{\alpha}a$ is a
molecule for every sign $\alpha$.
\item
If $a$ is an atom and $\dim a=p>1$, then
$d_{p-2}^{\beta}d_{p-1}^{\alpha}a=d_{p-2}^{\beta}a$ for every pair
of signs $\alpha$ and $\beta$.
\end{enumerate}
\end{theorem}

\begin{example}
Let $w^{J}$ be the atomic complex with atoms
$w^{J}[k,\varepsilon]$ ($k=0,1,\cdots$ and $\varepsilon=\pm$) such
that $\dim\,w^{J}[k,\varepsilon]=k$ and
$d_{k-1}^{\gamma}w^{J}[k,\varepsilon]=w^{J}[k-1, (-)^{J}\gamma]$
for $k>0$. It is clear that $w^{J}$ satisfies conditions in
Theorem \ref{atomic_omega}. Thus it is an $\omega$-complex. It is
also easy to see that the $\omega$-complex $w^{J}$ is equivalent
to infinite dimensional globe $w$ under an obvious equivalence of
$\omega$-complexes sending $w^{J}[k,(-)^{J}\varepsilon]$ to
$w[k,\varepsilon]$.
\end{example}

\begin{lemma} \label{interior}
In an $\omega$-complex, if x is a subcomplex, then
$d_{n}^{\alpha}x$ can be written as a union of interior of atoms.
\end{lemma}
\begin{proof}
Suppose that $\xi\in\dna x$. Let $a_{\xi}$ be the atom in $x$ such
that $\xi\in\Int a_{\xi}$. Then $\dim a_{\xi}\leq p$. We claim
$\Int a_{\xi}\subset\dna x$. Indeed, for every $\eta\in\Int
a_{\xi}$, we have $\eta\in a_{\xi}$. Moreover, suppose that
$\eta\in b$ for an atom $b\subset x$ with $\dim b=p+1$, then
$\xi\in a_{\xi}\subset b$. Hence $\xi\in\dna b$ by Definition
\ref{atomic_dpg}. Since $\dna b$ is a molecule, we have
$a_{\xi}\subset\dna b$. Therefore $\eta\in\dna b$. It follows from
Definition \ref{atomic_dpg} that $\eta\in\dna x$ for every
$\eta\in\Int a_{\xi}$. Therefore $\Int a_{\xi}\subset\dna x$.

Now it is evident that $$\dna x=\bigcup_{\xi\in\dna x} \Int
a_{\xi}$$ which shows that $\dna x$ is a union of interiors of
atoms.

This completes the proof.
\end{proof}

\begin{lemma}\label{lm5.2}
In an $\omega$-complex, let $x$ be a subcomplex and $a$ be an
atom. Then $\Int a\subset \dna x$ if and only if
\begin{enumerate}
\item
$a\subset x$ and $\dim a\leq p$;
\item
If $a\subset b\subset x$ for an atom $b\subset x$ with $\dim
b=p+1$, then $a\subset\dna b$.
\end{enumerate}
\end{lemma}
\begin{proof}
This is an direct consequence of Lemma \ref{interior}.
\end{proof}

\section{Products of $\omega$-complexes}\label{prm3}
In this section, we give the product construction for
$\omega$-complexes. Some treatments are different from (but of
course equivalent to) those in paper \cite{rjs}. This applies in
particular to structures of products of two infinite dimensional
globes.

\begin{prop}
Let $K$ and $L$ be atomic complexes. Then the product  $K\times L$
of sets is made into an atomic complex as follows. The atom set
$\A(K\times L)$ is given by $$\A(K\times L)=\{a\times b:\,\,
a\in\A(K)\text{ and } b\in\A(L)\};$$ the structure functions are
given by
\begin{itemize}
\item $\dim (a\times b)=\dim a+\dim b$;
\item $\partial^{\gamma}(a\times b)=(\partial^{\gamma}a\times b)\cup
(a\times \partial^{(-)^{\dim a}\gamma}b)$.
\end{itemize}
\end{prop}

The atomic complex $K\times L$ is called the {\em product} of $K$
and $L$.

\begin{theorem}\label{complex_product}
Let $K$ and $L$ be $\omega$-complexes. Then $K\times L$ is an
$\omega$-complex.
\end{theorem}

\begin{example}
Let u and v be infinite dimensional globes. Then the product
$u\times v$ of sets is made into an atomic complex as follows: the
atoms are of the form $u_{i}^{\alpha}\times v_{j}^{\beta}$ with
$i$ and $j$ run over all non-negative integers, and $\alpha$ and
$\beta$ run over all the signs; the structure functions are given
by
\begin{itemize}
\item $\dim (u_{i}^{\alpha}\times v_{j}^{\beta})=i+j$;
\item $\partial^{\gamma}(u_{i}^{\alpha}\times v_{j}^{\beta})
=(u_{i-1}^{\gamma}\times v_{j}^{\beta})\cup (u_{i}^{\alpha}\times
v_{j-1}^{(-)^{i}\gamma})$.
\end{itemize}
By Theorem \ref{complex_product}, $u\times v$ is actually an
$\omega$-complex.
\end{example}

It is straightforward to verify that the product construction is
associative. Therefore we can write a product of three
$\omega$-complexes $K$, $L$ and $M$ as $K\times L\times M$. By
Theorem \ref{complex_product}, The product $K\times L\times M$ is
still an $\omega$-complex. In particular, the atom set $\A(K\times
L\times M)$ is given by $$\A(K\times L\times M)=\{a\times b\times
c:\,\, a\in\A(K)\text{ and } b\in\A(L)\text{ and }b\in\A(M)\};$$
the structure functions are given by
\begin{itemize}
\item $\dim (a\times b\times c)=\dim a+\dim b+\dim c$;
\item $\partial^{\gamma}(a\times b\times c)=(\partial^{\gamma}a\times b\times c)\cup
(a\times \partial^{(-)^{\dim a}\gamma}b\times c)\cup (a\times
b\times \partial^{(-)^{\dim a+\dim b}\gamma}c)$.
\end{itemize}

\begin{example}
We now consider the product of three $1$-dimensional globes
$u_{1}\times v_{1}\times w_{1}$. Since the $1$-dimensional globe
is represented by the closed interval, the product $u_{1}\times
v_{1}\times w_{1}$ is a cube. Recall that the $1$-dimensional
globe consists of $3$ atoms. So the product $u_{1}\times
v_{1}\times w_{1}$ consists of $27$ atoms. The following figure
illustrates the source boundary $$\partial^{-}(u_{1}\times
v_{1}\times w_{1}) =(u_{1}\times v_{1}\times w_{0}^{-}) \cup
(u_{0}^{-}\times v_{1}\times w_{1}) \cup (u_{1}\times
v_{0}^{+}\times w_{1}) $$ of the cube, where $A=u_{1}\times
v_{0}^{+}\times w_{1}$, $B=u_{1}\times v_{1}\times w_{0}^{-}$ and
$C=u_{0}^{-}\times v_{1}\times w_{1}$. $$ \xymatrixrowsep{3pc}
\xymatrixcolsep{3pc} 
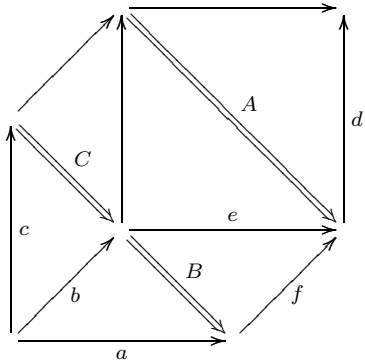
\diagram &\relax \rrtwocell<0>{\omit} \ar
@{=>} [ddrr] ^{A} &&\relax
\\
\ar @{=>} [dr] ^{C} \urtwocell<0>{\omit} &&
\\
&\uutwocell<0>{\omit} \ar @{=>} [dr] ^{B} \rrtwocell<0>^{e}{\omit}
&&\uutwocell<0>_{d}{\omit}
\\
\uutwocell<0>_{c}{\omit} \urtwocell<0>_{b}{\omit}
\rrtwocell<0>_{a}{\omit} &&\relax \urtwocell<0>_{f}{\omit}
\enddiagram
$$ We can identify edges and vertices. For example, the edge
$$b=B\cap C= (u_{1}\times v_{1}\times w_{0}^{-})\cap
(u_{0}^{-}\times v_{1}\times w_{1}) =u_{0}^{-}\times v_{1}\times
w_{0}^{-}$$ and the vertex $$a\cap b = (u_{1}\times
v_{0}^{-}\times w_{0}^{-})\cap (u_{0}^{-}\times v_{1}\times
w_{0}^{-})=u_{0}^{-}\times v_{0}^{-}\times w_{0}^{-}.$$ We can
then check that the directions of the edges and vertices are as
shown in the figure. For example, since
$$\partial^{-}b=\partial^{-}(u_{0}^{-}\times v_{1}\times
w_{0}^{-}) =u_{0}^{-}\times v_{0}^{-}\times w_{0}^{-}=a\cap b$$
and $$\partial^{+}b=\partial^{+}(u_{0}^{-}\times v_{1}\times
w_{0}^{-}) =u_{0}^{-}\times v_{0}^{+}\times w_{0}^{-}=b\cap e,$$
the direction of $b$ is as shown in the figure. Similarly, since
$$\partial^{-}B=\partial^{-}(u_{1}\times v_{1}\times w_{0}^{-})
=(u_{0}^{-}\times v_{1}\times w_{0}^{-}) \cup (u_{1}\times
v_{0}^{+}\times w_{0}^{-}) =b\cup e $$ and
$$\partial^{+}B=\partial^{+}(u_{1}\times v_{1}\times w_{0}^{-})
=(u_{0}^{+}\times v_{1}\times w_{0}^{-}) \cup (u_{1}\times
v_{0}^{-}\times w_{0}^{-}) =a\cup f $$ Thus the direction of
$B=u_{1}\times v_{1}\times w_{0}^{-}$ is as shown in the figure.

From the diagram of $u_{1}\times v_{1}\times w_{1}$, we can see
that all the subcomplexes in the following list are molecules.

1. $u_{1}\times v_{1}\times w_{1}$,

2. $u_{1}\times v_{1}\times w_{0}^{-}$,

3. $u_{1}\times v_{1}\times w_{0}^{+}$,

4. $u_{1}\times v_{0}^{-}\times w_{1}$,

5. $u_{1}\times v_{0}^{+}\times w_{1}$,

6. $u_{0}^{-}\times v_{1}\times w_{1}$,

7. $u_{0}^{+}\times v_{1}\times w_{1}$,

8. $u_{1}\times v_{0}^{-}\times w_{0}^{-}$,

9. $u_{1}\times v_{0}^{-}\times w_{0}^{+}$,

10. $u_{1}\times v_{0}^{+}\times w_{0}^{-}$,

11. $u_{1}\times v_{0}^{+}\times w_{0}^{+}$,

12. $u_{0}^{-}\times v_{1}\times w_{0}^{-}$,

13. $u_{0}^{-}\times v_{1}\times w_{0}^{+}$,

14. $u_{0}^{+}\times v_{1}\times w_{0}^{-}$,

15. $u_{0}^{+}\times v_{1}\times w_{0}^{+}$,

16. $u_{0}^{-}\times v_{0}^{-}\times w_{1}$,

17. $u_{0}^{-}\times v_{0}^{+}\times w_{1}$,

18. $u_{0}^{+}\times v_{0}^{-}\times w_{1}$,

19. $u_{0}^{+}\times v_{0}^{+}\times w_{1}$,

20. $u_{0}^{-}\times v_{0}^{-}\times w_{0}^{-}$,

21. $u_{0}^{-}\times v_{0}^{-}\times w_{0}^{+}$,

22. $u_{0}^{-}\times v_{0}^{+}\times w_{0}^{-}$,

23. $u_{0}^{-}\times v_{0}^{+}\times w_{0}^{+}$,

24. $u_{0}^{+}\times v_{0}^{-}\times w_{0}^{-}$,

25. $u_{0}^{+}\times v_{0}^{-}\times w_{0}^{+}$,

26. $u_{0}^{+}\times v_{0}^{+}\times w_{0}^{-}$,

27. $u_{0}^{+}\times v_{0}^{+}\times w_{0}^{+}$,

28. $(u_{1}\times v_{1}\times w_{0}^{+}) \cup (u_{0}^{+}\times
v_{1}\times w_{1}) \cup (u_{1}\times v_{0}^{-}\times w_{1})$,

29. $(u_{1}\times v_{1}\times w_{0}^{-}) \cup (u_{0}^{-}\times
v_{1}\times w_{1}) \cup (u_{1}\times v_{0}^{+}\times w_{1})$,

30. $(u_{1}\times v_{1}\times w_{0}^{+}) \cup (u_{1}\times
v_{0}^{-}\times w_{1})$,

31. $(u_{1}\times v_{1}\times w_{0}^{-}) \cup (u_{1}\times
v_{0}^{+}\times w_{1})$,

32. $(u_{1}\times v_{1}\times w_{0}^{+}) \cup (u_{0}^{-}\times
v_{0}^{-}\times w_{1})$,

33. $(u_{1}\times v_{1}\times w_{0}^{-}) \cup (u_{0}^{+}\times
v_{0}^{+}\times w_{1})$,

34. $(u_{0}^{+}\times v_{1}\times w_{1}) \cup (u_{1}\times
v_{0}^{-}\times w_{1})$,

35. $(u_{0}^{-}\times v_{1}\times w_{1}) \cup (u_{1}\times
v_{0}^{+}\times w_{1})$,

36. $(u_{0}^{+}\times v_{1}\times w_{1}) \cup (u_{1}\times
v_{0}^{-}\times w_{0}^{-})$,

37. $(u_{0}^{-}\times v_{1}\times w_{1}) \cup (u_{1}\times
v_{0}^{+}\times w_{0}^{+})$,

38. $(u_{0}^{+}\times v_{1}\times w_{0}^{+}) \cup (u_{1}\times
v_{0}^{-}\times w_{1})$,

39. $(u_{0}^{-}\times v_{1}\times w_{0}^{-}) \cup (u_{1}\times
v_{0}^{+}\times w_{1})$,

40. $(u_{0}^{+}\times v_{1}\times w_{0}^{+}) \cup (u_{1}\times
v_{0}^{-}\times w_{0}^{+})$,

41. $(u_{0}^{-}\times v_{1}\times w_{0}^{+}) \cup (u_{1}\times
v_{0}^{+}\times w_{0}^{+})$,

42. $(u_{0}^{+}\times v_{1}\times w_{0}^{-}) \cup (u_{1}\times
v_{0}^{-}\times w_{0}^{-})$,

43. $(u_{0}^{-}\times v_{1}\times w_{0}^{-}) \cup (u_{1}\times
v_{0}^{+}\times w_{0}^{-})$,

44. $(u_{0}^{-}\times v_{1}\times w_{0}^{+}) \cup (u_{0}^{-}\times
v_{0}^{-}\times w_{1})$,

45. $(u_{0}^{-}\times v_{1}\times w_{0}^{-}) \cup (u_{0}^{-}\times
v_{0}^{+}\times w_{1})$,

46. $(u_{0}^{+}\times v_{1}\times w_{0}^{+}) \cup (u_{0}^{+}\times
v_{0}^{-}\times w_{1})$,

47. $(u_{0}^{+}\times v_{1}\times w_{0}^{-}) \cup (u_{0}^{+}\times
v_{0}^{+}\times w_{1})$,

48. $(u_{0}^{+}\times v_{1}\times w_{0}^{-}) \cup (u_{1}\times
v_{0}^{-}\times w_{0}^{-}) \cup (u_{0}^{+}\times v_{0}^{+}\times
w_{1})$,

49. $(u_{0}^{+}\times v_{1}\times w_{0}^{+}) \cup (u_{1}\times
v_{0}^{-}\times w_{0}^{-}) \cup (u_{0}^{+}\times v_{0}^{-}\times
w_{1})$,

50. $(u_{0}^{-}\times v_{1}\times w_{0}^{-}) \cup (u_{1}\times
v_{0}^{+}\times w_{0}^{+}) \cup (u_{0}^{-}\times v_{0}^{+}\times
w_{1})$,

51. $(u_{0}^{-}\times v_{1}\times w_{0}^{+}) \cup (u_{1}\times
v_{0}^{+}\times w_{0}^{+}) \cup (u_{0}^{-}\times v_{0}^{-}\times
w_{1})$,

52. $(u_{0}^{+}\times v_{1}\times w_{0}^{+}) \cup (u_{1}\times
v_{0}^{-}\times w_{0}^{+}) \cup (u_{0}^{-}\times v_{0}^{-}\times
w_{1})$,

53. $(u_{0}^{-}\times v_{1}\times w_{0}^{-}) \cup (u_{1}\times
v_{0}^{+}\times w_{0}^{-}) \cup (u_{0}^{+}\times v_{0}^{+}\times
w_{1})$,

54. $(u_{1}\times v_{0}^{-}\times w_{0}^{+}) \cup (u_{0}^{-}\times
v_{0}^{-}\times w_{1})$,

55. $(u_{1}\times v_{0}^{-}\times w_{0}^{-}) \cup (u_{0}^{+}\times
v_{0}^{-}\times w_{1})$,

56. $(u_{1}\times v_{0}^{+}\times w_{0}^{+}) \cup (u_{0}^{-}\times
v_{0}^{+}\times w_{1})$,

57. $(u_{1}\times v_{0}^{+}\times w_{0}^{-}) \cup (u_{0}^{+}\times
v_{0}^{+}\times w_{1})$,

For example, from the figure, one can see that the 31st subcomplex
$A\cup B=(u_{1}\times v_{1}\times w_{0}^{-}) \cup (u_{1}\times
v_{0}^{+}\times w_{1})$ in the list can be decomposed into atoms
as $(b\#_{0}A)\#_{1}(B\#_{0}d)=[(u_{0}^{-}\times v_{1}\times
w_{0}^{-})\#_{0}(u_{1}\times v_{0}^{+}\times
w_{1})]\#_{1}[(u_{1}\times v_{1}\times w_{0}^{-})
\#_{0}(u_{0}^{+}\times v_{0}^{+}\times w_{1})]$.

One can show that every composite of molecules in the list is
still a molecule in the list. So we have a complete list of the
molecules in $u_{1}\times v_{1}\times w_{1}$.

In chapter 3, we will show how this list is compiled and how to
compile such lists for the molecules in the products of any three
finite dimensional globes.

\end{example}

\begin{theorem}\label{dt_product}\label{dpg_product}
Let $K$ and $L$ be $\omega$-complexes. Let $x$ and $y$ be
molecules in $K$ and $L$ respectively. Then $x\times y$ is a
molecule in $K\times L$ and $$d_{n}^{\gamma}(x\times y)=
(d_{n}^{\gamma}x\times d_{0}^{(-)^{n}\gamma}y)\cup
(d_{n-1}^{\gamma}x\times d_{1}^{(-)^{n-1}\gamma}y)\cup\cdots\cup
(d_{0}^{\gamma}x\times d_{n}^{\gamma}y).$$
\end{theorem}

The following Theorem is implicit in Paper \cite{al_rjs1}. To
avoid introducing more concepts, we give an independent proof.

\begin{theorem}\label{natural_homo}
Let $K_{i}$ and $L_{i}$ be $\omega$-complexes such that
$\M(K_{i})$ and $\M(L_{i})$ are $\omega$-categories for $1\leq
i\leq r$. Let $f_{i}:\M(K_{i})\to \M(L_{i})$ be homomorphisms of
partial $\omega$-categories for $1\leq i\leq r$. If
$\M(K_{1}\times \cdots\times K_{r})$ and
$\M(L_{1}\times\cdots\times L_{r})$ are $\omega$-categories, then
there is a {\em natural homomorphism}  $f: \M(K_{1}\times
\cdots\times K_{r})\to \M(L_{1}\times\cdots\times L_{r})$ of
partial $\omega$-categories such that $f(a_{1}\times \cdots\times
a_{r}) =f_{1}(a_{1})\times\cdots\times f_{r}(a_{r})$ for all atoms
$a_{1}\times\cdots\times a_{r}$ in $K_{1}\times\cdots\times
K_{r}$.
\end{theorem}
\begin{proof}
The arguments for different choices of $r$ are similar. We give
the proof for $r=2$.

Let $F:\C(K_{1}\times K_{2})\to \C(L_{1}\times L_{2})$ be the
union-preserving map such that $F(a_{1}\times
a_{2})=f_{1}(a_{1})\times f_{2}(a_{2})$ for all atoms $a_{1}$ and
$a_{2}$. To prove the theorem, it suffices to show that
$F(\M(K_{1}\times K_{2}))\subset\M(L_{1}\times L_{2})$ and
$F|_{\M(K_{1}\times K_{2})}:\M(K_{1}\times K_{2})\to M(L_{1}\times
L_{2})$ is a homomorphism of partial $\omega$-categories.

Firstly, we verify inductively that $F(x)$ is a molecule and
$F(d_{n}^{\gamma}x)=d_{n}^{\gamma}F(x)$ for all non-negative
integers $n$, all signs $\gamma$ and all molecules $x$ in
$K_{1}\times K_{2}$.

To begin the induction, let $a_{1}\times a_{2}$ be an atom in
$K_{1}\times K_{2}$. Then $F(a_{1}\times a_{2})=f_{1}(a_{1})\times
f_{2}(a_{2})$ by the definition of $F$. Since
$f_{i}(a_{i})\subset\M(K_{i})$, we have $F(a_{1}\times
a_{2})\subset\M(K_{1}\times K_{2})$ by Theorem \ref{dt_product}.
Moreover, by Theorem \ref{dt_product}, we have $$
\begin{array}{rl}
&F(d_{n}^{\gamma}(a_{1}\times a_{2}))\\
=&F((d_{n}^{\gamma}a_{1}\times d_{0}^{(-)^{n}\gamma}a_{2})\cup
(d_{n-1}^{\gamma}a_{1}\times
d_{1}^{(-)^{n-1}\gamma}a_{2})\cup\cdots\cup
(d_{0}^{\gamma}a_{1}\times d_{n}^{\gamma}a_{2}))\\
=&(f_{1}(d_{n}^{\gamma}a_{1})\times
f_{2}(d_{0}^{(-)^{n}\gamma}a_{2}))\cup
(f_{1}(d_{n-1}^{\gamma}a_{1})\times
f_{2}(d_{1}^{(-)^{n-1}\gamma}a_{2}))\cup\cdots\cup
(f_{1}(d_{0}^{\gamma}a_{1})\times f_{2}(d_{n}^{\gamma}a_{2}))\\
=&(d_{n}^{\gamma}f_{1}(a_{1})\times
d_{0}^{(-)^{n}\gamma}f_{2}(a_{2}))\cup
(d_{n-1}^{\gamma}f_{1}(a_{1})\times
d_{1}^{(-)^{n-1}\gamma}f_{2}(a_{2}))\cup\cdots\cup
(d_{0}^{\gamma}f_{1}(a_{1})\times d_{n}^{\gamma}f_{2}(a_{2}))\\
=&d_{n}^{\gamma}(f_{1}(a_{1})\times f_{2}(a_{2}))\\
=&d_{n}^{\gamma}F(a_{1}\times a_{2}).
\end{array}
$$ Thus $F(d_{n}^{\gamma}x)=d_{n}^{\gamma}F(x)$ holds when $x$ is
an atom in $K_{1}\times K_{2}$.

Next, suppose that $F(x')$ is a molecule and that
$F(d_{n}^{\gamma}x')=d_{n}^{\gamma}F(x')$ for every  molecule $x'$
in $K_{1}\times K_{2}$ which can be written as a composite of less
than $q$ atoms. Suppose also that $x$ is a molecule in
$K_{1}\times K_{2}$ which can be written as a composite of $q$
atoms. We verify that $F(x)$ is a molecule in $L_{1}\times L_{2}$
and that $F(d_{n}^{\gamma}x)=d_{n}^{\gamma}F(x)$. It is evident
that $x$ can be decomposed into molecules $x=y\#_{p}z$ such that
$y$ and $z$ can be written as less than $q$ atoms.

We first verify that $F(x)$ is a molecule. Indeed, since
$d_{p}^{+}F(y)=F(d_{p}^{+}y)=F(d_{p}^{-}z)=d_{p}^{-}F(z)$, the
composite $F(y)\#_{p}F(z)$ is defined. Moreover, since $F$
preserves unions of atoms, we have $F(x)=F(y\cup z)=F(y)\cup
F(z)=F(y)\#_{p}F(z)$. This shows that $F(x)$ is a molecule in
$L_{1}\times L_{2}$.

We next verify that $F(d_{n}^{\gamma}x)=d_{n}^{\gamma}F(x)$.
Indeed, if $n=p$, then
$$F(d_{n}^{-}x)=F(d_{n}^{-}(y\#_{n}z))=F(d_{n}^{-}y)=d_{n}^{-}F(y)
=d_{n}^{-}F(x)$$ and, similarly, $F(d_{n}^{+}x)=d_{n}^{+}F(x)$. If
$n>p$, then $$
\begin{array}{rl}
&F(d_{n}^{\gamma}x)\\ =&F(d_{n}^{\gamma}(y\#_{p}z))\\
=&F(d_{n}^{\gamma}y\#_{p} d_{n}^{\gamma}z)\\
=&F(d_{n}^{\gamma}y\cup d_{n}^{\gamma}z)\\
=&F(d_{n}^{\gamma}y)\cup F(d_{n}^{\gamma}z)\\
=&d_{n}^{\gamma}F(y)\cup d_{n}^{\gamma}F(z)
\end{array}
$$ and $$d_{p}^{+}d_{n}^{\gamma}F(y)=d_{p}^{+}F(y)=d_{p}^{-}F(z)
=d_{p}^{-}d_{n}^{\gamma}F(z); $$ thus
$$F(d_{n}^{\gamma}x)=d_{n}^{\gamma}F(y)\#_{p}
d_{n}^{\gamma}F(z)=d_{n}^{\gamma}(F(y)\#_{p}F(z))
=d_{n}^{\gamma}F(x).$$ If $n<p$, then
$$F(d_{n}^{\gamma}x)=F(d_{n}^{\gamma}(y\#_{p}z))=F(d_{n}^{\gamma}y)
=d_{n}^{\gamma}F(y)=d_{n}^{\gamma}F(x).$$ Therefore, $F(x)$ is a
molecule and $F(d_{n}^{\gamma}x)=d_{n}^{\gamma}F(x)$.

This shows that $F(x)$ is a molecule and
$F(d_{n}^{\gamma}x)=d_{n}^{\gamma}F(x)$ for all molecules $x$ in
$K_{1}\times K_{2}$ by induction.

Finally, by arguments similar to that in the proof of $F(x)$ being
a molecule above, we can see that $F|_{\M(K_{1}\times K_{2})}$
preserves composites of molecules. This completes the proof that
$F|_{\M(K_{1}\times K_{2})}$ is a homomorphism of partial
$\omega$-categories, as required.

\end{proof}

Let $x$ be a subcomplex in an $\omega$-complex. An atom $a$ is a
{\em maximal} atom in $x$ if $a\subset x$ and $a\subset b\subset
x$ implies $a=b$ for every atom $b$.

The following result characterises molecules in the product of two
infinite dimensional globes.

\begin{theorem} \label{2globe}
Let $u$ and $v$ be infinite dimensional globes. Then a subcomplex
$\Lambda$ of $u\times v$ is a molecule if and only if the
following conditions hold.

\begin{itemize}
\item
There are no distinct maximal atoms $u[i_{1},\alpha_{1}]\times
v[j_{1},\beta_{1}]$ and $u[i_{2},\alpha_{2}]\times
v[j_{2},\beta_{2}]$ such that $i_{1}\leq i_{2}$ and $j_{1}\leq
j_{2}$, so that the maximal atom in $\Lambda$ can be listed as
$\lambda_{1}$, $\dots$, $\lambda_{S}$ with
$\lambda_{s}=u[i_{s},\alpha_{s}]\times v[j_{s},\beta_{s}]$ for
$1\leq s\leq S$ such that $i_{1}>\cdots>i_{S}$ and
$j_{1}<\cdots<j_{S}$.

\item
If $\lambda_{s-1}$ and $\lambda_{s}$ are a pair of consecutive
maximal atoms in the above list, then
$\beta_{s-1}=-(-)^{i_{s}}\alpha_{s}$.
\end{itemize}
\end{theorem}

Now we give the construction of $d_{p}^{\gamma}\Lambda$ for a
molecule $\Lambda$ in $u\times v$.

\begin{theorem}\label{dp2}
Let $\Lambda$ be a molecule in $u\times v$. Then the set of
maximal atoms in $d_{p}^{\gamma}\Lambda$ consists of all the
maximal atoms $u[i',\alpha']\times v[j',\beta']$ in $\Lambda$ with
$i'+j'<p$ and all the atoms $u[i,\alpha]\times v[j,\beta]$ with
$i+j=p$ such that $i\leq i''$ and $j\leq j''$ for some maximal
atom $u[i'',\alpha'']\times v[j'',\beta'']$ of $\Lambda$, where
the signs $\alpha$ and $\beta$ are determined as follows:

\begin{enumerate}
\item
If  $u[i'',\alpha'']\times v[j'',\beta'']$ can be chosen such that
$i''>i$, then $\alpha=\gamma$; otherwise, $\alpha=\alpha^{''}$.

\item
If  $u[i'',\alpha'']\times v[j'',\beta'']$ can be chosen such that
$j''>j$, then $\beta=(-)^{i}\gamma$; otherwise,
$\beta=\beta^{''}$.

\end{enumerate}
\end{theorem}

The composition of molecules in $u\times v$ is characterised as
follows.

\begin{theorem}\label{comp2}
Let $\Lambda^{-}$ and $\Lambda^{+}$ be molecules of $u\times v$.
If $d_{p}^{+}\Lambda^{-}=d_{p}^{-}\Lambda^{+}$, then the composite
$\Lambda^{-}\#_{p}\Lambda^{+}$ is defined and the maximal atoms in
$\Lambda^{-}\#_{p}\Lambda^{+}$ consists of all the $q$-dimensional
common maximal atoms of $\Lambda^{-}$ and $\Lambda^{+}$ with
$q\leq p$ together with all the $r$-dimensional atoms in either
$\Lambda^{-}$ and $\Lambda^{+}$ with $r>p$.
\end{theorem}

\begin{coro}\label{molecule_omega}
The molecules in $u\times v$ form an $\omega$-category.
\end{coro}

\begin{example}
Let the subcomplexes $$
\begin{array}{rl}
\Lambda^{-}= &u_{5}^{-}\times v_{0}^{+}\\ \cup&u_{4}^{-}\times
v_{2}^{+}\\ \cup&u_{2}^{-}\times v_{3}^{-}\\ \cup&u_{1}^{-}\times
v_{4}^{+}\\ \cup&u_{0}^{-}\times v_{5}^{+}
\end{array}
$$ and $$
\begin{array}{rl}
  \Lambda^{+}=&u_{6}^{+}\times v_{0}^{-}\\
\cup&u_{5}^{-}\times v_{1}^{+}\\ \cup&u_{3}^{+}\times v_{2}^{+}\\
\cup&u_{2}^{-}\times v_{4}^{+}\\ \cup&u_{0}^{-}\times v_{5}^{+}.
\end{array}
$$ By Theorem \ref{2globe}, it is easy to see that $\Lambda^{-}$
and $\Lambda^{+}$ are molecules of $u\times v$. Moreover, by
Theorem \ref{dp2}, we have $d_{5}^{+}\Lambda^{-}$ and
$d_{5}^{-}\Lambda^{+}$ are both equal to the molecule
$$\begin{array}{rl}
    &u_{5}^{-}\times v_{0}^{+}\\
\cup&u_{4}^{-}\times v_{1}^{+}\\ \cup&u_{3}^{+}\times v_{2}^{+}\\
\cup&u_{2}^{-}\times v_{3}^{-}\\ \cup&u_{1}^{-}\times v_{4}^{+}\\
\cup&u_{0}^{-}\times v_{5}^{+}.
\end{array}
$$ Therefore, by Theorem \ref{comp2}, the composite
$\Lambda^{-}\#_{5}\Lambda^{+}$ is defined and the composite is the
following molecule. $$\begin{array}{rl}
    &u_{6}^{+}\times v_{0}^{-}\\
\cup&u_{5}^{-}\times v_{1}^{+}\\ \cup&u_{4}^{-}\times v_{2}^{+}\\
\cup&u_{2}^{-}\times v_{4}^{+}\\ \cup&u_{0}^{-}\times v_{5}^{+}.
\end{array}
$$
\end{example}

\begin{example}\label{2by1globe}
Let $u_{2}$ be the $2$-dimensional globe and let $v_{1}$ be a
$1$-dimensional globe. Geometrically, $u_{2}$ is a closed disk and
$v_{1}$ is a closed interval. Therefore the product $u_{2}\times
v_{1}$ is a cylinder. Since $u_{2}$ has $5$ atoms and $v_{1}$ has
$3$ atoms, the product $u_{2}\times v_{1}$ has $15$ atoms.

We calculate the boundaries $\partial^{-}(u_{2}\times v_{1})$ and
$\partial^{+}(u_{2}\times v_{1})$. By definition, we have
$$\partial^{-}(u_{2}\times v_{1}) =(u_{1}^{-}\times v_{1})\cup
(u_{2}\times v_{0}^{-})$$ and $$\partial^{+}(u_{2}\times v_{1})
=(u_{1}^{+}\times v_{1})\cup (u_{2}\times v_{0}^{+}).$$ If we put
the disk $u_{2}$ in a horizontal plane and put the interval
$u_{1}$ in a vertical line and represent $d_{0}^{-}$ and
$d_{0}^{+}$ on $v_{1}$ by bottom and top respectively, then
$\partial^{-}(u_{2}\times v_{1})$ is the union of the bottom disk
and half of the curved part of the boundary of the product
$u_{1}^{-}\times v_{1}$, as shown in the following figure, where
$a=u_{1}^{-}\times v_{0}^{-}$, $b=u_{0}^{+}\times v_{1}$,
$A=u_{1}^{-}\times v_{1}$ and $B=u_{2}\times v_{0}^{-}$.
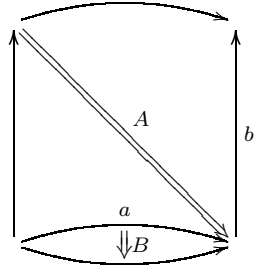
$$\xymatrixrowsep{3pc} \xymatrixcolsep{3pc} \diagram \ar @{=>}
[ddrr] ^{A} \relax \rrlowertwocell<3>^{}{\omit} &&\relax
\\
\\
\uutwocell<0>^{}{\omit} \relax\rrtwocell<3>^{a}_{}{B} &&\relax
\uutwocell<0>_{b}{\omit}
\enddiagram$$
Since $$\partial^{-}(u_{1}^{-}\times v_{1})=u_{0}^{-}\times
v_{1}\cup u_{1}^{-}\times v_{0}^{+}$$ and
$$\partial^{+}(u_{1}^{-}\times v_{1})=u_{0}^{+}\times v_{1}\cup
u_{1}^{-}\times v_{0}^{-},$$ one can easily see that the direction
of $A=u_{1}^{-}\times v_{1}$ is as indicated in the figure.
Similarly, we can get the direction for $B=u_{2}\times v_{0}^{-}$.
Moreover, it is easy to check that $$\partial^{-}(u_{2}\times
v_{1}) =(u_{1}^{-}\times v_{1})\cup (u_{2}\times v_{0}^{-})
=(u_{1}^{-}\times v_{1})\#_{1} [(u_{2}\times v_{0}^{-})\#_{0}
(u_{0}^{+}\times v_{1})].$$ One can also see this graphically from
the figure.

Similarly, $\partial^{+}(u_{2}\times v_{1})$ is the union of the
top disk and the other half of the curved part of the boundary of
the product $u_{1}^{+}\times v_{1}$, and we have
$$\partial^{+}(u_{2}\times v_{1}) =(u_{1}^{+}\times v_{1})\cup
(u_{2}\times v_{0}^{+}) =[(u_{0}^{-}\times v_{1})\#_{0}
(u_{2}\times v_{0}^{+})]\#_{1}(u_{1}^{+}\times v_{1}).$$ Therefore
$\partial^{-}(u_{2}\times v_{1})$ and $\partial^{+}(u_{2}\times
v_{1})$ are indeed molecules.

We can similarly workout the boundaries of other atoms.
\end{example}

\begin{example}
Let $u$  be a infinite dimensional globes. Let $u_{p}$ be a  $p$
dimensional globe. Recall that there is a homomorphism
$f_{p}^{u}:\M(u)\to \M(u_{p})$ of $\omega$-categories such that
for all atom $u_{i}^{\alpha}\in\M(u)$ $$f_{p}^{u}(u_{i}^{\alpha})=
  \begin{cases}
    u_{i}^{\alpha} & \text{when $i<p$}, \\
    u_{p} & \text{when $i\geq p$}.
  \end{cases}
$$ It follows from the Theorem \ref{natural_homo} and Corollary
\ref{molecule_omega} that there is a natural homomorphism
$f:\M(u\times v)\to\M(u_{p}\times v_{q})$ of $\omega$-categories
such that $f(u_{i}^{\alpha}\times
v_{j}^{\beta})=f_{p}^{u}(u_{i}^{\alpha})\times
f_{q}^{v}(v_{j}^{\beta})$ for all atoms $u_{i}^{\alpha}\times
v_{j}^{\beta}$ in $u\times v$.

\end{example}

\section{Decomposition of Molecules in Loop-Free
$\omega$-Complexes}\label{section1.4}

In this section, we  prove a decomposition theorem for molecules
in an loop-free $\omega$-complex. This theorem will be used later
in the thesis.

Firstly, we need to generalise some concepts and results from
loop-free directed complexes in paper \cite{rjs2} to loop-free
$\omega$-complexes.

\begin{definition}
A {\em directed precomplex} is a set $K$ together with functions
$\dim$, $\partial^{-}$ and $\partial^{+}$ on $K$ satisfying the
following conditions.

\begin{enumerate}
\item
If $\sigma\in K$, then $\dim\sigma$ is an non-negative integer,
called {\em dimension of $\sigma$}.
\item
If $\sigma\in K$ and $\dim\sigma>0$, then $\partial^{-}\sigma$ and
$\partial^{+}\sigma$ are subsets of $K$ consisting of
$\dim\sigma-1$ dimensional elements of $K$.
\end{enumerate}
\end{definition}

Let $K$ be a directed precomplex. A subset $x$ of $K$ is {\em
closed} if $\partial^{\alpha}\sigma\subset x$ for every $\sigma\in
x$ with $\dim\sigma>0$ and every sign $\alpha$. For a subset $y$
of $K$, the closure $\Cl(y)$ of $y$ is the smallest closed subset
of $K$ containing $y$. The closure $\Cl\{\sigma\}$ of a singleton
$\{\sigma\}$, denoted by $\bar{\sigma}$, is called an {\em atom}.

\begin{definition}
Let $K$ be a directed precomplex.

If $x\subset K$ and $\alpha=\pm$, then $$d_{n}^{\alpha}x
=\bigcup\{\sigma:\, \sigma\in x \text{ and } \dim \sigma\leq
n\}\setminus\bigcup\{\bar{\tau}\setminus
\Cl(\partial^{\alpha}\tau):\,\tau\in x \text{ and } \dim
\tau=n+1\}.$$

 If $x$ and $y$  are closed subsets of $K$, then the {\em composite} $x\#_{n}y$
is defined if and only if $x\cap y=d_{n}^{+}x=d_{n}^{-}y$; If
$x\#_{n}y$ is defined, then $x\#_{n}y=x\cup y$.

 A {\em molecule} is a subset generated from atoms by finitely applying the
composition operations $\#_{n}$ ($n=0,1,\dots$).

\end{definition}

\begin{definition}\label{directed}
A {\em directed complex} is a directed precomplex satisfying the
following conditions.
\begin{enumerate}
\item
If $\bar{\sigma}$ is an atom with $\dim\sigma=p>0$, then
$d_{p-1}^{\alpha}\bar{\sigma}$ is a molecule for $\alpha=\pm$.
\item
If $\bar{\sigma}$ is an atom with $\dim\sigma=p>1$, then
$d_{p-2}^{\beta}d_{p-1}^{\alpha}\bar{\sigma}=d_{p-2}^{\beta}\bar{\sigma}$
for $\alpha=\pm$ and $\beta=\pm$.
\end{enumerate}
\end{definition}

\begin{definition}
Let $K$ be a directed complex and $n$ be a non-negative integer.
Let $a$ and $b$ be elements in $K$.

\begin{itemize}
\item
An {\em $n$-path} of length $k$ from $a$ to $b$ is a sequence
$a=a_{0},\dots,a_{k}=b$ of elements in $K$ such that for $1\leq
i\leq k$ either
\begin{center}
$\dim a_{i-1}\leq n$ and $\dim a_{i}>n$ and $a_{i-1}\in
d_{n}^{-}\bar{a}_{i}\setminus (d_{n-1}^{-}\bar{a}_{i}\cup
d_{n-1}^{+}\bar{a}_{i})$
\end{center}
or
\begin{center}
$\dim a_{i-1}>n$ and $\dim a_{i}\leq n$ and $a_{i}\in
d_{n}^{+}\bar{a}_{i-1}\setminus (d_{n-1}^{-}\bar{a}_{i-1}\cup
d_{n-1}^{+}\bar{a}_{i-1})$.
\end{center}

\item
A {\em total path} of length $k$ from $a$ to $b$ is a sequence
$a=a_{0},\dots,a_{k}=b$ of elements in $K$ such that for $1\leq
i\leq k$ either $a_{i-1}\in\partial^{-}a_{i}$ or
$a_{i}\in\partial^{+}a_{i-1}$.

\item
An {\em $n$-loop} is an $n$-path of positive length from some
element of $K$ to itself; A {\em total loop} is a total path of
positive length from some element of $K$ to itself.

\item
A subset of $K$ is {\em loop-free} if it does not contain
$n$-loops for any $n$; A subset of $K$ is {\em total loop-free} if
it does not contain total loops.

\end{itemize}

\end{definition}

We can now generalise the concept of loop-freeness to
$\omega$-complexes, as follows.

\begin{definition}
Let $K$ be an $\omega$-complex and $n$ be a non-negative integer.
Let $a$ and $b$ be atoms.

\begin{itemize}
\item
An {\em $n$-path} of length $k$ from $a$ to $b$ is a sequence
$a=a_{0},\dots,a_{k}=b$ of atoms such that for $1\leq i\leq k$
either
\begin{center}
$\dim a_{i-1}\leq n$ and $\dim a_{i}>n$ and $\Int a_{i-1}\subset
d_{n}^{-}a_{i}\setminus (d_{n-1}^{-}a_{i}\cup d_{n-1}^{+}a_{i})$
\end{center}
or
\begin{center}
$\dim a_{i-1}>n$ and $\dim a_{i}\leq n$ and $\Int a_{i}\subset
d_{n}^{+}a_{i-1}\setminus (d_{n-1}^{-}a_{i-1}\cup
d_{n-1}^{+}a_{i-1})$.
\end{center}

\item
A {\em total path} of length $k$ from $a$ to $b$ is a sequence
$a=a_{0},\dots,a_{k}=b$ of atoms such that for $1\leq i\leq k$
either $\dim a_{i-1}=\dim a_{i}-1$ and
$a_{i-1}\subset\partial^{-}a_{i}$, or $\dim a_{i}=\dim a_{i-1}-1$
and $a_{i}\subset\partial^{+}a_{i-1}$.

\item
An {\em $n$-loop} is an $n$-path of positive length from some atom
to itself; A {\em total loop} is a total path of positive length
from some atom to itself.

\item
A subcomplex of $K$ is {\em loop-free} if it does not contain
$n$-loops for any $n$; A subcomplex of $K$ is {\em total
loop-free} if it does not contain total loops.

\end{itemize}

\end{definition}

For example, the $0$-path and total path in a 1-dimensional
$\omega$-complex is a directed path, regarded as a sequence of
alternate vertices and edges. In Example \ref{first_ex}, the
sequence $u$, $a$, $v$, $b$, $x$, $c$ is a total path; the
sequence $b$, $x$, $c$, $y$, $d$ is a $1$-path.

\begin{lemma}\label{p_minus_1}
Let $K$ be a loop-free $\omega$-complex.  Let $a$ be an atom in
$K$ with $\dim a=p>0$. Then $\partial^{\gamma} a$ is a union of
its $p-1$ dimensional atoms.
\end{lemma}

\begin{proof}
Suppose otherwise that there is a $p$-dimensional atom $a$ of $K$
such that $\partial^{\gamma} a$ is not a union of its $p-1$
dimensional atoms for some $\gamma$. Then $\partial^{\gamma} a $
has a maximal atom $b$ of dimension less than $p-1$.  Let $q=\dim
b$. By Lemma \ref{atomic_dpg}, we can see that $b$ is an maximal
atom in both $d_{q}^{-}\partial^{\gamma}a=d_{q}^{-}a$ and
$d_{q}^{+}\partial^{\gamma}a=d_{q}^{+}a$ which implies that $a$,
$b$, $a$ is a $q$-loop. This contradicts the assumption of the
loop-freeness of $K$.
\end{proof}

Let $K$ be a directed complex. According to the definition of
directed complexes and Theorem \ref{atomic_omega}, set $K$ is also
an $\omega$-complex such that the atoms in the $\omega$-complex
$K$ are exactly those in the directed complex $K$. Note that $\Int
\bar{\sigma}$ is a singleton $\{\sigma\}$ for every atom
$\bar{\sigma}$. By Theorem \ref{atomic_dpg}, Definition
\ref{equivalent} and Lemma \ref{p_minus_1}, it is easy to see that
every loop-free $\omega$-complex is equivalent to an
$\omega$-complex associated with a loop-free directed complex.
Thus all results for loop-free directed complexes can be
generalised to loop-free $\omega$-complexes. In particular, we
have the following definitions and theorems.

\begin{definition}
Let $x$ be a non-empty finite subcomplex of an $\omega$-complex
which is not an atom. Then the non-negative integer $$\max\{\dim
(a\cap b): a \text{ and } b \text{ are distinct maximal atoms in
$x$}\}$$ is called {\em frame dimension} of $x$, denoted by
$\text{fr }\dim x$.
\end{definition}

\begin{definition}
A  molecule $x$ in an $\omega$-complex is {\em split} if the
following conditions hold.
\begin{itemize}
\item
Let $a$ be a $p$-dimensional  atom in $x$. If $b$ is a $p-1$
dimensional atom in $\partial^{-} a$ and if $c$ is a $p-1$
dimensional atom in $\partial^{+}a$, then $b$ and $c$ are
distinct.

\item
If $y$ is a factor in some expression of $x$ as an iterated
composite, then there exists an expression of $y$ as an iterated
composite of atoms using the operations $\#_{n}$ only for $n\leq
\text{fr }\dim y$.
\end{itemize}
\end{definition}

\begin{prop}
If a subcomplex of an $\omega$-complex is total loop-free, then it
is loop-free.
\end{prop}

\begin{theorem}\label{split}
In a loop-free $\omega$-complex, all molecules  are split.
\end{theorem}

\begin{theorem}
If the atoms in an $\omega$-complex are all total loop-free, then
the molecules are all total loop-free, so that all the molecules
are split.
\end{theorem}

\begin{theorem}
Let $K$ and $L$ be $\omega$-complexes. If both $K$ and $L$ are
total loop-free, then so is $K\times L$.
\end{theorem}

According to this theorem, the products of infinite-dimensional
globes are total loop-free. Hence all molecules in products of
infinite-dimensional globes are split.

Now we can state the main theorem in this section.

\begin{theorem}\label{decom_mole}
Let $x$ be a molecule in a loop-free $\omega$-complex and $p=\fdim
x$. Let $q$ be an integer with $q\geq p$. If there is a maximal
atom $a_{1}$ in $x$ with $\dim a_{1}>p$ such that $a_{1}\cap
a'\subset d_{q}^{+}a_{1}\cap d_{q}^{-}a'$ for every other maximal
atom $a'$ in $x$ with $\dim a'>q$, then $x$ can be decomposed into
molecules $$x=x^{-}\#_{q}x^{+},$$ where $x^{-}=d_{q}^{-}x\cup
a_{1}$ and $x^{+}=d_{q}^{+}x\cup \bigcup\{a'': a'' \text{ is a
maximal atom in } x \text{ with } a''\neq a_{1}\}$.
\end{theorem}

The decomposition for $\partial^{-}(u_{2}\times v_{1})$ and
$\partial^{-}(u_{2}\times v_{1})$ in Example \ref{2by1globe} is
actually obtained by using this theorem.

The proof is separated into several lemmas.

\begin{lemma}
Let $x$ be a subcomplex and $y=y_{1}\cup\cdots\cup y_{n}$ be  a
union of subcomplexes. If $x\subset y$ and $x\cap y_{i}\subset
d_{p}^{\gamma}y_{i}$ for all $1\leq i\leq n$, then $x\subset
d_{p}^{\gamma}y$.
\end{lemma}

\begin{proof}
We give the proof only for $n=2$. The general case can be shown by
induction.

By Proposition \ref{dpgunion}, we have $d_{p}^{\gamma}y=(\dpg
y_{1}\cap\dpg y_{2})\cup(\dpg y_{1}\setminus y_{2})\cup (\dpg
y_{2}\setminus y_{1})$.

Suppose that $\xi\in x$. Then $\xi\in y_{1}$ or $\xi\in y_{2}$. If
$\xi\in y_{1}$ and $\xi\in y_{2}$, then $\xi\in \dpg y_{1}\cap\dpg
y_{2}\subset \dpg y$. If $\xi\in y_{1}$ but $\xi\not\in y_{2}$,
then $\xi\in\dpg y_{1}$ but $\xi\not\in y_{2}$; this implies that
$\xi\in\dpg y$. If $\xi\in y_{2}$ but $\xi\not\in y_{1}$, then
$\xi\in\dpg y_{2}$ but $\xi\not\in y_{1}$; this implies that
$\xi\in\dpg y$. This completes the proof that $x\subset\dpg y$.

\end{proof}

\begin{coro}\label{sub_dpg}
Let $x=x_{1}\cup\cdots x_{m}$ and $y=y_{1}\cup\dots\cup y_{n}$ be
a union of subcomplexes. If $x\subset y$ and $x_{i}\cap
y_{j}\subset d_{p}^{\gamma}y_{j}$ for all $1\leq i\leq m$ and
$1\leq j\leq n$, then $x\subset d_{p}^{\gamma}y$.
\end{coro}

\begin{lemma}\label{maximal_dpg}
Let $x$ be a subcomplex of an $\omega$-complex and $a$ be a
maximal atom with $\dim a\leq n$. If $d_{n}^{\alpha}x$ is a
subcomplex, then $a$ is a maximal atom in $d_{n}^{\alpha}x$ for
every sign $\alpha$.
\end{lemma}
\begin{proof}
This follows straightforwardly from Definition \ref{atomic_dpg}.
\end{proof}

\begin{lemma}\label{xi_dpg}
Let $x$ be a subcomplex of an $\omega$-complex and $\xi\in x$.
Then $\xi\in\dpg x$ if and only if $\xi\in\dpg a$ for every
maximal atom $a$ in $x$ with $\xi\in a$.
\end{lemma}
\begin{proof}
Suppose that $\xi\in\dpg x$. According to Proposition
\ref{dpg_sub}, we have $\xi\in a\cap\dpg x\subset\dpg a$ for every
maximal atom $a$ in $x$ with $\xi\in a$.

Conversely, suppose that $\xi\in\dpg a$ for every maximal atom $a$
in $x$ with $\xi\in a$. Then it is evident that $\xi\in a'$ for
some atom $a'$ in $x$ with $\dim a'\leq q$. Moreover, for every
$p+1$-dimensional  atom $b\subset x$ with $\xi\in b$, we have
$\xi\in b\cap\dpg b'\subset\dpg b$ by Proposition \ref{dpg_sub},
where $b'$ is a maximal atom containing $b$ . It follows from
Definition \ref{atomic_dpg} that $\xi\in\dpg x$, as required.

This completes the proof.
\end{proof}

\begin{lemma}
Let $x$, $x^{-}$ and $x^{+}$ be as in the statement of Theorem
\ref{decom_mole}. Then $x=x^{-}\#_{q}x^{+}$.
\end{lemma}
\begin{proof}
Recall that $$x^{-}=d_{q}^{-}x\cup a_{1}$$ and
$$x^{+}=d_{q}^{+}x\cup \bigcup\{a'': a'' \text{ is a maximal atom
in } x \text{ with } a''\neq a_{1}\}.$$ Thus $x=x^{-}\cup x^{+}$.

Since $x$ is a molecule, the set $\{a'': a'' \text{ is a maximal
atom in } x \text{ with } a''\neq a_{1}\}$ is finite. Moreover, it
is evident that $x^{-}$ and $x^{+}$ are subcomplexes.

Now we trivially have $x^{+}\cap d_{q}^{-}x\subset
d_{q}^{+}d_{q}^{-}x$; since $a_{1}\subset x$, we have
$d_{q}^{+}x\cap a_{1}\subset d_{q}^{+}a_{1}$ by Proposition
\ref{dpg_sub}; according to the assumption, we also have
$a_{1}\cap a''\subset d_{q}^{+}a_{1}\cap d_{q}^{-}a''\subset
d_{q}^{+}a_{1}$ for every maximal atom $a''$ in $x$ with $a''\neq
a_{1}$. It follows from Corollary \ref{sub_dpg} that $x^{-}\cap
x^{+}\subset d_{q}^{+}x^{-}$.

On the other hand, suppose that $\xi\in d_{q}^{+}x^{-}$ and
$\xi\not\in a''$ for every atom $a''$ distinct from $a_{1}$ such
that $\dim a''>q$. We claim that $\xi\in d_{q}^{+}x$ so that
$\xi\in x^{+}$ and hence $d_{q}^{+}x^{-}\subset x^{-}\cap x^{+}$.
Indeed, let $a$ be a maximal atom in $x$ with $\xi\in a$. If
$a=a_{1}$, then, by Proposition \ref{dpg_sub}, $\xi\in a_{1}\cap
d_{q}^{+}x^{-}\subset d_{q}^{+}a_{1}$; if $a\neq a_{1}$, then
$\dim a\leq q$ by the assumption; hence $\xi\in a=d_{q}^{+}a$. It
follows from Lemma \ref{xi_dpg} that $\xi\in d_{q}^{+}x$, as
required.

We have now shown that $x^{-}\cap x^{+}=d_{q}^{+}x^{-}$. By a
similar argument, we can also get $x^{-}\cap
x^{+}=d_{q}^{-}x^{+}$. This implies that $x^{-}\#_{q}x^{+}$ is
defined and hence $x=x^{-}\#_{q}x^{+}$, as required.

This completes the proof.

\end{proof}

\begin{lemma}
Let $x^{-}$ and $x^{+}$ be as in the statement of Theorem
\ref{decom_mole}. Then $x^{-}$ and $x^{+}$ are molecules.
\end{lemma}

\begin{proof}
Since $x$ is a molecule in a loop-free $\omega$ complex, it is
split by Theorem \ref{split}. Hence $x^{-}$ and $x^{+}$ are
molecules, as required.
\end{proof}

We have now completed the proof of Theorem \ref{decom_mole}.


\chapter{Molecules in the Product of Three Infinite-Dimensional
Globes}\label{ch2}
In this chapter, we study molecules in the product of three
infinite dimensional globes. We are going to give two equivalent
descriptions for the molecules in the product of three infinite
dimensional globes.

Throughout this chapter, infinite dimensional globes are denoted
by $u$, $v$ or $w$. An atom $u_{i}^{\alpha}$ is  denoted by
$\ui{}$. All subcomplexes refer to finite and non-empty
subcomplexes in the $\omega$-complex $u\times v\times w$.

\section{The Definition of Pairwise Molecular Subcomplexes}
In this section, we first define `projection maps' and prove some
of their basic properties. Then we state one of the main results
in this chapter which says that a subcomplex in products of three
infinite dimensional globes is a molecule if and only if it is
`projected' to molecules in (twisted) products of two infinite
dimensional globes together with a natural requirement. This leads
to the definition of pairwise molecular subcomplexes of $u\times
v\times w$.

Let $w^{J}$ be the atomic complex with atoms
$w^{J}[k,\varepsilon]$ ($k=0,1,\cdots$ and $\varepsilon=\pm$) such
that $\dim\,w^{J}[k,\varepsilon]=k$ and
$d_{k-1}^{\gamma}w^{J}[k,\varepsilon]=w^{J}[k-1, (-)^{J}\gamma]$
for $k>0$. It is clear that $w^{J}$ satisfies conditions in
Theorem \ref{atomic_omega}. Thus it is an $\omega$-complex. It is
also easy to see that the $\omega$-complex $w^{J}$ is equivalent
to infinite dimensional globe $w$ under an obvious equivalence of
$\omega$-complexes sending $w^{J}[k,(-)^{J}\varepsilon]$ to
$w[k,\varepsilon]$. Moreover, it is evident that this induces a
equivalence of $\omega$-complexes $u\times w$ and $u\times w^{J}$
sending $\ui{}\times\wk{}$ to $\ui{}\times w^{J}[k,\varepsilon]$.
By this equivalence, all the results for products of two infinite
dimensional globes in Section \ref{prm3} can be generalised to
$u\times w^{J}$. In particular, we have $ d_{p}^{\gamma}
(\ui{}\times w^{J}[k,\varepsilon])=\ui{}\times
w^{J}[k,\varepsilon]$ if $i+k\leq p$, while, if $i+k>p$, the
maximal atoms in $ d_{p}^{\gamma} (\ui{}\times
w^{J}[k,\varepsilon])$ consists of $\ul{}\times w^{J}[n,\omega]$
such that $l\leq i$, $n\leq k$ and $l+n=p$; the signs $\sigma$ and
$\omega$ are determined as follows:
\begin{enumerate}
\item
if $l=i$, then $\sigma=\alpha$; if $l<i$, then $\sigma=\gamma$;
\item
if $n=k$, then $\omega=\varepsilon$; if $n<k$, then
$\omega=(-)^{l+J}\gamma$.
\end{enumerate}

For an atom $\uvws{}$ in $u\times v\times w$, let $$
F_{J}^{v}(\Int\lambda)=
  \begin{cases}
    \Int(\ui{}\times w^{J}[k,\varepsilon]), & \text{when $j\geq J$}; \\
    \emptyset, & \text{when $j<J$}.
  \end{cases}
$$ This gives a map sending interiors of atoms in $u\times v\times
w$ to interiors of atoms in $u\times w^{J}$ or the empty set.

Since interiors of atoms are disjoint, it is clear that the map
$F_{J}^{v}$ can be extended uniquely to a map sending unions of
interiors of atoms in $u\times v\times w$ to unions of interiors
of atoms in $u\times w^{J}$ by requiring it preserves unions.

We can similarly define a map $F_{I}^{u}$ sending unions of
interiors of atoms in $u\times v\times w$ to unions of interiors
of atoms in $v^{I}\times w^{I}$ and a map $F_{K}^{w}$ sending
unions of interiors of atoms in $u\times v\times w$ to unions of
interiors of atoms in $u\times v$.

It is easy to see that every atom can be written as a union of
interiors of atoms. It follows that $F_{I}^{u}$, $F_{J}^{v}$ and
$F_{k}^{w}$ are defined on subcomplexes of $u\times v\times w$ and
preserve unions.

We next prove that $F_{I}^{u}$, $F_{J}^{v}$ and $F_{k}^{w}$ send
atoms to atoms or the empty set so that they send subcomplexes to
subcomplexes. We need two preliminary results.

\begin{lemma}\label{single}
$$ d_{p}^{\gamma} (\uvws{})
=
\bigcup\{d_{l}^{\gamma}\ui{}\times
d_{m}^{(-)^{l}\gamma}\vj{}\times d_{n}^{(-)^{l+m}\gamma}\wk{}
\,:\, l+m+n=p\} $$
\end{lemma}

\begin{proof}
According to Proposition \ref{dpg_product}, $$
\begin{array}{rl}
&d_{p}^{\gamma} (\ui{}\times \vj{}\times \wk{})\\ =&
d_{p}^{\gamma} [(\ui{}\times \vj{})\times \wk{}]\\ =&
\bigcup\{d_{s}^{\gamma}(\ui{}\times \vj{})\times
d_{t}^{(-)^{s}\gamma}\wk{} \,:\,s+t=p\}.
\end{array}
$$ Then the result follows easily by applying Proposition
\ref{dpg_product} again.
\end{proof}

\begin{prop}\label{dtatom}
Let $\lambda=\ui{}\times\vj{}\times\wk{}$ be an atom.
\begin{itemize}
\item
If $i+j+k\leq p$, then $d_{p}^{\gamma}\lambda=\lambda$.
\item
If $i+j+k>p$, then the set of maximal atoms in
$d_{p}^{\gamma}\lambda$ consists of all the atoms
$u_{l}^{\sigma}\times v_{m}^{\tau}\times w_{n}^{\omega}$ such that
$l+m+n=p$ and $l\leq i$, $m\leq j$ and $n\leq k$, where the signs
$\sigma$, $\tau$ and $\omega$ are determined as follows:

\begin{enumerate}
\item
If $l=i$, then $\sigma=\alpha$; if $l<i$, then $\sigma=\gamma$.
\item
If $m=j$, then $\tau=\beta$; if $m<j$, then $\tau=(-)^{l}\gamma$.
\item
If $n=k$, then $\omega=\varepsilon$; if $l<i$, then
$\omega=(-)^{l+m}\gamma$.

\end{enumerate}
\end{itemize}
\end{prop}

\begin{proof}
It is evident that $d_{p}^{\gamma}\lambda=\lambda$ when $i+j+k\leq
p$. We may assume in the following proof that $i+j+k>p$.

Let $\Lambda_{1}$ denote the union of the atoms described in this
lemma. We must show that $d_{p}^{\gamma}\lambda=\Lambda_{1}$.

By the formation of $\Lambda_{1}$, it is easy to see that every
maximal atom $\mu=\lmns{}$ in $\Lambda_{1}$ can be expressed as
$\mu=d_{l}^{\gamma}\ui{}\times d_{m}^{(-)^{l}\gamma}\vj{} \times
d_{n}^{(-)^{l+m}\gamma}\wk{}$. By Lemma \ref{single}, we can see
that $\mu\subset d_{p}^{\gamma}\lambda$, and hence
$\Lambda_{1}\subset d_{p}^{\gamma}\lambda$.

To prove the reverse inclusion, by Lemma \ref{single}, it suffices
to prove that $d_{l}^{\gamma}\ui{}\times
d_{m}^{(-)^{l}\gamma}\vj{} \times
d_{n}^{(-)^{l+m}\gamma}\wk{}\subset \Lambda_{1}$ for every triple
$(l,m,n)$ with $l+m+n=p$. By the formation of $\Lambda_{1}$, this
inclusion is obvious when $l\leq i$, $m\leq j$ and $n\leq k$. So
it suffices to prove that $d_{l}^{\gamma}\ui{}\times
d_{m}^{(-)^{l}\gamma}\vj{} \times d_{n}^{(-)^{l+m}\gamma}\wk{}$ is
not a  maximal atom in $d_{p}^{\gamma}\lambda$ when  $l>i$, $m>j$
or $n>k$.

Suppose that $l>i$. Then $m<j$ or $n<k$. If $m<j$, then
$d_{l}^{\gamma}\ui{}\times d_{m}^{(-)^{l}\gamma}\vj{} \times
d_{n}^{(-)^{l+m}\gamma}\wk{}\subsetneqq
d_{l-1}^{\gamma}\ui{}\times d_{m+1}^{(-)^{l-1}\gamma}\vj{} \times
d_{n}^{(-)^{l+m}\gamma}\wk{}\subset d_{p}^{\gamma}\lambda$. If
$m\geq j$, then $n<k$, so $d_{l}^{\gamma}\ui{}\times
d_{m}^{(-)^{l}\gamma}\vj{} \times
d_{n}^{(-)^{l+m}\gamma}\wk{}\subsetneqq
d_{l-1}^{\gamma}\ui{}\times d_{m}^{(-)^{l-1}\gamma}\vj{} \times
d_{n+1}^{(-)^{l+m-1}\gamma}\wk{}\subset d_{p}^{\gamma}\lambda$.
This shows that, in both cases, $d_{l}^{\gamma}\ui{}\times
d_{m}^{(-)^{l}\gamma}\vj{} \times d_{n}^{(-)^{l+m}\gamma}\wk{}$ is
not a maximal atom in $d_{p}^{\gamma}\lambda$. Similarly, the
above statement is true if $m>j$ or $n>k$.

This completes the proof of the lemma.
\end{proof}

\begin{prop}\label{pro_atom}
Let $\uvws{}$ be an atom  in $u\times v\times w$. Then
\begin{enumerate}
\item
$
F_{I}^{u}(\uvws{})=
  \begin{cases}
    \vIj{}\times \wIk{}, & \text{when $i\geq I$}; \\
    \emptyset, & \text{when $i<I$};
  \end{cases}
$
\item
$
F_{J}^{v}(\uvws{})=
  \begin{cases}
    \ui{}\times \wJk{}, & \text{when $j\geq J$}; \\
    \emptyset, & \text{when $j<J$};
  \end{cases}
$
\item
$
F_{K}^{w}(\uvws{})=
  \begin{cases}
    \ui{}\times\vj{}, & \text{when $k\geq K$}; \\
    \emptyset, & \text{when $k<K$}.
  \end{cases}
$
\end{enumerate}
In particular, $F_{I}^{u}$, $F_{J}^{v}$ and $F_{K}^{w}$ send atoms
to atoms or the empty set  so that they send subcomplexes to
subcomplexes.
\end{prop}

\begin{proof}
The argument for the three cases are similar. We only prove the
second one. The proof is given by induction on dimension of atoms.

For an atom $\lambda=\uvws{}$ in $u\times v\times w$, if
$\dim\lambda=0$, then $i=j=k=0$; hence $$
\begin{array}{rl}
&F_{J}^{v}(\lambda)\\ =&F_{J}^{v}(\Int\lambda)\\ =&
\begin{cases}
    \Int(\ui{}\times \wJk{}), & \text{when $J=0$}; \\
    \emptyset, & \text{when $J>0$}
  \end{cases} \\
=&
\begin{cases}
    \ui{}\times \wJk{}, & \text{when $J=0$}; \\
    \emptyset, & \text{when $J>0$},
  \end{cases}
\end{array}
$$ as required.

Suppose that the proposition holds for every atom of dimension
less then $p$. Suppose also that $\lambda=\uvws{}$ is a
$p$-dimensional atom. If $j<J$, then it is easy to see that
$F_{J}^{v}(\lambda)=\emptyset$, as required. If $j>J$, then we
have $$
\begin{array}{rl}
&F_{J}^{v}(\lambda)\\
=&F_{J}^{v}(\Int\lambda\cup\partial^{-}\lambda\cup\partial^{+}\lambda)\\
\supset & F_{J}^{v}(\partial^{+}\lambda)\\ \supset &
F_{J}^{v}(\ui{}\times v[j-1,(-)^{i}]\times \wk{})\\ =&\ui{}\times
\wJk{}
\end{array}
$$ since $\ui{}\times v[j-1,(-)^{i}]\times \wk{}$ is an atom of
dimension $p-1$; the reverse inclusion holds automatically; so
$F_{J}^{v}(\lambda)=\ui{}\times \wJk{}$, as required. Now suppose
that $j=J$. Then, by Lemma \ref{dtatom},
$\partial^{\gamma}\lambda$ is the union of atoms
$\ui{'}\times\vj{'}\times\wk{'}$ with $i'+j'+k'=p-1$ such that
\begin{enumerate}
\item
if $i'=i$, then $\alpha'=\alpha$; if $i'=i-1$, then
$\alpha'=\gamma$;
\item
if $j'=j$, then $\beta'=\beta$; if $j'=j-1$, then
$\beta'=(-)^{i}\gamma$;
\item
if $k'=k$, then $\varepsilon'=\varepsilon$; if $k'=k-1$, then
$\varepsilon'=(-)^{i+J}\gamma$.
\end{enumerate}
It follows easily from the induction hypothesis  that
$F_{J}^{v}(\partial^{\gamma}\lambda)=
\partial^{\gamma}(\ui{}\times \wJk{})$ for every sign $\gamma$.
Therefore $$
\begin{array}{rl}
&F_{J}^{v}(\lambda)\\ =&F_{J}^{v}(\Int\lambda)\cup
F_{J}^{v}(\partial^{-}\lambda)\cup
F_{J}^{v}(\partial^{-}\lambda)\\ =& \Int(\ui{}\times \wJk{})\cup
\partial^{-}(\ui{}\times \wJk{}) \cup \partial^{+}(\ui{}\times
\wJk{})\\ =&\ui{}\times \wJk{},
\end{array}
$$ as required.

This completes the proof of the proposition.
\end{proof}

Now we can state one of the main results in this chapter which
says that a subcomplex in $u\times v\times w$ is a molecule if and
only if it is {\em pairwise molecular,} i.e., it is 'projected` to
molecules in (twisted) products of two infinite dimensional globes
together with a natural condition (condition 1).

\begin{definition}\label{pairwise_def}
Let $\Lambda$ be a subcomplex in $u\times v\times w$. Then
$\Lambda$ is {\em pairwise molecular} if the following conditions
hold:
\begin{enumerate}
\item
there are no distinct maximal atoms $\uvws{}$ and \\ $\uvws{'}$ in
$\Lambda$ such that $i\leq i'$, $j\leq j'$ and $k\leq k'$;
\item
$F_{I}^{u}(\Lambda)$ is a molecule in $v^{I}\times w^{I}$ or the
empty set for every integer $I$;
\item
$F_{J}^{v}(\Lambda)$ is a molecule in $u\times w^{J}$ or the empty
set for every integer $J$;
\item
$F_{K}^{w}(\Lambda)$ is a molecule in $u\times v$ or the empty set
for every integer $K$.
\end{enumerate}
\end{definition}

\begin{example}\label{pairwise_example}
It is easy to check that the following subcomplex of $\uvwt$ is
pairwise molecular.

$$\begin{array}{rl} &u_{8}^{+}\times v_{2}^{+}\times w_{1}^{-}\\
\cup&u_{5}^{-}\times v_{2}^{+}\times w_{5}^{-}\\
\cup&u_{1}^{-}\times v_{2}^{+}\times w_{8}^{+}\\
\cup&u_{9}^{+}\times v_{1}^{-}\times w_{2}^{+}\\
\cup&u_{4}^{-}\times v_{1}^{-}\times w_{6}^{+}\\
\cup&u_{0}^{+}\times v_{1}^{+}\times w_{9}^{+}\\
\cup&u_{8}^{-}\times v_{0}^{-}\times w_{5}^{-}\\
\cup&u_{5}^{-}\times v_{0}^{+}\times w_{6}^{+}\\
\cup&u_{4}^{-}\times v_{0}^{-}\times w_{7}^{+}\\
\cup&u_{2}^{-}\times v_{0}^{-}\times w_{9}^{+}\\
\end{array}
$$

\end{example}

\begin{theorem}\label{3glb_main}
A subcomplex of $\uvwt$ is a molecule if and only if it is
pairwise molecular.
\end{theorem}

We end this section with a property of `projection' maps which is
used later in the thesis.

\begin{prop}\label{eq}
Let $\Lambda$ and $\Lambda'$ be subcomplexes of $\uvwt$ satisfying
condition 1 for pairwise molecular subcomplexes. If
$F_{I}^{u}(\Lambda)=F_{I}^{u}(\Lambda')$,
$F_{J}^{v}(\Lambda)=F_{J}^{v}(\Lambda')$ and
$F_{K}^{w}(\Lambda)=F_{K}^{w}(\Lambda')$ for all $I$, $J$ and $K$,
then $\Lambda=\Lambda'$.
\end{prop}
\begin{proof}
It suffices to prove that $\Lambda$ and $\Lambda'$ consists of the
same maximal atoms.

Let $\uvws{}$ be a maximal atom in $\Lambda$. It is easy to see
that $v^{i}[j,\beta]\times w^{i}[k,\varepsilon]$ is a maximal atom
in $F_{i}^{u}(\Lambda)=F_{i}^{u}(\Lambda')$. Thus $\Lambda'$ has a
maximal atom $\ui{'}\times\vj{}\times\wk{}$ with $i'\geq i$. Since
$v^{i}[j,\beta]\times w^{i}[k,\varepsilon]\not\subset
F_{i+1}^{u}(\Lambda)=F_{i+1}^{u}(\Lambda')$, we have $i'=i$.  One
can similarly get a maximal atom $\ui{}\times
v[j,\beta']\times\wk{}$ in $\Lambda'$. It follows from condition 1
for pairwise molecular subcomplexes that $\alpha'=\alpha$ and
$\beta'=\beta$. This shows that $\uvws{}$ is a maximal atom in
$\Lambda'$.

Symmetrically, we can see that every maximal atom in $\Lambda'$ is
a maximal atom in $\Lambda$.

This completes the proof that $\Lambda=\Lambda'$.
\end{proof}

\begin{remark}
The above proposition does not holds without Condition 1 for
pairwise molecular subcomplexes. This can be seen from the
following subcomplexes of $\uvwt$: $$
\begin{array}{rl}
\Lambda=&u_{1}^{+}\times v_{1}^{+}\times w_{1}^{+}\\
\cup&u_{1}^{+}\times v_{1}^{-}\times w_{1}^{-}\\
\cup&u_{1}^{-}\times v_{1}^{+}\times w_{1}^{-}\\
\cup&u_{1}^{-}\times v_{1}^{-}\times w_{1}^{+}
\end{array}
$$ and $$
\begin{array}{rl}
\Lambda'=&u_{1}^{+}\times v_{1}^{-}\times w_{1}^{+}\\
\cup&u_{1}^{+}\times v_{1}^{+}\times w_{1}^{-}\\
\cup&u_{1}^{-}\times v_{1}^{-}\times w_{1}^{-}\\
\cup&u_{1}^{-}\times v_{1}^{+}\times w_{1}^{+}.
\end{array}
$$

\end{remark}

\section{Molecules Are Pairwise Molecular}
In this section, we prove that molecules in $\uvwt$ are pairwise
molecular.

\begin{prop}\label{mlc1}
Let $\Lambda$ be a molecule in $\uvwt$. Then there are no distinct
maximal atoms $\uvws{_{1}}$ and $\uvws{_{2}}$ in $\Lambda$ such
that $i_{1}\leq i_{2}$, $j_{1}\leq j_{2}$ and $k_{1}\leq k_{2}$.
\end{prop}

\begin{proof}
Suppose otherwise that there are distinct maximal atoms
$\lambda_{1}=\uvws{_{1}}$ and  $\lambda_{2}=\uvws{_{2}}$ in
$\Lambda$ such that $i_{1}\leq i_{2}$, $j_{1}\leq j_{2}$ and
$k_{1}\leq k_{2}$. Then we have $\ui{_{2}}=u[i_{1},-\alpha_{1}]$,
$\vj{_{2}}=v[j_{1},-\beta_{1}]$ or
$\wk{_{2}}=w[k_{1},-\varepsilon_{1}]$. The arguments for various
cases are similar, we only give the proof for the case
$\ui{_{2}}=u[i_{1},-\alpha_{1}]$, $j_{1}<j_{2}$ and $k_{1}<k_{2}$.
In this case, it is easy to see that there is a natural
homomorphism $f:\M(u\times v\times w)\to \M(u_{i_{1}}\times
v\times w)$ of $\omega$-categories such that $f(\uvws{})=\uvws{}$
for $i<i_{1}$ and $f(\uvws{})=u_{i_{1}}\times \vj{}\times\wk{}$
for $i\geq i_{1}$. We are going to use this homomorphism to get a
contradiction.

Since $\lambda_{1}$ and $\lambda_{2}$ are maximal in the molecule
$\Lambda$, it is easy to see that there is a composite of
molecules $\Lambda_{1}\#_{n}\Lambda_{2}$ or
$\Lambda_{2}\#_{n}\Lambda_{1}$ such that $\lambda_{1}$ is a
maximal atom in $\Lambda_{1}$ and $\lambda_{2}$ is a maximal atom
in $\Lambda_{2}$ and $\lambda_{1}\not\subset\Lambda_{2}$ and
$\lambda_{2}\not\subset\Lambda_{1}$. We may assume that
$\Lambda_{1}\#_{n}\Lambda_{2}$ is defined. In this case, we have
$d_{n}^{+}\Lambda_{1}=d_{n}^{-}\Lambda_{2}=\Lambda_{1}\cap\Lambda_{2}$.
It follows from Lemma \ref{maximal_dpg} that
$n<\dim\lambda_{1}=i_{1}+j_{1}+k_{1}$. On the other hand, since
$f:\M(u\times v\times w)\to \M(u_{i_{1}}\times v\times w)$ is a
homomorphism, the composite $f(\Lambda_{1})\#_{n} f(\Lambda_{2})$
is defined; since $f(\lambda_{1})\subset f(\Lambda_{1})$ and
$f(\lambda_{1})=u_{i_{1}}\times \vj{_{1}}\times\wk{_{1}}\subset
u_{i_{1}}\times \vj{_{2}}\times\wk{_{2}}=f(\lambda_{2})\subset
f(\Lambda_{2})$, we have $f(\lambda_{1})\subset f(\Lambda_{1})\cap
f(\Lambda_{2})$; this implies that $n\geq\dim
d_{n}^{+}f(\Lambda_{1})=\dim(f(\Lambda_{1})\cap
f(\Lambda_{2}))\geq\dim f(\lambda_{1})=i_{1}+i_{2}+i_{3}$, a
contradiction.

This completes the proof.

\end{proof}

We have now proved that a molecule satisfies condition 1 for
pairwise molecular subcomplexes. We next prove that $F^{u}_{I}$,
$F^{v}_{J}$ and $F^{w}_{K}$ send molecules to molecules or
$\emptyset$. The arguments for the three maps are similar. We only
give the proof that $F^{v}_{J}$ sends molecules in $\uvwt$ to
molecules in $u\times w^{J}$ or the empty set.

Let $v_{J}$ be a $J$-dimensional globe. For $\Lambda$ a subcomplex
in $u\times v_{J}\times w$, let
$g_{J}^{v}(\Lambda)=\pr[\Lambda\cap (u\times \{\eta\}\times w)]$,
where $\eta\in \Int(v_{J})$ and $\pr$ is projection onto the first
and third factors. Then $g_{J}^{v}(\Lambda)\subset u\times w^{J}$.
We are going to show that $g_{J}^{v}(\Lambda)$ is a molecule in
$M(u\times w^{J})$ or the empty set for every molecule $\Lambda$.

We first investigate the image of $\dpg\lambda$ for an atom
$\lambda$ in $\uvwt$ under the map $g_{J}^{v}$.

\begin{lemma} \label{imatom}
Let $\lambda=u[i,\alpha]\times   v_{J}[j,\beta]\times
w[k,\varepsilon]$ be an atom in the $\omega$-complex $u\times
v_{J}\times w$ and $\Lambda,\Lambda'\in\M(u\times   v_{J}\times
w)$. Then
\begin{enumerate}

\item
$g_{J}^{v}(\lambda)\in\A(u\times w^{J})\cup\{\emptyset\}$;

\item
If  $\Lambda\#_{n}\Lambda'$ is defined, then
$g_{J}^{v}(\Lambda\#_{n}\Lambda')=g_{J}^{v}(\Lambda)\cup
g_{J}^{v}(\Lambda')$;

\item
$g_{J}^{v}(\Lambda)\neq\emptyset$ if and only if there is a
maximal atom $\uvws{}$ in $\Lambda$ such that $j=J$;

\item
$g_{J}^{v}(\dpg \lambda)=
\begin{cases}
d_{p-J}^{\gamma}g_{J}^{v}(\lambda) & \text{when $p\geq J$ and
$j=J$,}\\ \emptyset                & \text{when $p<J$ or $j<J$;}
\end{cases}
$

\end{enumerate}
\end{lemma}
\begin{proof}
The proofs of the first three conditions are  trivial verification
from the definition of $g_{J}^{v}$. we now verify condition 4.

If $p<J$ or $j<J$, then it is evident that
$g_{J}^{v}(d_{p}^{\gamma}\lambda)=\emptyset$ by the definition of
$g_{J}^{v}$.

Now, suppose that $p\geq J$ and $j=J$.  Then
$g_{J}^{v}(\lambda)=u[i,\alpha]\times w^{J}[k,\varepsilon]$. The
set of all maximal atoms in $\dpg \lambda$ consists of all
$\lmns{}$ with $l\leq i$, $m\leq J$ and $n\leq k$ by proposition
\ref{dtatom}, where the signs $\sigma$, $\tau$ and $\omega$ are
determined as follows:
\begin{enumerate}
\item
if $l=i$, then $\sigma=\alpha$; if $l<i$, then $\sigma=\gamma$;
\item
if $m=j$, then $\tau=\beta$; if $m<j$, then $\tau=(-)^{l}\gamma$;
\item
if $n=k$, then $\omega=\varepsilon$; if $n<k$, then
$\omega=(-)^{l+m}\gamma$.
\end{enumerate}
From this description and the formation of
$d_{p-J}^{\gamma}(u[i,\alpha]\times w^{J}[k,\varepsilon])$ in
$u\times w^{J}$, it is easy to see that $g_{J}^{v}(\dpg
\lambda)=d_{p-J}^{\gamma}g_{J}^{v}(\lambda)$, as required.

\end{proof}

Now we can prove that $g_{J}^{v}$ sends molecules to molecules or
the empty set.

\begin{theorem}
Let $g_{J}^{v}:\M(u\times   v_{J}\times w)\to \PS(u\times w^{J})$
be the map as above. Then
\begin{enumerate}
\item
$g_{J}^{v}(\M(u\times   v_{J}\times w))\subset\M(u\times
w^{J})\cup\{\emptyset\}$;
\item
For every molecule $\Lambda$ in $u\times   v_{J}\times w$, we have
$$ g_{J}^{v}(\dpg \Lambda)=
\begin{cases}
d_{p-J}^{\gamma}g_{J}^{v}(\Lambda) & \text{when $p\geq J$ and
                                 $g_{J}^{v}(\Lambda)\neq\emptyset$},\\
\emptyset                & \text{when $p<J$ or
                            $g_{J}^{v}(\Lambda)=\emptyset$}.
\end{cases}
$$

\item
If $\Lambda\#_{n} \Lambda'$ is defined, then
$$g_{J}^{v}(\Lambda\#_{n} \Lambda')=
\begin{cases}
g_{J}^{v}(\Lambda)\#_{n-J} g_{J}^{v}(\Lambda') & \text{when
$g_{J}^{v}(\Lambda)\neq\emptyset$ and
$g_{J}^{v}(\Lambda')\neq\emptyset$,}\\ g_{J}^{v}(\Lambda')      &
\text{when $g_{J}^{v}(\Lambda)=\emptyset$,} \\ g_{J}^{v}(\Lambda)
& \text{when $g_{J}^{v}(\Lambda')=\emptyset$.}
\end{cases}
$$
\end{enumerate}
\end{theorem}
\begin{proof}
 We are going to prove the first two conditions by induction and then
prove the third condition.

By Lemma \ref{imatom}, it is evident that the first two conditions
hold when $\Lambda$ is an atom.

Now suppose that $q>1$ and the first two conditions hold for every
molecule which can be written as a composite of less than $q$
atoms. Suppose also that $\Lambda$ is a molecule which can be
written as a composite of $q$ atoms. Since $q>1$, we have a proper
decomposition $\Lambda=\Lambda'\#_{n} \Lambda''$ such that
$\Lambda'$ and $\Lambda''$ are molecules. According to the
induction hypothesis, we know that the first two conditions hold
for $\Lambda'$ and $\Lambda''$. We must show that the first two
conditions in the proposition hold for $\Lambda$. There are two
cases, as follows.

1. Suppose that $g_{J}^{v}(\Lambda')=\emptyset$ or
$g_{J}^{v}(\Lambda'')=\emptyset$. We may assume that
$g_{J}^{v}(\Lambda')=\emptyset$. We have
$g_{J}^{v}(\Lambda)=g_{J}^{v}(\Lambda'')$. Thus
$g_{J}^{v}(\Lambda)\in\M(u\times w^{J})\cup\{\emptyset\}$ as
required by the first condition. Moreover, if $p\neq n$, then $$
\begin{array}{rl}
 &g_{J}^{v}(d_{p}^{\gamma} \Lambda)\\
=&g_{J}^{v}(d_{p}^{\gamma}\Lambda'\#_{n} d_{p}^{\gamma}
\Lambda'')\\ =&g_{J}^{v}(d_{p}^{\gamma} \Lambda')\cup
g_{J}^{v}(d_{p}^{\gamma} \Lambda'')\\ =&
\begin{cases}
d_{p-J}^{\gamma}g_{J}^{v}(\Lambda'') & \text{when $p\geq J$ and
$g_{J}^{v}(\Lambda'')\neq\emptyset$,}\\ \emptyset & \text{when
$g_{J}^{v}(\Lambda'')=\emptyset$ or $p<J$},
\end{cases}\\
=&
\begin{cases}
d_{p-J}^{\gamma}g_{J}^{v}(\Lambda) & \text{when $p\geq J$ and
$g_{J}^{v}(\Lambda)\neq\emptyset$,}\\ \emptyset & \text{when
$g_{J}^{v}(\Lambda)=\emptyset$ or $p<J$},
\end{cases}
\end{array}
$$ as required by the second condition. Suppose that $p=n\geq J$.
Then $g_{J}^{v}(d_{p}^{+}\Lambda')=\emptyset$. So
$g_{J}^{v}(d_{p}^{-}\Lambda'')=\emptyset$. Hence, by the
hypothesis, one gets $g_{J}^{v}(\Lambda'')=\emptyset$. Therefore
$g_{J}^{v}(\Lambda)=\emptyset$ and
$g_{J}^{v}(d_{p}^{\gamma}\Lambda)=\emptyset$, as required by the
second condition.

2. Suppose that $g_{J}^{v}(\Lambda')\neq\emptyset$ and
$g_{J}^{v}(\Lambda'')\neq\emptyset$. Then there is a maximal atom
$\lambda'=\uvws{'}$ in $\Lambda'$ and a maximal atom
$\lambda''=\uvws{''}$ in  $\Lambda''$ such that $j'=j''=J$. We
claim that $n\geq J$. There are two cases, as follows:

a. Suppose that both $\lambda'$ and $\lambda''$ are maximal in
$\Lambda$. By Proposition \ref{mlc1}, we have $i'\neq i''$ and
$k'\neq k''$. So
$\lambda'\cap\lambda''\subset\Lambda'\cap\Lambda''
=d_{n}^{+}\Lambda'=d_{n}^{-}\Lambda''$; Since
$\lambda'\cap\lambda''\neq\emptyset$ and $j=j'=J$ and $\dim
(d_{n}^{+}\Lambda')\leq n$, we can see that $J\leq n$, as
required.

b. Suppose that $\lambda'$ is not maximal in $\Lambda$. Then
$\Lambda$ has a maximal atom $\lambda_{1}'=\uvws{_{1}'}$ distinct
from $\lambda'$ with $\lambda'\subset\lambda_{1}'$. Hence
$j_{1}'=J$. It is easy to see that $\lambda_{1}'$ is  maximal in
$\Lambda''$. So we have $\lambda'\subset\Lambda'\cap\Lambda''
=d_{n}^{+}\Lambda'=d_{n}^{-}\Lambda''$. Since $\dim
(d_{n}^{+}\Lambda')\leq n$, we have $J\leq n$, as required.

Now since $g_{J}^{v}(\Lambda)=g_{J}^{v}(\Lambda')\cup
g_{J}^{v}(\Lambda'')$, and $$
\begin{array}{rl}
 &d_{n-J}^{+}g_{J}^{v}(\Lambda')\\
=&g_{J}^{v}(d_{n}^{+} \Lambda')\\ =&g_{J}^{v}(d_{n}^{-}
\Lambda'')\\ =&d_{n-J}^{-}g_{J}^{v}(\Lambda''),
\end{array}
$$ we can see that $g_{J}^{v}(\Lambda')\#_{n-J}
g_{J}^{v}(\Lambda'')$ is defined,  and
$g_{J}^{v}(\Lambda)=g_{J}^{v}(\Lambda')\#_{n-J}
g_{J}^{v}(\Lambda'')$. So $g_{J}^{v}(\Lambda)$ is a molecule, as
required by the first condition. We now verify that $\Lambda$
satisfies the second condition. If $p<J$, then
$g_{J}^{v}(\dpg\Lambda)=\emptyset$, as required. If $p=n\geq J$,
then $$
\begin{array}{rl}
 &g_{J}^{v}(d_{p}^{-}\Lambda)\\
=&g_{J}^{v}(d_{p}^{-}\Lambda')\\
=&d_{p-J}^{-}g_{J}^{v}(\Lambda')\\
=&d_{p-J}^{-}g_{J}^{v}(\Lambda);
\end{array}
$$ and similarly we have $g_{J}^{v}(d_{p}^{+}\Lambda)
=d_{p-J}^{+}g_{J}^{v}(\Lambda)$. If $J\leq p<n$, then $$
\begin{array}{rl}
 &g_{J}^{v}(\dpg \Lambda)\\
=&g_{J}^{v}(\dpg \Lambda')\\
=&d_{p-J}^{\gamma}g_{J}^{v}(\Lambda')\\
=&d_{p-J}^{\gamma}g_{J}^{v}(\Lambda)
\end{array}
$$ If $p\geq J$ and $p>n$, then $$
\begin{array}{rl}
 &g_{J}^{v}(\dpg \Lambda)\\
=&g_{J}^{v}(\dpg \Lambda'\#_{n}\dpg \Lambda'')\\ =&g_{J}^{v}(\dpg
\Lambda')\cup g_{J}^{v}(\dpg \Lambda'')\\
=&d_{p-J}^{\gamma}g_{J}^{v}(\Lambda')\cup
d_{p-J}^{\gamma}g_{J}^{v}(\Lambda'').
\end{array}
$$ and $$
\begin{array}{rl}
&d_{n-J}^{+}d_{p-J}^{\gamma}g_{J}^{v}(\Lambda')\\
=&d_{n-J}^{+}g_{J}^{v}(\Lambda')\\
=&d_{n-J}^{-}g_{J}^{v}(\Lambda'')\\
=&d_{n-J}^{-}d_{p-J}^{\gamma}g_{J}^{v}(\Lambda''),
\end{array}
$$ thus $d_{p-J}^{\gamma}g_{J}^{v}(\Lambda')\#_{n-J}
d_{p-J}^{\gamma}g_{J}^{v}(\Lambda'')$ is defined and $$
\begin{array}{rl}
 & g_{J}^{v}(\dpg \Lambda)\\
=&d_{p-J}^{\gamma}g_{J}^{v}(\Lambda')\#_{n-J}
d_{p-J}^{\gamma}g_{J}^{v}(\Lambda'')\\
=&d_{p-J}^{\gamma}[g_{J}^{v}(\Lambda')\#_{n-J}g_{J}^{v}(\Lambda'')]\\
=&d_{p-J}^{\gamma}g_{J}^{v}(\Lambda).
\end{array}
$$ Therefore $\Lambda$ satisfies the second condition.

Finally,  condition 3 can be easily verified by using condition 2
and the fact that $g_{J}^{v}$ preserves unions.

This completes the proof.

\end{proof}

Recall that there is a natural homomorphism $f_{J}^{v}:\M(u\times
v\times w)\to\M(u\times v_{J}\times w)$ of $\omega$-categories
sending every atom $\uvws{}$ to $\ui{}\times
  v_{J}[j',\beta']\times w[k,\varepsilon]$, where
$  v_{J}[j',\beta']=  v_{J}[j,\beta]$ whenever $j<J$, and $
v_{J}[j',\beta']=  v_{J}$ whenever $j\geq J$. According to the
definitions of $g_{J}^{v}$ and $f_{J}^{v}$, it is easy to see that
$F_{J}^{v}=g_{J}^{v}\circ f_{J}^{v}$. Thus $F_{J}^{v}$ sends
molecules in $\uvwt$ to molecules in $u\times w^{J}$ or the empty
set.

We can similarly define maps $g_{I}^{u}: \M(u_{I}\times v\times
w)\to v^{I}\times w^{I}$ and $g_{K}^{w}:\M(u\times v\times
w_{K})\to u\times v $ which send molecules to molecules or the
empty set. Moreover, we have natural homomorphisms
$f_{I}^{u}:\M(u_{I}\times v\times w)$ and $f_{K}^{w}:\M(u\times
v\times w_{K})$ of $\omega$-categories and we can see that
$F_{I}^{u}=g_{I}^{u}\circ f_{I}^{u}$ and $F_{K}^{w}=g_{K}^{w}\circ
f_{K}^{w}$. Therefore $F_{I}^{u}$ and $F_{K}^{w}$ sends molecules
to molecules or the empty set.

We have now proved the following theorem

\begin{theorem}
Molecules in $\uvwt$ are pairwise molecular.
\end{theorem}

\section{Properties of Pairwise Molecular Subcomplexes}

In this section, we prove some basic properties of pairwise
molecular subcomplexes. In the next section, we are going to show
that some of these properties characterise pairwise molecular
subcomplexes in $\uvwt$.

\begin{lemma} \label{mlc3}
Let $\Lambda$ be a pairwise molecular subcomplex of $u\times
v\times w$ and let $\uvws{_1}$ and $\uvws{_{2}}$ be distinct
maximal atoms in $\Lambda$.
\begin{enumerate}
\item
If $i_{1}=i_{2}$ and $\alpha_{1}=-\alpha_{2}$, then $\Lambda$ has
a maximal atom $\uvws{}$ with $i>i_{1}=i_{2}$,
$\vj{}\supset\vj{_{1}}\cap\vj{_{2}}$ and
$k\geq\min\{k_{1},k_{2}\}$.

\item
If $i_{1}=i_{2}$ and $\alpha_{1}=-\alpha_{2}$, then $\Lambda$ has
a maximal atom $\uvws{}$ with $i>i_{1}=i_{2}$,
$j\geq\min\{j_{1},j_{2}\}$ and
$\wk{}\supset\wk{_{1}}\cap\wk{_{2}}$.

\item
If $j_{1}=j_{2}$ and $\beta_{1}=-\beta_{2}$, then $\Lambda$ has a
maximal atom $\uvws{}$ with $j>j_{1}=j_{2}$,
$\ui{}\supset\ui{_{1}}\cap\ui{_{2}}$ and
$k\geq\min\{k_{1},k_{2}\}$.

\item
If $j_{1}=j_{2}$ and $\beta_{1}=-\beta_{2}$, then $\Lambda$ has a
maximal atom $\uvws{}$ with $j>j_{1}=j_{2}$,
$i\geq\min\{i_{1},i_{2}\}$ and
$\wk{}\supset\wk{_{1}}\cap\wk{_{2}}$.

\item
If $k_{1}=k_{2}$ and $\varepsilon_{1}=-\varepsilon_{2}$, then
$\Lambda$ has a maximal atom $\uvws{}$ with $k>k_{1}=k_{2}$,
$\ui{}\supset\ui{_{1}}\cap\ui{_{2}}$ and
$j\geq\min\{j_{1},j_{2}\}$.

\item
If $k_{1}=k_{2}$ and $\varepsilon_{1}=-\varepsilon_{2}$, then
$\Lambda$ has a maximal atom $\uvws{}$ with $k>k_{1}=k_{2}$,
$i\geq\min\{i_{1},i_{2}\}$ and
$\vj{}\supset\vj{_{1}}\cap\vj{_{2}}$.

\end{enumerate}
\end{lemma}
\begin{proof}
The proof of these conditions are similar, we only prove the
second one. Suppose that $i_{1}=i_{2}$ and
$\alpha_{1}=-\alpha_{2}$. Let $\lambda_{1}=\uvws{_{1}}$ and
$\lambda_{2}=\uvws{_{2}}$. Let $J=\min\{j_{1},j_{2}\}$. It is
evident that $F_{J}^{v}(\lambda_{1})=\ui{_{1}}\times
w^{J}[k_{1},\varepsilon_{1}]\subset F_{J}^{v}(\Lambda)$ and
$F_{J}^{v}(\lambda_{2})=\ui{_{2}}\times
w^{J}[k_{2},\varepsilon_{2}]\subset F_{J}^{v}(\Lambda)$. Since
$F_{J}^{v}(\Lambda)$ is a molecule in $u\times w^{J}$, it follows
from the formation of maximal atoms in $F_{J}^{v}(\Lambda)$ that
$F_{J}^{v}(\lambda_{1})$ or $F_{J}^{v}(\lambda_{2})$ is not
maximal in $F_{J}^{v}(\Lambda)$, and $F_{J}^{v}(\Lambda)$ has a
maximal atom $\mu=u[l,\sigma]\times w^{J}[n,\omega]$ with $l>i$
and $w[n,\omega]\supset\wk{_{1}}$ or
$w[n,\omega]\supset\wk{_{2}}$. By the definition of $F_{J}^{v}$,
it is easy to see that every maximal atom in $F_{J}^{v}(\Lambda)$
is an image of a maximal atom in $\Lambda$. Therefore $\Lambda$
has a maximal atom $\lambda=\uvws{}$ with $\ui{}=u[l,\sigma]$,
$j\geq J$ and $\wk{}=w[n,\omega]$, as required.

This completes the proof.
\end{proof}

The next property says that certain signs in a pair of `adjacent'
maximal atoms of a pairwise molecular subcomplexes are related.
Before we prove this property, we need to give the precise
definition of adjacency of a pair of maximal atoms.

\begin{definition}
Let $\Lambda$ be a subcomplex. A pair of distinct maximal atoms
$\lambda_{1}=\uvws{_{1}}$ and $\lambda_{2}=\uvws{_{2}}$ in
$\Lambda$ is {\em adjacent} if, for every maximal atom
$\lambda=\uvws{}$ in $\Lambda$ with $i\geq\min\{i_{1},i_{2}\}$,
$j\geq\min\{j_{1},j_{2}\}$ and $k\geq\min\{k_{1},k_{2}\}$, one has
$$
\begin{array}{rl}
&\min\{i_{1},i\}+\min\{j_{1},j\}+\min\{k_{1},k\}\\
=&\min\{i_{1},i_{2}\}+\min\{j_{1},j_{2}\}+\min\{k_{1},k_{2}\}
\end{array}
$$ or $$
\begin{array}{rl}
&\min\{i_{2},i\}+\min\{j_{2},j\}+\min\{k_{2},k\}\\ = &
\min\{i_{1},i_{2}\}+\min\{j_{1},j_{2}\}+\min\{k_{1},k_{2}\}.
\end{array}
$$
\end{definition}

The following proposition may be helpful to understand the concept
of adjacency.

\begin{prop}\label{adjacent_eq}
Let $\Lambda$ be a subcomplex satisfying condition 1 for pairwise
molecular subcomplexes. A pair of distinct maximal atoms
$\lambda_{1}=\uvws{_{1}}$ and $\lambda_{2}=\uvws{_{2}}$ in
$\Lambda$ is adjacent if and only if the following conditions
hold.

\begin{itemize}
\item
If $i_{1}=i_{2}$, then there is no maximal atom $\uvws{}$ such
that $i\geq i_{1}=i_{2}$, $j>\min\{j_{2},j_{2}\}$ and
$k>\min\{k_{1},k_{2}\}$.

\item
If $j_{1}=j_{2}$, then there is no maximal atom $\uvws{}$ such
that $j\geq j_{1}=j_{2}$, $i>\min\{i_{2},i_{2}\}$ and
$k>\min\{k_{1},k_{2}\}$.

\item
If $k_{1}=k_{2}$, then there is no maximal atom $\uvws{}$ such
that $k\geq k_{1}=k_{2}$, $i>\min\{i_{2},i_{2}\}$ and
$j>\min\{j_{2},j_{2}\}$.

\item
If $i_{1}>i_{2}$ and $j_{1}>j_{2}$, then there is no maximal atom
$\uvws{}$ such that $i>i_{2}$, $j\geq j_{2}$ and $k>k_{1}$; and
there is no maximal atom $\uvws{}$  such that $i\geq i_{2}$,
$j>j_{2}$ and $k>k_{1}$.

\item
If $i_{1}>i_{2}$ and $k_{1}>k_{2}$, then there is no maximal atom
$\uvws{}$ such that $i>i_{2}$, $j>j_{1}$ and $k\geq k_{2}$; and
there is no maximal atom $\uvws{}$  such that $i\geq i_{2}$,
$j>j_{1}$ and $k>k_{2}$.
\item
If $j_{1}>j_{2}$ and $k_{1}>k_{2}$, then there is no maximal atom
$\uvws{}$ such that $i>i_{1}$, $j>j_{2}$ and $k\geq k_{2}$; and
there is no maximal atom $\uvws{}$  such that $i>i_{1}$, $j\geq
j_{2}$ and $k>k_{2}$.
\end{itemize}
\end{prop}
\begin{proof}
The proof is a straightforward verification from the definition of
adjacency and condition 1 for pairwise molecular subcomplexes.
\end{proof}

\begin{example}
For the subcomplex in Example \ref{pairwise_example}, all the
adjacent pairs of maximal atoms are

\begin{center}
$u_{8}^{+}\times v_{2}^{+}\times w_{1}^{-}$ and $u_{5}^{-}\times
v_{2}^{+}\times w_{5}^{-}$;
\end{center}

\begin{center}
$u_{8}^{+}\times v_{2}^{+}\times w_{1}^{-}$ and $u_{9}^{+}\times
v_{1}^{-}\times w_{2}^{+}$;
\end{center}

\begin{center}
$u_{5}^{-}\times v_{2}^{+}\times w_{5}^{-}$ and $u_{1}^{-}\times
v_{2}^{+}\times w_{8}^{+}$;
\end{center}

\begin{center}
$u_{5}^{-}\times v_{2}^{+}\times w_{5}^{-}$ and $u_{9}^{+}\times
v_{1}^{-}\times w_{2}^{+}$;
\end{center}

\begin{center}
$u_{5}^{-}\times v_{2}^{+}\times w_{5}^{-}$ and $u_{4}^{-}\times
v_{1}^{-}\times w_{6}^{+}$;
\end{center}

\begin{center}
$u_{5}^{-}\times v_{2}^{+}\times w_{5}^{-}$ and $u_{8}^{-}\times
v_{0}^{-}\times w_{5}^{-}$;
\end{center}

\begin{center}
$u_{5}^{-}\times v_{2}^{+}\times w_{5}^{-}$ and $u_{5}^{-}\times
v_{0}^{+}\times w_{6}^{+}$;
\end{center}

\begin{center}
$u_{1}^{-}\times v_{2}^{+}\times w_{8}^{+}$ and $u_{9}^{+}\times
v_{1}^{-}\times w_{2}^{+}$;
\end{center}

\begin{center}
$u_{1}^{-}\times v_{2}^{+}\times w_{8}^{+}$ and $u_{4}^{-}\times
v_{1}^{-}\times w_{6}^{+}$;
\end{center}

\begin{center}
$u_{1}^{-}\times v_{2}^{+}\times w_{8}^{+}$ and $u_{0}^{+}\times
v_{1}^{+}\times w_{9}^{+}$;
\end{center}

\begin{center}\label{>}
$u_{1}^{-}\times v_{2}^{+}\times w_{8}^{+}$ and $u_{2}^{-}\times
v_{0}^{-}\times w_{9}^{+}$;
\end{center}

\begin{center}
$u_{9}^{+}\times v_{1}^{-}\times w_{2}^{+}$ and $u_{8}^{-}\times
v_{0}^{-}\times w_{5}^{-}$;
\end{center}

\begin{center}
$u_{4}^{-}\times v_{1}^{-}\times w_{6}^{+}$ and $u_{4}^{-}\times
v_{0}^{-}\times w_{7}^{+}$;
\end{center}

\begin{center}
$u_{0}^{+}\times v_{1}^{+}\times w_{9}^{+}$ and $u_{2}^{-}\times
v_{0}^{-}\times w_{9}^{+}$;
\end{center}

\begin{center}
$u_{8}^{-}\times v_{0}^{-}\times w_{5}^{-}$ and $u_{5}^{-}\times
v_{0}^{+}\times w_{6}^{+}$;
\end{center}

\begin{center}
$u_{5}^{-}\times v_{0}^{+}\times w_{6}^{+}$ and $u_{4}^{-}\times
v_{0}^{-}\times w_{7}^{+}$;
\end{center}

\begin{center}
$u_{4}^{-}\times v_{0}^{-}\times w_{7}^{+}$ and $u_{2}^{-}\times
v_{0}^{-}\times w_{9}^{+}$.
\end{center}

\end{example}

Let $\Lambda$ be a  subcomplex of $\uvwt$ satisfying condition 1
for pairwise molecular subcomplexes. Let $J$ be a fixed
non-negative integer. A maximal atom $\uvws{}$ in $\Lambda$ is
{\em $(v,J)$-projection maximal} if $j\geq J$ and there is no
maximal atom $\uvws{'}$ with $i'\geq i$, $J\leq j'<j$ and $k'\geq
k$.

Similarly, we can define  a maximal atom to be {\em
$(u,I)$-projection maximal} and {\em $(w,K)$-projection maximal.}

It is evident that a maximal atom $\lambda$ in $\Lambda$ is
$(v,J)$-projection maximal implies that $F_{J}^{v}(\lambda)$ is
maximal in $F_{J}^{v}(\Lambda)$. Conversely, for every maximal
atom $\mu$ in $F_{J}^{v}(\Lambda)$, there is a maximal atom $\mu'$
in $\Lambda$ such that $F_{J}^{v}(\mu')=\mu$. The following
proposition implies that $\mu'$ is actually $(v,J)$-projection
maximal.

\begin{prop}\label{projection_reason}
Let $\Lambda$ be a pairwise molecular subcomplex of $\uvwt$ and
$\lambda$ be a maximal atom in $\Lambda$. Then
\begin{enumerate}
\item
$\lambda$ is $(u,I)$-projection maximal if and only if
$F_{I}^{u}(\lambda)$ is maximal in $F_{J}^{v}(\Lambda)$.
\item
$\lambda$ is $(v,J)$-projection maximal if and only if
$F_{J}^{v}(\lambda)$ is maximal in $F_{J}^{v}(\Lambda)$.
\item
$\lambda$ is $(w,K)$-projection maximal if and only if
$F_{K}^{w}(\lambda)$ is maximal in $F_{K}^{w}(\Lambda)$.
\end{enumerate}
\end{prop}
\begin{proof}
The arguments for the three cases are similar. We only give the
proof for the second one.

Suppose that $\lambda$ is not $(v,J)$-projection maximal. Let
$\lambda=\uvws{}$. Then there is a maximal atom
$\lambda'=\uvws{'}$ in $\Lambda$ such that $J\leq j'<j$ $i'\geq i$
and $k'\geq k$. By condition 1 for pairwise molecular
subcomplexes, we have $i'>i$ or $k'>k$. If $\ui{'}\supset\ui{}$
and $\wk{'}\supset\wk{}$, then it is evident that
$F_{J}^{v}(\lambda)\subsetneqq F_{J}^{v}(\mu)$ so that
$F_{J}^{v}(\lambda)$ is not maximal in $F_{J}^{v}(\Lambda)$. Now
suppose that $\ui{'}\not\supset\ui{}$ or $\wk{'}\not\supset\wk{}$.
$\ui{'}=u[i,-\alpha]$ or $\wk{'}=w[k,-\varepsilon]$. Thus we can
get a maximal atom $\lambda''=\uvws{''}$ such that $J\leq j''<j$
and $\ui{''}\supset\ui{}$ and $\wk{''}\supset\wk{}$ by  applying
Lemma \ref{mlc3}. It follows that $F_{J}^{v}(\lambda)$ is not
maximal in $F_{J}^{v}(\Lambda)$.

Conversely, suppose that $F_{J}^{v}(\lambda)$ is not maximal in
$F_{J}^{v}(\lambda)$. It follows evidently from the definition
that $\lambda$ is not $(v,J)$-projection maximal.

This completes the proof.

\end{proof}

\begin{example}
For the subcomplex in Example \ref{pairwise_example}, there is no
$(v,J)$-projection maximal atoms for $J>2$. The $(v,2)$-projection
maximal atoms are $$u_{8}^{+}\times v_{2}^{+}\times w_{1}^{-},$$
$$u_{5}^{-}\times v_{2}^{+}\times w_{5}^{-},$$ $$u_{1}^{-}\times
v_{2}^{+}\times w_{8}^{+}.$$ The $(v,1)$-projection maximal atoms
are $$u_{9}^{+}\times v_{1}^{-}\times w_{2}^{+},$$
$$u_{5}^{-}\times v_{2}^{+}\times w_{5}^{-},$$ $$u_{4}^{-}\times
v_{1}^{-}\times w_{6}^{+},$$ $$u_{1}^{-}\times v_{2}^{+}\times
w_{8}^{+},$$ $$u_{0}^{+}\times v_{1}^{+}\times w_{9}^{+}.$$ The
$(v,0)$-projection maximal atoms are $$u_{9}^{+}\times
v_{1}^{-}\times w_{2}^{+},$$ $$u_{8}^{-}\times v_{0}^{-}\times
w_{5}^{-},$$ $$u_{5}^{-}\times v_{0}^{+}\times w_{6}^{+},$$
$$u_{4}^{-}\times v_{0}^{-}\times w_{7}^{+},$$ $$u_{2}^{-}\times
v_{0}^{-}\times w_{9}^{+}.$$ One can similarly work out all the
$(u,I)$-projection maximal atoms and all the $(w,K)$-projection
maximal atoms.

\end{example}

\begin{lemma}\label{adlowest}
Let $\Lambda$ be a pairwise molecular subcomplex of $u\times
v\times w$. Let $\lambda_{1}=\uvws{_1}$ and
$\lambda_{2}=\uvws{_2}$ be a pair of adjacent maximal atoms in
$\Lambda$.

\begin{enumerate}
\item
If $i_{1}>i_{2}$ and $j_{1}<j_{2}$, then there is a pair of
adjacent $(w,K)$-projection maximal atoms
$\lambda_{1}'=\uvws{_{1}'}$ and $\lambda_{2}'=\uvws{_{2}'}$ with
$K=\min\{k_{1},k_{2}\}$ such that $\ui{_{2}'}=\ui{_{2}}$,
$\vj{_{1}'}=\vj{_{1}}$ and $\min\{k_{1}',k_{2}'\}=K$.

\item
If $i_{1}>i_{2}$ and $k_{1}<k_{2}$, then there is a pair of
adjacent $(v,J)$-projection maximal atoms
$\lambda_{1}'=\uvws{_{1}'}$ and $\lambda_{2}'=\uvws{_{2}'}$ with
$J=\min\{j_{1},j_{2}\}$ such that $\ui{_{2}'}=\ui{_{2}}$,
$\wk{_{1}'}=\wk{_{1}}$ and $\min\{j_{1}',j_{2}'\}=J$.

\item
If $j_{1}>j_{2}$ and $k_{1}<k_{2}$, then there is a pair of
adjacent $(u,I)$-projection maximal atoms
$\lambda_{1}'=\uvws{_{1}'}$ and $\lambda_{2}'=\uvws{_{2}'}$ with
$I=\min\{i_{1},i_{2}\}$ such that $\vj{_{2}'}=\vj{_{2}}$,
$\wk{_{1}'}=\wk{_{1}}$ and $\min\{i_{1}',i_{2}'\}=I$.
\end{enumerate}
\end{lemma}
\begin{proof}
The arguments for these three cases are similar. We only give the
proof for the second case.

Let $\lambda_{1}'=\uvws{_{1}'}$ and $\lambda_{2}'=\uvws{_{2}'}$ be
the $(v,J)$-projection maximal atoms such that $i'_{t}\geq i_{t}$
and $k'_{t}\geq k_{t}$ for $t=1,2$. It follows from Lemma
\ref{mlc3} and the adjacency of $\lambda_{1}$ and $\lambda_{2}$
that $\ui{_{2}'}=\ui{_{2}}$, $\wk{_{1}'}=\wk{_{1}}$ and
$\min\{j_{1}',j_{2}'\}=\min\{j_{1},j_{2}\}$, and $\lambda_{1}'$
and $\lambda_{2}'$ are adjacent, as required.

This completes the proof.
\end{proof}

Now we can prove the sign conditions for pairwise molecular
subcomplexes.

\begin{prop} \label{mlc2}
Let $\Lambda$ be a pairwise molecular subcomplex of $u\times
v\times w$. Then the following sign conditions hold.

Sign conditions: for a pair of adjacent maximal atoms
$\lambda_{1}=\uvws{_{1}}$ and $\lambda_{2}=\uvws{_{2}}$ in
$\Lambda$, let $i=\min\{i_{1},i_{2}\}$, $j=\min\{j_{1},j_{2}\}$
and $k=\min\{k_{1},k_{2}\}$.

\begin{enumerate}

\item
If $i=i_{1}<i_{2}$ and $j=j_{2}<j_{1}$, then
$\beta_{2}=-(-)^{i}\alpha_{1}$;
\item
if $i=i_{1}<i_{2}$ and $k=k_{2}<k_{1}$, then
$\varepsilon_{2}=-(-)^{i+j}\alpha_{1}$;
\item
if $j=j_{1}<j_{2}$ and $k=k_{2}<k_{1}$, then
$\varepsilon_{2}=-(-)^{j}\beta_{1}$.
\end{enumerate}
\end{prop}
\begin{proof}
Suppose that $i_{1}>i_{2}$ and $k_{1}<k_{2}$. Let
$J=\min\{j_{1},j_{1}\}$. We must prove
$\varepsilon_{1}=-(-)^{i_{2}+J}\alpha_{2}$.

According to Lemma \ref{adlowest}, we may assume that
$\lambda_{1}$ and $\lambda_{2}$ are $(v,J)$-projection maximal. It
is evident that $F_{J}^{v}(\lambda_{1})=\ui{_1}\times
w^{J}[k_{1},\varepsilon_{1}]$ and
$F_{J}^{v}(\lambda_{2})=\ui{_2}\times
w^{J}[k_{2},\varepsilon_{2}]$, and they are maximal atoms in the
molecule $F_{J}^{v}(\Lambda)$. Moreover, by the adjacency of
$\lambda_{1}$ and $\lambda_{2}$, we can see that
$F_{J}^{v}(\lambda_{1})$ and $F_{J}^{v}(\lambda_{2})$ are adjacent
maximal atoms in $F_{J}^{v}(\Lambda)$. Since $F_{J}^{v}(\Lambda)$
is a molecule in $u\times w^{J}$, we have
$\varepsilon_{1}=-(-)^{i_{2}+J}\alpha_{2}$, as required.

The other cases can be proved similarly.

This completes the proof.
\end{proof}

Compared with the properties for molecules in the product of two
globes, there is a new feature caused by the middle factors, as
follows.

\begin{prop}\label{mlc4}
Let $\Lambda$ be a pairwise molecular subcomplex of $\uvwt$. Let
$\lambda_{1}=\uvws{_{1}}$ and $\lambda_{2}=\uvws{_{2}}$ be a pair
of adjacent maximal atoms in $\Lambda$. If $i_{1}>i_{2}$,
$k_{1}<k_{2}$ and $\min\{j_{1},j_{2}\}>0$, then there is a maximal
atom $\lambda=\uvws{}$ such that $j=\min\{j_{1},j_{2}\}-1$,
$i>i_{2}$ and $k>k_{1}$.
\end{prop}
\begin{proof}
Let $J=\min\{j_{1},j_{2}\}$. According to Lemma \ref{adlowest}, we
may assume that $\lambda_{1}$ and $\lambda_{2}$ are
$(v,J)$-projection maximal. There are several cases, as follows.

1. Suppose that both $\lambda_{1}$ and $\lambda_{2}$ are
$(v,J-1)$-projection maximal. Then $F_{J-1}(\lambda_{1})$ and
$F_{J-1}(\lambda_{2})$ are maximal atoms in the molecule
$F_{J-1}(\Lambda)$ of the $\omega$-complex $u\times w^{J-1}$. It
is evident that $F_{J-1}(\lambda_{t})=\ui{_{t}}\times
w^{J-1}[k_{t},\varepsilon_{t}]$ for $t=1,2$, and
$\varepsilon_{1}=-(-)^{i_{2}+J}\alpha_{2}$ by Proposition
\ref{mlc2}. Hence, according to the formation of molecules in
$u\times w^{J-1}$, we can see that $F_{J-1}(\lambda_{1})$ and
$F_{J-1}(\lambda_{2})$ are not adjacent in $u\times w^{J-1}$. So
$F_{J-1}(\Lambda)$ has a maximal atom $\mu=\ui{}\times
w^{J-1}[k,\varepsilon]$ with $i>i_{2}$ and $k>k_{1}$. It follows
that there is a maximal atom $\lambda=\uvws{}$ such that
$F_{J-1}(\lambda)=\mu$ and hence $i>i_{2}$ and $k>k_{1}$. By the
adjacency of $\lambda_{1}$ and $\lambda_{2}$, we must have
$j=J-1$. Therefore $\lambda$ is as required.

2. Suppose that $\lambda_{1}$ is not $(v,J-1)$-projection maximal.
Then there is a maximal atom $\lambda_{1}'=\uvws{_{1}'}$ with
$j_{1}'\geq J-1$ such that $i_{1}'\geq i_{1}$, $j_{1}'<j_{1}$ and
$k_{1}'\geq k_{1}$. It is evident that $j_{1}'=J-1$. If
$k_{1}'>k_{1}$, then $\lambda_{1}'$ is as required. Suppose that
$k_{1}'=k_{1}$ and $\varepsilon_{1}'=-\varepsilon_{1}$. By
applying Lemma \ref{mlc3} to $\lambda_{1}$ and $\lambda_{1}'$, one
can get a maximal atom as required. The argument is similar if
$\lambda_{2}$ is not $(v,J-1)$-projection maximal and
$i_{2}'>i_{2}$, or if $\lambda_{2}$ is not $(v,J-1)$-projection
maximal, $i_{2}'=i_{2}$ and $\alpha_{2}'=-\alpha_{2}$, where
$\lambda_{2}'=\uvws{_{2}'}$ is a maximal atom with $j_{2}'=J-1$
such that $i_{2}'\geq i_{2}$, $j_{2}'<j_{2}$ and $k_{2}'\geq
k_{2}$. There remains the case that $\ui{_{2}'}=\ui{_{2}}$ and
$\wk{_{1}'}=\wk{_1}$. In this case, since $\lambda_{1}'$ and
$\lambda_{2}'$ are maximal atoms with  $j_{1}'=j_{2}'=J-1$ and
$\varepsilon_{1}'=-(-)^{i_{2}'+J}\alpha_{2}'$, it follows from
Proposition \ref{mlc2} that $\lambda_{1}'$ and $\lambda_{2}'$ are
not adjacent in $\Lambda$. Therefore $\Lambda$ has a maximal atom
as required.

This completes the proof.
\end{proof}

\begin{prop}
Let $\Lambda$ be a pairwise molecular subcomplex of $\uvwt$. If
$\Lambda$ has three pairwise adjacent maximal atoms
$\lambda_{1}=u[i_{1},\alpha_{1}]\times v[j,\beta_{1}]\times
w[k,\varepsilon_{1}]$ $\lambda_{2}=u[i,\alpha_{2}]\times
v[j_{2},\beta_{2}]\times w[k,\varepsilon_{2}]$ and
$\lambda_{3}=u[i,\alpha_{3}]\times v[j,\beta_{3}]\times
w[k_{3},\varepsilon_{2}]$ with $i_{1}>i$, $j_{2}>j$ and $k_{3}>k$,
then $\alpha_{2}=\alpha_{3}$ or $\beta_{1}=\beta_{3}$ or
$\varepsilon_{1}=\varepsilon_{2}$.
\end{prop}
\begin{proof}
Suppose otherwise that $\alpha_{2}=-\alpha_{3}$ and
$\beta_{1}=-\beta_{3}$ and $\varepsilon_{1}=-\varepsilon_{2}$.
Applying Lemma \ref{mlc3} to $\lambda_{1}$ and $\lambda_{2}$, one
can get a maximal atom $\lambda'=\uvws{'}$ with $k'>k$,
$\ui{'}\supset u[i,\alpha_{2}]$ and $j'\geq j$. Since
$\lambda_{1}$ and $\lambda_{3}$ are adjacent, we must have $i'=i$
and $\alpha'=\alpha_{2}=-\alpha_{3}$. Since $\lambda_{2}$ and
$\lambda_{3}$ are adjacent, we also have $j'=j$. Note that
$\lambda'$ and $\lambda_{3}$ are distinct, we get a contradiction
to the first condition for pairwise molecular subcomplexes.

This completes the proof.
\end{proof}

\begin{prop}\label{sign_4_4}
Let $\Lambda$ be a pairwise molecular subcomplex in $\uvwt$. Let
$\lambda_{1}=\uvws{_{1}}$ and $\lambda_{2}=\uvws{_{2}}$ be maximal
atoms in $\Lambda$.
\begin{enumerate}
\item\label{sign_4_41}
If $i_{1}<i_{2}$ and $j_{1}>j_{2}$, and  if there is no maximal
atom $\uvws{}$ such that $i>i_{1}$, $j>j_{2}$ and
$k\geq\min\{k_{1},k_{2}\}$, then
$\beta_{2}=-(-)^{i_{1}}\alpha_{1}$;
\item \label{sign_4_42}
if $i_{1}<i_{2}$ and $k_{1}>k_{2}$, and  if there is no maximal
atom $\uvws{}$ such that $i>i_{1}$, $j\geq\min\{j_{1},j_{2}\}$ and
$k>k_{2}$, then
$\varepsilon_{2}=-(-)^{i_{1}+\min\{j_{1},j_{2}\}}\alpha_{1}$;
\item\label{sign_4_43}
if $j_{1}<j_{2}$ and $k_{1}>k_{2}$, and  if there is no maximal
atom $\uvws{}$ such that $i\geq\min\{i_{1},i_{2}\}$, $j>j_{1}$ and
$k>k_{2}$, then $\varepsilon_{2}=-(-)^{j_{1}}\beta_{1}$.
\end{enumerate}
\end{prop}

\begin{note}
We some times say a pair of maximal atoms as in condition
\ref{sign_4_41} to be {\em $(1,2)$-adjacent}; a pair of maximal
atoms as in condition \ref{sign_4_42} to be {\em
$(1,3)$-adjacent}; and a pair of maximal atoms as in condition
\ref{sign_4_43} to be {\em $(2,3)$-adjacent}. It is evident that
two maximal atoms are $(r,s)$-adjacent ($1\leq r<s\leq 3$) if they
are adjacent. However, in general, the reverse is not true. For
example, in the pairwise molecular subcomplex in Example
\ref{pairwise_example}, the maximal atoms $u_{5}^{-}\times
v_{2}^{+}\times w_{5}^{-}$ and $u_{4}^{-}\times v_{0}^{-}\times
w_{7}^{+}$ are $(1,2)$-adjacent, but they are not adjacent.

\end{note}

\begin{proof}
The arguments for the first and the third cases are similar. We
give the proof for the first and the second case.

In the first case, let $\lambda_{1}'=\uvws{_{1}'}$ be the maximal
atom in $\Lambda$ with $i_{1}'\geq i_{1}$, $j_{1}'>j_{2}$ and
$k_{1}'\geq \min\{k_{1},k_{2}\}$ such that  $j_{1}'$ is minimal;
let $\lambda_{2}'=\uvws{_{2}'}$ be the maximal atom in $\Lambda$
with $i_{2}'>i_{1}$, $j_{2}'\geq j_{2}$ and $k_{2}'\geq
\min\{k_{1},k_{2}\}$ such that  $i_{2}'$ is minimal. According to
the assumption and Lemma \ref{mlc3}, we have
$\ui{_{1}'}=\ui{_{1}}$ and $\vj{_{2}'}=\vj{_{2}}$. It is evident
that $\lambda_{1}'$ and $\lambda_{2}'$ are adjacent. It follows
from the sign condition for $\lambda_{1}'$ and $\lambda_{2}'$ that
$\beta_{2}=-(-)^{i_{1}}\alpha_{1}$, as required.

In the second case, we claim that $\lambda_{1}$ is adjacent to
$\lambda_{2}$ so that
$\varepsilon_{2}=-(-)^{i_{1}+\min\{j_{1},j_{2}\}}\alpha_{1}$, as
required. In fact, suppose otherwise that $\lambda_{1}$ and
$\lambda_{2}$ are not adjacent. Then $j_{1}\neq j_{2}$. We may
assume that $j_{1}<j_{2}$. In this case,  there exists a maximal
atom $\lambda_{1}'=\uvws{_{1}'}$ such that $i_{1}'\geq i_{1}$,
$j_{1}'\geq j_{1}$ and $k_{1}'>k_{2}$. By the assumption, we have
$i_{1}'=i_{1}$. By condition 1 for pairwise molecular
subcomplexes, we have $j_{1}'>j_{1}$ and $k_{2}<k_{1}'<k_{1}$. It
follows from Lemma \ref{mlc4} that there is a maximal atom
$\mu=\lmns{}$ such that $l>i_{1}'=i_{1}$,
$m\geq\min\{j_{1}',j_{2}\}-1\geq j_{1}$ and $n>k_{2}$. This
contradicts the assumption.

This completes the proof.

\end{proof}

\section{An Alternative Description for Pairwise Molecular
Subcomplexes}

In this section, we give an alternative description for pairwise
molecular subcomplexes of $\uvwt$, as follows.

\begin{theorem}\label{pairwise_eq}
Let $\Lambda$ be a subcomplex of $\uvwt$. Then $\Lambda$ is
pairwise molecular if and only if the following conditions hold.

\begin{enumerate}
\item \label{m1}
There are no distinct maximal atoms $\uvws{}$ and $\uvws{'}$ such
that $i\leq i'$, $j\leq j'$ and $k\leq k'$.

\item \label{m2}
Sign conditions: for a pair of adjacent maximal atoms
$\lambda_{1}=\uvws{_{1}}$ and $\lambda_{2}=\uvws{_{2}}$ in
$\Lambda$, let $i=\min\{i_{1},i_{2}\}$, $j=\min\{j_{1},j_{2}\}$
and $k=\min\{k_{1},k_{2}\}$. If $i=i_{1}<i_{2}$ and
$j=j_{2}<j_{1}$, then $\beta_{2}=-(-)^{i}\alpha_{1}$; if
$i=i_{1}<i_{2}$ and $k=k_{2}<k_{1}$, then
$\varepsilon_{2}=-(-)^{i+j}\alpha_{1}$; if $j=j_{1}<j_{2}$ and
$k=k_{2}<k_{1}$, then $\varepsilon_{2}=-(-)^{j}\beta_{1}$.

\item \label{m3}
Let $\uvws{_{1}}$ and $\uvws{_{2}}$ be a pair of maximal atoms in
$\Lambda$. If $i_{1}=i_{2}$ and $\alpha_{1}=-\alpha_{2}$, then
$\Lambda$ has a maximal atom $\uvws{}$ with $i>i_{1}=i_{2}$,
$j\geq\min\{j_{1},j_{2}\}$ and $k\geq\min\{k_{1},k_{2}\}$; if
$j_{1}=j_{2}$ and $\beta_{1}=-\beta_{2}$, then $\Lambda$ has a
maximal atom $\uvws{}$ with $j>j_{1}=j_{2}$,
$i\geq\min\{i_{1},i_{2}\}$ and $k\geq\min\{k_{1},k_{2}\}$; if
$k_{1}=k_{2}$ and $\varepsilon_{1}=-\varepsilon_{2}$, then
$\Lambda$ has a maximal atom $\uvws{}$ with $k>k_{1}=k_{2}$,
$i\geq\min\{i_{1},i_{2}\}$ and $j\geq\min\{j_{1},j_{2}\}$.

\item \label{m4}
If $\Lambda$ has a pair of adjacent maximal atoms
$\lambda_{1}=\uvws{_{1}}$ and $\lambda_{2}=\uvws{_{2}}$ with
$i_{2}<i_{1}$, $k_{1}<k_{2}$ and $\min\{j_{1},j_{2}\}>0$, then
$\Lambda$ has a maximal atom $\uvws{}$ with $i>i_{2}$,
$j=\min\{j_{1},j_{2}\}-1$ and $k>k_{1}$.

\item \label{m5}
If $\Lambda$ has three pairwise adjacent maximal atoms
$\lambda_{1}=u[i_{1},\alpha_{1}]\times v[j,\beta_{1}]\times
w[k,\varepsilon_{1}]$ $\lambda_{2}=u[i,\alpha_{2}]\times
v[j_{2},\beta_{2}]\times w[k,\varepsilon_{2}]$ and
$\lambda_{3}=u[i,\alpha_{3}]\times v[j,\beta_{3}]\times
w[k_{3},\varepsilon_{3}]$ with $i_{1}>i$, $j_{2}>j$ and $k_{3}>k$,
then $\alpha_{2}=\alpha_{3}$ or $\beta_{1}=\beta_{3}$ or
$\varepsilon_{1}=\varepsilon_{2}$.

\end{enumerate}

\end{theorem}

\begin{note}\label{m4_note}
In condition \ref{m4}, it is easy to see that
$\beta=-(-)^{i_{2}}\alpha_{2}$ and $\varepsilon_{1}=-(-)^{j}\beta$
by sign conditions and condition \ref{m3}.
\end{note}

In the last section, we have proved that the five conditions in
the theorem are necessary for a pairwise molecular subcomplexes.
The sufficiency is implied by the following Proposition \ref{lvp2}
and the comments after the proposition.

Some of the following lemmas are preliminaries for the proof of
Proposition \ref{lvp2}, while some of them are designed for better
understanding the five conditions in Theorem \ref{pairwise_eq}.

\begin{lemma}
Let $\Lambda$ be a subcomplex of $\uvwt$ satisfying the five
conditions in Theorem \ref{pairwise_eq}. For a pair of adjacent
maximal atoms $\lambda_{1}=\uvws{_{1}}$ and
$\lambda_{2}=\uvws{_{2}}$ in $\Lambda$, let
$i=\min\{i_{1},i_{2}\}$, $j=\min\{j_{1},j_{2}\}$ and
$k=\min\{k_{1},k_{2}\}$.
\begin{enumerate}
\item
If $i=i_{1}<i_{2}$ and $j=j_{1}<j_{2}$, then
$\beta_{1}=(-)^{i}\alpha_{1}$;
\item
if $i=i_{1}<i_{2}$ and $k=k_{1}<k_{2}$, then
$\varepsilon_{1}=(-)^{i+j}\alpha_{1}$;
\item
if $j=j_{1}<j_{2}$ and $k=k_{1}<k_{2}$, then
$\varepsilon_{2}=(-)^{j}\beta_{1}$.
\end{enumerate}
\end{lemma}
\begin{proof}
Suppose that $i=i_{1}<i_{2}$ and $j=j_{1}<j_{2}$. By condition 1,
we have $k_{1}>k_{2}$. It follows from sign conditions that
$\varepsilon_{2}=-(-)^{i+j}\alpha_{1}$ and
$\varepsilon_{2}=-(-)^{j}\beta_{1}$. Thus
$\beta_{1}=(-)^{i}\alpha_{1}$, as required.

The other cases can be argued similarly.
\end{proof}

\begin{lemma} \label{lm3}
Let $\Lambda$ be a subcomplex satisfying the five conditions in
Theorem \ref{pairwise_eq}. Let $\uvws{_{1}}$ and $\uvws{_{2}}$  be
a pair of maximal atoms in  $\Lambda$.
\begin{enumerate}
\item
If $i_{1}=i_{2}$, $\alpha_{1}=-\alpha_{2}$, $j_{1}<j_{2}$ and
$k_{1}>k_{2}$, then $\Lambda$ has a maximal atom $\uvws{'}$ such
that $i'>i_{1}=i_{2}$, $\vj{'}\supset\vj{_{1}}$ and
$\wk{'}\supset\wk{_{2}}$;
\item
if $j_{1}=j_{2}$, $\beta_{1}=-\beta_{2}$, $i_{1}<i_{2}$ and
$k_{1}>k_{2}$, then $\Lambda$ has a maximal atom $\uvws{'}$ such
that $j'>j_{1}=j_{2}$, $\ui{'}\supset\ui{_{1}}$ and
$\wk{'}\supset\wk{_{2}}$;
\item
if $k_{1}=k_{2}$, $\varepsilon_{1}=-\varepsilon_{2}$,
$i_{1}<i_{2}$ and $j_{1}>j_{2}$, then $\Lambda$ has a maximal atom
$\uvws{'}$ with $k'>k_{1}=k_{2}$, $\ui{'}\supset\ui{_{1}}$ and
$\vj{'}\supset\vj{_{2}}$.
\end{enumerate}
\end{lemma}
\begin{proof}
The arguments for the above three cases are similar. We prove only
the first case.

Let $\lambda_{1}=\uvws{_{1}}$ and $\lambda_{2}=\uvws{_{2}}$. Let
$i=i_{1}=i_{2}$. Suppose that $\lambda_{1}$ and $\lambda_{2}$ are
not adjacent. Then, by the definition of adjacency, $\Lambda$ has
a maximal atom $\lambda_{1}'=\uvws{_{1}'}$ with $i_{1}'\geq i$,
$j_{1}<j_{1}'<j_{2}$ $k_{2}<k_{1}'<k_{1}$. If $i_{1}'>i$, then
$\lambda_{1}'$ is as required by the lemma. If $i_{1}'=i$, then
$\alpha_{1}'=-\alpha_{1}$ or $\alpha_{1}'=-\alpha_{2}$. By
repeating this process, we can  get either a maximal atom as
required or a pair of adjacent maximal atoms
$\lambda_{1}''=\uvws{_{1}''}$ and $\lambda_{2}''=\uvws{_{2}''}$
with $i_{1}''=i_{2}''=i_{1}=i_{2}$, $\alpha_{1}''=-\alpha_{2}''$,
$\vj{_{1}}\subset \vj{_{1}''}\cap\vj{_{2}''}$ and
$\wk{_{2}}\subset \wk{_{1}''}\cap\wk{_{2}''}$. In the following
proof,  we may assume that $\lambda_{1}=\uvws{_{1}}$ and
$\lambda_{2}=\uvws{_{2}}$ are adjacent.

Let $\alpha_{1}=-\gamma$, $j=j_{1}<j_{2}$, $\beta=\beta_{1}$,
$k=k_{2}$ and $\varepsilon=\varepsilon_{2}$. Thus
$\alpha_{2}=\gamma$, $k=k_{2}<k_{1}$ and
$\varepsilon=-(-)^{j}\beta$. By condition 3, $\Lambda$ has a
maximal atom $\lambda'=\uvws{'}$ with $i'>i$, $j'\geq j$ and
$k'\geq k$. We choose $\lambda'$ such that $i'$ is minimal. By
condition 1, we have $j'<j_{2}$ and $k'<k_{1}$. Since
$\lambda_{1}$ and $\lambda_{2}$ are adjacent, we have $j'=j$ or
$k'=k$. Now there are two cases, as follows.

1. If $j'=j$ and $k'>k$, we claim that $\beta'=\beta$ which means
that $\lambda'$ is as required.

Indeed, suppose otherwise that $\beta'=-\beta$, then, by condition
3 in Theorem \ref{pairwise_eq}, there is a maximal atom
$\lambda''=\uvws{''}$ in $\Lambda$ with $j''> j$, $i''\geq i$ and
$k''\geq k'>k$. This contradicts the adjacency of $\lambda_{1}$
and $\lambda_{2}$.

The argument for the case  $j'>j$ and $k'=k$ is similar.

2. Suppose that $j'=j$ and $k'=k$.  By the choice of $\lambda'$,
it is easy to see that $\lambda'$ is adjacent to both
$\lambda_{1}$ and $\lambda_{2}$. So $\beta'=-(-)^{i}\gamma$ and
$\varepsilon'=(-)^{i+j}\gamma$. Thus
$\varepsilon'=-(-)^{j}\beta'$. Since $\lambda_{1}$ is adjacent to
$\lambda_{2}$, we can see that
$\varepsilon_{2}=-(-)^{j}\beta_{1}$. By condition 5, one has
$\beta'=\beta_{1}$ or $\varepsilon'=\varepsilon_{2}$. Therefore
$\beta'=\beta_{1}$ and $\varepsilon'=\varepsilon_{2}$ which means
that $\lambda'$ is as required.

This completes the proof of the lemma.
\end{proof}

\begin{lemma}\label{cond5}
Let $\Lambda$ be a subcomplex of $\uvwt$ satisfying conditions 1
and 2 in Theorem \ref{pairwise_eq} . Then $\Lambda$ satisfies
condition 5 if and only if for any triple of pairwise adjacent
maximal atoms $\lambda_{1}=u[i_{1},\alpha_{1}]\times
v[j,\beta_{1}]\times w[k,\varepsilon_{1}]$,
$\lambda_{2}=u[i,\alpha_{2}]\times v[j_{2},\beta_{2}]\times
w[k,\varepsilon_{2}]$ and $\lambda_{3}=u[i,\alpha_{3}]\times
v[j,\beta_{3}]\times w[k_{3},\varepsilon_{2}]$ with $i_{1}>i$,
$j_{2}>j$ and $k_{3}>k$, there is a maximal atom in $\Lambda$
containing $u[i,\gamma]\times v[j,(-)^{i}\gamma]\times
w[k,(-)^{i+j}\gamma]$ for $\gamma=+$ or $\gamma=-$.
\end{lemma}
\begin{proof}
Suppose that $\Lambda$ satisfies condition 5 in Theorem
\ref{pairwise_eq}. Then $\alpha_{2}=\alpha_{3}$ or
$\beta_{1}=\beta_{3}$ or $\varepsilon_{1}=\varepsilon_{2}$.

Suppose that $\alpha_{2}=\alpha_{3}$ and let
$\gamma=\alpha_{1}=\alpha_{2}$. Then $\beta_{1}=-(-)^{i}\gamma$
and $\varepsilon_{1}=-(-)^{i+j}\gamma$ by the sign conditions. If
$\beta_{3}=(-)^{i}\gamma$, then $u[i,\gamma]\times
v[j,(-)^{i}\gamma]\times w[k,(-)^{i+j}\gamma]\subset \lambda_{3}$
and $u[i,-\gamma]\times v[j,-(-)^{i}\gamma]\times
w[k,-(-)^{i+j}\gamma]\subset \lambda_{1}$, as required. If
$\beta_{3}=-(-)^{i}\gamma$, then $\varepsilon_{2}=(-)^{i+j}\gamma$
by the sign condition for $\lambda_{2}$ and $\lambda_{3}$.
Therefore $u[i,\gamma]\times v[j,(-)^{i}\gamma]\times
w[k,(-)^{i+j}\gamma]\subset \lambda_{2}$ and $u[i,-\gamma]\times
v[j,-(-)^{i}\gamma]\times w[k,-(-)^{i+j}\gamma]\subset
\lambda_{1}$, as required.

The other cases can be argued similarly.

Conversely, suppose that $\Lambda$ has a maximal atom
$\lambda'=\uvws{'}$ containing $u[i,\gamma]\times
v[j,(-)^{i}\gamma]\times w[k,(-)^{i+j}\gamma]$ for $\gamma=+$ or
$\gamma=-$. By the pairwise adjacency of $\lambda_{1}$,
$\lambda_{2}$ and $\lambda_{3}$, it is easy to see that $\lambda'$
must be $\lambda_{1}$, $\lambda_{2}$ or $\lambda_{3}$. If
$\lambda'=\lambda_{1}$, then $\beta_{1}=(-)^{i}\gamma$ and
$\varepsilon_{1}=(-)^{i+j}\gamma$. It follows from the sign
condition in Theorem \ref{pairwise_eq} that
$\alpha_{2}=\alpha_{3}=-\gamma$, as required by condition 5 in
Theorem \ref{pairwise_eq}.

The other cases can be argued similarly.

This completes the proof.
\end{proof}

\begin{lemma} \label{cd3''}
Let $\Lambda$ be a subcomplex satisfying conditions 1, 2, 3 and 5
in Theorem \ref{pairwise_eq}. Then $\Lambda$ satisfies condition 4
if and only if, for any pair of adjacent maximal atoms
$\lambda_{1}=\ui{_{1}}\times \vj{_{1}} \times \wk{}$ and
$\lambda_{2}=\ui{_{}}\times \vj{_{2}} \times \wk{_2}$ with
$i<i_{1}$, $k<k_{2}$ and $j=\min\{j_{1},j_{2}\}>0$, there is a
maximal atom containing
$
u[i,-\alpha]\times v[j-1,-(-)^{i}\alpha]\times
w[k,(-)^{i+j}\alpha].
$
\end{lemma}
\begin{proof}
The necessity is obvious. We now prove the sufficiency. We can
assume that $j=j_{1}\leq j_{2}$, and hence
$\varepsilon=-(-)^{i+j}\alpha$.

By the assumption, there is a maximal atom $\lambda'=\uvws{'}$
containing
$
u[i,-\alpha]\times v[j-1,-(-)^{i}\alpha]\times
w[k,(-)^{i+j}\alpha].
$
Thus $i'\geq i$, $j'\geq j-1$ and $k'\geq k$. We claim that
$j'=j-1$ and hence $\beta'=-(-)^{i}\alpha$.

Indeed, suppose otherwise that $j'>j-1$. Then $i'=i$ by the
adjacency of $\lambda_{1}$ and $\lambda_{2}$. Hence
$\alpha'=-\alpha$. Note that the proof of Lemma \ref{lm3} does not
use condition 4. So, by applying Lemma \ref{lm3} to $\lambda_{2}$
and $\lambda'$, one can get a maximal atom
$\lambda_{3}=\uvws{_{3}}$ with $i_{3}>i$, $j_{3}\geq j$ and
$\wk{_{3}}\supset \wk{_{2}}\cap\wk{'}$. Since $k_{2}>k$,
$\wk{'}\supset w[k,(-)^{i+j}\alpha]=w[k,-\varepsilon]$ and
$k_{2}\neq k'$, we can see that $\lambda_{3}$ is distinct from
$\lambda_{1}$ and $\lambda_{2}$. This contradicts the adjacency of
$\lambda_{1}$ and $\lambda_{2}$.

Now, if $i'>i$ and $k'>k$, then $\lambda'$ is as required. Suppose
that $i'=i$. Then $\alpha'=-\alpha$ and $k'>k$. Thus, by Lemma
\ref{lm3}, $\Lambda$ has a maximal atom $\lambda''=\uvws{''}$ with
$i''>i$, $j''=j-1$ (by the adjacency of $\lambda_{1}$ and
$\lambda_{2}$), $\beta''=\beta'=-(-)^{i}\alpha$ and $k''\geq
k'>k$, with the required property. The argument for the case
$k'=k$ is similar.

This completes the proof of the lemma.
\end{proof}

\begin{lemma}\label{lvp1}
Let $\Lambda$ be a subcomplex of $\uvwt$ satisfying the five
conditions in Theorem \ref{pairwise_eq} . Then
\begin{enumerate}
\item\label{lv1}
Every maximal atom $\uvws{}$ with $j=J$ is $(v,J)$-projection
maximal.
\item\label{lv2}
For every maximal atom $\uvws{}$ with $j\geq J$, there is a
$(v,J)$-projection maximal atom $\uvws{'}$ such that
$\ui{}\subset\ui{'}$ and $\wk{}\subset\wk{'}$.

\item\label{lv3}
All the $(v,J)$-projection maximal atoms, if exist, can be listed
as $\lambda_{1}$, $\cdots$, $\lambda_{S}$ with
$\lambda_{s}=\uvws{_{s}}$ such that $i_{1}>\cdots>i_{S}$ and
$k_{1}<\cdots<k_{S}$ and
$\varepsilon_{s-1}=(-)^{i_{s}+J}\alpha_{s}$ for $1<s\leq S$.

\item\label{lv4}
For two consecutive $(v,J)$-projection maximal atoms
$\lambda_{s-1}$ and $\lambda_{s}$ in the above list, either
$j_{s-1}=J$ or $j_{s}=J$ for $1<s\leq S$.
\end{enumerate}
\end{lemma}
\begin{proof}
In this proof, all the maximal atoms refer to maximal atoms with
dimension of second factors not less than $J$.

Condition \ref{lv1} follows from the definition of projection
maximal.

To prove condition \ref{lv2}, suppose that $\lambda=\uvws{}$ is
not $(v,J)$-projection maximal. Then there is a $(v,J)$-projection
maximal atom $\lambda_{1}=\uvws{_{1}}$ such that $i_{1}\geq i$,
$j_{1}<j$ and $k_{1}\geq k$. If $\ui{_{1}}\supset\ui{}$ and
$\wk{_{1}}\supset\wk{}$, then $\lambda_{1}$ is as required.
Suppose that $\ui{_{1}}\not\supset\ui{}$. Then $i_{1}=i$ and
$\alpha_{1}=-\alpha$. Moreover, we have $k_{1}>k$ by condition 1
in Theorem \ref{pairwise_eq}. Hence, by Lemma \ref{lm3}, there is
a maximal atom $\lambda_{2}=\uvws{_{2}}$ in $\Lambda$ such that
$i_{2}>i$, $j_{2}<j$, $\vj{_{1}}\subset\vj{_{2}}$ and
$\wk{}\subset\wk{_{2}}$. This shows that $\ui{_{2}}\supset\ui{}$
and $\wk{_{2}}\supset\wk{}$ and $j_{2}<j$. Therefore condition
\ref{lv2} holds by induction. The argument for the case
$\wk{_{1}}\not\supset\wk{}$ is similar.

Condition \ref{lv4} follows easily from condition 4 in Theorem
\ref{pairwise_eq}, while Condition \ref{lv3} follows easily from
the definition of projection maximal and condition 1 and condition
2 (sign conditions) in Theorem \ref{pairwise_eq}.

This completes the proof.
\end{proof}

\begin{prop}\label{lvp2}
Let $\Lambda$ be a subcomplex of $\uvwt$ satisfying the five
conditions in Theorem \ref{pairwise_eq}. Let the
$(v,J)$-projection maximal atoms in $\Lambda$ be listed as
$\lambda_{1}$, $\cdots$, $\lambda_{S}$ with
$\lambda_{s}=\uvws{_{s}}$ for $1\leq s\leq S$ such that
$i_{1}>\cdots>i_{S}$ and $k_{1}<\cdots<k_{S}$. Then
$F_{J}^{v}(\Lambda)=\ui{_{1}}\times \wJk{_{1}}\cup\cdots\cup
\ui{_{S}}\times\wJk{_{S}}$.
\end{prop}
\begin{proof}
This is a direct consequence of Proposition \ref{pro_atom} for
$F_{J}^{v}$ and Lemma \ref{lvp1}.
\end{proof}

\begin{coro}
Let $\Lambda$ be a subcomplex of $\uvwt$ satisfying the five
conditions in Theorem \ref{pairwise_eq}. Then $F_{J}^{v}(\Lambda)$
is a molecule in $u\times w^{J}$ or the empty set for every
non-negative integer $J$.
\end{coro}

We can similarly show that $F_{I}^{u}(\Lambda)$ and
$F_{K}^{w}(\Lambda)$ are molecules or the empty set for a
subcomplex $\Lambda$ of $\uvwt$ satisfying the five conditions in
Theorem \ref{pairwise_eq}. This completes the proof of Theorem
\ref{pairwise_eq}.

\section{Sources and Targets of Pairwise Molecular Subcomplexes}
In this section, we study source and target operators
$d_{p}^{\gamma}$ on pairwise molecular subcomplexes in $\uvwt$.
The main result in this section is that $d_{p}^{\gamma}\Lambda$ is
pairwise molecular for every pairwise molecular subcomplex
$\Lambda$ of $\uvwt$.

Recall that $d_{p}^{\gamma}\Lambda$ is a union of interiors of
atoms. We first prove that $d_{p}^{\gamma}\Lambda$ is a subcomplex
of $\Lambda$.

\begin{lemma}\label{lm5.3}
Let $\Lambda$ be a subcomplex of $u\times v\times w$ and
$\lambda=\uvws{}$ be a $p$-dimensional atom in $\Lambda$ with
$\Int\lambda\subset\dpg\Lambda$.
\begin{enumerate}
\item
If there is an atom $\lambda'=\uvws{'}$ in $\Lambda$ such that
$\lambda\subset\lambda'$ and $i'>i$, then $\alpha=\gamma$.
\item
If there is an atom $\lambda'=\uvws{'}$ in $\Lambda$ such that
$\lambda\subset\lambda'$ and $j'>j$, then $\beta=(-)^{i}\gamma$.
\item
If there is an atom $\lambda'=\uvws{'}$ in $\Lambda$ such that
$\lambda\subset\lambda'$ and $k'>k$, then
$\varepsilon=(-)^{i+j}\gamma$.
\end{enumerate}
\end{lemma}
\begin{proof}
Suppose that there is an atom $\lambda'=\uvws{'}$ in $\Lambda$
such that $\lambda\subset\lambda'$ and $i'>i$. Then
$\lambda\subset
u[i+1,\alpha']\times\vj{}\times\wk{}\subset\Lambda$ and $\dim
(u[i+1,\alpha']\times\vj{}\times\wk{})=p+1$. Since
$\Int\lambda\subset\dpg\Lambda$, we have $\lambda\subset\dpg
(u[i+1,\alpha']\times\vj{}\times\wk{})$. It follows easily from
Lemma \ref{dtatom} that $\alpha=\gamma$, as required.

The arguments for other cases are similar.
\end{proof}

\begin{prop}\label{int_sign}
Let $\Lambda$ be a pairwise molecular subcomplex of $\uvwt$. Let
$\lambda=\uvws{}$ be a $p$-dimensional atom such that
$\Int\lambda\subset\dpg\Lambda$.
\begin{enumerate}
\item
If there is a maximal atom $\lambda'$ in $\Lambda$ such that
$i'>i$, $j'\geq j$ and $k'\geq k$, then $\alpha=\gamma$;
\item
if there is a maximal atom $\lambda'$ in $\Lambda$ such that
$i'\geq i$, $j'>j$ and $k'\geq k$, then $\beta=(-)^{i}\gamma$;
\item
if there is a maximal atom $\lambda'$ in $\Lambda$ such that
$i'\geq i$, $j'\geq j$ and $k'>k$, then
$\varepsilon=(-)^{i+j}\gamma$.
\end{enumerate}
\end{prop}

\begin{proof}
The arguments for the three cases are similar. We give the proof
for the first case.

Since $\Int\lambda\subset\Lambda$, there is a maximal atom
$\mu=\lmns{}$ such that $\lambda\subset\mu$. If $\mu$ can be
chosen such that $l>i$, then we have $\alpha=\gamma$ by Lemma
\ref{lm5.3} , as required. In the following proof, we may assume
that $\mu$ cannot be chosen such that $l>i$ so that $\ul{}=\ui{}$.

Suppose that there is a maximal atom $\lambda'=\uvws{'}$ such that
$i'>i$ $j'\geq j$ and $k'\geq k$. Then we have
$\vj{'}=v[j,-\beta]$ or $\wk{'}=w[k,-\varepsilon]$. By applying
Lemma \ref{lm3}, we may assume that $\vj{'}=v[j,-\beta]$ and
$m>j$, or assume that $\wk{'}=w[k,-\varepsilon]$ and $n>k$.

Suppose that $\wk{'}=w[k,-\varepsilon]$ and $n>k$. Then
$\varepsilon=(-)^{i+j}\gamma$ by Lemma \ref{lm5.3}. If
$\min\{j',m\}=j$, and if $\lambda'$ is $(1,3)$-adjacent to $\mu$,
then $\alpha=\sigma=-(-)^{i+j}\varepsilon'=\gamma$ by Proposition
\ref{sign_4_4}, as required. Otherwise, by the definition of
adjacency or condition \ref{m4} in Theorem \ref{pairwise_eq}, we
may chose $\lambda'$ and $\mu$ such that $\min\{k',n\}>k$ so that
$\vj{'}=v[j,-\beta]$; according Lemma \ref{lm3}, we may also
assume that $m>j$; thus $\beta=(-)^{i}\gamma$. In this case,
according to the assumptions, $\lambda'$ must be $(1,2)$-adjacent
to $\mu$. It follows from Proposition \ref{sign_4_4} that
$\alpha=\sigma=-(-)^{i}\beta'=\gamma$, as required.

Suppose that $\vj{'}=v[j,-\beta]$ and $m>j$. Then we can get
$\alpha=\gamma$, as required, by a similar argument.

This completes the proof.

\end{proof}

\begin{lemma}\label{include_p-1}
Let $\Lambda$ be a pairwise molecular subcomplex of $\uvwt$. Let
$\lambda=\uvws{}$ be a $p-1$ dimensional atom  such that
$\Int\lambda\subset\dpg\Lambda$.
\begin{enumerate}
\item
If there is a maximal atom $\lambda'$ in $\Lambda$ with
$\lambda'\supset\lambda$ such that $i'>i$ and $j'>j$, then
$\alpha=\gamma$ or $\beta=-(-)^{i}\gamma$;
\item
if there is a maximal atom $\lambda'$ in $\Lambda$ with
$\lambda'\supset\lambda$ such that $i'>i$ and $k'>k$, then
$\alpha=\gamma$ or $\varepsilon=-(-)^{i+j}\gamma$;
\item
if there is a maximal atom $\lambda'$ in $\Lambda$ with
$\lambda'\supset\lambda$ such that $j'>j$ and $k'>k$, then
$\beta=(-)^{i}\gamma$ or $\varepsilon=-(-)^{i+j}\gamma$.
\end{enumerate}
\end{lemma}
\begin{proof}
The arguments for the three cases are similar. We give the proof
for the first one.

Suppose that there is a maximal atom $\lambda'$ in $\Lambda$ with
$\lambda'\supset\lambda$ such that $i'>i$ and $j'>j$. Then
$\lambda\subset u[i+1,\alpha']\times
v[j+1,\beta']\times\wk{}\subset\Lambda$. Since
$\Int\lambda\subset\dpg\Lambda$, we have
$\lambda\subset\dpg(u[i+1,\alpha']\times
v[j+1,\beta']\times\wk{})$. It follows easily from Lemma
\ref{dtatom} that $\alpha=\gamma$ or $\beta=-(-)^{i}\gamma$, as
required.

This completes the proof.

\end{proof}

\begin{prop}\label{sign_p-1}
Let $\Lambda$ be a pairwise molecular subcomplex of $\uvwt$. Let
$\lambda=\uvws{}$ be a $p-1$ dimensional atom  such that
$\Int\lambda\subset\dpg\Lambda$.
\begin{enumerate}
\item
If there is a maximal atom $\lambda'$ in $\Lambda$ such that
$i'>i$, $j'>j$ and $k'\geq k$, then $\alpha=\gamma$ or
$\beta=-(-)^{i}\gamma$;
\item
if there is a maximal atom $\lambda'$ in $\Lambda$ such that
$i'>i$, $j'\geq j$ and $k'>k$, then $\alpha=\gamma$ or
$\varepsilon=-(-)^{i+j}\gamma$;
\item
if there is a maximal atom $\lambda'$ in $\Lambda$ such that
$i'\geq i$, $j'>j$ and $k'>k$, then $\beta=(-)^{i}\gamma$ or
$\varepsilon=-(-)^{i+j}\gamma$.
\end{enumerate}
\end{prop}
\begin{proof}
The arguments for case 1 and case 3 are similar. We give the
proofs for case 1 and case 2.

1. Suppose that there is a maximal atom $\lambda'=\uvws{'}$ in
$\Lambda$ such that $i'>i$, $j'>j$ and $k'\geq k$. If $\lambda'$
can be chosen such that $\lambda'\supset\lambda$, then we have
$\alpha=\gamma$ or $\beta=-(-)^{i}\gamma$, as required, by Lemma
\ref{include_p-1}. In the following, we assume that $\lambda'$
cannot be chosen such that $\lambda'\supset\lambda$ so that
$\wk{'}=w[k,-\varepsilon]$. Let $\lambda_{1}=\uvws{_{1}}$ be a
maximal atom in $\Lambda$ such that $\lambda\subset\lambda_{1}$.
Then $\ui{_{1}}=\ui{}$ or $\vj{_{1}}=\vj{}$ by the assumption.
According to Lemma \ref{lm3}, we may assume that $k_{1}>k$. Now
there are several cases, as follows.

Suppose that $\lambda_{1}$ cannot be chosen such that $i_{1}>i$ or
$j_{1}>j$. According to Lemma \ref{lm3}, it is easy to see that
$\lambda_{1}$ and $\lambda'$ are adjacent. Thus $\alpha=\gamma$
(when $\varepsilon=-(-)^{i+j}\gamma$)  or $\beta=-(-)^{i}\gamma$
(when $\varepsilon=(-)^{i+j}\gamma$) by sign conditions, as
required.

Suppose that $\lambda_{1}$ can be chosen such that $i_{1}>i$.
Suppose also that $\alpha=-\gamma$. Then $\vj{_{1}}=\vj{}$ by the
assumptions. According to Lemma \ref{include_p-1}, it is easy to
see that $\varepsilon=-(-)^{i+j}\gamma$, hence
$\varepsilon'=(-)^{i+j}\gamma$. It is evident that $\lambda_{1}$
and $\lambda'$ are $(2,3)$-adjacent. It follows from Proposition
\ref{sign_4_4} that $\beta=\beta_{1}=-(-)^{i}\gamma$, as required.

Suppose that $\lambda_{1}$ can be chosen such that $j_{1}>j$ and
that $\lambda_{1}$ cannot be chosen such that $i_{1}>i$. Suppose
also that $\beta=(-)^{i}\gamma$. According to Lemma \ref{lm3},
condition \ref{m4} in Theorem \ref{pairwise_eq} and the
assumptions, it is easy to see that $\lambda_{1}$ and $\lambda'$
are adjacent and $\min\{j',j_{1}\}=j+1$. It follows from
condition \ref{m4} in Theorem \ref{pairwise_eq} that there is a
maximal atom $\lambda''=\uvws{''}$ such that $i''>i$, $j''=j$ and
$k''>k$. Moreover, we have $\beta''=-(-)^{i}\alpha$ by Note
\ref{m4_note}. By the assumptions, we have
$\beta''=-\beta=-(-)^{i}\gamma$. It follows that $\alpha=\gamma$,
as required.

This completes the proof for case 1.

2. Suppose that there is a maximal atom $\lambda'$ in $\Lambda$
such that $i'>i$, $j'\geq j$ and $k'>k$. If $\lambda'$ can be
chosen such that $\lambda'\supset\lambda$, then we have
$\alpha=\gamma$ or $\varepsilon=-(-)^{i}\gamma$, as required, by
Lemma \ref{include_p-1}.  In the following, we assume that
$\lambda'$ cannot be chosen such that $\lambda'\supset\lambda$ so
that $\vj{'}=v[j,-\beta]$. Let $\lambda_{1}=\uvws{_{1}}$ be a
maximal atom in $\Lambda$ such that $\lambda\subset\lambda_{1}$.
Then $\ui{_{1}}=\ui{}$ or $\wk{_{1}}=\wk{}$ by the assumption.
According to Lemma \ref{lm3}, we may assume that $j_{1}>j$. Now
there are several cases, as follows.

Suppose that $\lambda_{1}$ cannot be chosen such that $i_{1}>i$ or
$k_{1}>k$. Then it is easy to see that $\lambda_{1}$ adjacent to
$\lambda'$. It follows from sign conditions that $\alpha=\gamma$
(when $\beta'=-(-)^{i}\gamma$) or $\varepsilon=-(-)^{i+j}\gamma$
(when $\beta'=(-)^{i}\gamma$), as required.

Suppose that $\lambda_{1}$ can be chosen such that $i_{1}>i$.
Suppose also that $\alpha=-\gamma$. Then $\beta=-(-)^{i}\gamma$ by
Lemma \ref{include_p-1}, and hence $\beta'=-\beta=(-)^{i}\gamma$.
Moreover, we can see that $\lambda'$ and $\lambda_{1}$ are
$(2,3)$-adjacent. It follows from Lemma \ref{sign_4_4} that
$\varepsilon=-(-)^{i+j}\gamma$, as required.

Suppose that $\lambda_{1}$ can be chosen such that $k_{1}>k$. By a
similar argument as in the above case, we can get $\alpha=\gamma$
or $\varepsilon=-(-)^{i+j}\gamma$, as required.

This completes the proof.

\end{proof}

\begin{lemma}\label{lm5.6}
Let $x$ be a union of interiors of atoms in an $\omega$-complex.
Then $x$ is a subcomplex if and only if for every atom $a$ in $x$
with $\Int a\subset x$ and every atom $b$ with $b\subset a$, one
has $\Int b\subset x$.
\end{lemma}
\begin{proof} The necessity is evident. To prove the sufficiency, it
suffices to prove that for every atom $a$ with $\Int a\subset x$
we have $a\subset x$. Note that $a$ can be written as a union of
interiors of atoms $b$ with $b\subset a$. The sufficiency follows.
\end{proof}

\begin{prop}\label{prop5.5}
Let $\Lambda$ be a pairwise molecular subcomplex of $\uvwt$. Then
$\dpg\Lambda$ is a subcomplex.
\end{prop}
\begin{proof}
From Lemma \ref{interior}, we have already known that
$\dpg\Lambda$ is a union of interiors of atoms. By Lemma
\ref{lm5.6}, it suffices to prove that for every atom $\lambda$
with $\Int\lambda\subset\dpg\Lambda$ and every atom $\lambda_{1}$
with $\lambda_{1}\subset\lambda$, one has
$\Int\lambda_{1}\subset\dpg\Lambda$. It is evident that there is a
sequence $\lambda\supset\lambda_{1}^{1}\supset\lambda^{2}_{1}
\cdots\supset\lambda_{1}$ such that the difference of the
dimensions of any pair of consecutive atoms is $1$. We may assume
that $\dim\lambda_{1}=\dim\lambda-1$.

Let $\lambda=\uvws{}$. Since
$\Int\lambda\subset\dpg\Lambda\subset\Lambda$ and $\Lambda$ is a
subcomplex, we have $\lambda_{1}\subset\lambda\subset\Lambda$ and
$\dim\lambda_{1}\leq\dim\lambda\leq p$. Suppose that $\mu=\lmns{}$
is an atom with $\dim\mu=p+1$ and
$\lambda_{1}\subset\mu\subset\Lambda$. We must prove
$\lambda_{1}\subset\dpg\mu$.

If $\lambda\subset\mu$, then
$\lambda_{1}\subset\lambda\subset\dpg\mu$ since
$\lambda\subset\dpg\Lambda$. If $l>i+1$ or $m>j+1$ or $n>k+1$,
then we evidently have $\lambda_{1}\subset\dpg\mu$ by Lemma
\ref{dtatom}. In the following, we may further assume that
$\lambda\not\subset\mu$ and that $l\leq i+1$ and $m\leq j+1$ and
$n\leq k+1$. Thus $\ui{}\not\subset\ul{}$ or
$\vj{}\not\subset\vm{}$ or $\wk{}\not\subset\wn{}$; moreover, if
$\ui{}\not\subset\ul{}$, then we have $\ul{}=u[i,-\alpha]$ or
$\ul{}=u[i-1,\sigma]$, we also have $\vj{}\subset\vm{}$ and
$\wk{}\subset\wn{}$; if $\vj{}\not\subset\vm{}$, then we have
$\vm{}=v[j,-\beta]$  or $\vm{}=v[j-1,\tau]$, we also have
$\ui{}\subset\ul{}$ and $\wk{}\subset\wn{}$; if
$\wk{}\not\subset\wn{}$, then we have $\wn{}=w[k,-\varepsilon]$ or
$\wn{}=u[k-1,\omega]$, we also have $\ui{}\subset\ul{}$ and
$\vj{}\subset\vm{}$;  Note that $\dim\mu=p+1$ and $\dim\lambda\leq
p$, we now have 3 cases, as follows.

1. Suppose that $\ul{}=u[i,-\alpha]$ or $\vm{}=v[j,-\beta]$ or
$\wn{}=w[k,-\varepsilon]$; suppose also that $\dim\lambda=p$. Then
only one of the equations $l=i+1$, $m=j+1$ and $n=k+1$ holds. The
arguments for the three cases are similar, we only give the proof
for the case $\vm{}=v[j,-\beta]$ and $\dim\lambda=p$. In this
case, we must have $\mu=u[i+1,\sigma]\times
v[j,-\beta]\times\wk{}$ or $\mu=\ui{}\times v[j,-\beta]\times
w[k+1,\omega]$. Hence $\lambda_{1}$ is of the form
$\lambda_{1}=\ui{}\times v[j-1,\tilde{\beta}]\times\wk{}$.

Suppose that $\mu=u[i+1,\sigma]\times v[j,-\beta]\times\wk{}$.
Then there is a maximal atom $\mu'=\lmns{'}$ such that $l'>i$,
$m'\geq j$ and $n'\geq k$. It follows from Proposition
\ref{int_sign} that $\alpha=\gamma$. This implies
$\lambda_{1}\subset\dpg\mu$, as required.

Suppose that $\mu=\ui{}\times v[j,-\beta]\times w[k+1,\omega]$.
Then there is a maximal atom $\mu'=\lmns{'}$ such that $l'\geq i$,
$m'\geq j$ and $n'>k$. It follows from Proposition \ref{int_sign}
that $\varepsilon=(-)^{i+j}\gamma$. This implies
$\lambda_{1}\subset\dpg\mu$, as required.

2. Suppose that $l=i-1$ or $m=j-1$ or $n=k-1$; suppose also that
$\dim\lambda=p$. The arguments for these three cases are similar.
We only give the proof for the case $m=j-1$ and $\dim\lambda=p$.
In this case, we have $l=i+1$ and $n=k+1$ because $\dim\mu=p+1$;
we also have $\lambda_{1}=\ui{}\times v[j-1,\tau]\times\wk{}$. To
get $\lambda_{1}\subset\dpg\mu$, by Lemma \ref{dtatom}, it
suffices to prove that $\alpha=\gamma$ or
$\varepsilon=(-)^{i+j}\gamma$.

Let $\lambda'=\uvws{'}$ be a maximal atom in $\Lambda$ such that
$\lambda\subset\lambda'$. Let $\mu'=\lmns{'}$ be a maximal atom in
$\Lambda$ such that $\mu\subset\mu'$. If $\lambda'$ can be chosen
such that $i'>i$ or $k'>k$, then we have $\alpha=\gamma$ or
$\varepsilon=(-)^{i+j}\gamma$ which implies that
$\lambda_{1}\subset\dpg\mu$, as required. If there is a maximal
atom $\mu''=\lmns{''}$ with $\mu''\supset\lambda_{1}$ and $m''>m$
such that $l''>i$ or $n''>k$, then, by Proposition \ref{int_sign},
we have $\alpha=\gamma$ or $\varepsilon=(-)^{i+j}\gamma$ which
implies that $\lambda_{1}\subset\dpg\mu$, as required. Now suppose
that  $\lambda'$ cannot be chosen  such that $i'>i$ or $k'>k$.
Suppose also that there is no maximal atom $\mu''=\lmns{''}$ with
$\mu''\supset\lambda_{1}$ and $m''>m$ such that $l''>i$ or
$n''>k$. Then $\ui{'}=\ui{}$, $\wk{'}=\wk{}$ and $\vm{'}=\vm{}$.
Moreover, it is easy to see that $\lambda'$ and $\mu'$ are
adjacent. It follows from the sign condition for $\lambda'$ and
$\mu'$ that $\alpha=\gamma$ (when $\tau=-(-)^{i}\gamma$) or
$\varepsilon=(-)^{i+j}\gamma$ (when $\tau=-(-)^{i}\gamma$). This
implies that $\lambda_{1}\subset\dpg\mu$, as required.

3. Suppose that $\ul{}=u[i,-\alpha]$ or $\vm{}=v[j,-\beta]$ or
$\wn{}=w[k,-\varepsilon]$; suppose also that $\dim\lambda=p-1$.
The arguments for these three cases are similar. We only give the
proof for the case $\vm{}=v[j,-\beta]$ and $\dim\lambda=p-1$. In
this case, we have $l=i+1$ and $n=k+1$.  Moreover, we can see that
$\lambda_{1}$ is of the form $\lambda_{1}=\ui{}\times
v[j-1,\tilde{\beta}]\times \wk{}$. According to Lemma
\ref{sign_p-1}, we have $\alpha=\gamma$ or
$\varepsilon=-(-)^{i+j}\gamma$. This implies that
$\lambda_{1}\subset\dpg\mu$, as required.

This completes the proof.
\end{proof}

We can now start to prove that $d_{p}^{\gamma}\Lambda$ is pairwise
molecular for a molecular subcomplex $\Lambda$ in $\uvwt$ by
verifying conditions in Definition \ref{pairwise_def}.

By Lemma \ref{interior}, the maps $F_{I}^{u}$, $F_{J}^{v}$ and
$F_{K}^{w}$ are defined on $\dpg\Lambda$ for every subcomplex
$\Lambda$ of $u\times v\times w$.

\begin{prop}\label{prop5.6}
Let $\Lambda$ be a pairwise molecular subcomplex of $u\times
v\times w$. If $p\geq J$ and $F_{J}^{v}(\Lambda)\neq\emptyset$,
then $F_{J}^{v}(\dpg\Lambda)=d_{p-J}^{\gamma}F_{J}^{v}(\Lambda)$;
therefore $F_{J}^{v}(\dpg\Lambda)$ is a molecule in $u\times
w^{J}$.
\end{prop}
\begin{proof} Firstly, we prove that
$d_{p-J}^{\gamma}F_{J}^{v}(\Lambda)\subset
F_{J}^{v}(\dpg\Lambda)$.

Let $\ui{}\times\wJk{}$ be an atom in $u\times w^{J}$ such that
$\Int(\ui{}\times\wJk{})\subset
d_{p-J}^{\gamma}F_{J}^{v}(\Lambda)$. We must show that
$\Int(\ui{}\times\wJk{})\subset F_{J}^{v}(\dpg\Lambda)$. Clearly,
we have $\ui{}\times\wJk{}\subset F_{J}^{v}(\Lambda)$. So it is
easy to see that $\ui{}\times v[J,\beta]\times \wk{}\subset
\Lambda$ for some sign $\beta$. We are going to prove
$\Int(\ui{}\times v[J,\beta]\times \wk{})\subset \dpg\Lambda$ by
verifying conditions in Lemma \ref{lm5.2}. It is evident that
$\dim(\ui{}\times v[J,\beta]\times \wk{})\leq p$. To verify the
other conditions, we consider two cases, as follows.

1. Suppose that $\beta$ can be chosen such that
$\beta=(-)^{i}\gamma$. Suppose also that there is an atom
$\ui{'}\times \vj{'}\times \wk{'}\subset\Lambda$  such that
$\uvws{}\subset\uvws{'}$. Then
$\ui{}\times\wJk{}\subset\ui{'}\times\wJk{'}$ in $u\times w^{J}$.
Therefore $\ui{}\times\wJk{}\subset
d_{p-J}^{\gamma}(\ui{'}\times\wJk{'})$. It follows easily that
$\ui{}\times v[J,\beta]\times\wk{}\subset\dpg (\uvws{'})$, as
required by the second condition of Lemma \ref{lm5.2}.

2. Suppose that $\beta$ cannot be chosen such that
$\beta=(-)^{i}\gamma$. Suppose also that there is an atom
$\ui{'}\times \vj{'}\times \wk{'}\subset\Lambda$  such that
$\uvws{}\subset\uvws{'}$. Then $j'=J$ and
$\beta'=\beta=-(-)^{i}\gamma$ from Lemma \ref{lm5.3}. By an
argument similar to the above case, it is easy to see that
$\ui{}\times v[J,\beta]\times\wk{}\subset\dpg(\uvws{'})$, as
required by the second condition of Lemma \ref{lm5.2}.

We have now shown that $\Int(\ui{}\times
v[J,\beta]\times\wk{})\subset\dpg\Lambda$. It follows that
$\Int(\ui{}\times\wJk{})=F_{J}^{v}(\Int(\ui{}\times
v[J,\beta]\times\wk{}))\subset F_{J}^{v}(\dpg\Lambda)$. This
completes the proof that
$d_{p-J}^{\gamma}F_{J}^{v}(\Lambda)\subset
F_{J}^{v}(\dpg\Lambda)$.

Conversely, let $\lambda=\ui{}\times\wJk{}$ be an atom such that
$\Int\lambda\subset F_{J}^{v}(\dpg\Lambda)$. We must show that
$\Int\lambda\subset d_{p-J}^{\gamma}F_{J}^{v}(\Lambda)$. It is
easy to see that there is an atom $\ui{}\times\vj{}\times\wk{}$ in
$\Lambda$ such that
$\Int(\ui{}\times\vj{}\times\wk{})\subset\dpg\Lambda$ and $j\geq
J$. Since $\dpg\Lambda$ is a subcomplex of $\uvwt{}$, we can see
that $\ui{}\times v[J,\beta']\times\wk{}\subset\dpg\Lambda$ for
some sign $\beta'$. It follows that $\dim\lambda\leq p-J$.
Clearly, we have $\lambda\subset F_{J}^{v}(\Lambda)$. This shows
that the first condition of Lemma \ref{lm5.2} is satisfied. To
verify the other condition of \ref{lm5.2}, let $\mu=\ul{}\times
w^{J}[n,\omega]$ be an atom in $F_{J}^{v}(\Lambda)$ such that
$\lambda\subset\mu$ and $\dim\mu=p-J+1$. We must prove that
$\lambda\subset d_{p-J}^{\gamma}\mu$. It is evident that there is
an atom $\ul{}\times v[J,\tau']\times\wn{}$ in $\Lambda$ for some
sign $\tau'$. If $l>i+1$ or $n>k+1$, then it is evident that
$\lambda\subset\dpg\mu$, as required. In the following proof, we
may assume that $l\leq i+1$ and $n\leq k+1$ so that
$\dim\lambda=p-J$ or $\dim\lambda=p-J-1$. Now there are various
cases, as follows.

Suppose that $\beta'$ and $\tau'$ can be chosen such that
$\beta'=\tau'$. Then $\ui{}\times
v[J,\beta']\times\wk{}\subset(\dpg\Lambda\cap (\ul{}\times
v[J,\tau']\times\wn{}) \subset\dpg(\ul{}\times
v[J,\tau']\times\wn{})$ by Proposition \ref{dpg_sub}. It follows
easily that $\lambda\subset d_{p-J}^{\gamma}\mu$, as required.

Suppose that $\beta'$ and $\tau'$ cannot be chosen such that
$\beta'=\tau'$. Suppose also that $J>0$. Since $\dpg\Lambda$ is a
subcomplex, we know that $\ui{}\times
v[J-1,\pm]\times\wk{}\subset\dpg\Lambda$. This implies that
$\ui{}\times v[J-1,\pm]\times\wk{}\subset \dpg(\ul{}\times
v[J,\tau']\times\wn{})$. It follows easily that $\lambda\subset
d_{p-J}^{\gamma}\mu$, as required.

There remain the case that $J=0$ and $\beta'$ and $\tau'$ cannot
be chosen such that $\beta'=\tau'$. If $\dim\lambda=p$, by
Proposition \ref{int_sign}, we can get $\alpha=\gamma$ when $l>i$,
while $\varepsilon=(-)^{i}\gamma$ when $n>k$; thus
$\lambda\subset\dpg\mu$, as required. If $\dim\lambda=p-1$, then
$l=i+1$ and $n=k+1$; by Proposition \ref{sign_p-1}, we can get
$\alpha=\gamma$ or $\varepsilon=-(-)^{i}\gamma$; thus
$\lambda\subset\dpg\mu$, as required.

This completes the proof.

\end{proof}

We can prove the following two results by similar arguments.

\begin{prop}\label{prop5.7}
Let $\Lambda$ be a pairwise molecular subcomplex of $u\times
v\times w$. If $p\geq I$ and $F_{I}^{u}(\Lambda)\neq\emptyset$,
then $F_{I}^{u}(\dpg\Lambda)=d_{p-I}^{\gamma}F_{I}^{u}(\Lambda)$.
\end{prop}

\begin{prop}\label{prop5.8}
Let $\Lambda$ be a pairwise molecular subcomplex of $u\times
v\times w$. If $p\geq K$ and $F_{K}^{w}(\Lambda)\neq\emptyset$,
then $F_{K}^{w}(\dpg\Lambda)=d_{p-K}^{\gamma}F_{K}^{w}(\Lambda)$.
\end{prop}

We also need to show that $\dpg\Lambda$ satisfies condition 1 for
pairwise molecular subcomplexes for a pairwise molecular
subcomplex $\Lambda$.

\begin{lemma}
Let $\Lambda$ be a pairwise molecular subcomplex. Then there are
no distinct maximal atoms  $\lambda=\uvws{}$ and
$\lambda'=\uvws{'}$ in $\dpg\Lambda$ such that $i\leq i'$, $j\leq
j'$ and $k\leq k'$.
\end{lemma}

\begin{proof}
Let $\lambda=\uvws{}$ and $\lambda'=\uvws{'}$ be a pair of maximal
atom in $\dpg\Lambda$ such that $i\leq i'$, $j\leq j'$ and $k\leq
k'$.  We must prove that $\lambda=\lambda'$. Suppose that
$\dim\lambda<p$ or $\dim\lambda'<p$. By Lemma \ref{maximal_dpg},
we can see that $\lambda$ is a maximal atom in $\Lambda$ when
$\dim\lambda<p$ and $\lambda'$ is a maximal atom in $\Lambda$ when
$\dim\lambda'<p$. According to condition 1 for pairwise molecular
subcomplex $\Lambda$ of $\uvwt$, it is evident that
$\lambda=\lambda'$, as required. In the following argument, we may
assume that $\dim\lambda=p$ and $\dim\lambda=p$ so that $i=i'$,
$j=j'$ and $k=k'$.

Now suppose otherwise that $\lambda\neq\lambda'$. Then
$\alpha'=-\alpha$ or $\beta'=-\beta$ or
$\varepsilon'=-\varepsilon$. We may assume that $\alpha'=-\alpha$.
In this case, we have $F_{j}^{v}(\lambda)\subset
F_{j}^{v}(\dpg\Lambda) =d_{p-j}^{\gamma}F_{j}^{v}(\Lambda)$ and
similarly $F_{j}^{v}(\lambda') \subset
d_{p-j}^{\gamma}F_{j}^{v}(\Lambda)$ by Proposition \ref{prop5.6}.
Since $\dim F_{j}^{v}(\lambda)=\dim F_{j}^{v}(\lambda')=p-j$ and
$\dim (d_{p-j}^{\gamma}F_{j}^{v}(\Lambda))\leq p-j$, we can see
that $F_{j}^{v}(\lambda)$ and $F_{j}^{v}(\lambda')$ are maximal
atoms in the molecule $d_{p-j}^{\gamma}F_{j}^{v}(\Lambda)$. Note
that $F_{j}^{v}(\lambda)=\ui{}\times w^{j}[k,\varepsilon]$ and
$F_{j}^{v}(\lambda')=u[i,-\alpha]\times w^{j}[k,\varepsilon']$. We
get a contradiction to condition 1 in Theorem \ref{2globe}.

The arguments for the case $\beta'=-\beta$ or
$\varepsilon'=-\varepsilon$ are similar.

This completes the proof.

\end{proof}

Now we can prove the main result in this section.

\begin{prop}
Let $\Lambda$ be a pairwise molecular subcomplex. Then so is
$\dpg\Lambda$.
\end{prop}
\begin{proof} We have shown  that $\dpg\Lambda$ satisfies condition 1 for
pairwise molecular subcomplexes. Moreover, by Proposition
\ref{prop5.6}, \ref{prop5.7} and \ref{prop5.8}, we have
$F_{I}^{u}(\dpg\Lambda)=d_{p-I}^{\gamma}F_{I}^{u}(\Lambda)$,
$F_{J}^{v}(\dpg\Lambda)=d_{p-J}^{\gamma}F_{J}^{v}(\Lambda)$ and
$F_{K}^{w}(\dpg\Lambda)=d_{p-K}^{\gamma}F_{K}^{w}(\Lambda)$ for
all $I\geq p$, $J\geq p$ and $K\geq p$. Since
$F_{I}^{u}(\Lambda)$, $F_{J}^{v}(\Lambda)$ and
$F_{K}^{w}(\Lambda)$ are molecules or the empty set for all $I$,
$J$ and $K$, we can see that $F_{I}^{u}(\dpg\Lambda)$,
$F_{J}^{v}(\dpg\Lambda)$ and $F_{K}^{w}(\dpg\Lambda)$ are
molecules or the empty set for all $I$, $J$ and $K$. It follows
that $\dpg\Lambda$ is pairwise molecular.

This completes the proof.

\end{proof}

The following theorem gives the algorithm of constructing
$d_{p}^{\gamma}\Lambda$ for a pairwise molecular subcomplex
$\Lambda$ in $\uvwt$.

\begin{theorem} \label{dt}
Let $\Lambda$ be a pairwise molecular subcomplex. Then the
dimension of every maximal atom in $d_{p}^{\gamma}\Lambda$ is not
greater than $p$. Moreover, an atom of dimension less than $p$ is
a maximal atom in $d_{p}^{\gamma}\Lambda$ if and only if it is a
maximal atom in $\Lambda$; an atom $\uvws{}$ of dimension $p$ is a
maximal atom in $d_{p}^{\gamma}\Lambda$ if and only if there is a
maximal atom $\uvws{''}$ in $\Lambda$ such that $i''\geq i$,
$j''\geq j$ and $k''\geq k$, and the signs  $\alpha$, $\beta$ and
$\gamma$ satisfy the following conditions:

\begin{enumerate}
\item
if $\uvws{''}$ can be chosen such that $i''>i$, then
$\alpha=\gamma$; otherwise $\alpha=\alpha''$;

\item
if $\uvws{''}$ can be chosen such that $j''>j$, then
$\beta=(-)^{i}\gamma$; otherwise $\beta=\beta''$;

\item
if $\uvws{''}$ can be chosen such that $k''>k$, then
$\varepsilon=(-)^{i+j}\gamma$; otherwise
$\varepsilon=\varepsilon''$.

\end{enumerate}
\end{theorem}

\begin{note}
It follows easily from condition \ref{main4_3} in Theorem
\ref{main4} that $\alpha$, $\beta$ and $\gamma$ are well defined.
\end{note}

\begin{proof}
Evidently, the dimension of every maximal atom in
$d_{p}^{\gamma}\Lambda$ is not greater than $p$. Let $\Lambda_{1}$
be the union of the atoms as described in this theorem. It is easy
to see that $\Lambda_{1}$ satisfies condition 1 for pairwise
molecular subcomplexes. To prove the theorem, by Proposition
\ref{eq}, it suffices to prove that
$F_{I}^{u}(\Lambda_{1})=F_{I}^{u}(\dpg\Lambda)$,
$F_{J}^{v}(\Lambda_{1})=F_{J}^{v}(\dpg\Lambda)$ and
$F_{K}^{w}(\Lambda_{1})=F_{K}^{w}(\dpg\Lambda)$ for all $I$, $J$
and $K$. The arguments for the three equations are similar, we
prove only the second one.  If $J>p$, then it is easy to see that
$F_{J}^{v}(\Lambda_{1})=\emptyset=F_{J}^{v}(\dpg\Lambda)$, as
required. In the remaining proof, we may assume that $J\leq p$. We
have known that
$F_{J}^{v}(\dpg\Lambda)=d_{p-J}^{\gamma}F_{J}^{v}(\Lambda)$. we
need only to prove that
$F_{J}^{v}(\Lambda_{1})=d_{p-J}^{\gamma}F_{J}^{v}(\Lambda)$.

By the definition of $F_{J}^{v}$, it is easy to see that
$F_{J}^{v}(\Lambda_{1})$ and $d_{p-J}^{\gamma}F(\Lambda)$ are
subcomplexes of $u\times w^{J}$. We are going to prove that
$F_{J}^{v}(\Lambda_{1})$ and $d_{p-J}^{\gamma}F_{J}^{v}(\Lambda)$
consist of the same maximal atoms so that they are equal.

Let $\mu=\ui{}\times\wJk{}$ be a maximal atom in
$F_{J}^{v}(\Lambda_{1})$. Then $\Lambda_{1}$ has a
$(v,J)$-projection maximal atom $\lambda$  of the form
$\lambda=\uvws{}$. Hence $\Lambda$ has a maximal atom
$\lambda'=\uvws{'}$ with $i\leq i'$, $j\leq j'$ and $k\leq k'$.

Suppose that $j=J$ and $i+j+k=p$.  Since $\ui{'}\times\wJk{'}$ is
an atom in $F_{J}^{v}(\Lambda)$ and $i+k=p-J$, we know that
$d_{p-J}^{\gamma}F_{J}^{v}(\Lambda)$ has a maximal atom of the
form $u[i,\alpha'']\times w^{J}[k,\varepsilon'']$. Moreover, we
can see that there is a maximal atom $\lmns{}$ in $\Lambda$ such
that $l>i$, $m\geq j$ and $n\geq k$ if and only if there is a
maximal atom $\ul{}\times w^{J}[n,\omega]$ $F_{J}^{v}(\Lambda)$
such that $l>i$ and $n\geq k$; and we can also see that there is a
maximal atom $\lmns{}$ in $\Lambda$ such that $l\geq i$, $m\geq j$
and $n>k$ if and only if there is a maximal atom $\ul{}\times
w^{J}[n,\omega]$ in $F_{J}^{v}(\Lambda)$ such that $l\geq i$ and
$n>k$. It follows from \ref{2globe} that $\alpha=\alpha''$ and
$\varepsilon=\varepsilon''$. This implies that $\mu$ is a maximal
atom in $d_{p-J}^{\gamma}F(\Lambda)$.

Suppose that $j=J$ and $i+j+k<p$. Then $\lambda$ is also a maximal
atom in $\Lambda$. Therefore $\mu=F_{J}^{v}(\lambda)$ is a maximal
atom in $F_{J}^{v}(\Lambda)$. Since $i+k<p-J$, we know that $\mu$
is a maximal atom in $d_{p-J}^{\gamma}F(\Lambda)$.

There remains the case that $j>J$. In this case, there are no
maximal atom $\lmns{}$ in $\Lambda$ with $l\geq i$ and $m\geq J$
and $n\geq k$  such that $l>i$ or $n>k$. So $i=i'$,
$\alpha=\alpha'$, $k=k'$ and $\varepsilon=\varepsilon'$. On the
other hand, since
$\mu=\ui{}\times\wJk{}=\ui{'}\times\wJk{'}=F_{J}^{v}(\lambda')$,
we see that $\mu$ is a maximal atom in $F_{J}^{v}(\Lambda)$.
Because $i'+k'=i+k\leq p-j<p-J$, we know that $\mu$ is a maximal
atom in $d_{p-J}^{\gamma}F_{J}^{v}(\Lambda)$.

This shows  that every maximal atom in $F_{J}^{v}(\Lambda_{1})$ is
a maximal atom in $d_{p-J}^{\gamma}F_{J}^{v}(\Lambda)$.

Conversely, let $\mu=\ui{}\times\wJk{}$ be a maximal atom in
$d_{p-J}^{\gamma}F_{J}^{v}(\Lambda)$. Then $F_{J}^{v}(\Lambda)$
has a maximal atom $\mu'=\ui{'}\times\wJk{'}$ with $i\leq i'$ and
$k\leq k'$. Therefore $\Lambda$ has a $(v,J)$-projection maximal
atom of the form $\lambda'=\uvws{'}$.

Suppose that $i+k=p-J$. Then $\Lambda_{1}$ has a
$(v,J)$-projection maximal atom of the form
$\lambda=u[i,\alpha'']\times v[J,\beta'']\times
w[k,\varepsilon'']$. We can see that there is a maximal atom
$\lmns{}$ in $\Lambda$ such that $l>i$, $m\geq J$ and $n\geq k$ if
and only if there is a maximal atom $\ul{}\times w^{J}[n,\omega]$
$F_{J}^{v}(\Lambda)$ such that $l>i$ and $n\geq k$; and we can
also see that there is a maximal atom $\lmns{}$ in $\Lambda$ such
that $l\geq i$, $m\geq J$ and $n>k$ if and only if there is a
maximal atom $\ul{}\times w^{J}[n,\omega]$ in $F_{J}^{v}(\Lambda)$
such that $l\geq i$ and $n>k$. So $\alpha''=\alpha$ and
$\varepsilon''=\varepsilon$. Since
$F_{J}^{v}(\lambda)=u[i,\alpha'']\times
w^{J}[k,\varepsilon'']=\mu$, we can see that $\mu$ is a maximal
atom in $F_{J}^{v}(\Lambda_{1})$.

Suppose that $i+k<p-J$. Then $\mu=\ui{}\times\wJk{}$ is also a
maximal atom in $F_{J}^{v}(\Lambda)$. So $\Lambda$ has a
$(v,J)$-projection maximal atom $\lambda'=\ui{}\times\vj{'}\times
\wk{}$.  Now, if $j'=J$, then $i+j'+k<p$; hence $\lambda'$ is also
a maximal atom in $\Lambda_{1}$ and
$F_{J}^{v}(\lambda)=\ui{}\times w^{J}[k,\varepsilon]$ is a maximal
atom in $F_{J}^{v}(\Lambda_{1})$. Suppose that $j'>J$. Then it is
easy to see that there is no maximal atom $\lmns{}$ in $\Lambda$
with $l\geq i$ and $m\geq J$ and $n\geq k$  such that $l>i$ or
$n>k$. Hence $\Lambda_{1}$ has a $(v,J)$-projection maximal atom
of the form $\lambda''=\ui{}\times v[j'',\beta'']\times\wk{}$.
Therefore we see that $F_{J}^{v}(\lambda'')=\ui{}\times\wJk{}=\mu$
is a maximal atom in $F_{J}^{v}(\Lambda_{1})$.

This shows that every maximal atom in
$d_{p-J}^{\gamma}F_{J}^{v}(\Lambda)$ is a maximal atom in
$F_{J}^{v}(\Lambda_{1})$.

This completes the proof.

\end{proof}

\section{Composition of Pairwise Molecular Subcomplexes}

In this section, we consider composition of pairwise molecular
subcomplexes in $\uvwt$. We first give the construction of
composites of pairwise molecular subcomplexes. Then we show that
composites of pairwise molecular subcomplexes are pairwise
molecular.

\begin{lemma}\label{lessp}
Let $\Lambda^{-}$ and $\Lambda^{+}$ be pairwise molecular
subcomplexes. If $d_{p}^{+}\Lambda^{-}=d_{p}^{-}\Lambda^{+}$, then
for every maximal atom $\lambda^{-}=\uvws{^{-}}$ in $\Lambda^{-}$
and every maximal atom $\lambda^{+}=\uvws{^{+}}$ in $\Lambda^{+}$
one has
$\min\{i^{-},i^{+}\}+\min\{j^{-},j^{+}\}+\min\{k^{-},k^{+}\}\leq
p$.
\end{lemma}
\begin{proof} Let $l=\min\{i^{-},i^{+}\}$, $m=\min\{j^{-},j^{+}\}$ and
$n=\min\{k^{-},k^{+}\}$. Suppose otherwise that $l+m+n>p$. Then
there is an ordered triple $\{i,j,k\}$ with $i\leq l$, $j\leq m$,
$k\leq n$ and $i+j+k=p$. Since $l+m+n>p$, we have $i<l$, $j<m$ or
$k<m$. If $i<l$, then $d_{p}^{+}\Lambda^{-}$ has a maximal atom of
the form $u[i,+]\times v[j,\beta]\times w[k,\varepsilon]$, while
$d_{p}^{-}\Lambda^{+}$ has a maximal atom of the form
$u[i,-]\times v[j,\beta']\times w[k,\varepsilon']$ by Theorem
\ref{dt}. This contradicts condition 1 for pairwise molecular
subcomplex $d_{p}^{+}\Lambda^{-}=d_{p}^{-}\Lambda^{+}$. The
arguments for the cases $j<m$ and $k<m$ are similar.
\end{proof}

\begin{lemma}
Let $\Lambda^{-}$ and $\Lambda^{+}$ be pairwise molecular
subcomplexes in $u\times v\times w$. If
$d_{p}^{+}\Lambda^{-}=d_{p}^{-}\Lambda^{+}$, then
$$F_{I}^{u}(\Lambda^{-})\cap F_{I}^{u}(\Lambda^{+})
=F_{I}^{u}(\Lambda^{-}\cap\Lambda^{+})
=F_{I}^{u}(d_{p}^{+}\Lambda^{-})
=F_{I}^{u}(d_{p}^{-}\Lambda^{+}),$$ $$F_{J}^{v}(\Lambda^{-})\cap
F_{J}^{v}(\Lambda^{+}) =F_{J}^{v}(\Lambda^{-}\cap\Lambda^{+})
=F_{J}^{v}(d_{p}^{+}\Lambda^{-})
=F_{J}^{v}(d_{p}^{-}\Lambda^{+})$$ and
$$F_{K}^{w}(\Lambda^{-})\cap F_{K}^{w}(\Lambda^{+})
=F_{K}^{w}(\Lambda^{-}\cap\Lambda^{+})
=F_{K}^{w}(d_{p}^{+}\Lambda^{-})
=F_{K}^{w}(d_{p}^{-}\Lambda^{+})$$ for all $I$, $J$ and $K$.
\end{lemma}

\begin{proof}
The arguments  for the three formulae  are similar.  We give the
proof for the second one. There are two cases, as follows.

1. Suppose that $J>p$. We claim that $F_{J}^{v}(\Lambda^{-})\cap
F_{J}^{v}(\Lambda^{+})=\emptyset$.

Indeed, suppose otherwise that $F_{J}^{v}(\Lambda^{-})\cap
F_{J}^{v}(\Lambda^{+})\neq\emptyset$. Then it is evident that
there are atoms $\mu^{-}=\lmns{^{-}}$ in $\Lambda^{-}$ and
$\mu^{+}=\lmns{^{+}}$ in $\Lambda^{+}$ such that $m^{-}\geq J>p$
and $m^{+}\geq J>p$. According to Theorem \ref{dt}, this implies
that there are maximal atoms  $u[0,\alpha_{1}']\times v[p,
+]\times w[0,\varepsilon']$ and $u[0,\alpha_{1}'']\times v[p,
-]\times w[0,\varepsilon'']\times$ in $d_{p}^{+}\Lambda^{-}$ and
$d_{p}^{-}\Lambda^{+}$ respectively. This contradicts the
condition 1 for pairwise molecular subcomplex
$d_{p}^{+}\Lambda^{-}=d_{p}^{-}\Lambda^{+}$.

Now we have $F_{J}^{v}(d_{p}^{+}\Lambda^{-})\subset
F_{J}^{v}(\Lambda^{-}\cap\Lambda^{+})\subset
F_{J}^{v}(\Lambda^{-})\cap F_{J}^{v}(\Lambda^{+})=\emptyset$.
Therefore $F_{J}^{v}(d_{p}^{+}\Lambda^{-})=
F_{J}^{v}(\Lambda^{-}\cap \Lambda^{+})=F_{J}^{v}(\Lambda^{-})\cap
F_{J}^{v}(\Lambda^{+})$, as required.

2. Suppose that $J\leq p$. Since
$d_{p}^{+}\Lambda^{-}=d_{p}^{-}\Lambda^{+}$, we have
$d_{p-J}^{+}F_{J}^{v}(\Lambda^{-})
=F_{J}^{v}(d_{p}^{+}\Lambda^{-}) =F_{J}^{v}(d_{p}^{-}\Lambda^{+})
=d_{p-J}^{-}F_{J}^{v}(\Lambda^{+}).
$
Because $F_{J}^{v}(\Lambda^{-})$ and $F_{J}^{v}(\Lambda^{+})$ are
molecules, we can see that $F_{J}^{v}(\Lambda^{-})\#_{p-J}
F_{J}^{v}(\Lambda^{+})$ is defined. Hence
$
 F_{J}^{v}(\Lambda^{-})\cap F_{J}^{v}(\Lambda^{+})
=d_{p-J}^{+}F_{J}^{v}(\Lambda^{-})=F_{J}^{v}(d_{p}^{+}\Lambda^{-})
\subset F_{J}^{v}(\Lambda^{-}\cap\Lambda^{+})$. Since we
automatically have $F_{J}^{v}(\Lambda^{-}\cap\Lambda^{+})\subset
F_{J}^{v}(\Lambda^{-})\cap F_{J}^{v}(\Lambda^{+})$, we get
$F_{J}^{v}(\Lambda^{-})\cap
F_{J}^{v}(\Lambda^{+})=F_{J}^{v}(\Lambda^{-}\cap\Lambda^{+})
=F_{J}^{v}(d_{p}^{+}\Lambda^{-})$, as required.

This completes the proof.

\end{proof}

\begin{prop}\label{dfd}
Let $\Lambda^{-}$ and $\Lambda^{+}$ be pairwise molecular
subcomplexes. If $d_{p}^{+}\Lambda^{-}=d_{p}^{-}\Lambda^{+}$, then
$\Lambda^{-}\cap\Lambda^{+}=
d_{p}^{+}\Lambda^{-}(=d_{p}^{-}\Lambda^{+})$; hence
$\Lambda^{-}\#_{p}\Lambda^{+}$ is defined.
\end{prop}

\begin{proof}
Let $M=d_{p}^{+}\Lambda^{-}=d_{p}^{-}\Lambda^{+}$. It is evident
that $M\subset\Lambda^{-}\cap\Lambda^{+}$. To prove the reverse
inclusion, it suffices to prove that every maximal atom in
$\Lambda^{-}\cap\Lambda^{+}$ is contained in $M$.

Suppose otherwise that there is a maximal atom $\lambda=\uvws{}$
in $\Lambda^{-}\cap\Lambda^{+}$ such that $\lambda\not\subset M$.
Since $\ui{}\times\vj{}=F_{k}^{w}(\lambda)\subset
F_{k}^{w}(\Lambda^{-}\cap\Lambda^{+})=F_{k}^{w}(M)$, we can see
that $M$ has a maximal atom $\lambda'=\uvws{'}$ such that
$\ui{}\subset\ui{'}$ and $\vj{}\subset\vj{'}$ and $k'\geq k$.
Because $\lambda=\uvws{}$ is maximal in
$\Lambda^{-}\cap\Lambda^{+}$ and $M\subset
\Lambda^{-}\cap\Lambda^{+}$, we have $k'=k$ and
$\varepsilon'=-\varepsilon$. Now we know that
$\lambda\cup\lambda'\subset\Lambda^{-}$ and
$\lambda\cup\lambda'\subset\Lambda^{+}$.  By applying Lemma
\ref{lm3} to $\Lambda^{-}$ and $\Lambda^{+}$, it is easy to see
that there are maximal atoms $\lambda^{-}=\uvws{^{-}}$ in
$\Lambda^{-}$ and $\lambda^{+}=\uvws{^{+}}$ in $\Lambda^{+}$ such
that $\ui{^{-}}\cap\ui{^{+}}\supset \ui{}$ and
$\vj{^{-}}\cap\vj{^{+}}\supset \vj{}$  and
$\min\{k^{-},k^{+}\}>k$. Since $\lambda$ is maximal atom in
$\Lambda^{-}\cap\Lambda^{+}$, we have $k^{-}=k^{+}=k+1$ and
$\varepsilon^{-}=-\varepsilon^{+}$.

Now, we have $\ui{}\times\vj{}\subset F_{k+1}^{w}(\Lambda^{-})\cap
F_{k+1}^{w}(\Lambda^{+})
=F_{k+1}^{w}(\Lambda^{-}\cap\Lambda^{+})$. Therefore
$\Lambda^{-}\cap\Lambda^{+}$ has a maximal atom
$\lambda''=\uvws{''}$ with $\ui{''}\supset\ui{}$ and
$\vj{''}\supset\vj{}$ and $k''>k$. This contradicts that $\lambda$
is a maximal atom in $\Lambda^{-}\cap\Lambda^{+}$.

This completes the proof.

\end{proof}

The following Proposition tells us how to construct the composite
of a pair of pairwise molecular subcomplexes of $\uvwt$.

\begin{prop}\label{comp}
Let $\Lambda^{-}$ and $\Lambda^{+}$ be pairwise molecular
subcomplexes of $\uvwt$. If
$d_{p}^{+}\Lambda^{-}=d_{p}^{-}\Lambda^{+}$, then the maximal
atoms in the composite $\Lambda^{-}\#_{p}\Lambda^{+}$ are the
$q$-dimensional common maximal atoms of $\Lambda^{-}$ and
$\Lambda^{+}$ with $q\leq p$ and  the $r$-dimensional atoms in
either $\Lambda^{-}$ and $\Lambda^{+}$ with $r>p$.
\end{prop}
\begin{proof} Let $\Lambda$ be the subcomplex of $\uvwt$ as described in the
proposition. We must prove that
$\Lambda=\Lambda^{-}\cup\Lambda^{+}$. Clearly, we have
$\Lambda\subset\Lambda^{-}\cup\Lambda^{+}$; it suffices to prove
that $\Lambda^{-}\cup\Lambda^{+}\subset\Lambda$. By the formation
of $\Lambda$, we must prove that, for each maximal atom
$\lambda=\uvws{}$ in either $\Lambda^{-}$ or $\Lambda^{+}$ with
$i+j+k\leq p$ and such that $\lambda$ is not a common maximal atom
in $\Lambda^{-}$ and $\Lambda^{+}$, $\lambda\subset \Lambda$. It
is easy to see that this  can only happen when $i+j+k=p$. Suppose
that $\lambda$ is a maximal atom in $\Lambda^{\gamma}$ which is
not a maximal atom in $\Lambda^{-\gamma}$. Then $\lambda$ must be
a maximal atom  in $d_{p}^{+}\Lambda^{-}=d_{p}^{-}\Lambda^{+}$
which implies that $\lambda\subset \lambda^{-\gamma}$ for some
maximal atom $\lambda^{-\gamma}=\uvws{^{-\gamma}}$ with
$i^{-\gamma}+j^{-\gamma}+k^{-\gamma}>p$. Thus $\lambda\subset
\Lambda$. Therefore, we have
$\Lambda^{-}\cup\Lambda^{+}\subset\Lambda$.

This completes the proof.

\end{proof}

Now we can show that the composites of pairwise molecular
subcomplexes in $\uvwt$ are pairwise molecular.

\begin{prop}
Let $\Lambda^{-}$ and $\Lambda^{+}$ be pairwise molecular
subcomplexes. If $d_{p}^{+}\Lambda^{-}=d_{p}^{-}\Lambda^{+}$, then
$\Lambda^{-}\#_{p}\Lambda^{+}$ is a pairwise molecular subcomplex
of $u\times v\times w$.
\end{prop}
\begin{proof} Let $\Lambda=\Lambda^{-}\#_{p}\Lambda^{+}$. According to Lemma
\ref{lessp}, it is easy to see that $\Lambda$ satisfies condition
1 for pairwise molecular subcomplexes. Moreover, we have
$F_{I}^{u}(\Lambda^{-}\#_{p}\Lambda^{+})
=F_{I}^{u}(\Lambda^{-}\cup\Lambda^{+}) =F_{I}^{u}(\Lambda^{-})\cup
F_{I}^{u}(\Lambda^{+})$.

Now suppose that $p\geq I$. We  have
$d_{p-I}^{+}F_{I}^{u}(\Lambda^{-})
=F_{I}^{u}(d_{p}^{+}\Lambda^{-}) =F_{I}^{u}(d_{p}^{-}\Lambda^{+})
=d_{p-I}^{-}F_{I}^{u}(\Lambda^{+})$. Thus
$F_{I}^{u}(\Lambda^{-}\#_{p}\Lambda^{+})
=F_{I}^{u}(\Lambda^{-})\#_{p-I} F_{I}^{u}(\Lambda^{+})$. Therefore
$F_{I}^{u}(\Lambda^{-}\#_{p}\Lambda^{+})$ is a molecule.

Suppose that $p<I$. Then it is easy to see that
$F_{I}^{u}(\Lambda^{-})=\emptyset$ or
$F_{I}^{u}(\Lambda^{+})=\emptyset$. (Otherwise, we have
$F_{I}^{u}(\Lambda^{-}\cap\Lambda^{+})\neq\emptyset$. This would
lead to a contradiction to Lemma \ref{lessp}.) Therefore
$F_{I}^{u}(\Lambda^{-}\#_{p}\Lambda^{+})$ is a molecule or the
empty set.

We have now proved that $F_{I}^{u}(\Lambda^{-}\#_{p}\Lambda^{+})$
is a molecule or the empty set for all $I$.

Similarly, we can see that
$F_{J}^{v}(\Lambda^{-}\#_{p}\Lambda^{+})$ and
$F_{K}^{w}(\Lambda^{-}\#_{p}\Lambda^{+})$ are molecules or the
empty set for all $J$ and $K$.

It follows from Definition \ref{pairwise_def} that $\Lambda$ is a
pairwise molecular subcomplex of $u\times v\times w$.

\end{proof}

\section{Decomposition of Pairwise Molecular Subcomplexes}

The aim of this section is to prove the main theorem in this
chapter.

\begin{theorem} \label{pairwise_mole}
If  $\Lambda$ is a pairwise molecular subcomplex of $\uvwt$, then
$\Lambda$ is a molecule.
\end{theorem}

It is trivial that the theorem holds when $\Lambda$ is an atom.
Thus we may assume that $\Lambda$ is a pairwise molecular
subcomplex in $\uvwt$ which is not an atom throughout this
section. We are going to show that $\Lambda$ is a molecule.

Let
\begin{center}
$p=\max\{\dim(\lambda\cap\mu)$:\,\,$\lambda$ and $\mu$ are
distinct maximal atoms in $\Lambda\}$.
\end{center}
Recall that $p$ is called frame dimension of $\Lambda$. It is
evident that there are at least two maximal atoms $\lambda$ and
$\mu$ in $\Lambda$ with $\dim\lambda>p$ and $\dim\mu>p$. By Lemma
\ref{lm3}, it is easy to see that $p$ is the maximal number among
the numbers
$\min\{i_{1},i_{2}\}+\min\{j_{1},j_{2}\}+\min\{k_{1},k_{2}\}$,
where $\lambda_{1}=\uvws{_{1}}$ and $\lambda_{2}=\uvws{_{2}}$ run
over all pairs of distinct maximal atoms in $\Lambda$.

\begin{lemma}\label{n_order_p}
Let $\lambda=\uvws{}$ and $\lambda'=\uvws{'}$ are maximal atoms in
$\Lambda$ with $\min\{i,i'\}+\min\{j,j'\}+\min\{k,k'\}=p$.
\begin{enumerate}
\item
If $i=i'$, $\alpha=-\alpha'$ and $j<j'$,  then
$\beta=(-)^{i}\alpha$;
\item
If $j=j'$, $\beta=-\beta'$ and $k<k'$, then
$\varepsilon=(-)^{j}\beta$;
\item
If $k=k'$, $\varepsilon=-\varepsilon'$ and $j<j'$, then
$\varepsilon=(-)^{j}\beta$;

\end{enumerate}

\end{lemma}
\begin{proof}
The arguments for the three cases are similar, we prove only for
the first case.

Suppose that $i=i'$, $\alpha=-\alpha'$ and $j<j'$. According to
Lemma \ref{lm3}, we can get a maximal atom $\lambda''=\uvws{''}$
with $i''>i$, $\vj{''}\supset\vj{}$ and $\wk{''}\supset\wk{'}$.
Since $\min\{i,i'\}+\min\{j,j'\}+\min\{k,k'\}=p$, we have $j''=j$
and $k''=k'$. Hence $\vj{''}=\vj{}$ and $\wk{''}=\wk{'}$.
Moreover, it is easy to see that $\lambda$, $\lambda'$ and
$\lambda''$ are pairwise adjacent by the choice of $p$. It follows
easily from the sign conditions that $\beta=(-)^{i}\alpha$, as
required.

This completes the proof.

\end{proof}

We are going to prove that a pairwise molecular subcomplex
$\Lambda$ in $\uvwt$ is a molecule by showing that $\Lambda$ can
be properly decomposed into pairwise molecular subcomplexes. This
decomposition depends essentially on the following total order on
the set of  maximal atoms in $\Lambda$.

For a pair of atoms $\lambda=\uvws{}$ and $\lambda'=\uvws{'}$ in
$\Lambda$, we write $\lambda<\lambda'$ if one of the following
holds:
\begin{itemize}
\item
$\alpha=\alpha'=-$ and $i<i'$;
\item
$\alpha=\alpha'=+$ and $i>i'$;
\item
$\alpha=-$ and $\alpha'=+$;
\item
$i=i'$ are even, $\alpha=\alpha'$, $\beta=\beta'=-$ and $j<j'$;
\item
$i=i'$ are even, $\alpha=\alpha'$, $\beta=\beta'=+$ and $j'<j$;
\item
$i=i'$ are even, $\alpha=\alpha'$, $\beta=-$ and $\beta'=+$
\item
$i=i'$ are odd, $\alpha=\alpha'$, $\beta=\beta'=+$ and $j<j'$;
\item
$i=i'$ are odd, $\alpha=\alpha'$, $\beta=\beta'=-$ and $j'<j$;
\item
$i=i'$ are odd, $\alpha=\alpha'$, $\beta=+$ and $\beta'=-$.
\end{itemize}
It is evident that the relation $<$ is a total order on the set of
maximal atoms in $\Lambda$.

\begin{lemma}
For any pair of maximal atoms $\lambda$ and $\lambda'$ in
$\Lambda$ with $\dim\lambda>p$ and $\dim\lambda'>p$, if
$\lambda<\lambda'$, then $\lambda\cap\lambda'\subset
d_{p}^{+}\lambda\cap d_{p}^{-}\lambda'$.
\end{lemma}

\begin{proof}
Let $\lambda=\uvws{}$ and $\lambda'=\uvws{'}$. According to the
choice of $p$, it is evident that
$\min\{i,i'\}+\min\{j,j'\}+\min\{k,k'\}\leq p$. We now consider
five cases, as follows.

1. Suppose that $\min\{i,i'\}+\min\{j,j'\}+\min\{k,k'\}=p$. Then
$\lambda$ and $\lambda'$ are adjacent by the choice of $p$.
According to Lemma \ref{n_order_p} and sign conditions for
pairwise molecular subcomplexes, it is easy to see that
$\lambda\cap\lambda'\subset d_{p}^{+}\lambda\cap
d_{p}^{-}\lambda'$, as required.

2. Suppose that $\min\{i,i'\}+\min\{j,j'\}+\min\{k,k'\}<p-1$. Then
it is easy to see that $\lambda\cap\lambda'\subset
d_{p}^{+}\lambda\cap d_{p}^{-}\lambda'$, as required.

3. Suppose that $\min\{i,i'\}+\min\{j,j'\}+\min\{k,k'\}=p-1$ and
that $\lambda$ and $\lambda'$ are adjacent. There are two case, as
follows: (1) $i=i'$; (2) $i\neq i'$ . In  case (1), it is evident
that $\lambda\cap\lambda'\subset d_{p}^{+}\lambda\cap
d_{p}^{-}\lambda'$, as required. In case (2), it follows easily
from the sign conditions that $\lambda\cap\lambda'\subset
d_{p-1}^{+}\lambda\cap d_{p-1}^{-}\lambda'$; thus
$\lambda\cap\lambda'\subset d_{p}^{+}\lambda\cap
d_{p}^{-}\lambda'$, as required.

4. Suppose that $\min\{i,i'\}+\min\{j,j'\}+\min\{k,k'\}=p-1$ and
that $\lambda$ and $\lambda'$ are not adjacent.  Suppose also that
$i=i'$ or $j=j'$ or $k=k'$. Then it is easy to see that
$\lambda\cap\lambda'\subset d_{p}^{+}\lambda\cap
d_{p}^{-}\lambda'$, as required.

5. Suppose that $\min\{i,i'\}+\min\{j,j'\}+\min\{k,k'\}=p-1$ and
that $\lambda$ and $\lambda'$ are not adjacent. Suppose also that
$i\neq i'$ and $j\neq j'$ and $k\neq k'$. Then there are several
cases, as follows. (1) $i<i'$ and $j<j'$, or $i<i'$ and $k<k'$;
(2) $i<i'$ and $j>j'$ and $k>k'$; (3)  $i>i'$ and  $j>j'$,  or
$i>i'$ and $k>k'$; (4) $i>i'$ and $j<j'$ and $k<k'$. In case (1),
we have $\alpha=-$; it follows easily that
$\lambda\cap\lambda'\subset d_{p}^{+}\lambda\cap
d_{p}^{-}\lambda'$, as required. Similarly, in case (3), we have
$\alpha'=+$; this also implies that $\lambda\cap\lambda'\subset
d_{p}^{+}\lambda\cap d_{p}^{-}\lambda'$, as required. There remain
case (2) and case (4).

To give the proof for case (2), suppose that
$\min\{i,i'\}+\min\{j,j'\}+\min\{k,k'\}=p-1$ and that $\lambda$
and $\lambda'$ are not adjacent; suppose also that $i<i'$ and
$j>j'$ and $k>k'$. Then $\alpha=-$ and there is a maximal atom
$\lambda''=\uvws{''}$ in $\Lambda$ distinct from $\lambda'$ such
that $i''>i$, $j''\geq j'$ and $k''\geq k'$. By the choice of $p$,
we can see that $\lambda''$ is adjacent to both $\lambda$ and
$\lambda'$, and we have $i''=i+1$. According to condition 1 for
pairwise molecular subcomplexes, we have $j''>j'$ or $k''>k'$.

In case (2), suppose that $j''>j$. Then $\min\{j'',j\}=j'+1$ and
$k''=k'$ by the choice of $p$. Hence
$\varepsilon''=-[-(-)^{i+j'+1}]=-(-)^{i+j'}$. If
$\varepsilon'=\varepsilon''=-(-)^{i+j'}$, then it is easy to see
that $\lambda\cap\lambda'\subset d_{p}^{+}\lambda\cap
d_{p}^{-}\lambda'$, as required. If
$\varepsilon'=-\varepsilon''=(-)^{i+j'}$, then we can get
$\varepsilon'=(-)^{j'}\beta'$, i.e., $(-)^{j'}\beta'=(-)^{i+j'}$;
thus $\beta'=(-)^{i}$; this implies that
$\lambda\cap\lambda'\subset d_{p}^{+}\lambda\cap
d_{p}^{-}\lambda'$, as required.

In case (2), suppose that $k''>k$. Then $j''=j'$ by the choice of
$p$. We can also have $\beta''=-(-)^{i}\alpha=(-)^{i}$ by the sign
conditions. If $\beta'=\beta''=(-)^{i}$, then it is easy to see
that $\lambda\cap\lambda'\subset d_{p}^{+}\lambda\cap
d_{p}^{-}\lambda'$, as required. If $\beta'=-\beta''=-(-)^{i}$,
then we can get $\varepsilon'=(-)^{j'}\beta'=-(-)^{i+j'}$; this
implies that $\lambda\cap\lambda'\subset d_{p}^{+}\lambda\cap
d_{p}^{-}\lambda'$, as required.

This completes the proof for case (2).

To give the proof for case (4), suppose that
$\min\{i,i'\}+\min\{j,j'\}+\min\{k,k'\}=p-1$ and that $\lambda$
and $\lambda'$ are not adjacent; suppose also that $i>i'$ and
$j<j'$ and $k<k'$. Then $\alpha=\alpha'=+$ and there is a maximal
atom $\lambda''=\uvws{''}$ in $\Lambda$ distinct from $\lambda'$
such that $i''>i'$, $j''\geq j$ and $k''\geq k$. By the choice of
$p$, we can see that $\lambda''$ is adjacent to both $\lambda$ and
$\lambda'$, and we have $i''=i'+1$. According to condition 1 for
pairwise molecular subcomplexes, we have $j''>j$ or $k''>k$.

In case (4),  suppose that $j''>j$. Then $\min\{j'',j'\}=j+1$ and
$k''=k$ by the choice of $p$. Hence
$\varepsilon''=[-(-)^{i'+j+1}]=(-)^{i'+j}$. If
$\varepsilon=\varepsilon''=(-)^{i'+j}$, then it is easy to see
that $\lambda\cap\lambda'\subset d_{p}^{+}\lambda\cap
d_{p}^{-}\lambda'$, as required. If
$\varepsilon=-\varepsilon''=-(-)^{i'+j}$, then we can get
$\varepsilon=(-)^{j}\beta$, i.e., $(-)^{j}\beta=-(-)^{i'+j}$; thus
$\beta'=-(-)^{i'}$; this implies that $\lambda\cap\lambda'\subset
d_{p}^{+}\lambda\cap d_{p}^{-}\lambda'$, as required.

In case (4), suppose that $k''>k$. Then $j''=j$ by the choice of
$p$. We can also have $\beta''=-(-)^{i'}\alpha=-(-)^{i'}$ by the
sign conditions. If $\beta=\beta''=-(-)^{i'}$, then it is easy to
see that $\lambda\cap\lambda'\subset d_{p}^{+}\lambda\cap
d_{p}^{-}\lambda'$, as required. If $\beta=-\beta''=(-)^{i'}$,
then we can get $\varepsilon=(-)^{j}\beta=(-)^{i'+j}$; this
implies that $\lambda\cap\lambda'\subset d_{p}^{+}\lambda\cap
d_{p}^{-}\lambda'$, as required.

This completes the proof for case (4), thus completes the proof of
the lemma.

\end{proof}

By this lemma, we can arrange all the maximal atoms in $\Lambda$
with dimension greater than $p$ as $$\lambda_{1},
\lambda_{2},\cdots,\lambda_{n}$$ such $\lambda_{i}\cap
\lambda_{j}\subset d_{p}^{+}\lambda_{i}\cap d_{p}^{-}\lambda_{j}$
for $i<j$.

Let $\Lambda^{-}=d_{p}^{-}\Lambda\cup\lambda_{1}$ and
$\Lambda^{+}=d_{p}^{+}\Lambda\cup\lambda_{2}\cdots\lambda_{n}$. We
are going to prove that $\Lambda^{-}$ and $\Lambda^{+}$ are
pairwise molecular subcomplexes and $\Lambda$ can be decomposed
into $\Lambda^{-}$ and $\Lambda^{+}$.

\begin{lemma} \label{L-c1}
$\Lambda^{-}$ satisfies condition 1 for pairwise molecular
subcomplexes.
\end{lemma}
\begin{proof}
We first prove that $d_{p}^{-}\lambda_{1}\subset
d_{p}^{-}\Lambda$. Suppose that $\xi\in d_{p}^{-}\lambda_{1}$.
Then, for every maximal atom $\lambda'$ in $\Lambda$ with
$\xi\in\lambda'$, if $\lambda'=\lambda_{t}$ for some $t>1$, then
$\xi\in \lambda_{1}\cap\lambda_{t}\subset
d_{p}^{-}\lambda_{t}=d_{p}^{-}\lambda'$; if $\dim\lambda'\leq p$,
then we automatically have $\xi\in d_{p}^{-}\lambda'$.  It follows
from Lemma \ref{xi_dpg} that $d_{p}^{-}\lambda_{1}\subset
d_{p}^{-}\Lambda$, as required.

We now verify that $\Lambda^{-}$ satisfies condition 1 for
pairwise molecular subcomplexes. It suffices to  prove that any
maximal atom $\lambda=\uvws{}$ in $d_{p}^{-}\Lambda$ with $i\leq
i_{1}$, $j\leq j_{1}$ and $k\leq k_{1}$ is contained in
$\lambda_{1}$.  By the formation of $d_{p}^{-}\lambda_{1}$ and
$d_{p}^{-}\Lambda$, it is easy to see that $\lambda$ is a maximal
atom in $d_{p}^{-}\lambda_{1}$, and hence
$\lambda\subset\lambda_{1}$, as required.
\end{proof}

\begin{lemma}\label{L+c1}
$\Lambda^{+}$ satisfies condition 1 for pairwise molecular
subcomplexes.
\end{lemma}
\begin{proof} It suffices to prove that any maximal atom $\lambda=\uvws{}$ in
$d_{p}^{+}\Lambda$ with $i\leq i_{t}$, $j\leq j_{t}$ and $k\leq
k_{t}$ for some $2\leq t\leq n$ is contained in some $\lambda_{s}$
for $2\leq s\leq n$. It is evident that $i+j+k=p$.

Let $r$ be the maximal integer between $2$ and $n$ such that
$i\leq i_{r}$, $j\leq j_{r}$ and $k\leq k_{r}$. Then
$d_{p}^{+}\lambda_{r}$ has a maximal atom
$\lambda'=u[i,\alpha']\times v[j,\beta']\times w[k,\varepsilon']$.
By the choice of $r$, it is evident that
$\Int\lambda'\cap\lambda_{t}=\emptyset$ for any $t>r$.  Moreover,
for any $1\leq s<r$, we have
$\lambda'\cap\lambda_{s}\subset\lambda_{r}\cap\lambda_{s}\subset
d_{p}^{+}\lambda_{s}$. By Lemma \ref{xi_dpg}, it is easy to see
that $\Int\lambda'\subset d_{p}^{+}\Lambda$ and hence
$\lambda'\subset d_{p}^{+}\Lambda$. So, by condition 1 for the
pairwise molecular subcomplex $d_{p}^{+}\Lambda$, we can see that
$\lambda=\lambda'\subset\lambda_{r}$, as required.

This completes the proof.
\end{proof}

\begin{lemma}\label{Fiu}
Let $p\geq I$ and let $\lambda_{1}$ be  $(u,I)$-projection
maximal. Then
\begin{enumerate}
\item
$F_{I}^{u}(\Lambda^{-})$ and $F_{I}^{u}(\Lambda^{+})$ are
molecules in $v^{I}\times w^{I}$.
\item
$d_{p-I}^{+}F_{I}^{u}(\Lambda^{-})=d_{p-I}^{-}F_{I}^{u}(\Lambda^{+})$,
hence $F_{I}^{u}(\Lambda^{-})\#_{p-I}F_{I}^{u}(\Lambda^{+})$ is
defined.
\item
$F_{I}^{u}(\Lambda)=F_{I}^{u}(\Lambda^{-})\#_{p-I}F_{I}^{u}(\Lambda^{+})$.
\end{enumerate}
\end{lemma}

\begin{proof}
Since $F_{I}^{u}$ preserves unions, we have
$F_{I}^{u}(\Lambda^{-})
=F_{I}^{u}(d_{p}^{-}\Lambda\cup\lambda_{1})
=F_{I}^{u}(d_{p}^{-}\Lambda)\cup F_{I}^{u}(\lambda_{1})$ and
$F_{I}^{u}(\Lambda^{+})
=F_{I}^{u}(d_{p}^{+}\Lambda\cup\lambda_{2}\cup\cdots\cup\lambda_{n})
=F_{I}^{u}(d_{p}^{+}\Lambda)\cup
F_{I}^{u}(\lambda_{2})\cup\cdots\cup F_{I}^{u}(\lambda_{n})$.  If
$\dim F_{I}^{u}(\lambda_{1})=j_{1}+k_{1}\leq p-I$, then
$F_{I}^{u}(\lambda_{1})$ is a maximal atom in
$d_{p-I}^{-}F_{I}^{u}(\Lambda)$ by Theorem \ref{maximal_dpg};
hence $F_{I}^{u}(\Lambda^{-})=F_{I}^{u}(d_{p}^{-}\Lambda)
=d_{p-I}^{-}F_{I}^{u}(\Lambda)$ and similarly
$F_{I}^{u}(\Lambda^{+})=F_{I}^{u}(\Lambda)$; it follows easily
that $F_{I}^{u}(\Lambda^{-})$ and $F_{I}^{u}(\Lambda^{+})$ are
molecules and
$d_{p-I}^{+}F_{I}^{u}(\Lambda^{-})=d_{p-I}^{-}F_{I}^{u}(\Lambda^{+})$,
as required. If $F_{I}^{u}(\Lambda)=F_{I}^{u}(\lambda_{1})$, then
$\lambda_{s}$ are not $(u,I)$-projection maximal  for
$s=2,\dots,n$; thus
$F_{I}^{u}(\lambda_{s})=F_{I}^{u}(\lambda_{1}\cap\lambda_{s})\subset
F_{I}^{u}(d_{p}^{+}\lambda_{1})=d_{p-I}^{+}F_{I}^{u}(\lambda_{1})=
d_{p-I}^{+}F_{I}^{u}(\Lambda)$; it follows easily that
$F_{I}^{u}(\Lambda^{+})=F_{I}^{u}(d_{p}^{+}\Lambda)
=d_{p-I}^{+}F_{I}^{u}(\Lambda)$; it is also evident that
$F_{I}^{u}(\Lambda^{-})=F_{I}^{u}(\lambda_{1})=F_{I}^{u}(\Lambda)$;
therefore $F_{I}^{u}(\Lambda^{-})$ and $F_{I}^{u}(\Lambda^{+})$
are molecules and
$d_{p-I}^{+}F_{I}^{u}(\Lambda^{-})=d_{p-I}^{-}F_{I}^{u}(\Lambda^{+})$,
as required. In the following proof, we may assume that $\dim
F_{I}^{u}(\lambda_{1})>p-I$ and $F_{I}^{u}(\Lambda)$ has at least
two distinct maximal atoms.

Let $$q=\max\{\dim(\mu\cap\mu'):\text{$\mu$ and $\mu'$ are
distinct maximal atoms in $F_{I}^{u}(\Lambda)$}\}.$$ It is clear
that $q\leq p-I$ by the choice of $p$. Let
$\mu=v^{I}[m,\tau]\times w^{I}[n,\omega]$ be a maximal atom in
$F_{I}^{u}(\Lambda)$ distinct from $F_{I}^{u}(\lambda_{1})$. If
$\dim (F_{I}^{u}(\lambda_{1})\cap\mu)<p-I$, then it is easy to see
that $F_{I}^{u}(\lambda_{1})\cap\mu\subset
d_{p-I}^{+}F_{I}^{u}(\lambda_{1})\cap d_{p-I}^{-}\mu$ by the
construction of molecule $F_{I}^{u}(\Lambda)$ in $v^{I}\times
w^{I}$ (Theorem \ref{2globe}). Suppose that $\dim
(F_{I}^{u}(\lambda_{1})\cap\mu)=p-I$. Then there is a maximal atom
$\lambda'=\uvws{'}$ in $\Lambda$ such that
$F_{I}^{u}(\lambda')=\mu$. If $i_{1}\leq i'$ and $j_{1}<j'$, then
$i_{1}=I$ by the choice of $p$, $\alpha_{1}=-$ when $i_{1}<i'$ by
the definition of natural order and $k'<k_{1}$ by condition 1 for
pairwise molecular subcomplexes; hence $\beta_{1}=-(-)^{I}$ by the
sign conditions for $\lambda_{1}$ and $\lambda'$ or by the
definition of natural order; thus
$\omega=\varepsilon'=(-)^{I+j_{1}}$ by the sign condition for
$\lambda_{1}$ and $\lambda'$; it follows easily that
$F_{I}^{u}(\lambda_{1})\cap\mu \subset
d_{p-I}^{+}F_{I}^{u}(\lambda_{1})\cap d_{p-I}^{-}\mu$.  If
$i_{1}\leq i'$ and $j_{1}>j'$, then it is easy to see that
$i_{1}=I$, $\tau=\beta'=(-)^{I}$ and $\varepsilon_{1}=-(-)^{I+m}$
by the sign condition for pairwise molecular subcomplexes  or the
definition of the natural order; it follows easily that
$F_{I}^{u}(\lambda_{1})\cap\mu \subset
d_{p-I}^{+}F_{I}^{u}(\lambda_{1})\cap d_{p-I}^{-}\mu$.  If
$i_{1}>i'$, then, by an similar argument, one can get
$F_{I}^{u}(\lambda_{1})\cap \mu\subset
d_{p-I}^{+}F_{I}^{u}(\lambda_{1})\cap d_{p-I}^{-}\mu$.  We have
now shown that $F_{I}^{u}(\lambda_{1})\cap \mu\subset
d_{p-I}^{+}F_{I}^{u}(\lambda_{1})\cap d_{p-I}^{-}\mu$ for every
maximal atom $\mu$ in $F_{I}^{u}(\Lambda)$.

Moreover, we have
$F_{I}^{u}(\Lambda^{-})=d_{p-I}^{-}F_{I}^{u}(\Lambda)\cup
F_{I}^{u}(\lambda_{1})$ and
$$F_{I}^{u}(\Lambda^{+})=d_{p-I}^{+}F_{I}^{u}(\Lambda)\cup\bigcup\{\mu:
\mu\text{ is a maximal atom in }F_{I}^{u}(\Lambda)\text{ with }
\mu\neq F_{I}^{u}(\lambda_{1})\}$$ (Notice that it is possible
that  $F_{I}^{u}(\Lambda^{+})=d_{p-I}^{+}F_{I}^{u}(\Lambda)$). It
follows from Theorem \ref{decom_mole} that
$F_{I}^{u}(\Lambda^{-})$ and $F_{I}^{u}(\Lambda^{+})$ are
molecules in $v^{I}\times w^{I}$,
$d_{p-I}^{+}F_{I}^{u}(\Lambda^{-})=
d_{p-I}^{-}F_{I}^{u}(\Lambda^{+})$ and
$F_{I}^{u}(\Lambda)=F_{I}^{u}(\Lambda^{-})\#_{p-I}F_{I}^{u}(\Lambda^{+})$,
as required.

This completes the proof.

\end{proof}

\begin{lemma}\label{FJv}
Let $p\geq J$ and let $\lambda_{1}$ be a $(v,J)$-projection
maximal atom. Then
\begin{enumerate}
\item
$F_{J}^{v}(\Lambda^{-})$ and $F_{J}^{v}(\Lambda^{+})$ are
molecules in $u\times w^{J}$.
\item
$d_{p-J}^{+}F_{J}^{v}(\Lambda^{-})=d_{p-J}^{-}F_{J}^{v}(\Lambda^{+})$,
hence $F_{J}^{v}(\Lambda^{-})\#_{p-J}F_{J}^{v}(\Lambda^{+})$ is
defined.
\item
$F_{J}^{v}(\Lambda)=F_{J}^{v}(\Lambda^{-})\#_{p-J}F_{J}^{v}(\Lambda^{+})$.
\end{enumerate}
\end{lemma}
\begin{proof}
Since $F_{J}^{v}$ preserves unions, we have
$F_{J}^{v}(\Lambda^{-})
=F_{J}^{v}(d_{p}^{-}\Lambda\cup\lambda_{1})
=F_{J}^{v}(d_{p}^{-}\Lambda)\cup F_{J}^{v}(\lambda_{1})$ and
$F_{J}^{v}(\Lambda^{+})
=F_{J}^{v}(d_{p}^{+}\Lambda\cup\lambda_{2}\cup\cdots\cup\lambda_{n})
=F_{J}^{v}(d_{p}^{+}\Lambda)\cup
F_{J}^{v}(\lambda_{2})\cup\cdots\cup F_{J}^{v}(\lambda_{n})$. If
$\dim F_{J}^{v}(\lambda_{1})=i_{1}+k_{1}\leq p-J$, then it is
evident that $F_{J}^{v}(\Lambda^{-})=F_{J}^{v}(d_{p}^{-}\Lambda)
=d_{p-J}^{-}F_{J}^{v}(\Lambda)$ and
$F_{J}^{v}(\Lambda^{+})=F_{J}^{v}(\Lambda)$; it follows easily
that $F_{J}^{v}(\Lambda^{-})$ and $F_{J}^{v}(\Lambda^{+})$ are
molecules and
$d_{p-J}^{+}F_{J}^{v}(\Lambda^{-})=d_{p-J}^{-}F_{J}^{v}(\Lambda^{+})$,
as required. If $F_{J}^{v}(\Lambda)=F_{J}^{v}(\lambda_{1})$, then
$\lambda_{s}$ are not $(v,J)$-projection maximal for $s\neq 1$;
thus
$F_{J}^{v}(\lambda_{s})=F_{J}^{v}(\lambda_{1}\cap\lambda_{s})\subset
F_{J}^{v}(d_{p}^{+}\lambda_{1})=d_{p-J}^{+}F_{J}^{v}(\lambda_{1})=
d_{p-J}^{+}F_{J}^{v}(\Lambda)$; it follows easily that
$F_{J}^{v}(\Lambda^{+})=F_{J}^{v}(d_{p}^{+}\Lambda)
=d_{p-J}^{+}F_{J}^{v}(\Lambda)$; it is also evident that
$F_{J}^{v}(\Lambda^{-})=F_{J}^{v}(\lambda_{1})=F_{J}^{v}(\Lambda)$;
therefore $F_{J}^{v}(\Lambda^{-})$ and $F_{J}^{v}(\Lambda^{+})$
are molecules and
$d_{p-J}^{+}F_{J}^{v}(\Lambda^{-})=d_{p-J}^{-}F_{J}^{v}(\Lambda^{+})$,
as required. In the following proof, we may assume that $\dim
F_{J}^{v}(\lambda_{1})>p-J$ and $F_{J}^{v}(\Lambda)$ has at least
two distinct maximal atoms.

Let $$q=\max\{\dim(\mu\cap\mu'):\text{$\mu$ and $\mu'$ are
distinct maximal atoms in $F_{J}^{v}(\Lambda)$}\}.$$ It is clear
that $q\leq p-J$ by the choice of $p$. Let $\mu=\ul{}\times
w^{J}[n,\omega]$ be a maximal atom in $F_{J}^{v}(\Lambda)$
distinct from $F_{J}^{v}(\lambda_{1})$. If $\dim
(F_{J}^{v}(\lambda_{1})\cap\mu)<p-J$, then it is easy to see that
$F_{J}^{v}(\lambda_{1})\cap\mu\subset
d_{p-J}^{+}F_{J}^{v}(\lambda_{1})\cap d_{p-J}^{-}\mu$ by the
construction of molecule $F_{J}^{v}(\Lambda)$ in $u\times w^{J}$
(Theorem \ref{2globe}). Now suppose that $\dim
(F_{J}^{v}(\lambda_{1})\cap\mu)=p-J$. Let
$\lambda'=\ul{}\times\vj{'}\times w[n,\omega]$ be the
$(v,J)$-projection maximal atom in $\Lambda$ such that
$F_{J}^{v}(\lambda')=\mu$. Then $\dim\lambda'>p$. We can also see
that $\min\{j_{1},j'\}=J$ by the choice of $p$ and $\lambda$ is
adjacent to $\lambda'$. Since $F_{J}^{v}(\Lambda)$ is a molecule
in $u\times w^{J}$, we have $i_{1}\neq l$ and $k_{1}\neq n$. If
$i_{1}<l$, then $\alpha_{1}=-$ and $k_{1}>n$; it follows from the
sign condition for $\lambda_{1}$ and $\lambda'$ that
$\omega=(-)^{i_{1}+J}$ which implies that
$F_{J}^{v}(\lambda_{1})\cap\mu\subset
d_{p-J}^{+}F_{J}^{v}(\lambda_{1})\cap d_{p-J}^{-}\mu$. Similarly,
if $i_{1}>l$, then $\sigma_{1}=+$ and $k_{1}<n$; it follows from
the sign condition for $\lambda_{1}$ and $\lambda'$ that
$\varepsilon_{1}=-(-)^{l+J}$ which implies that
$F_{J}^{v}(\lambda_{1})\cap\mu\subset
d_{p-J}^{+}F_{J}^{v}(\lambda_{1})\cap d_{p-J}^{-}\mu$.

Moreover, we have
$F_{J}^{v}(\Lambda^{-})=d_{p-J}^{-}F_{J}^{v}(\Lambda)\cup
F_{J}^{v}(\lambda_{1})$ and
$$F_{J}^{v}(\Lambda^{+})=d_{p-J}^{+}F_{J}^{v}(\Lambda)\cup\bigcup\{\mu:
\mu\text{ is a maximal atom in }F_{J}^{v}(\Lambda)\text{ with }
\mu\neq F_{J}^{v}(\lambda_{1})\}.$$ According to Proposition
\ref{decom_mole}, we can see that $F_{J}^{v}(\Lambda^{-})$ and
$F_{J}^{v}(\Lambda^{+})$ are molecules in $u\times w^{J}$,
$d_{p-J}^{+}F_{J}^{v}(\Lambda^{-})=
d_{p-J}^{-}F_{J}^{v}(\Lambda^{+})$ and
$F_{J}^{v}(\Lambda)=F_{J}^{v}(\Lambda^{-})\#_{p-J}F_{J}^{v}(\Lambda^{+})$,
as required.

This completes the proof.

\end{proof}

\begin{lemma}\label{Fkw}
Let $p\geq K$ and $\lambda_{1}$ be a $(w,K)$-projection maximal
atom. Then
\begin{enumerate}
\item
$F_{K}^{w}(\Lambda^{-})$ and $F_{K}^{w}(\Lambda^{+})$ are
molecules in $u\times v$.
\item
$d_{p-K}^{+}F_{K}^{w}(\Lambda^{-})=d_{p-K}^{-}F_{K}^{w}(\Lambda^{+})$,
hence $F_{K}^{w}(\Lambda^{-})\#_{p-K}F_{K}^{w}(\Lambda^{+})$ is
defined.
\item
$F_{K}^{w}(\Lambda)=F_{K}^{w}(\Lambda^{-})\#_{p-K}F_{K}^{w}(\Lambda^{+})$.
\end{enumerate}
\end{lemma}
\begin{proof}
The argument is similar to  the proof of Lemma \ref{FJv}.

\end{proof}

\begin{prop}
Let $\Lambda$ be a pairwise molecular subcomplex. Then
\begin{enumerate}
\item
$\Lambda^{-}$ and $\Lambda^{+}$ are pairwise molecular
subcomplexes.
\item
$d_{p}^{+}\Lambda^{-}=d_{p}^{-}\Lambda^{+}$, hence  the composite
$\Lambda^{-}\#_{p}\Lambda^{+}$ is defined.
\item
$\Lambda=\Lambda^{-}\#_{p}\Lambda_{+}$.
\end{enumerate}
\end{prop}
\begin{proof}
If  $\lambda_{1}$ is not a $(v,J)$-projection maximal atom in
$\Lambda$, then it is easy to see that
$F_{J}^{v}(\Lambda^{-})=F_{J}^{v}(d_{p}^{-}\Lambda)$ and
$F_{J}^{v}(\Lambda^{+})=F_{J}^{v}(\Lambda)$ by the choice of $p$
and Lemmas \ref{L-c1} and \ref{L+c1}; hence
$F_{J}^{v}(\Lambda^{-})$ and $F_{J}^{v}(\Lambda^{+})$ are the
empty set or  molecules in $u\times w^{J}$. Similarly, if
$\lambda_{1}$ is not $(u,I)$-projection maximal atom in $\Lambda$,
then $F_{I}^{u}(\Lambda^{-})$ and $F_{I}^{u}(\Lambda^{+})$ are
molecules in $v^{I}\times w^{I}$ or the empty set; if
$\lambda_{1}$ is not $(w,K)$-projection maximal atom in $\Lambda$,
then $F_{K}^{w}(\Lambda^{-})$ and $F_{K}^{w}(\Lambda^{+})$ are
molecules in $u\times v$ or the empty set.

According to the above argument and Lemmas \ref{Fiu} to \ref{Fkw},
we can see that $F_{I}^{u}(\Lambda^{-})$,
$F_{I}^{u}(\Lambda^{+})$, $F_{J}^{v}(\Lambda^{-})$,
$F_{J}^{v}(\Lambda^{+})$, $F_{K}^{w}(\Lambda^{-})$ and
$F_{K}^{w}(\Lambda^{+})$ are molecules in the corresponding
$\omega$-complexes or the empty set for all $I$, $J$ and $K$. Thus
$\Lambda^{-}$ and $\Lambda^{+}$ are pairwise molecular.

Now, if $p\geq J$ and $\lambda_{1}$ is not $(v,J)$-projection
maximal, then $F_{J}^{v}(d_{p}^{+}\Lambda^{-})
=d_{p-J}^{+}F_{J}^{v}(\Lambda^{-})
=d_{p-J}^{+}F_{J}^{v}(d_{p}^{-}\Lambda))
=d_{p-J}^{-}F_{J}^{v}(\Lambda) =d_{p-J}^{-}F_{J}^{v}(\Lambda^{+})
=F_{J}^{v}(d_{p}^{-}\Lambda^{+})$; if $p<J$, then
$F_{J}^{v}(d_{p}^{+}\Lambda^{-})=\emptyset=F_{J}^{v}(d_{p}^{-}\Lambda^{+})$.
It follows from Lemmas \ref{Fiu} to \ref{Fkw} and Propositions
\ref{prop5.6} to \ref{prop5.8} that
$F_{I}^{u}(d_{p}^{+}\Lambda^{-})=F_{I}^{u}(d_{p}^{-}\Lambda^{+})$,
$F_{J}^{v}(d_{p}^{+}\Lambda^{-})=F_{J}^{v}(d_{p}^{-}\Lambda^{+})$
and
$F_{K}^{w}(d_{p}^{+}\Lambda^{-})=F_{K}^{w}(d_{p}^{-}\Lambda^{+})$
for all $I$, $J$ and $K$. By  Lemma \ref{eq}, we can see that
$d_{p}^{+}\Lambda^{-}=d_{p}^{-}\Lambda^{+}$. Hence
$\Lambda^{-}\#_{p}\Lambda^{+}$ is defined. Clearly, we have
$\Lambda=\Lambda^{-}\cup\Lambda^{+}$. Therefore
$\Lambda=\Lambda^{-}\#_{p}\Lambda^{+}$.

This completes the proof.
\end{proof}

We have now proved that a pairwise molecular subcomplex $\Lambda$
in $\uvwt$ can be decomposed into pairwise molecular subcomplexes
$\Lambda=\Lambda^{-}\#_{p}\Lambda^{+}$. It is evident that this is
a proper decomposition. By induction, we can see that $\Lambda$
can be eventually decomposed into atoms. Thus $\Lambda$ is a
molecule. So we get the proof for Theorem \ref{pairwise_mole}.

\chapter{Construction of Molecules in the Product of Three Infinite-Dimensional
Globes}\label{3glb}
According to Proposition \ref{mlc1}, the maximal atoms in a
molecule of $\uvwt$ can be listed as $\lambda_{1}$, $\lambda_{2}$,
$\dots$, $\lambda_{R}$ with $\lambda_{r}=\uvw{_{r}}$ such that
$j_{1}\geq\cdots\geq j_{R}$ and such that $i_{r}>i_{r+1}$ when
$1\leq r<R$ and $j_{r}=j_{r+1}$.

In this chapter, we aim to construct molecules by listing their
maximal atoms as described above. The point in this chapter is
that this is easily achieved inductively. In more detail, let
maximal atoms $\lambda_{1}$, $\dots$, $\lambda_{r}$ be an initial
segment of the list. One can easily determine whether
$\lambda_{1}\cup\cdots\cup\lambda_{r}$ is already a molecule and
determine the set of possible next maximal atoms $\lambda_{r+1}$.

Throughout this chapter, the $(v,J+1)$-projection maximal atoms in
a subcomplex of $\uvwt$ are called the {\em lowest} maximal atoms
above level $J$.  An atom with dimension of second factor equal to
$J$ is said to be {\em at level $J$}, while an atom with dimension
of second factor great than $J$ is said to be {\em above level
$J$}. For the convenience of the statement, we allow $J$ to be
$-1$.

\section{Another Description of Molecules}

In this section, we give another description of molecules in terms
of the second factor on which the construction of the molecules is
based.

\begin{prop} \label{bridge}
Let $\Lambda$ be a subcomplex. Suppose that all the maximal atoms
above level $J$ satisfy all the conditions in Theorem
\ref{pairwise_eq}. Suppose also that all the maximal atoms at
level $J$ together with all the lowest maximal atoms above level
$J$ satisfy all the conditions in Theorem \ref{pairwise_eq}. Then
all the maximal atoms above level $J-1$ satisfy all the conditions
in Theorem \ref{pairwise_eq}.
\end{prop}
\begin{proof}

Let $\lambda=\uvws{}$ be a maximal atom above level $J$. Suppose
that $\lambda$ is not lowest above level $J$. Then $j>J+1$ and
there is a lowest maximal atom $\lambda'=\uvws{'}$ above level $J$
such that $i'\geq i$, $J<j'<j$ and $k'\geq k$. Let
$\mu=\ul{}\times v[J,\tau]\times\wn{}$ be a maximal atom at level
$J$. Note that there are no three pairwise adjacent maximal atoms
as in the hypothesis of the condition 5 such that two of them are
at level $J$ and one of them is above level $J$ and not lowest,
hence the condition 5 is automatically satisfied by maximal atoms
above level $J-1$. Now, it suffices to prove that $\lambda$ and
$\mu$ satisfies the conditions 1 to 4.

The condition 1 for $\lambda$ and $\mu$ follows easily from the
condition 1 for $\lambda'$ and $\mu$.

To verify the conditions 2, 3 and 4 for $\lambda$ and $\mu$,
suppose that $\lambda$ and $\mu$ are adjacent. Then $i\geq l$ or
$k\geq n$. The arguments for these two cases are similar. We only
give the proof for the case $i\geq l$.

Suppose that $i\geq l$. Then $k<n$ by condition 1 for $\lambda$
and $\mu$ and $k'=k$ by the adjacency of $\lambda$ and $\mu$.
Hence $\varepsilon'=\varepsilon$ by condition 3 in Theorem
\ref{pairwise_eq} for $\lambda$ and $\lambda'$ and the adjacency
of $\lambda$ and $\mu$. Moreover, we can see that $i'>i$ by
condition 1 for $\lambda$ and $\lambda'$. Thus the condition 3 for
$\lambda$ and $\mu$ is automatically satisfied (whenever $i=l$ and
$\alpha=-\sigma$). Finally, we can see that the conditions 2 and 4
for $\lambda$ and $\mu$ follow from the corresponding conditions
for $\lambda'$ and $\mu$.

This completes the proof.

\end{proof}

The following proposition also characterises molecules in $\uvwt$.

\begin{prop}\label{lv}
Let $\Lambda$ be a subcomplex. Then $\Lambda$ is a molecule if and
only if the following conditions hold for every non-negative
integer $J$:
\begin{enumerate}
\item \label{lvl1}
For every non-negative integer $J$, all the maximal atoms
$\uvws{}$ at level $J$, if there are any, can be listed by
decreasing $i$ and increasing $k$.

\item \label{lvl2}
Suppose that $\uvws{}$ is a lowest maximal atom above level $J$
and $\lmns{}$ is a maximal atom at level $J$. If $l\leq i$, then
$n>k$.

\item\label{lvl3}
Let all the lowest maximal atoms $\lambda_{s}=\uvws{_{s}}$ above
level $J$, if there are any,  be listed as $\lambda_{1}$,
$\cdots$, $\lambda_{S}$ by decreasing $i_{s}$ and increasing
$k_{s}$; let all the maximal atoms $\mu_{t}=\lmns{_{t}}$ at level
$J$, if there are any, be listed as $\mu_{1}$, $\cdots$, $\mu_{T}$
by decreasing $l_{t}$ and increasing $n_{t}$.
\begin{enumerate}

\item\label{lvl31}
If $1<s\leq S$, then there exists $\mu_{t}$ such that
$l_{t}>i_{s}$ and $n_{t}>k_{s-1}$.

\item \label{lvl35}
If $l_{t}>i_{s}$ and $n_{t}>k_{s-1}$ ($1<s\leq S$), then
$\tau_{t}=-(-)^{i_{s}}\alpha_{s}=-(-)^{J}\varepsilon_{s-1}$;

if $l_{t}>i_{1}$, then $\tau_{t}=-(-)^{i_{1}}\alpha_{1}$;

if $n_{t}>k_{S}$, then $\tau_{t}=-(-)^{J}\varepsilon_{S}$;

if $J$ is the greatest dimension of second factors of maximal
atoms in $\Lambda$, then $\tau_{1}=\cdots=\tau_{T}$.

\item\label{lvl32}
If $1<t\leq T$ and if there is no $\lambda_{s}$ such that
$i_{s}>l_{t}$ and $k_{s}>n_{t-1}$, then
$\omega_{t-1}=-(-)^{i_{t}+J}\sigma_{t}$.

\item \label{lvl33}
Suppose that $n_{t}<k_{s}$. If  $l_{t+1}\leq i_{s}$ ($1\leq t<T$),
or if $s=S$ and $t=T$, then $\omega_{t}=-(-)^{i_{s}+J}\alpha_{s}$.

\item \label{lvl34}
Suppose that $l_{t}<i_{s}$. If $n_{t-1}\leq k_{s}$ ($1<t\leq T$)
or if $s=t=1$, then $\sigma_{t}=-(-)^{l_{t}+J}\varepsilon_{s}$.

\item \label{lvl36}
Suppose that $i_{s}=l_{t}$. If $k_{s}>n_{t-1}$ ($1<t\leq T$), or
if $s=t=1$, then $\alpha_{s}=\sigma_{t}$.

\item \label{lvl37}
Suppose that  $k_{s}=n_{t}$. If $i_{s}>l_{t+1}$ ($1\leq t<T$), or
if $s=S$ and $t=T$, then $\varepsilon_{s}=\omega_{t}$.

\item \label{lvl38}
If $1\leq t<T$ and $i_{s}=l_{t+1}$ and  $k_{s}=n_{t}$, then
$\alpha_{s}=\sigma_{t+1}$  or $\varepsilon_{s}=\omega_{t}$.

\end{enumerate}
\end{enumerate}
\end{prop}

\bigskip

\noindent Remark. 1. By induction, it follows easily from
condition \ref{lvl1} and \ref{lvl2} in the proposition that, for
every integer $J$ less than the greatest dimension of second
factors of maximal atoms of $\Lambda$, all the lowest maximal
atoms $\uvws{_{s}}$ above level $J$ can be listed by decreasing
$i_{s}$ and increasing $k_{s}$, as required by the assumption in
Condition \ref{lvl3}.

2. By  condition \ref{lvl35} in the proposition, if $l_{t}>i_{s}$
and $n_{t}\geq k_{s}$, then $\tau_{t}=-(-)^{i_{s}}\alpha_{s}$; if
$l_{t}\geq i_{s}$ and $n_{t}>k_{s}$, then
$\tau_{t}=-(-)^{J}\varepsilon_{s}$. (Hence $l_{t}>i_{s}$ and
$n_{t}>k_{s}$ cannot hold simultaneously unless
$\varepsilon_{s}=(-)^{i_{s}+J}\alpha_{s}$.)

3. It follows from the first part of  condition \ref{lvl35} that
$\varepsilon_{s-1}=(-)^{i_{s}+J}\alpha_{s}$ which we have known
from earlier part of construction.

\bigskip
\begin{proof}
Suppose that $\Lambda$ is a molecule. Then $\Lambda$ satisfies all
the conditions in Theorem \ref{pairwise_eq}. We are going to
verify all the conditions in this proposition.

Firstly, it follows easily from condition 1 in Theorem
\ref{pairwise_eq} that, for every integer $J$, all the maximal
atoms $\uvws{}$ at level $J$, if there are any, can be listed by
decreasing $i$ and increasing $k$, as required.

Next, suppose that $\uvws{}$ is a lowest maximal atom above level
$J$ and $\lmns{}$ is a maximal atom at level $J$. If $l\leq i$,
then it follows easily from condition 1 in Theorem
\ref{pairwise_eq} that $n>k$.

Finally, let all the lowest maximal atoms
$\lambda_{s}=\uvws{_{s}}$, above level $J$, if there are any,  be
listed as $\lambda_{1}$, $\cdots$, $\lambda_{S}$ by decreasing
$i_{s}$ and increasing $k_{s}$; let all the maximal atoms
$\mu_{t}=\lmns{_{t}}$ at level $J$, if there are any, be listed as
$\mu_{1}$, $\cdots$, $\mu_{T}$ by decreasing $l_{t}$ and
increasing $n_{t}$. (These can be done by condition 1 in Theorem
\ref{pairwise_eq}.) We must verify conditions \ref{lvl31} to
\ref{lvl38}. By the definition of lowest, it is easy to see that
every pair of consecutive maximal atoms in the list $\lambda_{1}$,
$\cdots$, $\lambda_{S}$ are adjacent.

\begin{enumerate}

\item[\ref{lvl31}]
Condition \ref{lvl31} follows from condition 4 in Theorem
\ref{pairwise_eq}.

\item[\ref{lvl35}]
Condition \ref{lvl35} follows from  conditions 2 and 3 in Theorem
\ref{pairwise_eq}.

\item[\ref{lvl32}]
Condition \ref{lvl32} follows from condition 2 in Theorem
\ref{pairwise_eq} since $\mu_{t-1}$ and $\mu_{t}$ are adjacent
under the hypothesis of condition \ref{lvl32}.

\item[\ref{lvl33}]
Condition \ref{lvl33} follows from condition 2 in Theorem
\ref{pairwise_eq} since $\lambda_{s}$ and $\mu_{t}$ are adjacent
under the hypothesis of condition \ref{lvl33}.

\item[\ref{lvl34}]
Condition \ref{lvl34} holds by an argument similar to the proof of
condition \ref{lvl33}.

\item[\ref{lvl36}]
Condition \ref{lvl36} follows from condition 3 in Theorem
\ref{pairwise_eq}.

\item[\ref{lvl37}]
Condition \ref{lvl37} also follows from condition 3 in Theorem
\ref{pairwise_eq}.

\item[\ref{lvl38}]
Condition \ref{lvl38} follows from condition 5 in Theorem
\ref{pairwise_eq}.

\end{enumerate}

To prove the sufficiency, suppose that $\Lambda$ satisfies all the
conditions in the proposition. It is evident that the maximal
atoms at the highest level satisfy conditions 1 to 5 in Theorem
\ref{pairwise_eq}. Suppose that $J$ less than the highest level
and all maximal atoms above level $J$ satisfy conditions 1 to 5 in
Theorem \ref{pairwise_eq}. By induction and the proposition
\ref{bridge}, it suffices to prove that all the maximal atoms at
level $J$ together with all the lowest maximal atoms above level
$J$ satisfy conditions 1 to 5 in Theorem \ref{pairwise_eq}.

Condition 1. By the conditions \ref{lvl1} and \ref{lvl2}  in the
proposition, condition 1 in Theorem \ref{pairwise_eq} is satisfied
by all the maximal atoms at level $J$ together with all the lowest
maximal atoms above level $J$.

Condition 2. By condition \ref{lvl32} in the proposition, a pair
of adjacent maximal atoms at level $J$ satisfies condition 2 in
Theorem \ref{pairwise_eq}. Let $\lambda_{s}$ be a lowest maximal
atom  above level $J$ and let $\mu_{t}$ be a maximal atom at level
$J$. Suppose that $\lambda_{s}$ and $\mu_{t}$ are adjacent.

Case 1.  If $l_{t}\geq i_{s}$ and $n_{t}\geq k_{s}$, then
condition 2 for $\lambda_{s}$ and $\mu_{t}$ is satisfied by
remark 2 after the proposition.

Case 2. Suppose that $n_{t}<k_{s}$. Then $l_{t}>i_{s}$ and, by the
adjacency of $\lambda_{s}$ and $\mu_{t}$, we have $l_{t}>i_{s}$
and $l_{t+1}\leq i_{s}$ whenever $t<T$. Hence $\lambda_{s}$ and
$\mu_{t}$ satisfy condition 2 in Theorem \ref{pairwise_eq} by
conditions \ref{lvl31}, \ref{lvl35} and  \ref{lvl33} in this
proposition.

Case 3. Suppose that $l_{t}<i_{s}$. The argument is similar to the
above case.

This completes the proof that all the maximal atoms at level $J$
together with all the lowest maximal atoms above level $J$ satisfy
conditions 2 in Theorem \ref{pairwise_eq}.

Condition 3. Suppose that $\mu_{t}$ and $\mu_{t+1}$ are a pair of
adjacent maximal atoms at level $J$. Suppose also that there is no
maximal atom $\lambda=\uvws{}$ above level $J$ with $i\geq
l_{t+1}$ and $k\geq n_{t}$, then conditions \ref{lvl35} implies
that $\tau_{t}=\tau_{t+1}$, as required by condition 3 in Theorem
\ref{pairwise_eq}.

Indeed, if $l_{t+1}>i_{1}$ or $n_{t}>k_{S}$, then it follows
easily from condition \ref{lvl35} in the proposition that
$\tau_{t}=\tau_{t+1}$, as required. If $l_{t+1}\leq i_{1}$ and
$n_{t}\leq k_{S}$, then $k_{1}<n_{t}$ and $i_{S}<l_{t+1}$ by the
hypothesis, i.e., $i_{S}<l_{t+1}\leq i_{1}$ and $k_{1}<n_{t}\leq
k_{S}$. Now let $\lambda_{s}$ be such that $i_{s+1}<l_{t+1}\leq
i_{s}$. Then $k_{s}<n_{t}$ by the hypothesis. So we have
$l_{t}>l_{t+1}>i_{s+1}$ and $n_{t+1}>n_{t}>k_{s}$. So by condition
\ref{lvl35} in the proposition, it is easy to see that
$\tau_{t}=\tau_{t+1}$, as required by condition 3 in Theorem
\ref{pairwise_eq}.

To finish the proof of condition 3, let $\lambda_{s}$ be a lowest
maximal atom above level $J$ and $\mu_{t}$ be a maximal atom at
level $J$. If $i_{s}=l_{t}$, or if $k_{s}=k_{t}$, then
$\lambda_{s}$ and $\mu_{t}$ are adjacent.   Therefore, by
conditions \ref{lvl36} and \ref{lvl37} in the proposition, it is
evident that condition 3 in Theorem \ref{pairwise_eq} hold for
$\lambda_{s}$ and $\mu_{t}$.

Condition 4. By conditions \ref{lvl31} in the proposition, it is
evident that condition 4 in Theorem \ref{pairwise_eq} is satisfied
by a pair of adjacent maximal atoms at level $J$ since they are
consecutive in the list of lowest maximal atoms above level $J-1$.
Now if $\lambda_{s}$ is a lowest maximal atom above level $J-1$,
and if $\lambda_{s}$ and $\mu_{t}$ are adjacent, then
$\lambda_{s}$ and $\mu_{t}$ are consecutive in the list of lowest
maximal atoms above level $J-1$. So, similar to the above case,
the condition 4 in Theorem \ref{pairwise_eq} holds for
$\lambda_{s}$ and $\mu_{t}$. Suppose that $\lambda_{s}$ is not the
lowest maximal atom above level $J-1$. Suppose also that
$l_{t}<i_{s}$. Then $n_{t}>k_{s}$. In this case, there must be a
maximal atom $\mu'=\lmns{'}$ at level $J$ such that $n'=k_{s}$.
Hence $l'>i_{s}$. It is evident that $\mu'$ and $\mu_{t}$ are
adjacent. Since we have known that condition 4 in Theorem
\ref{pairwise_eq} holds for $\mu'$ and $\mu_{t}$, we can see that
condition 4 in Theorem \ref{pairwise_eq} hold for $\lambda_{s}$
and $\mu_{t}$. If $n_{t}<k_{s}$, then we can see that condition 4
in Theorem \ref{pairwise_eq} holds for $\lambda_{s}$ and $\mu_{t}$
by a similar argument.

Condition 5. By condition \ref{lvl38} in the proposition,
condition 5 in Theorem \ref{pairwise_eq} is satisfied by all the
maximal atoms at level $J$ together with all the lowest maximal
atoms above level $J$.

This completes the proof.

\end{proof}

\bigskip

We can now characterise the sets of  maximal atoms in molecules of
$\uvwt$.

Let $\cal A$ be a finite and non-empty set of atoms in $u\times
v\times w$. For a fixed integer $J$,  An atom $\uvws{}$ above
level $J$ in $\cal A$ is {\em lowest above level $J$} if there is
no atom $\uvws{'}$ in $\cal A$ with $i'\geq i$, $J<j'<j$ and
$k'\geq k$.

Suppose that there are no distinct atoms $\uvws{}$ and $\uvws{'}$
in $\cal A$ such that $i\leq i'$, $j\leq j'$ and $k\leq k'$. Let
$\Lambda$ be the union of atoms in $\cal A$. Then it is evident
that the maximal atoms in $\Lambda$ are exactly the atoms in $\cal
A$. Moreover, it is easy to see that, for every integer $J$ with
$J\geq -1$,  an maximal atom in $\Lambda$ is lowest above level
$J$ in $\Lambda$ if and only if it is lowest above level $J$ in
$\cal A$.

\bigskip

\begin{prop}\label{le}
Let $\cal A$ be a finite and non-empty set of atoms in $u\times
v\times w$. Then $\cal A$ is the set of maximal atoms in a
molecule if and only if the following conditions hold for every
non-negative integer $J$:
\begin{enumerate}
\item \label{le1}
For every non-negative integer $J$, all the atoms $\uvws{}$ at
level $J$ in $\cal A$, if there are any, can be listed by
decreasing $i$ and increasing $k$.

\item \label{le2}
Suppose that $\uvws{}$ is a lowest atom above level $J$ in $\cal
A$ and $\lmns{}$ is an atom at level $J$ in $\cal A$. If $l\leq
i$, then $n>k$.

\item \label{le3}
Let all the lowest atoms $\lambda_{s}=\uvws{_{s}}$ above level $J$
in $\cal A$, if there are any,  be listed as $\lambda_{1}$,
$\cdots$, $\lambda_{S}$ by decreasing $i_{s}$ and increasing
$k_{s}$; let all the atoms $\mu_{t}=\lmns{_{t}}$ at level $J$ in
$\cal A$, if there are any, be listed as $\mu_{1}$, $\cdots$,
$\mu_{T}$ by decreasing $l_{t}$ and increasing $n_{t}$.
\begin{enumerate}

\item\label{le31}
For $1<s\leq S$, there exists $\mu_{t}$ such that $l_{t}>i_{s}$
and $n_{t}>k_{s-1}$.

\item \label{le35}
If $l_{t}>i_{s}$ and $n_{t}>k_{s-1}$ ($1<s\leq S$), then
$\tau_{t}=-(-)^{i_{s}}\alpha_{s}=-(-)^{J}\varepsilon_{s-1}$;

if $l_{t}>i_{1}$, then $\tau_{t}=-(-)^{i_{1}}\alpha_{1}$;

if $n_{t}>k_{S}$, then $\tau_{t}=-(-)^{J}\varepsilon_{S}$;

if $J$ is the greatest dimension of second factors of atoms in
$\cal A$, then $\tau_{1}=\cdots=\tau_{T}$.

\item\label{le32}
If $1<t\leq T$ and if there is no $\lambda_{s}$ such that
$i_{s}>l_{t}$ and $k_{s}>n_{t-1}$, then
$\omega_{t-1}=-(-)^{i_{t}+J}\sigma_{t}$.

\item \label{le33}
Suppose that $n_{t}<k_{s}$. If  $l_{t+1}\leq i_{s}$ ($1\leq t<T$),
or if $s=S$ and $t=T$, then $\omega_{t}=-(-)^{i_{s}+J}\alpha_{s}$.

\item \label{le34}
Suppose that $l_{t}<i_{s}$. If $n_{t-1}\leq k_{s}$ ($1<t\leq T$)
or if $s=t=1$, then $\sigma_{t}=-(-)^{l_{t}+J}\varepsilon_{s}$.

\item \label{le36}
Suppose that $i_{s}=l_{t}$. If $k_{s}>n_{t-1}$ ($1<t\leq T$), or
if $s=t=1$, then $\alpha_{s}=\sigma_{t}$.

\item \label{le37}
Suppose that  $k_{s}=n_{t}$. If $i_{s}>l_{t+1}$ ($1\leq t<T$), or
if $s=S$ and $t=T$, then $\varepsilon_{s}=\omega_{t}$.

\item \label{le38}
If $1\leq t<T$ and $i_{s}=l_{t+1}$ and  $k_{s}=n_{t}$, then
$\alpha_{s}=\sigma_{t+1}$ or $\varepsilon_{s}=\omega_{t}$.

\end{enumerate}
\end{enumerate}
\end{prop}

\noindent Note: By induction, it follows easily from condition
\ref{le1} and condition \ref{le2} in the proposition that, for
every integer $J$ less than the greatest dimension of second
factors of  atoms in $\cal A$, all the lowest atoms $\uvws{_{s}}$
above level $J$ in $\cal A$ can be listed by decreasing $i_{s}$
and increasing $k_{s}$, as required by the assumption in condition
\ref{le3} of the proposition.

\begin{proof}
Suppose that $\cal A$ is the set of maximal atoms in a molecule
$\Lambda$. Then an atom $\uvws{}$ in $\cal A$ is at level $J$ in
$\cal A$ if and only if it is at level $J$ in $\Lambda$; $\uvws{}$
is above level $J$ in $\cal A$ if and only if it is above level
$J$ in $\Lambda$; while $\uvws{}$ is lowest above level $J$ in
$\cal A$ if and only if it is lowest above level $J$ in $\Lambda$.
So the necessity follows from the necessity part of Proposition
\ref{lv}.

Conversely, suppose that a finite and non-empty set $\cal A$
satisfy conditions \ref{le1} to \ref{le3}. It follows from
condition \ref{le1} and \ref{le2} that $\cal A$ is the set of
maximal atoms in a subcomplex $\Lambda$. As in the proof of the
necessity, an atom $\uvws{}$ in $\cal A$ is at level $J$ in $\cal
A$ if and only if it is at level $J$ in $\Lambda$; $\uvws{}$ is
above level $J$ in $\cal A$ if and only if it is above level $J$
in $\Lambda$; while $\uvws{}$ is lowest above level $J$ in $\cal
A$ if and only if it is lowest above level $J$ in $\Lambda$.
Therefore the sufficiency follows from the sufficiency part of
Proposition \ref{lv}.

This completes the proof
\end{proof}

\section{Construction of Molecules}

In this section, we propose an approach of constructing all the
molecules in $\uvwt$ based on Proposition \ref{le}. The
justification will be given in the next section.

We start at the top level and go down.

First choose top level $\bar{J}$ and a fixed sign $\bar{\beta}$
associated with the top level; then choose a list of atoms of the
form $u[\bar{i}_{s},\bar{\alpha}_{s}]\times
v[\bar{J},\bar{\beta}]\times w[\bar{k}_{s},\bar{\varepsilon}_{s}]$
for $1\leq s\leq\bar{S}$, where $\bar{S}\geq 1$, such that
$\bar{i}_{1}>\cdots>\bar{i}_{\bar{S}}$ and
$\bar{k}_{1}<\cdots<\bar{k}_{\bar{S}}$ and
$\bar{\varepsilon}_{s-1}=-(-)^{\bar{i}_{s}+\bar{J}}\bar{\alpha}_{s}$
for $\bar{s}>1$.

For an integer $J$ with $0\leq J<\bar{J}$, suppose that the
 atoms above level $J$ are already constructed. Suppose
also that the lowest  atoms above level $J$ are $\uvws{_{s}}$ with
$1\leq s\leq S$ such that $i_{1}>\cdots>i_{S}$ and
$k_{1}<\cdots<k_{S}$. By condition 1 in Theorem \ref{pairwise_eq},
the  atoms at level $J$, if there are any, can be listed as a
sequence $\ul{_{t}}\times v[J,\tau_{t}]\times\wn{_{t}}$ with
$1\leq t\leq T$, where $T\geq 1$, such that $l_{1}>\cdots>l_{T}$
and $n_{1}<\cdots<n_{T}$.

We are going to give all possibilities for the sequence of atoms
at level $J$.

We first determine the possibilities for the sequence $(l_{1},
n_{1},\cdots,l_{T},n_{T})$ working from left to right.

\begin{enumerate}

\item
We now determine all the possibilities for $l_{1}$ and $n_{1}$.

We determine $l_{1}$ as follows.
\begin{enumerate}
\item
If $S=1$, then there may or may not be  atoms at level $J$; if
there is at least one  atom at level $J$, then $l_{1}\geq 0$.
\item
If $S>1$, then there must be at least one  atom at level $J$ and
$l_{1}>i_{2}$.
\end{enumerate}

For a fixed $l_{1}$, we determine $n_{1}$  as follows.

\begin{enumerate}
\item
If $l_{1}>i_{1}$ and $\varepsilon_{s}=(-)^{i_{s}+J}\alpha_{s}$ for
every $s$, then $n_{1}\geq 0$.

\item
If $l_{1}>i_{1}$ and if there exists $s$ such that
$\varepsilon_{s}=-(-)^{i_{s}+J}\alpha_{s}$, then $0\leq n_{1}\leq
k_{s_{1}}$, where $s_{1}$ is the least $s$ with
$\varepsilon_{s}=-(-)^{i_{s}+J}\alpha_{s}$.

\item
If $l_{1}\leq i_{1}$ and $\varepsilon_{s}=(-)^{i_{s}+J}\alpha_{s}$
for every $s$ with $s>1$, then $n_{1}>k_{1}$.

\item
If $l_{1}\leq i_{1}$ and there exists $s$ with $s>1$ such that
$\varepsilon_{s}=-(-)^{i_{s}+J}\alpha_{s}$, then $k_{1}<n_{1}\leq
k_{s_{2}}$, where $s_{2}$ is the least $s$ with $s>1$ and
$\varepsilon_{s}=-(-)^{i_{s}+J}\alpha_{s}$.
\end{enumerate}

\item
Suppose that $t_{0}>1$ and that $l_{t}$ and $n_{t}$ for all
$t<t_{0}$ are already constructed. We are going to determine all
the possibilities for $l_{t_{0}}$ and $n_{t_{0}}$.

We determine $l_{t_{0}}$ as follows. There are various cases.

\begin{enumerate}
\item
If $n_{t_{0}-1}>k_{S}$ and $l_{t_{0}-1}=0$, then there are no more
 atoms at level $J$.

\item
If $n_{t_{0}-1}>k_{S}$ and $l_{t_{0}-1}>0$, then there may or may
not be another  atom at level $J$; if there is another atom at
level $J$, then $0\leq l_{t_{0}}<l_{t_{0}-1}$.

\item
Suppose that $n_{t_{0}-1}=k_{S}$. Then there may or may not be
another  atom at level $J$. Suppose also that  there is another
 atom at level $J$. If
$\varepsilon_{S}=(-)^{i_{S}+J}\alpha_{S}$, then $0\leq l_{t_{0}}<
l_{t_{0}-1}$; if $\varepsilon_{S}=-(-)^{i_{S}+J}\alpha_{S}$, then
$0\leq l_{t_{0}}\leq i_{S}$.

\item
If $S>1$ and $k_{S-1}<n_{t_{0}-1}<k_{S}$, then there may or may
not be another  atom at level $J$; if there is another atom at
level $J$, then $0\leq l_{t_{0}}<l_{t_{0}-1}$.

\item
If $S=1$ and $n_{t_{0}-1}<k_{1}$, then there may or may not be
another  atom at level $J$; if there is another  atom at level
$J$, then $0\leq l_{t_{0}}<l_{t_{0}-1}$.

\item
Suppose that $1\leq s<S$ and $n_{t_{0}-1}=k_{s}$.  Then there must
be another  atom at level $J$. Moreover, if
$\varepsilon_{s}=(-)^{i_{s}+J}\alpha_{s}$, then
$i_{s+1}<l_{t_{0}}<l_{t_{0}-1}$; if
$\varepsilon_{s}=-(-)^{i_{s}+J}\alpha_{s}$, then
$i_{s+1}<l_{t_{0}}\leq i_{s}$.

\item
If $1<s<S$ and $k_{s-1}<n_{t_{0}-1}<k_{s}$, then there must be
another  atom at level $J$ and $i_{s+1}<l_{t_{0}}<l_{t_{0}-1}$.

\item
If $S>1$ and $n_{t_{0}-1}<k_{1}$ , then there must be another
 atom at level $J$ and  $i_{2}<l_{t_{0}}<l_{t_{0}-1}$.

\end{enumerate}

For a fixed $l_{t_{0}}$, we can determine $n_{t_{0}}$ as follows.

\begin{enumerate}
\item
If $l_{t_{0}}>i_{1}$ and $\varepsilon_{s}=(-)^{i_{s}+J}\alpha_{s}$
for every $s$, then $n_{t_{0}}>n_{t_{0}-1}$.

\item
If $l_{t_{0}}>i_{1}$ and there is $s$ such that
$\varepsilon_{s}=-(-)^{i_{s}+J}\alpha_{s}$, then
$n_{t_{0}-1}<n_{t_{0}}\leq k_{s_{3}}$, where $s_{3}$ be the least
$s$ with $\varepsilon_{s}=-(-)^{i_{s}+J}\alpha_{s}$.

\item
If $l_{t_{0}}\leq i_{S}$, then
$n_{t_{0}}>\max\{n_{t_{0}-1},k_{S}\}$.

\item
If $S>1$ and $i_{s_{4}}<l_{t_{0}}\leq i_{s_{4}-1}$ for some
$s_{4}$, and if $\varepsilon_{s}=(-)^{i_{s}+J}\alpha_{s}$ for
every $s$ with $s\geq s_{4}$, then
$n_{t_{0}}>\max\{k_{s_{4}-1},n_{t_{0}-1}\}$.

\item
If $S>1$ and $i_{s_{4}}<l_{t_{0}}\leq i_{s_{4}-1}$ for some
$s_{4}$, and if there is $s$ such that $s\geq s_{4}$ and
$\varepsilon_{s}=-(-)^{i_{s}+J}\alpha_{s}$, then
$\max\{k_{s_{4}-1},n_{t_{0}-1}\}<n_{t_{0}}\leq k_{s_{5}}$, where
$s_{5}$ is the least $s$ with $s\geq s_{4}$ and
$\varepsilon_{s}=-(-)^{i_{s}+J}\alpha_{s}$.

\end{enumerate}
\end{enumerate}

This completes the construction of the sequence
$(l_{1},n_{1},\cdots,l_{T},n_{T})$.

We now determine the signs $\sigma_{t}$, $\tau_{t}$ and
$\omega_{t}$ for each $t$.

We can determine $\tau_{t}$ for each $t$, as follows.

\begin{enumerate}
\item
If $l_{t}>i_{1}$,  then $\tau_{t}=-(-)^{i_{1}}\alpha_{1}$.

\item
If $s>1$ and $i_{s}<l_{t}\leq i_{s-1}$, then
$\tau_{t}=-(-)^{i_{s}}\alpha_{s}$.

\item
If $l_{t}\leq i_{S}$ (in this case, we have $n_{t}>k_{S}$), then
$\tau_{t}=-(-)^{J}\varepsilon_{S}$.
\end{enumerate}

We now determine signs $\sigma_{t}$ and $\omega_{t}$ for each $t$.

We first determine $\sigma_{1}$.

\begin{enumerate}
\item
If $l_{1}>i_{1}$, then $\sigma_{1}$ is arbitrary.

\item
If $l_{1}=i_{1}$, then $\sigma_{1}=\alpha_{1}$.

\item
If $l_{1}<i_{1}$, then $\sigma_{1}=-(-)^{l_{1}+J}\varepsilon_{1}$.
\end{enumerate}

We next determine $\omega_{t-1}$ and $\sigma_{t}$ for $1<t\leq T$.
Note that there can be at most one value of $s$ such that
$i_{s}\geq l_{t}$ and $k_{s}\geq n_{t-1}$ by the construction of
$l_{t}$. There are various cases, as follows.

\begin{enumerate}
\item
If there is no $s$ such that $i_{s}\geq l_{t}$ and $k_{s}\geq
n_{t-1}$, then $\omega_{t-1}$ is arbitrary and
$\sigma_{t}=-(-)^{l_{t}+J}\omega_{t-1}$ for a fixed
$\omega_{t-1}$.

\item
If there exists $s$ such that $i_{s}>l_{t}$ and $k_{s}>n_{t-1}$,
then $\omega_{t-1}=-(-)^{i_{s}+J}\alpha_{s}$ and
$\sigma_{t}=-(-)^{l_{t}+J}\varepsilon_{s}$.

\item
If there exists $s$ such that $i_{s}=l_{t}$ and $k_{s}>n_{t-1}$,
then $\omega_{t-1}=-(-)^{i_{s}+J}\alpha_{s}$ and
$\sigma_{t}=\alpha_{s}$.

\item
If there exists $s$ such that $i_{s}>l_{t}$ and $k_{s}=n_{t-1}$,
then $\omega_{t-1}=\varepsilon_{s}$ and
$\sigma_{t}=-(-)^{l_{t}+J}\varepsilon_{s}$.

\item
Suppose that there exists $s$ such that $i_{s}=l_{t}$ and
$k_{s}=n_{t-1}$. If $\varepsilon_{s}=(-)^{i_{s}+J}\alpha_{s}$,
then $\omega_{t-1}$ is arbitrary and
$\sigma_{t}=-(-)^{l_{t}+J}\omega_{t-1}$ for a fixed
$\omega_{t-1}$. If $\varepsilon_{s}=-(-)^{i_{s}+J}\alpha_{s}$,
then $\omega_{t-1}=\varepsilon_{s}$ and $\sigma_{t}=\alpha_{s}$.

\end{enumerate}

Finally, we determine $\omega_{T}$.

\begin{enumerate}
\item
If $n_{T}>k_{S}$, then $\omega_{T}$ is arbitrary.

\item
If $n_{T}=k_{S}$, then $\omega_{T}=\varepsilon_{T}$.

\item
If $n_{T}<k_{S}$, then $\omega_{T}=-(-)^{i_{S}+J}\alpha_{S}$.
\end{enumerate}

This completes the construction of all the possibilities for the
sequence of  atoms at level $J$. Therefore, by induction, we can
construct all the molecules in $\uvwt$.

\bigskip

 \begin{remark}
In a subcomplex as constructed in the last section, we verify that
the permitted value of  $l_{t}$ and $n_{t}$ form non-empty
intervals of integers for each $t$.

By the construction of  atoms at level $J$, it is evident that a
lowest  atom above level $J$ and an atom at level $J$ satisfy
condition 1 in Theorem \ref{pairwise_eq}.

\begin{enumerate}

\item
It is evident that the permitted value of  $l_{1}$ and $n_{1}$
form non-empty intervals of integers.

\item
In the construction of $l_{t_{0}}$, it is evident that the
permitted values of $l_{t_{0}}$ form a non-empty interval of
integers in (c) part two and (f) part two. If $n_{t_{0}-1}\leq
k_{S}$, then we have $l_{t_{0}-1}>i_{S}\geq 0$.  Therefore the
permitted values of $l_{t_{0}}$ in (b), (c) part one, (d) and (e)
form a non-empty interval of integers.  Finally, if $s<S$ and
$n_{t_{0}-1}\leq k_{s}$, then we have $l_{t_{0}-1}>i_{s}\geq
i_{s+1}+1$ by condition 1 for $\uvws{_{s}}$ and
$\ul{_{t_{0}-1}}\times v[J,\tau]\times\wn{_{t_{0}-1}}$. This
implies that the permitted values of $l_{t_{0}}$ form a non-empty
interval of integers in (f) part one, (g) and (h).

\item
In the construction of $n_{t_{0}}$, it is evident that the
permitted values of $n_{t_{0}}$ forms a non-empty interval of
integers in (a), (c) and (d).

Suppose that $l_{t_{0}}>i_{s_{4}}$ for some $s_{4}$. Suppose also
that there is $s$ with $s\geq s_{4}$ such that
$\varepsilon_{s}=-(-)^{i_{s}+J}\alpha_{s}$. Let $s_{5}$ be the
least $s$ with $s\geq s_{4}$ and
$\varepsilon_{s}=-(-)^{i_{s}+J}\alpha_{s}$. We claim that
$n_{t_{0}-1}<k_{s_{5}}$ which implies that the permitted values of
$n_{t_{0}}$ forms a non-empty interval of integers in (b) and (e).

Indeed, since $l_{t_{0}-1}>l_{t_{0}}>i_{s_{4}}\geq i_{s_{5}}$, we
have $n_{t_{0}-1}\leq k_{s_{5}}$ by the construction of
$l_{t_{0}-1}$ and $n_{t_{0}-1}$. If $n_{t_{0}-1}=k_{s_{5}}$, then
$l_{t_{0}}\leq i_{s_{5}}\leq i_{s_{4}}$ by the construction of
$l_{t_{0}}$; this contradicts the assumption on $l_{t_{0}}$.
Therefore we have $n_{t_{0}-1}<k_{s_{5}}$, as required.

\end{enumerate}
Therefore, the permitted value of  $l_{t}$ and $n_{t}$ form
non-empty intervals of integers for each $t$.
\end{remark}

\begin{example}

The the molecule in Example \ref{pairwise_example} is really
constructed by the approach in this section. The construction of
the example involves most of the above cases.

\end{example}

\section{Justification}

In this section, we prove that the construction in the last
section indeed gives molecules in $\uvwt$.

\begin{lemma}\label{llvJ}
In a subcomplex as constructed in the last section, for every
level $J$, the  atoms $\ul{_{t}}\times
v[J,\tau_{t}]\times\wn{_{t}}$ with $1\leq t\leq T$ at level $J$
satisfy $l_{1}>\cdots>l_{T}$ and $n_{1}<\cdots<n_{T}$.
\end{lemma}
\begin{proof}

By induction, it suffices to verify that $l_{t_{0}}<l_{t_{0}-1}$
and $n_{t_{0}}>n_{t_{0}-1}$ in the construction of $l_{t_{0}}$ and
$n_{t_{0}}$.

In the construction of $l_{t_{0}}$, we have already required that
$l_{t_{0}}<l_{t_{0}-1}$ except in (c) part two and (f) part two.
Now if $n_{t_{0}-1}=k_{s}$ for some $s$, then $i_{s}<l_{t_{0}-1}$
by the earlier part of construction (or, more precisely, by the
induction hypothesis); hence $l_{t_{0}}<l_{t_{0}-1}$ in (c) part
two and (f) part two, as required.

In the construction of $n_{t_{0}}$, we have already required that
$n_{t_{0}}>n_{t_{0}-1}$ in all cases.

Therefore, the  atoms $\ul{_{t}}\times
v[J,\tau_{t}]\times\wn{_{t}}$ at level $J$ as constructed can be
listed by decreasing $l_{t}$ and increasing $n_{t}$ for each level
$J$, as required.

This completes the proof.
\end{proof}

\begin{lemma}\label{llvJJ}
In a subcomplex as constructed in the last section, all the atoms
constructed satisfy condition 1 in Theorem \ref{pairwise_eq}.
Hence all the lowest atoms $\uvws{_{s}}$ above  level $J-1$ can be
listed such that $i_{1}>\cdots>i_{S}$ and $k_{1}<\cdots<k_{S}$.
\end{lemma}
\begin{proof}
We first show that all the atoms constructed satisfy condition 1
in Theorem \ref{pairwise_eq}.

It is evident that all the  atoms at the top level $\bar{J}$
satisfy condition 1 in Theorem \ref{pairwise_eq}.

Suppose that $J<\bar{J}$ and that all the   atoms above level $J$
satisfy condition 1 in Theorem \ref{pairwise_eq}. We are going to
show that all the  atoms above level $J-1$ satisfy condition 1 in
Theorem \ref{pairwise_eq}.

It follows from the Lemma \ref{llvJ} that a pair of atoms at level
$J$ satisfy condition 1 in Theorem \ref{pairwise_eq}. Let
$\uvws{}$ be an atom above level $J$ and $\ul{}\times
v[J,\omega]\times\wn{}$ be an atom at level $J$. If $\uvws{}$ is
lowest above level $J$, then it is evident that $\uvws{}$ and
$\ul{}\times v[J,\omega]\times\wn{}$ satisfy condition 1 in
Theorem \ref{pairwise_eq} by the construction of  atoms at level
$J$. If $\uvws{}$ is not lowest above level $J$, then there is a
lowest
 atom $\uvws{'}$ above level $J$ with $i'\geq i$, $j'<j$ and
$k'\geq k$; hence condition 1 for $\uvws{}$ and $\ul{}\times
v[J,\omega]\times\wn{}$ follows easily from condition 1 for
$\uvws{}$ and $\ul{}\times v[J,\omega]\times\wn{}$. Thus all the
 atoms above level $J-1$ satisfy condition 1 in Theorem \ref{pairwise_eq}.

Therefore, all the  atoms satisfy condition 1 in Theorem
\ref{pairwise_eq}.

Now let $\uvws{}$ and $\uvws{'}$ be a pair of lowest  atoms above
level $J$. By the definition of lowest, we have $i\neq i'$ and
$k\neq k'$. Moreover, it follows easily from condition 1 for
pairwise molecular subcomplexes for $\uvws{}$ and $\uvws{'}$ that
$i>i'$ if and only if $k<k'$. Hence all  the lowest  atoms
$\uvws{_{s}}$ above level $J$ can be listed by decreasing $i_{s}$
and increasing $k_{s}$, as required.
\end{proof}

\begin{lemma}\label{lvvJ}
In a subcomplex as constructed in the last section, let all the
lowest atoms $\huvws{_{s}}$ above level $J-1$  be listed such that
$\hat{i}_{1}>\cdots>\hat{i}_{\hat{S}}$ and
$\hat{k}_{1}<\cdots<\hat{k}_{\hat{S}}$. Then $\min\{\hat{j}_{s-1},
\hat{j}_{s}\}=J$ and
$\hat{\varepsilon}_{s-1}=-(-)^{\hat{i}_{s}+J}\hat{\alpha}_{s}$ for
every $1<s\leq \hat{S}$.
\end{lemma}
\begin{proof}

Let $\huvws{_{s-1}}$ and $\huvws{_{s}}$ be a pair of consecutive
lowest  atoms  above level $J-1$. We first show that
$\min\{\hat{j}_{s-1}, \hat{j}_{s}\}=J$.

Indeed, suppose otherwise that $\min\{\hat{j}_{s-1},
\hat{j}_{s}\}>J$. Then we can see that $\huvws{_{s-1}}$ and
$\huvws{_{s}}$ are lowest
 atoms above level $J$ and they are consecutive in the list
for lowest  atoms above level $J$.  It follows from the
construction that there is an atom $\ul{}\times
v[J,\hat{\beta}]\times\wn{}$ at level $J$ with $l>\hat{i}_{s}$ and
$n>\hat{k}_{s-1}$. Since $\huvws{_{s-1}}$ and $\huvws{_{s}}$ are
lowest atom  above level $J-1$, we have
$\hat{i}_{s-1}>l>\hat{i}_{s}$ and $\hat{k}_{s-1}<n<\hat{k}_{s}$.
This contradicts the assumption that $\huvws{_{s-1}}$ and
$\huvws{_{s}}$ are a pair of consecutive lowest  atoms above level
$J-1$.

Now we are going to show that
$\hat{\varepsilon}_{s-1}=-(-)^{\hat{i}_{s}+J}\hat{\alpha}_{s}$ for
every $1<s\leq \hat{S}$. Note that either $\huvws{_{s-1}}$ or
$\huvws{_{s}}$ is an atom at level $J$ by the first part of the
lemma. Now there are several cases, as follows.

If both $\huvws{_{s-1}}$ and $\huvws{_{s}}$ are  atoms at level
$J$, then, by the construction of the signs for  atoms at level
$J$, it is evident that
$\hat{\varepsilon}_{s-1}=-(-)^{\hat{i}_{s}+J}\hat{\alpha}_{s}$ for
every $1<s\leq \hat{S}$, as required.

Suppose that $\huvws{_{s-1}}$ is an atom above level $J$ and
$\huvws{_{s}}$ is an atom at level $J$. Then
$\huvws{_{s}}=\ul{_{t}}\times v[J,\tau_{t}]\times \wn{_{t}}$ for
some $t$ in the construction. If $t=1$, then we have
$\huvws{_{s-1}}=\uvws{_{1}}$ and
$\hat{\varepsilon}_{s-1}=-(-)^{\hat{i}_{s}+J}\hat{\alpha}_{s}$ by
the construction, as required. If $t>1$, then it is easy to see
that $s>2$ and $\huvws{_{s-2}}=\ul{_{t-1}}\times
v[J,\tau_{t-1}]\times\wn{_{t-1}}$ by the first part of this lemma;
thus we have
$\hat{\varepsilon}_{s-1}=-(-)^{\hat{i}_{s}+J}\hat{\alpha}_{s}$ by
the construction of signs, as required.

Suppose that $\huvws{_{s-1}}$ is an atom above level $J$ and
$\huvws{_{s}}$ is an atom at level $J$. By an argument similar to
the above case, one can also get
$\hat{\varepsilon}_{s-1}=-(-)^{\hat{i}_{s}+J}\hat{\alpha}_{s}$ for
every $1<s\leq \hat{S}$, as required.

This completes the proof.
\end{proof}

By Proposition \ref{le} and the remark after the statement of the
proposition, it is easy to see that every molecule can be
constructed as above. Now we are going to prove that every
subcomplex  of $\uvwt$ constructed as above is indeed a molecule.

\begin{prop}
Let $\Lambda$ be a subcomplex whose maximal atoms are as
constructed above. Then $\Lambda$ is a molecule.
\end{prop}
\begin{proof}
Let $\cal{A}$ be a set of  atoms as constructed above. It suffices
to show that $\cal{A}$ satisfies all the conditions in Proposition
\ref{le}.

By Lemmas \ref{llvJ} and \ref{llvJJ}, it is easy to see that
conditions \ref{le1} and \ref{le2} hold.

Now let all the lowest atoms $\lambda_{s}=\uvws{_{s}}$ with
dimension of second factors greater than $J$ in $\cal A$, if there
are any, be listed as $\lambda_{1}$, $\cdots$, $\lambda_{S}$ by
decreasing $i_{s}$ and increasing $k_{s}$; let all the atoms
$\mu_{t}=\lmns{_{t}}$ with dimension of second factors equal to
$J$ in $\cal A$, if there are any, be listed as $\mu_{1}$,
$\cdots$, $\mu_{T}$ by decreasing $l_{t}$ and increasing $n_{t}$.
By the construction of $l_{t}$ and $n_{t}$ for all $t$, we can see
that condition \ref{le31} hold. Moreover, by the construction of
signs $\sigma_{t}$ and $\omega_{t}$, it is easy to see that
conditions \ref{le32} to \ref{le38} hold. To complete the proof,
we need only to verify condition \ref{le35}.

Suppose that $l_{t}>i_{s}$ and $n_{t}>k_{s-1}$ ($1<s\leq S$). Let
$\hat{s}$ be such that $i_{\hat{s}}<l_{t}\leq i_{\hat{s}-1}$. By
the construction, we have
$\tau_{t}=-(-)^{i_{\hat{s}}}\alpha_{i_{\hat{s}}}$. If $s=\hat{s}$,
then
$\tau_{t}=-(-)^{i_{s}}\alpha_{s}=-(-)^{i_{s}}[-(-)^{i_{s}+J+1}]
\varepsilon_{s-1}= -(-)^{J}\varepsilon_{s-1}$, as required,  since
$\varepsilon_{s-1}=-(-)^{i_{s}+J+1}\alpha_{s}$ by Lemma
\ref{lvvJ}. Suppose that $\hat{s}<s$. Then
$l_{t}>i_{\hat{s}}>\cdots >i_{s}$ and $n_{t}>k_{s-1}>\cdots
>k_{\hat{s}}$. Hence
$\varepsilon_{\hat{s}}=(-)^{i_{\hat{s}}+J}\alpha_{\hat{s}}$,
$\cdots$, $\varepsilon_{s-1}=(-)^{i_{s-1}+J}\alpha_{s-1}$ by the
construction of signs. It follows that $$
\begin{array}{rl}
\tau_{t}&=-(-)^{i_{\hat{s}}}\alpha_{i_{\hat{s}}}\\
&=-(-)^{i_{\hat{s}}}(-)^{i_{\hat{s}}+J}\varepsilon_{\hat{s}}\\
&=-(-)^{J}\varepsilon_{\hat{s}}\\
&=-(-)^{J}[-(-)^{i_{\hat{s}+1}+J+1}\alpha_{i_{\hat{s}+1}}]\\
&=-(-)^{i_{\hat{s}+1}}\alpha_{i_{\hat{s}+1}}\\ &=\cdots\\
&=-(-)^{J}\varepsilon_{s-1}\\ &=-(-)^{i_{s}}\alpha_{s}
\end{array}
$$ as required by condition \ref{le35}.

The other parts of condition \ref{le35} can be seen easily from
the construction of the sign $\tau_{t}$ for each $t$.

This completes the proof.

\end{proof}


\chapter{Molecules in the Product of Four Infinite-Dimensional Globes}
In this chapter, we study molecules in the product of four
infinite dimensional globes. Similar to the results for molecules
in the product of three infinite dimensional globes, we are going
to give some equivalent descriptions for the molecules in the
product of four infinite dimensional globes. The discussion is in
parallel to that in chapter 2. There are some new features because
of the two `middle' factors.

In this chapter, all the subcomplexes refer to finite and
non-empty subcomplexes in the $\omega$-complex $\ubbb$, and all
the integers refer to non-negative integers.

Recall that the $\omega$-complex $u^{I}$ is equivalent to the
infinite dimensional globe $u$. It is easy to see that this
equivalence induces an equivalence $u\times v\times w$ to
$u^{I}\times v^{J}\times w^{K}$ of $\omega$-complexes sending
every atom $\uvw{}$ to $u^{I}[i,(-)^{I}\alpha]\times
v^{J}[j,(-)^{J}\beta]\times w^{K}[k,(-)^{K}\varepsilon]$. Thus all
the results for the molecules in the product of three globes can
be generalised to the molecules in the product of three `twisted'
infinite dimensional globes $u^{I}\times v^{J}\times w^{K}$. In
particular, a pairwise molecular subcomplex in $u^{I}\times
v^{J}\times w^{K}$ is defined as the image of a pairwise molecular
subcomplex in $\uvwt$ under the above equivalence of
$\omega$-complexes and a subcomplex of $u^{I}\times v^{J}\times
w^{K}$ is a molecule if and only if it is pairwise molecular. We
are not going to make any more comments of this kind.

\section{The Definition of Pairwise Molecular Subcomplexes}

In this section, we define projection maps and give the definition
of pairwise molecular subcomplexes of $\ubbb$. Some proofs are
omitted because the arguments are very similar to that in Chapter
\ref{3glb}.

For an atom $\lambda=\uaaa{i}{\alpha}{}$ in $\ubbb$, let $$
F_{I_{1}}^{u_{1}}(\Int\lambda)=
  \begin{cases}
    \Int(u^{I_{1}}_{2}[i_{2},\alpha_{2}]\times
         u^{I_{1}}_{3}[i_{3},\alpha_{3}]\times
         u^{I_{1}}_{4}[i_{4},\alpha_{4}]), & \text{when $i_{1}\geq I_{1}$};
         \\
    \emptyset, & \text{when $j<J$}.
  \end{cases}
$$ This gives a map sending interiors of atoms in $\ubbb$ to
interiors of atoms in $u_{2}^{I_{1}}\times u_{3}^{I_{1}}\times
u_{4}^{I_{1}}$ or the empty set.

Since interiors of atoms are disjoint, it is clear that the map
$F_{I_{1}}^{u_{1}}$ can be extended uniquely to a map sending
unions of interiors of atoms in $\uaaa{i}{\alpha}{}$ to unions of
interiors of atoms in $u_{2}^{I_{1}}\times u_{3}^{I_{1}}\times
u_{4}^{I_{1}}$ by requiring it union-preserving.

We can similarly define a map $F_{I_{2}}^{u_{2}}$ sending unions
of interiors of atoms in $\ubbb$ to unions of interiors of atoms
in $u_{1}\times u_{3}^{I_{2}}\times u_{4}^{I_{2}}$, a map
$F_{I_{3}}^{u_{3}}$ sending unions of interiors of atoms in
$\ubbb$ to unions of interiors of atoms in $u_{1}\times
u_{2}\times u_{4}^{I_{3}}$ and a map $F_{I_{4}}^{u_{4}}$ sending
unions of interiors of atoms in $\ubbb$ to unions of interiors of
atoms in $u_{1}\times u_{2}\times u_{3}$.

It is easy to see that every atom in $\ubbb$ can be written as a
union of interiors of atoms. It follows that $F_{I_{1}}^{u_{1}}$,
$F_{I_{2}}^{u_{2}}$, $F_{I_{3}}^{u_{3}}$ and $F_{I_{4}}^{u_{4}}$
are defined on subcomplexes of $\ubbb$ and preserve unions.

We shall prove that $F_{I_{s}}^{u_{s}}$ sends atoms to atoms or
the empty set so that it sends subcomplexes to subcomplexes for
every $s$. We need a preliminary result.

\begin{lemma}\label{dtatom4}
Let $\lambda=\uaaa{i}{\alpha}{}$ be an atom in $\ubbb$.
\begin{enumerate}
\item
If $i_{1}+i_{2}+i_{3}+i_{4}\leq p$, then $\dpg\lambda=\lambda$.
\item
If $i_{1}+i_{2}+i_{3}+i_{4}>p$, then the set of maximal atoms in
$\dpg\lambda$ consists of all the atoms $\uaaa{l}{\sigma}{}$ such
that $l_{1}\leq i_{1}$, $l_{2}\leq i_{2}$, $l_{3}\leq i_{3}$,
$l_{4}\leq i_{4}$, where the signs $\sigma_{1}$, $\sigma_{2}$,
$\sigma_{3}$ and $\sigma_{4}$ are determined as follows:
\begin{enumerate}
\item
if $l_{1}=i_{1}$, then $\sigma_{1}=\alpha_{1}$; If $l_{1}<i_{1}$,
then $\sigma_{1}=\gamma$;

\item
if $l_{2}=i_{2}$, then $\sigma_{2}=\alpha_{2}$; If $l_{2}<i_{2}$,
then $\sigma_{2}=(-)^{l_{1}}\gamma$;

\item
if $l_{3}=i_{3}$, then $\sigma_{3}=\alpha_{3}$; If $l_{3}<i_{3}$,
then $\sigma_{3}=(-)^{l_{1}+l_{2}}\gamma$;

\item
if $l_{4}=i_{4}$, then $\sigma_{4}=\alpha_{4}$; If $l_{4}<i_{4}$,
then $\sigma_{4}=(-)^{l_{1}+l_{2}+l_{3}}\gamma$.
\end{enumerate}
\end{enumerate}

\end{lemma}

\begin{prop}
Let $\lambda=\uaaa{i}{\alpha}{}$ be an atom  in $\ubbb$. Then
\begin{enumerate}
\item
$
F_{I_{1}}^{u_{1}}(\lambda)=
  \begin{cases}
    u^{I_{1}}_{2}[i_{2},\alpha_{2}]\times
         u^{I_{1}}_{3}[i_{3},\alpha_{3}]\times
         u^{I_{1}}_{4}[i_{4},\alpha_{4}], & \text{when $i_{1}\geq I_{1}$}; \\
    \emptyset, & \text{when $i_{1}<I_{1}$};
  \end{cases}
$
\item
$
F_{I_{2}}^{u_{2}}(\lambda)=
  \begin{cases}
    u_{1}[i_{1},\alpha_{1}]\times
         u^{I_{2}}_{3}[i_{3},\alpha_{3}]\times
         u^{I_{2}}_{4}[i_{4},\alpha_{4}], & \text{when $i_{2}\geq I_{2}$}; \\
    \emptyset, & \text{when $i_{2}<I_{2}$};
  \end{cases}
$

\item
$
F_{I_{3}}^{u_{3}}(\lambda)=
  \begin{cases}
    u_{1}[i_{1},\alpha_{1}]\times
         u_{2}[i_{2},\alpha_{2}]\times
         u^{I_{3}}_{4}[i_{4},\alpha_{4}], & \text{when $i_{3}\geq I_{3}$}; \\
    \emptyset, & \text{when $i_{3}<I_{3}$};
  \end{cases}
$
\item
$
F_{I_{4}}^{u_{4}}(\lambda)=
  \begin{cases}
    u_{1}[i_{1},\alpha_{1}]\times
         u_{2}[i_{2},\alpha_{2}]\times
         u_{3}[i_{3},\alpha_{3}], & \text{when $i_{4}\geq I_{4}$}; \\
    \emptyset, & \text{when $i_{4}<I_{4}$}.
  \end{cases}
$

\end{enumerate}

In particular, the maps $F_{I_{1}}^{u_{1}}$, $F_{I_{2}}^{u_{2}}$,
$F_{I_{3}}^{u_{3}}$ and $F_{I_{4}}^{u_{4}}$ send atoms to atoms or
the empty set.
\end{prop}

\begin{proof}
The arguments for the four cases are similar. We only prove the
second one. The proof is given by induction on dimension of atoms.

For an atom $\lambda=\uaaa{i}{\alpha}{}$ in $\ubbb$, if
$\dim\lambda=0$, then $i_{1}=i_{2}=i_{3}=i_{4}=0$; hence $$
\begin{array}{rl}
&F_{I_{2}}^{u_{2}}(\lambda)\\ =&F_{I_{2}}^{u_{2}}(\Int\lambda)\\
=&
\begin{cases}
    \Int(\ua{1}\times u^{I_{2}}_{3}[i_{3},\alpha_{3}]\times
    u^{I_{2}}_{4}[i_{4},\alpha_{4}]), & \text{when $I_{2}=0$}; \\
    \emptyset, & \text{when $I_{2}>0$}
  \end{cases} \\
=&
\begin{cases}
\ua{1}\times u^{I_{2}}_{3}[i_{3},\alpha_{3}]\times
    u^{I_{2}}_{4}[i_{4},\alpha_{4}], & \text{when $I_{2}=0$}; \\
    \emptyset, & \text{when $I_{2}>0$}
  \end{cases}
\end{array}
$$ as required.

Suppose that $p>0$ and that the proposition holds for every atom
of dimension less then $p$. Suppose also that
$\lambda=\uaaa{i}{\alpha}{}$ is a $p$-dimensional atom. If
$i_{2}<I_{2}$, then it is easy to see that
$F_{I_{2}}^{u_{2}}(\lambda)=\emptyset$, as required. If
$i_{2}>I_{2}$, then we have $$
\begin{array}{rl}
&F_{I_{2}}^{u_{2}}(\lambda)\\
=&F_{I_{2}}^{u_{2}}(\Int\lambda\cup\partial^{-}\lambda\cup\partial^{+}\lambda)\\
\supset & F_{I_{2}}^{u_{2}}(\partial^{+}\lambda)\\ \supset &
F_{I_{2}}^{u_{2}}(\ua{1}\times u_{2}[i_{2}-1,(-)^{i_{1}}]\times
\ua{3}\times\ua{4})\\ =&\ua{1}\times
u^{I_{2}}_{3}[i_{3},\alpha_{3}]\times
u^{I_{2}}_{4}[i_{4},\alpha_{4}]
\end{array}
$$ since $\ui{}\times v[j-1,(-)^{i}]\times \wk{}$ is an atom of
dimension $p-1$; the reverse inclusion holds automatically; so
$F_{I_{2}}^{u_{2}}(\lambda)=\ua{1}\times
u^{I_{2}}_{3}[i_{3},\alpha_{3}]\times
u^{I_{2}}_{4}[i_{4},\alpha_{4}]$, as required. Suppose that
$i_{2}=I_{2}$. Then $\partial^{\gamma}\lambda$ is the union of
atoms $\uaaa{l}{\sigma}{}$ with $l_{1}+l_{2}+l_{3}+l_{4}=p-1$ such
that
\begin{enumerate}
\item
if $l_{1}=i_{1}$, then $\sigma_{1}=\alpha_{1}$; if
$l_{1}=i_{1}-1$, then $\sigma_{1}=\gamma$;
\item
if $l_{2}=i_{2}$, then $\sigma_{2}=\alpha_{2}$; if
$l_{2}=i_{2}-1$, then $\sigma_{2}=(-)^{i_{1}}\gamma$;
\item
if $l_{3}=i_{3}$, then $\sigma_{3}=\alpha_{3}$; if
$l_{3}=i_{3}-1$, then $\sigma_{3}=(-)^{i_{1}+I_{2}}\gamma$.
\item
if $l_{4}=i_{4}$, then $\sigma_{4}=\alpha_{4}$; if
$l_{4}=i_{4}-1$, then $\sigma_{4}=(-)^{i_{1}+I_{2}+i_{3}}\gamma$.
\end{enumerate}
It follows easily from the inductive hypothesis and Theorem
\ref{dt} that $F_{I_{2}}^{u_{2}}(\partial^{\gamma}\lambda)=
\partial^{\gamma}(\ua{1}\times
u^{I_{2}}_{3}[i_{3},\alpha_{3}]\times
u^{I_{2}}_{4}[i_{4},\alpha_{4}])$ for every sign $\gamma$.
Therefore $$
\begin{array}{rl}
 &F_{I_{2}}^{u_{2}}(\lambda)\\
=&F_{I_{2}}^{u_{2}}(\Int\lambda)\cup
F_{I_{2}}^{u_{2}}(\partial^{-}\lambda)\cup
F_{I_{2}}^{u_{2}}(\partial^{+}\lambda)\\ =& \Int(\ua{1}\times
u^{I_{2}}_{3}[i_{3},\alpha_{3}]\times
u^{I_{2}}_{4}[i_{4},\alpha_{4}])\cup
\partial^{-}(\ua{1}\times
u^{I_{2}}_{3}[i_{3},\alpha_{3}]\times
u^{I_{2}}_{4}[i_{4},\alpha_{4}])\cup \\ &
\partial^{+}(\ua{1}\times u^{I_{2}}_{3}[i_{3},\alpha_{3}]\times
u^{I_{2}}_{4}[i_{4},\alpha_{4}])\\ =&\ua{1}\times
u^{I_{2}}_{3}[i_{3},\alpha_{3}]\times
u^{I_{2}}_{4}[i_{4},\alpha_{4}],
\end{array}
$$ as required.

This completes the proof of the proposition.
\end{proof}

We now define the concept of pairwise molecular subcomplexes as
follows.

\begin{definition}
Let $\Lambda$ be a subcomplex of $u_{1}\times u_{2}\times
u_{3}\times u_{4}$. Then $\Lambda$ is {\em pairwise molecular} if
\begin{enumerate}
\item
There are no distinct maximal atoms $u_{1}[i_{1},\alpha_{1}]\times
u[i_{2},\alpha_{2}]\times u[i_{3},\alpha_{3}]\times
u[i_{4},\alpha_{4}]$ and $u_{1}[i_{1}',\alpha_{1}']\times
u[i_{2}',\alpha_{2}']\times u[i_{3}',\alpha_{3}']\times
u[i_{4}',\alpha_{4}']$ in $\Lambda$ such that $i_{1}\leq i_{1}'$,
$i_{2}\leq i_{2}'$,  $i_{3}\leq i_{3}'$ and $i_{4}\leq i_{4}'$.
\item
$F_{I_{1}}^{u_{1}}(\Lambda)$ is a molecule in $u_{2}^{I_{1}}\times
u_{3}^{I_{1}}\times u_{4}^{I_{1}}$ or the empty set for every
integer $I_{1}$.
\item
$F_{I_{2}}^{u_{2}}(\Lambda)$ is a molecule in $u_{1}\times
u_{3}^{I_{2}}\times u_{4}^{I_{2}}$ or the empty set for every
integer $I_{2}$.
\item
$F_{I_{3}}^{u_{3}}(\Lambda)$ is a molecule in $u_{1}\times
u_{2}\times u_{4}^{I_{3}}$ or the empty set for every integer
$I_{3}$.
\item
$F_{I_{4}}^{u_{4}}(\Lambda)$ is a molecule in $u_{1}\times
u_{2}\times u_{3}$ or the empty set for every integer $I_{4}$.

\end{enumerate}
\end{definition}

Note. The reason that a subcomplex satisfying the above conditions
is said to be pairwise molecular is made clear in the following
Proposition \ref{pairwise_reason}.

One of the main result in this chapter is as follows.

\begin{theorem}\label{4glb_main}
A subcomplex in $\ubbb$ is a molecule if and only if it is
pairwise molecular.
\end{theorem}

\begin{prop}
Let $\Lambda$ be a  subcomplex of $\ubbb$. Then
$F_{I_{s}}^{u_{s}}[F_{I_{t}}^{u_{t}}(\Lambda)]
=F_{I_{t}}^{u_{t}}[F_{I_{s}}^{u_{s}}(\Lambda)]$ for all $s$ and
$t$ with $1\leq s, t\leq 4$ and $s\neq t$.
\end{prop}
\begin{proof}
Let $\lambda$ be an atom in $\ubbb$. It is evident from the
definition that $F_{I_{s}}^{u_{s}}[F_{I_{t}}^{u_{t}}(\lambda)]
=F_{I_{t}}^{u_{t}}[F_{I_{s}}^{u_{s}}(\lambda)]$ for all $s$ and
$t$ with $1\leq s, t\leq 4$ and $s\neq t$. Since
$F_{I_{s}}^{u_{s}}$ and $F_{I_{t}}^{u_{t}}$ preserve unions, we
can see that $F_{I_{s}}^{u_{s}}[F_{I_{t}}^{u_{t}}(\Lambda)]
=F_{I_{t}}^{u_{t}}[F_{I_{s}}^{u_{s}}(\Lambda)]$, as required.

\end{proof}

For every finite non-empty subcomplex $\Lambda$ of $\ubbb$, the
subcomplex $F_{I_{s}}^{u_{s}}[F_{I_{t}}^{u_{t}}(\Lambda)]
=F_{I_{t}}^{u_{t}}[F_{I_{s}}^{u_{s}}(\Lambda)]$ is denoted by
$F_{I_{s},I_{t}}^{u_{s},u_{t}}(\Lambda)$.

\begin{prop}\label{pairwise_reason}
Let $\Lambda$ be a subcomplex of $u_{1}\times u_{2}\times
u_{3}\times u_{4}$. Then $\Lambda$ is pairwise molecular if and
only if the following conditions hold.
\begin{enumerate}
\item
There are no distinct maximal atoms $u_{1}[i_{1},\alpha_{1}]\times
u_{2}[i_{2},\alpha_{2}]\times u_{3}[i_{3},\alpha_{3}]\times
u_{4}[i_{4},\alpha_{4}]$ and $u_{1}[i_{1}',\alpha_{1}']\times
u_{2}[i_{2}',\alpha_{2}']\times u_{3}[i_{3}',\alpha_{3}']\times
u_{4}[i_{4}',\alpha_{4}']$ in $\Lambda$ such that $i_{1}\leq
i_{1}'$, $i_{2}\leq i_{2}'$,  $i_{3}\leq i_{3}'$ and $i_{4}\leq
i_{4}'$.
\item
If $F_{I_{1}}^{u_{1}}(\Lambda)\neq\emptyset$, then
$F_{I_{1}}^{u_{1}}(\Lambda)$ satisfies condition 1 for pairwise
molecular subcomplexes in $u_{2}^{I_{1}}\times u_{3}^{I_{1}}\times
u_{4}^{I_{1}}$.

\item
If $F_{I_{2}}^{u_{2}}(\Lambda)\neq\emptyset$, then
$F_{I_{2}}^{u_{2}}(\Lambda)$ satisfies condition 1 for pairwise
molecular subcomplexes in $u_{1}\times u_{3}^{I_{2}}\times
u_{4}^{I_{2}}$.

\item
If $F_{I_{3}}^{u_{}}(\Lambda)\neq\emptyset$, then
$F_{I_{3}}^{u_{3}}(\Lambda)$ satisfies condition 1 for pairwise
molecular subcomplexes in $u_{1}\times u_{2}\times u_{4}^{I_{3}}$.

\item
If $F_{I_{4}}^{u_{4}}(\Lambda)\neq\emptyset$, then
$F_{I_{4}}^{u_{4}}(\Lambda)$ satisfies condition 1 for pairwise
molecular subcomplexes in $u_{1}\times u_{2}\times u_{3}$.

\item
$F_{I_{1},I_{2}}^{u_{1},u_{2}}(\Lambda)$ is a molecule in
$u_{3}^{I_{1}+I_{2}}\times u_{4}^{I_{1}+I_{2}}$ or the empty set
for every pair of integers $I_{1}$ and $I_{2}$.
\item
$F_{I_{1},I_{3}}^{u_{1},u_{3}}(\Lambda)$ is a molecule in
$u_{2}^{I_{1}}\times u_{4}^{I_{1}+I_{3}}$ or the empty set for
every pair of integers $I_{1}$ and  $I_{3}$.
\item
$F_{I_{1},I_{4}}^{u_{1},u_{4}}(\Lambda)$ is a molecule in
$u_{2}^{I_{1}}\times u_{3}^{I_{1}}$ or the empty set for every
pair of integers $I_{1}$ and  $I_{4}$.
\item
$F_{I_{2},I_{3}}^{u_{2},u_{3}}(\Lambda)$ is a molecule in
$u_{1}\times u_{4}^{I_{2}+I_{3}}$ or the empty set for every pair
of integers $I_{2}$ and $I_{3}$.
\item
$F_{I_{2},I_{4}}^{u_{2},u_{4}}(\Lambda)$ is a molecule in
$u_{1}\times u_{3}^{I_{2}}$ or the empty set for every pair of
integers $I_{2}$ and $I_{4}$.
\item
$F_{I_{3},I_{4}}^{u_{3},u_{4}}(\Lambda)$ is a molecule in
$u_{1}\times u_{2}$ or the empty set for every pair of integers
$I_{3}$ and $I_{4}$.
\end{enumerate}
\end{prop}

\begin{proof}

Suppose that $\Lambda$ is pairwise molecular. Then
$F_{I_{s}}^{u_{s}}(\Lambda)$ is a molecule or the empty set for
every $s$. It follows from definition of
$F_{I_{s},I_{t}}^{u_{s},u_{t}}$ and Theorem \ref{3glb_main} that
conditions 1 to 11 hold.

Conversely, suppose that $\Lambda$ satisfies condition 1 to 11. By
the definition of $F_{I_{s},I_{t}}^{u_{s},u_{t}}$ , we can see
that $F_{I_{1}}^{u_{1}}[F_{I_{4}}^{u_{4}}(\Lambda)]$,
$F_{I_{2}}^{u_{2}}[F_{I_{4}}^{u_{4}}(\Lambda)]$ and
$F_{I_{3}}^{u_{3}}[F_{I_{4}}^{u_{4}}(\Lambda)]$ are molecules in
the corresponding (twisted) products of two globes or the empty
set. Since $F_{I_{4}}^{u_{4}}(\Lambda)$ satisfies condition 1 for
pairwise molecular subcomplexes, it follows from Theorem
\ref{3glb_main} that $F_{I_{4}}^{u_{4}}(\Lambda)$ is a molecule in
$u_{1}\times u_{2}\times u_{3}$ or the empty set. Similarly, we
can prove that $F_{I_{4}}^{u_{4}}(\Lambda)$,
$F_{I_{4}}^{u_{4}}(\Lambda)$ and $F_{I_{4}}^{u_{4}}(\Lambda)$ are
molecules in the corresponding (twisted) product of three globes
or the empty set. This shows that $\Lambda$ is pairwise molecular,
as required.

This completes the proof.
\end{proof}

We end this section by a proposition which is used later in this
chapter.

\begin{prop}\label{eq4}
Let $\Lambda$ and $\Lambda'$ be subcomplexes of $\ubbb$ satisfying
condition 1 for pairwise molecular subcomplexes. If
$F_{I_{s}}^{u_{s}}(\Lambda)=F_{I_{s}}^{u_{s}}(\Lambda')$ for every
$s$ and every $I_{s}$ with $1\leq s\leq 4$, then
$\Lambda=\Lambda'$.
\end{prop}
\begin{proof}
It suffices to prove that $\Lambda$ and $\Lambda'$ consists of the
same maximal atoms.

Let $\uaaa{i}{\alpha}{}$ be a maximal atom in $\Lambda$. It is
easy to see that $\ua{1}\times\ua{2}\times\ua{3}$ is a maximal
atom in $F_{i_{4}}^{u_{4}}(\Lambda)=F_{i_{4}}^{u_{4}}(\Lambda')$.
Thus $\Lambda'$ has a maximal atom $u_{1}[i_{1},\alpha_{1}]\times
u_{2}[i_{2},\alpha_{2}]\times u_{3}[i_{3},\alpha_{3}]\times
u_{4}[i_{4}',\alpha_{4}']$ with $i_{4}'\geq i_{4}$. Since
$u_{1}[i_{1},\alpha_{1}]\times u_{2}[i_{2},\alpha_{2}]\times
u_{3}[i_{3},\alpha_{3}]\not\subset
F_{i_{4}+1}^{u_{4}}(\Lambda)=F_{i_{4}+1}^{u_{4}}(\Lambda')$, we
have $i_{4}'=i_{4}$. One can similarly get a maximal atom
$u_{1}[i_{1},\alpha_{1}]\times u_{2}[i_{2},\alpha_{2}]\times
u_{3}[i_{3},\alpha_{3}']\times u_{4}[i_{4},\alpha_{4}]$ of
$\Lambda'$. It follows from condition 1 for pairwise molecular
subcomplexes that $\alpha_{4}'=\alpha_{4}$ and
$\alpha_{3}'=\alpha_{3}$. This shows that $\uaaa{i}{\alpha}{}$ is
a maximal atom in $\Lambda'$.

Symmetrically, we can see that every maximal atom in $\Lambda'$ is
a maximal atom in $\Lambda$.

This completes the proof that $\Lambda=\Lambda'$.
\end{proof}

\section{Molecules Are Pairwise Molecular}
In this section, we prove that molecules in $\ubbb$ are pairwise
molecular.  The argument is different from that in section 2 of
chapter 2. We show that $F_{I_{s}}^{u_{s}}$ sends molecules to
molecules or the empty set for every value of $s$ without
introducing $g_{I_{s}}^{u_{s}}$ (see section 2 of chapter 2).

We first show that molecules satisfy condition 1 for pairwise
molecular subcomplexes.

\begin{prop}\label{mlc14}
In a molecule of $\ubbb$,  there are no distinct maximal atoms
$\uaaa{i}{\alpha}{}$ and $\uaaa{i}{\alpha}{'}$ such that
$i_{s}\leq i_{s}'$ for all  $1\leq s\leq 4$.
\end{prop}
\begin{proof}
Suppose otherwise that there are maximal atoms
$\lambda=\uaaa{i}{\alpha}{}$ and $\lambda'=\uaaa{i}{\alpha}{'}$ in
the molecule with $\lambda\neq\lambda'$ such that $i_{s}\leq
i_{s}'$ for all $1\leq s\leq 4$. By decomposing the given
molecule, one can get composite of molecules
$\Lambda\#_{n}\Lambda'$ or $\Lambda'\#_{n}\Lambda$ such that
$\lambda\subset\Lambda$, $\lambda\not\subset\Lambda'$,
$\lambda'\subset\Lambda'$ and $\lambda'\not\subset\Lambda$. We may
assume that the given molecule is decomposed into
$\Lambda\#_{n}\Lambda'$. We now consider two cases, as follows.

1. Suppose that $n\geq i_{1}+i_{2}+i_{3}+i_{4}$. Then, by Lemma
\ref{maximal_dpg}, we have $\lambda\subset
d_{n}^{+}\Lambda=d_{n}^{-}\Lambda'\subset\Lambda'$. This is a
contradiction.

2. Suppose that $n<i_{1}+i_{2}+i_{3}+i_{4}$. Consider the
(natural) homomorphism $F:\ubbb\to u_{i_{1}}\times u_{i_{2}}\times
u_{i_{3}}\times u_{i_{4}}$. Since $\Lambda\#_{n}\Lambda'$ exists,
we know that
$F(\Lambda\#_{n}\Lambda')=F(\Lambda)\#_{n}F(\Lambda')$ exists. On
the other hand, one can see that $u_{i_{1}}\times u_{i_{2}}\times
u_{i_{3}}\times u_{i_{4}}\subset F(\lambda)\cap F(\lambda')\subset
F(\Lambda)\cap F(\Lambda')$. Therefore $\dim[F(\Lambda)\cap
F(\Lambda')]\geq i_{1}+i_{2}+i_{3}+i_{4}>n$. Since
$\dim[d_{n}^{+}F(\Lambda)]\leq n$, one gets
$d_{n}^{+}F(\Lambda)\neq F(\Lambda)\cap F(\Lambda')$. This
contradicts that
$F(\Lambda\#_{n}\Lambda')=F(\Lambda)\#_{n}F(\Lambda')$ exists.

This completes the proof.

\end{proof}

We next show that $F_{I_{s}}^{u_{s}}$ sends molecules to molecules
or the empty set for every value of $s$. The arguments for
different values of $s$ are similar. We only give the proof for
$s=2$.

\begin{lemma} \label{imatom4}
Let $\lambda=\uaaa{i}{\alpha}{}$ be an atom in the
$\omega$-complex $\ubbb$ and $\Lambda,\Lambda'\in\M(\ubbb)$. Then
\begin{enumerate}

\item
$F_{I_{2}}^{u_{2}}(\lambda)\in\A(u_{1}\times u_{3}^{I_{2}}\times
u_{4}^{I_{2}})\cup\{\emptyset\}$;

\item
If  $\Lambda\#_{n}\Lambda'$ is defined, then
$F_{I_{2}}^{u_{2}}(\Lambda\#_{n}\Lambda')=F_{I_{2}}^{u_{2}}(\Lambda)\cup
F_{I_{2}}^{u_{2}}(\Lambda')$;

\item
$F_{I_{2}}^{u_{2}}(\Lambda)\neq\emptyset$ if and only if there is
a maximal atom $\uaaa{i}{\alpha}{}$ in $\Lambda$ such that
$i_{2}\geq I_{2}$;

\item
$F_{I_{2}}^{u_{2}}(\dpg \lambda)=
\begin{cases}
d_{p-I_{2}}^{\gamma}F_{I_{2}}^{u_{2}}(\lambda) & \text{when $p\geq
I_{2}$ and $i_{2}\geq I_{2}$,}\\ \emptyset                &
\text{when $p<I_{2}$ or $i_{2}<I_{2}$;}
\end{cases}
$

\end{enumerate}
\end{lemma}
\begin{proof}
The proof of the first three conditions is  a trivial verification
from the definition of $F_{I_{2}}^{u_{2}}$. We now verify
condition 4.

If $p<I_{2}$ or $i_{2}<I_{2}$, then it is evident that
$F_{I_{2}}^{u_{2}}(\dpg \lambda)=\emptyset$ by the definition of
$F_{I_{2}}^{u_{2}}$.

Now, suppose that $p\geq I_{2}$ and $i_{2}\geq I_{2}$.  Then
$F_{I_{2}}^{u_{2}}(\lambda)=u_{1}[i_{1},\alpha_{1}]\times
u_{3}^{I_{2}}[i_{3},\alpha_{3}]\times
u_{4}^{I_{2}}[i_{4},\alpha_{4}]$. By proposition \ref{dtatom4},
the set of all maximal atoms in $\dpg \lambda$ consists of all
$\uaaa{l}{\sigma}{}$ with $l_{s}\leq i_{s}$ for all $1\leq s\leq
4$ such that $l_{1}+l_{2}+l_{3}+l_{4}=p$, where the signs
$\sigma_{s}$ ($s=1, 2, 3, 4$) are determined as follows:
\begin{enumerate}
\item
If $l_{1}=i_{1}$, then $\sigma_{1}=\alpha_{1}$; if $l_{1}<i_{1}$,
then $\sigma_{1}=\gamma$.
\item
If $l_{2}=i_{2}$, then $\sigma_{2}=\alpha_{2}$; if $l_{2}<i_{2}$,
then $\sigma_{2}=(-)^{l_{1}}\gamma$.
\item
If $l_{3}=i_{3}$, then $\sigma_{3}=\alpha_{3}$; if $l_{3}<i_{3}$,
then $\sigma_{3}=(-)^{l_{1}+l_{2}}\gamma$.
\item
If $l_{4}=i_{4}$, then $\sigma_{4}=\alpha_{4}$; if $l_{4}<i_{4}$,
then $\sigma_{4}=(-)^{l_{1}+l_{2}+l_{3}}\gamma$.
\end{enumerate}
From this description and the formation of
$d_{p-I_{2}}^{\gamma}(u_{1}[i_{1},\alpha_{1}]\times
u_{3}^{I_{2}}[i_{3},\alpha_{3}]\times
u_{4}^{I_{2}}[i_{4},\alpha_{4}])$ in $u_{1}\times
u_{3}^{I_{2}}\times u_{4}^{I_{2}}$ (Theorem \ref{dt}), it is easy
to see that $F_{I_{2}}^{u_{2}}(\dpg \lambda)=
d_{p-I_{2}}^{\gamma}F_{I_{2}}^{u_{2}}(\lambda)$, as required.
\end{proof}

\begin{lemma} \label{i2geqn4}
Let $\Lambda$ be a molecule in $\ubbb$. If $\Lambda$ is decomposed
into $\Lambda=\Lambda'\#_{n}\Lambda''$ and if
$F_{I_{2}}^{u_{2}}(\Lambda')\neq\emptyset$ and
$F_{I_{2}}^{u_{2}}(\Lambda')\neq\emptyset$, then $n\geq I_{2}$.
\end{lemma}
\begin{proof}
Let $f_{I_{2}}^{u_{2}}:\M(\ubbb)\to\M(u_{1}\times u_{I_{2}}\times
u_{3}\times u_{4})$ be the natural homomorphism of
$\omega$-categories sending every maximal atom
$\uaaa{i}{\alpha}{}$ in $\ubbb$ with $i_{2}<I_{2}$ to
$\uaaa{i}{\alpha}{}$, and sending every maximal atom
$\uaaa{i}{\alpha}{}$ in $\ubbb$ with $i_{2}\geq I_{2}$ to
$\uaa{i}{\alpha}{}{1}\times u_{I_{2}}\times
\uaa{i}{\alpha}{}{3}\times \uaa{i}{\alpha}{}{4}$. Then
$f_{I_{2}}^{u_{2}}(\Lambda)=f_{I_{2}}^{u_{2}}(\Lambda')\#_{n}
f_{I_{2}}^{u_{2}}(\Lambda'')$ is defined. Thus
$d_{n}^{+}f_{I_{2}}^{u_{2}}(\Lambda')
=d_{n}^{-}f_{I_{2}}^{u_{2}}(\Lambda'')
=f_{I_{2}}^{u_{2}}(\Lambda')\cap f_{I_{2}}^{u_{2}}(\Lambda'')$.
Since $F_{I_{2}}^{u_{2}}(\Lambda')\neq\emptyset$ and
$F_{I_{2}}^{u_{2}}(\Lambda'')\neq\emptyset$, we know that there
are maximal atoms $\lambda'=\uaaa{i}{\alpha}{'}$ and
$\lambda''=\uaaa{i}{\alpha}{''}$ in $\Lambda'$ and $\Lambda''$
respectively with $i_{2}'\geq I_{2}$ and $i_{2}''\geq I_{2}$ such
that $f_{I_{2}}^{u_{2}}(\lambda')$ and
$f_{I_{2}}^{u_{2}}(\lambda'')$ are maximal atoms in
$f_{I_{2}}^{u_{2}}(\Lambda')$ and $f_{I_{2}}^{u_{2}}(\Lambda'')$
respectively. If $f_{I_{2}}^{u_{2}}(\lambda')$ is not maximal in
$f_{I_{2}}^{u_{2}}(\Lambda)$, then it is easy to see that there is
a maximal atom in $f_{I_{2}}^{u_{2}}(\Lambda'')$ containing
$f_{I_{2}}^{u_{2}}(\lambda')$; it follows that $n\geq\dim
d_{n}^{-}f_{I_{2}}^{u_{2}}(\Lambda'')=\dim(f_{I_{2}}^{u_{2}}(\Lambda')
\cap f_{I_{2}}^{u_{2}}(\Lambda''))\geq \dim
f_{I_{2}}^{u_{2}}(\lambda')\geq I_{2}$. Similarly, if
$f_{I_{2}}^{u_{2}}(\lambda'')$ is not maximal in
$f_{I_{2}}^{u_{2}}(\Lambda)$, then $n\geq I_{2}$, as required. In
the following proof, we may assume that both
$f_{I_{2}}^{u_{2}}(\lambda')$ and $f_{I_{2}}^{u_{2}}(\lambda'')$
are maximal in $f_{I_{2}}^{u_{2}}(\Lambda)$. Now there are two
cases as follows.

1. Suppose that $f_{I_{2}}^{u_{2}}(\lambda')\cap
f_{I_{2}}^{u_{2}}(\lambda'')\neq\emptyset$. Then it is easy to see
that $n\geq\dim d_{n}^{+}f_{I_{2}}^{u_{2}}(\Lambda')
=\dim(f_{I_{2}}^{u_{2}}(\Lambda')\cap
f_{I_{2}}^{u_{2}}(\Lambda''))\geq I_{2}$, as required.

2. Suppose that $f_{I_{2}}^{u_{2}}(\lambda')\cap
f_{I_{2}}^{u_{2}}(\lambda'')=\emptyset$, then we must have
$i_{1}'=i_{1}''=0$, $i_{3}'=i_{3}''=0$ or $i_{4}'=i_{4}''=0$. We
may assume that $i_{1}'=i_{1}''=0$. Thus
$\alpha_{1}'=-\alpha_{1}''$. In this case, consider the natural
homomorphism $f_{0,I_{2}}^{u_{1},u_{2}}:
\M(\ubbb)\to\M(u_{0}\times u_{I_{2}}\times u_{3}\times u_{4})$. It
is easy to see that $f_{0,I_{2}}^{u_{1},u_{2}}(\lambda')$ and
$f_{0,I_{2}}^{u_{1},u_{2}}(\lambda'')$ are maximal in
$f_{0,I_{2}}^{u_{1},u_{2}}(\Lambda)$. It follows that $i_{3}'\neq
i_{3}''$ and $i_{4}'\neq i_{4}''$ and hence
$f_{0,I_{2}}^{u_{1},u_{2}}(\Lambda')\cap
f_{0,I_{2}}^{u_{1},u_{2}}(\Lambda'')\neq\emptyset$. Since
$f_{0,I_{2}}^{u_{1},u_{2}}(\Lambda)
=f_{0,I_{2}}^{u_{1},u_{2}}(\Lambda')\#_{n}
f_{0,I_{2}}^{u_{1},u_{2}}(\Lambda'')$ we can see that $n\geq \dim
d_{n}^{+}f_{0,I_{2}}^{u_{1},u_{2}}(\Lambda')
=\dim(f_{0,I_{2}}^{u_{1},u_{2}}(\Lambda')\cap
f_{0,I_{2}}^{u_{1},u_{2}}(\Lambda''))\geq I_{2}$, as required.

This completes the proof.
\end{proof}

\begin{prop}
Let $F_{I_{2}}^{u_{2}}:\M(\ubbb)\to \C(u_{1}\times
u_{3}^{I_{2}}\times u_{4}^{I_{2}})$ be the map as above. Then
\begin{enumerate}
\item
$F_{I_{2}}^{u_{2}}(\M(\ubbb))\subset \M(u_{1}\times
u_{3}^{I_{2}}\times u_{4}^{I_{2}})\cup\{\emptyset\}$;
\item
For every molecule $\Lambda$ in $\ubbb$, we have $$
F_{I_{2}}^{u_{2}}(\dpg\Lambda)=
\begin{cases}
d_{p-I_{2}}^{\gamma}F_{I_{2}}^{u_{2}}(\Lambda) & \text{when $p\geq
I_{2}$ and
                                 $F_{I_{2}}^{u_{2}}(\Lambda)\neq\emptyset$},\\
\emptyset                & \text{when $p<I_{2}$ or
                            $F_{I_{2}}^{u_{2}}(\Lambda)=\emptyset$}.
\end{cases}
$$

\item
If $\Lambda\#_{n} \Lambda'$ is defined, then
$$F_{I_{2}}^{u_{2}}(\Lambda\#_{n} \Lambda')=
\begin{cases}
F_{I_{2}}^{u_{2}}(\Lambda)\#_{n-I_{2}} F_{I_{2}}^{u_{2}}(\Lambda')
& \text{when $F_{I_{2}}^{u_{2}}(\Lambda)\neq\emptyset$ and
$F_{I_{2}}^{u_{2}}(\Lambda')\neq\emptyset$,}\\
F_{I_{2}}^{u_{2}}(\Lambda')      & \text{when
$F_{I_{2}}^{u_{2}}(\Lambda)=\emptyset$,} \\
F_{I_{2}}^{u_{2}}(\Lambda) & \text{when
$F_{I_{2}}^{u_{2}}(\Lambda')=\emptyset$.}
\end{cases}
$$
\end{enumerate}
\end{prop}
\begin{proof}
We are going to prove the first two conditions by induction and
then prove the third condition.

By Lemma \ref{imatom4}, it is evident that the first two
conditions hold when $\Lambda$ is an atom.

Now suppose that $q>1$ and the first two conditions hold for
molecules which can be written as a composite of less than $q$
atoms. Suppose also that $\Lambda$ is a molecule which can be
written as a composite of $q$ atoms. Since $q>1$, we have a proper
decomposition $\Lambda=\Lambda'\#_{n} \Lambda''$ such that
$\Lambda'$ and $\Lambda''$ are molecules which can be written as
composites of less than $q$ atoms. According to the induction
hypothesis, we know that the first two conditions hold for
$\Lambda'$ and $\Lambda''$. We must show that the first two
conditions in the proposition hold for $\Lambda$. There are two
cases, as follows.

1. Suppose that $F_{I_{2}}^{u_{2}}(\Lambda')=\emptyset$ or
$F_{I_{2}}^{u_{2}}(\Lambda'')=\emptyset$. We may assume that
$F_{I_{2}}^{u_{2}}(\Lambda')=\emptyset$. In this case, we have
$F_{I_{2}}^{u_{2}}(\Lambda)=F_{I_{2}}^{u_{2}}(\Lambda'')$. Thus
$F_{I_{2}}^{u_{2}}(\Lambda)\in\M(u_{1}\times u_{3}^{I_{2}}\times
u_{4}^{I_{2}})\cup\{\emptyset\}$ as required by the first
condition. Moreover, if $p\neq n$, then $$
\begin{array}{rl}
 &F_{I_{2}}^{u_{2}}(d_{p}^{\gamma} \Lambda)\\
=&F_{I_{2}}^{u_{2}}(d_{p}^{\gamma}\Lambda'\#_{n} d_{p}^{\gamma}
\Lambda'')\\ =&F_{I_{2}}^{u_{2}}(d_{p}^{\gamma} \Lambda')\cup
F_{I_{2}}^{u_{2}}(d_{p}^{\gamma} \Lambda'')\\ =&
\begin{cases}
d_{p-I_{2}}^{\gamma}F_{I_{2}}^{u_{2}}(\Lambda'') & \text{when
$p\geq I_{2}$ and $F_{I_{2}}^{u_{2}}(\Lambda'')\neq\emptyset$,}\\
\emptyset & \text{when $F_{I_{2}}^{u_{2}}(\Lambda'')=\emptyset$ or
$p<I_{2}$},
\end{cases}\\
=&
\begin{cases}
d_{p-I_{2}}^{\gamma}F_{I_{2}}^{u_{2}}(\Lambda) & \text{when $p\geq
I_{2}$ and $F_{I_{2}}^{u_{2}}(\Lambda)\neq\emptyset$,}\\ \emptyset
& \text{when $F_{I_{2}}^{u_{2}}(\Lambda)=\emptyset$ or $p<I_{2}$},
\end{cases}
\end{array}
$$ as required by the second condition. Suppose that $p=n\geq
I_{2}$. Then $F_{I_{2}}^{u_{2}}(d_{p}^{+}\Lambda')=\emptyset$. So
$F_{I_{2}}^{u_{2}}(d_{p}^{-}\Lambda'')=\emptyset$. Hence, by the
hypothesis, one gets $F_{I_{2}}^{u_{2}}(\Lambda'')=\emptyset$.
Therefore $F_{I_{2}}^{u_{2}}(\Lambda)=\emptyset$ and
$F_{I_{2}}^{u_{2}}(d_{p}^{\gamma}\Lambda)=\emptyset$, as required
by the second condition.

2. Suppose that $F_{I_{2}}^{u_{2}}(\Lambda')\neq\emptyset$ and
$F_{I_{2}}^{u_{2}}(\Lambda'')\neq\emptyset$.  By Lemma
\ref{i2geqn4}, we have $n\geq I_{2}$.  Since
$F_{I_{2}}^{u_{2}}(\Lambda)=F_{I_{2}}^{u_{2}}(\Lambda')\cup
F_{I_{2}}^{u_{2}}(\Lambda'')$ and $$
\begin{array}{rl}
 &d_{n-I_{2}}^{+}F_{I_{2}}^{u_{2}}(\Lambda')\\
=&F_{I_{2}}^{u_{2}}(d_{n}^{+} \Lambda')\\
=&F_{I_{2}}^{u_{2}}(d_{n}^{-} \Lambda'')\\
=&d_{n-I_{2}}^{-}F_{I_{2}}^{u_{2}}(\Lambda''),
\end{array}
$$ we can see that $F_{I_{2}}^{u_{2}}(\Lambda')\#_{n-I_{2}}
F_{I_{2}}^{u_{2}}(\Lambda'')$ is defined,  and
$F_{I_{2}}^{u_{2}}(\Lambda)=F_{I_{2}}^{u_{2}}(\Lambda')\#_{n-I_{2}}
F_{I_{2}}^{u_{2}}(\Lambda'')$. So $F_{I_{2}}^{u_{2}}(\Lambda)$ is
a molecule, as required by the first condition. We now verify that
$\Lambda$ satisfy the second condition. If $p<I_{2}$, then
$F_{I_{2}}^{u_{2}}(\dpg\Lambda)=\emptyset$, as required. If
$p=n\geq I_{2}$, then $$
\begin{array}{rl}
 &F_{I_{2}}^{u_{2}}(d_{p}^{-}\Lambda)\\
=&F_{I_{2}}^{u_{2}}(d_{p}^{-}\Lambda')\\
=&d_{p-I_{2}}^{-}F_{I_{2}}^{u_{2}}(\Lambda')\\
=&d_{p-I_{2}}^{-}F_{I_{2}}^{u_{2}}(\Lambda);
\end{array}
$$ and similarly we have $F_{I_{2}}^{u_{2}}(d_{p}^{+}\Lambda)
=d_{p-I_{2}}^{+}F_{I_{2}}^{u_{2}}(\Lambda)$. If $I_{2}\leq p<n$,
then $$
\begin{array}{rl}
 &F_{I_{2}}^{u_{2}}(\dpg \Lambda)\\
=&F_{I_{2}}^{u_{2}}(\dpg \Lambda')\\
=&d_{p-I_{2}}^{\gamma}F_{I_{2}}^{u_{2}}(\Lambda')\\
=&d_{p-I_{2}}^{\gamma}F_{I_{2}}^{u_{2}}(\Lambda).
\end{array}
$$ If $p\geq I_{2}$ and $p>n$, then $$
\begin{array}{rl}
 &F_{I_{2}}^{u_{2}}(\dpg \Lambda)\\
=&F_{I_{2}}^{u_{2}}(\dpg \Lambda'\#_{n}\dpg \Lambda'')\\
=&F_{I_{2}}^{u_{2}}(\dpg \Lambda')\cup F_{I_{2}}^{u_{2}}(\dpg
\Lambda'')\\ =&d_{p-I_{2}}^{\gamma}F_{I_{2}}^{u_{2}}(\Lambda')\cup
d_{p-I_{2}}^{\gamma}F_{I_{2}}^{u_{2}}(\Lambda'')
\end{array}
$$ and $$
\begin{array}{rl}
&d_{n-I_{2}}^{+}d_{p-I_{2}}^{\gamma}F_{I_{2}}^{u_{2}}(\Lambda')\\
=&d_{n-I_{2}}^{+}F_{I_{2}}^{u_{2}}(\Lambda')\\
=&d_{n-I_{2}}^{-}F_{I_{2}}^{u_{2}}(\Lambda'')\\
=&d_{n-I_{2}}^{-}d_{p-I_{2}}^{\gamma}F_{I_{2}}^{u_{2}}(\Lambda''),
\end{array}
$$ thus
$d_{p-I_{2}}^{\gamma}F_{I_{2}}^{u_{2}}(\Lambda')\#_{n-I_{2}}
d_{p-I_{2}}^{\gamma}F_{I_{2}}^{u_{2}}(\Lambda'')$ is defined and
$$
\begin{array}{rl}
 & F_{I_{2}}^{u_{2}}(\dpg \Lambda)\\
=&d_{p-I_{2}}^{\gamma}F_{I_{2}}^{u_{2}}(\Lambda')\#_{n-I_{2}}
d_{p-I_{2}}^{\gamma}F_{I_{2}}^{u_{2}}(\Lambda'')\\
=&d_{p-I_{2}}^{\gamma}[F_{I_{2}}^{u_{2}}(\Lambda')\#_{n-I_{2}}F_{I_{2}}^{u_{2}}(\Lambda'')]\\
=&d_{p-I_{2}}^{\gamma}F_{I_{2}}^{u_{2}}(\Lambda).
\end{array}
$$ Therefore $\Lambda$ satisfies the second condition.

Finally,  condition 3 can be easily verified by using condition 2
and the fact that $F_{I_{2}}^{u_{2}}$ preserves unions.

This complete the proof.

\end{proof}

In particular, we have that $F_{I_{2}}^{u_{2}}$ sends molecules to
molecules in $u_{1}\times u_{3}^{I_{2}}\times u_{4}^{I_{2}}$ or
the empty set.

We can similarly prove that $F_{I_{s}}^{u_{s}}$ sends molecules to
molecules in the corresponding (twisted) product of three infinite
dimensional globes or the empty set for every value of $s$. Thus
we have proved the main theorem in this section.

\begin{theorem}
Molecules in $\ubbb$ are pairwise molecular.
\end{theorem}

We finish this section by a property of molecules in $\ubbb$. It
can be proved from the results later in this chapter. But the
proof here is also interesting.

\begin{prop}
Let $\Lambda$ be a molecule in $\ubbb$. Let
$\lambda_{1}=\uaa{i}{\alpha}{}{1}\times \uaa{i}{\alpha}{}{2}\times
\uaa{i}{\alpha}{}{3}\times \uaa{i}{\alpha}{}{4}$ and
$\lambda_{2}=u_{1}[i_{1},-\alpha_{1}]\times
\uaa{i}{\alpha}{}{2}\times \uaa{i}{\alpha}{}{3}\times
\uaa{i}{\alpha}{}{4}$. If $\lambda_{1}\subset\Lambda$ and
$\lambda_{2}\subset\Lambda$, then there is an atom
$\lambda'=\uaaa{i}{\alpha}{'}$ in $\Lambda$ with $i'>i_{1}=i_{2}$
such that $\lambda'\supset\lambda_{1}$ and
$\lambda'\supset\lambda_{2}$.

\end{prop}
\begin{proof}
If $\Lambda$ is an atom, then the required property holds
automatically.

Suppose that the required property hold for every molecule which
can be written as a composite of less than $q$ atoms. Suppose also
that $\Lambda$ can be written as a composite of $q$ atoms.

It is easy to see that there is a composite of molecules
$\Lambda_{1}\#_{n}\Lambda_{2}$ with $\lambda_{1}\subset
\Lambda_{1}$ and $\lambda_{2}\subset\Lambda_{2}$ such that both
$\Lambda_{1}$ and $\Lambda_{2}$ can be written as composites of
less than $q$ atoms. By applying the natural homomorphism
$f_{i_{1}}^{u_{1}}:\M(\ubbb)\to\M((u_{i_{1}}\times u_{2}\times
u_{3}\times u_{4})$ of $\omega$-categories, one gets $$
f_{i_{1}}^{u_{1}}(\Lambda_{1}\#_{n}\Lambda_{2})=
f_{i_{1}}^{u_{1}}(\Lambda_{1})\#_{n}
f_{i_{1}}^{u_{1}}(\Lambda_{2})). $$ This implies that $$
f_{i_{1}}^{u_{1}}(\Lambda_{1}\cap\Lambda_{2})=
f_{i_{1}}^{u_{1}}(\Lambda_{1})\cap
f_{i_{1}}^{u_{1}}(\Lambda_{2})). $$ Since
$f_{i_{1}}^{u_{1}}(\lambda_{1})=f_{i_{1}}^{u_{1}}(\lambda_{2})\subset
f_{i_{1}}^{u_{1}}(\Lambda_{1})\cap
f_{i_{1}}^{u_{1}}(\Lambda_{2}))$, we have
$f_{i_{1}}^{u_{1}}(\lambda_{1})=f_{i_{1}}^{u_{1}}(\lambda_{2})\subset
f_{i_{1}}^{u_{1}}(\Lambda_{1}\cap \Lambda_{2})$. It follows that
$\lambda_{3}=\uaa{i}{\beta}{}{1}\times \uaa{i}{\alpha}{}{2}\times
\uaa{i}{\alpha}{}{3}\times \uaa{i}{\alpha}{}{4}\subset\Lambda_{1}
\cap\Lambda_{2}$ for some sign $\beta$. Now if
$\beta_{1}=-\alpha$, then one can get an atom in $\Lambda_{1}$ as
required by applying the induction hypothesis on $\lambda_{1}$ and
$\lambda_{3}$ in $\Lambda_{1}$; if $\beta_{1}=\alpha$, then one
can get an atom in $\Lambda_{2}$ as required by applying the
induction hypothesis on $\lambda_{2}$ and $\lambda_{3}$ in
$\Lambda_{2}$.

This completes the proof.
\end{proof}

\section{Properties of Pairwise Molecular Subcomplexes}
In this section, we study some properties of pairwise molecular
subcomplexes in $\ubbb$. In the next section, we are going to
prove that some of these conditions are sufficient for a
subcomplex of $\ubbb$ to be pairwise molecular.

\begin{lemma} \label{mlc34}
Let $\Lambda$ be a pairwise molecular subcomplex of $\ubbb$.  Let
$\uaaa{i}{\alpha}{}$ and $\uaaa{i}{\alpha}{'}$ be a pair of
distinct maximal atoms in $\Lambda$. If $i_{s}=i_{s}'$ and
$\alpha_{s}=-\alpha_{s}'$ for some $1\leq s\leq 4$,  then for
every $t$ with $t\neq s$, there is a maximal atom
$\uaaa{l}{\sigma}{}$ in $\Lambda$ such that $l_{s}>i_{s}=i_{s}'$,
$l_{t}\geq\min\{i_{t},i_{t}'\}$ and
$\uaa{l}{\sigma}{}{r}\supset\uaa{i}{\alpha}{}{r}\cap\uaa{i}{\alpha}{'}{r}$
for all $r\in\{1,2,3,4\}\setminus\{s,t\}$.
\end{lemma}

\begin{proof}
Let $\lambda=\uaaa{i}{\alpha}{}$ and
$\lambda'=\uaaa{i}{\alpha}{'}$. The arguments for various cases
are similar, we give the proof for $s=1$ and $t=4$. Let
$I_{4}=\min\{i_{4},i_{4}'\}$. Since $\Lambda$ is pairwise
molecular, we can see that $F_{I_{4}}^{u_{4}}(\Lambda)$ is a
molecule in $u_{1}\times u_{2}\times u_{3}$. Since
$F_{I_{4}}^{u_{4}}(\lambda)=\uaa{i}{\alpha}{}{1}\times
\uaa{i}{\alpha}{}{2}\times \uaa{i}{\alpha}{}{3}$ and
$F_{I_{4}}^{u_{4}}(\lambda)=\uaa{i}{\alpha}{'}{1}\times
\uaa{i}{\alpha}{'}{2}\times \uaa{i}{\alpha}{'}{3}$, it follows
easily from Lemma \ref{lm3} that there is a maximal atom
$\uaa{l}{\sigma}{}{1}\times \uaa{l}{\sigma}{}{2}\times
\uaa{l}{\sigma}{}{3}$ such that $l_{1}>i_{1}=i_{1}'$,
$\uaa{l}{\sigma}{}{2}\supset
\uaa{i}{\alpha}{}{2}\cap\uaa{i}{\alpha}{'}{2}$ and
$\uaa{l}{\sigma}{}{3}\supset \uaa{i}{\alpha}{}{3}\cap
\uaa{i}{\alpha}{'}{3}$. Therefore $\Lambda$ has a maximal atom
$\uaaa{l}{\sigma}{}$ such that $l_{1}>i_{1}=i_{1}'$,
$\uaa{l}{\sigma}{}{2}\supset
\uaa{i}{\alpha}{}{2}\cap\uaa{i}{\alpha}{'}{2}$,
$\uaa{l}{\sigma}{}{3}\supset \uaa{i}{\alpha}{}{3}\cap
\uaa{i}{\alpha}{'}{3}$ and $l_{4}\geq I_{4}$, as required.

This completes the proof.
\end{proof}

The following definition of {\em adjacency} for a pair of maximal
atoms in a subcomplex of $\ubbb$ is inspired by Propositions
\ref{adjacent_eq} and \ref{sign_4_4}.

\begin{definition}
Let $\Lambda$ be a subcomplex of $\ubbb$. Let $1\leq s<t\leq 4$. A
pair of maximal atoms $\uaaa{i}{\alpha}{}$ and
$\uaaa{i}{\alpha}{'}$ are {\em $(s,t)$-adjacent} if
$\max\{i_{s},i_{s}'\}+\max\{i_{t},i_{t}'\}>
\max\{i_{s}+i_{t},i_{s}'+i_{t}'\}$ and if there is no maximal atom
$\uaaa{j}{\beta}{}$ with $j_{r}\geq\min\{i_{r},i_{r}'\}$ for all
$1\leq r\leq 4$ such that
$\min\{i_{s},j_{s}\}+\min\{i_{t},j_{t}\}>
\min\{i_{s},i_{s}'\}+\min\{i_{t},i_{t}'\}$ and
$\min\{i_{s}',j_{s}\}+\min\{i_{t}',j_{t}\}>
\min\{i_{s},i_{s}'\}+\min\{i_{t},i_{t}'\}$. A pair of maximal
atoms $\uaaa{i}{\alpha}{}$ and $\uaaa{i}{\alpha}{'}$ are {\em
adjacent} if they are $(s,t)$-adjacent for all $1\leq s<t\leq 4$
such that $\max\{i_{s},i_{s}'\}+\max\{i_{t},i_{t}'\}>
\max\{i_{s}+i_{t},i_{s}'+i_{t}'\}$.
\end{definition}

\begin{example}
Suppose that $\lambda=u_{1}[5,\alpha_{1}]\times
u_{2}[0,\alpha_{2}]\times u_{3}[1,\alpha_{3}]\times
u_{4}[1,\alpha_{4}]$ and $\mu =u_{1}[0,\beta_{1}]\times
u_{2}[5,\beta_{2}]\times u_{3}[1,\beta_{3}]\times
u_{4}[2,\beta_{4}]$ are a pair of maximal atoms in a subcomplex.
If $\Lambda$ has a maximal atom $\nu
=u_{1}[1,\varepsilon_{1}]\times u_{2}[1,\varepsilon_{2}]\times
u_{3}[2,\varepsilon_{3}]\times u_{4}[1,\varepsilon_{4}]$, then
$\lambda$ and $\mu$ are not $(1,2)$-adjacent.
\end{example}

The following proposition shows that the definition of adjacency
is in consistent with that in Chapter \ref{ch2}.

\begin{prop} \label{adjacent_eq4}
Let $\Lambda$ be a subcomplex of $\ubbb$ satisfying condition 1
for pairwise molecular subcomplexes. Then a pair of distinct
maximal atoms $\uaaa{i}{\alpha}{}$ and $\uaaa{i}{\alpha}{'}$ are
adjacent if and only if there is no maximal atom
$\uaaa{j}{\beta}{}$ with $j_{r}\geq\min\{i_{r},i_{r}'\}$ for all
$1\leq r\leq 4$ such that
$$\sum_{r=1}^{4}\min\{i_{r},j_{r}\}>\sum_{r=1}^{4}\min\{i_{r},i_{r}'\}$$
and
$$\sum_{r=1}^{4}\min\{i_{r}',j_{r}\}>\sum_{r=1}^{4}\min\{i_{r},i_{r}'\}.$$
\end{prop}
\begin{proof}
The proof is a straightforward verification from the definition of
adjacency.
\end{proof}

The concept of projection maximal can be defined in the similar
way as that in chapter \ref{ch2}.

Let $\Lambda$ be a subcomplex of $\ubbb$ satisfying condition 1
for pairwise molecular subcomplex. Let $I_{r}$ be a fixed
non-negative integer. A maximal atom $\uaaa{i}{\alpha}{}$ in
$\Lambda$ is {\em $(u_{r},I_{r})$-projection maximal} if
$i_{r}\geq I_{r}$ and there is no maximal atom
$\uaaa{i}{\alpha}{'}$ such that  $I_{r}\leq i_{r}'<i_{r}$ and
$i_{s}'\geq i_{s}$ for all $s\in\{1,2,3,4\}\setminus \{r\}$.

Evidently, if a maximal atom $\lambda$ in $\Lambda$ is
$(u_{r},I_{r})$-projection maximal, then
$F_{I_{r}}^{u_{r}}(\lambda)$ is maximal in
$F_{I_{r}}^{u_{r}}(\Lambda)$. Conversely, for every maximal atom
$\mu$ in $F_{I_{r}}^{u_{r}}(\Lambda)$, there is a maximal atom
$\mu'$ in $\Lambda$ such that $F_{I_{r}}^{u_{r}}(\mu')=\mu$. The
following proposition implies that $\mu'$ is actually
$(u_{r},I_{r})$-projection maximal.

\begin{prop}\label{projection_reason4}
Let $\Lambda$ be a pairwise molecular subcomplex of $\ubbb$ and
$\lambda$ be a maximal atom in $\Lambda$. Let $1\leq r\leq 4$.
Then $\lambda$ is $(u_{r},I_{r})$-projection maximal if and only
if $F_{I_{r}}^{u_{r}}(\lambda)$ is maximal in
$F_{I_{r}}^{u_{r}}(\Lambda)$.
\end{prop}
\begin{proof}
The proof is similar to that in Proposition
\ref{projection_reason}.
\end{proof}

\begin{lemma}\label{adlowest42}
Let $\Lambda$ be a pairwise molecular subcomplex in $\ubbb$. Let
$\lambda=\uaaa{i}{\alpha}{}$ and $\mu=\uaaa{j}{\beta}{}$ be a pair
of $(s,t)$-adjacent maximal atoms in $\Lambda$ for some $1\leq
s<t\leq 4$. If $i_{s}>j_{s}$ and $i_{t}<j_{t}$, then for
$r\in\{1,2,3,4\}\setminus\{s,t\}$ there is a pair of
$(s,t)$-adjacent and $(u_{r},I_{r})$-projection maximal atoms
$\uaaa{i}{\alpha}{'}$ and $\uaaa{j}{\beta}{'}$  such that
$\uaa{j}{\beta}{'}{s}=\uaa{j}{\beta}{}{s}$,
$\uaa{i}{\alpha}{'}{t}=\uaa{i}{\alpha}{}{t}$ and
$\min\{i_{r}',j_{r}'\}=\min\{i_{r},j_{r}\}$, where
$I_{r}=\min\{i_{r},j_{r}\}$. Moreover, for
$\bar{r}\in\{1,2,3,4\}\setminus\{r,s,t\}$, if $s<\bar{r}<t$, then
$\min\{i_{\bar{r}}',j_{\bar{r}}'\}=\min\{i_{\bar{r}},j_{\bar{r}}\}$.
\end{lemma}
\begin{proof}
Let $\lambda'=\uaaa{i}{\alpha}{'}$ be the
$(u_{r},I_{r})$-projection maximal atom such that $i_{s}'\geq
i_{s}$, $i_{t}'\geq i_{t}$ and $i_{\bar{r}}'\geq i_{\bar{r}}$. Let
$\mu'=\uaaa{j}{\beta}{'}$ be the $(u_{r},I_{r})$-projection
maximal atom such that $j_{s}'\geq j_{s}$, $j_{t}'\geq j_{t}$ and
$j_{\bar{r}}'\geq j_{\bar{r}}'$. It follows easily from Lemma
\ref{mlc34} and the adjacency of $\lambda$ and $\mu$ that
$\lambda'$ and $\mu'$ are $(s,t)$-adjacent and that
$\uaa{j}{\beta}{'}{s}=\uaa{j}{\beta}{}{s}$,
$\uaa{i}{\alpha}{'}{t}=\uaa{i}{\alpha}{}{t}$,
$\min\{i_{r}',j_{r}'\}= I_{r}$ and
$\min\{i_{\bar{r}}',j_{\bar{r}}'\}\geq
\min\{i_{\bar{r}},j_{\bar{r}}\}$.

Moreover, if $s\leq\bar{r}\leq t$,  we show that
$\min\{i_{\bar{r}}',j_{\bar{r}}'\}=\min\{i_{\bar{r}},j_{\bar{r}}\}$
by contradiction. The arguments for various choices of $r$, $s$
and $t$ are similar. We give the proof for $s=1$, $t=4$ and $r=3$.
In this case, we have $\bar{r}=2$. Suppose otherwise that
$\min\{i_{2}',j_{2}'\}>\min\{i_{2},j_{2}\}$. Then
$F_{I_{3}}^{u_{3}}(\lambda')$ and $F_{I_{3}}^{u_{3}}(\mu')$ are
maximal atoms in $F_{I_{3}}^{u_{3}}(\Lambda)$. Note that
$F_{I_{3}}^{u_{3}}(\lambda')=
\uaa{i}{\alpha}{'}{1}\times\uaa{i}{\alpha}{'}{2}\times
u_{4}^{I_{3}}[i_{4},\alpha_{4}]$ and $F_{I_{3}}^{u_{3}}(\mu')=
\uaa{j}{\beta}{}{1}\times\uaa{j}{\beta}{'}{2}\times
u_{4}^{I_{3}}[j_{4}',\beta_{4}']$. Since
$F_{I_{3}}^{u_{3}}(\Lambda)$ is a molecule in $u_{1}\times
u_{2}\times u_{4}^{I_{3}}$, it follows from condition \ref{m4} in
Theorem \ref{pairwise_eq} that there is a maximal atom
$\uaa{l}{\sigma}{}{1}\times\uaa{l}{\sigma}{}{2}\times
u_{4}^{I_{3}}[l_{4},\sigma_{4}]$ in $F_{I_{3}}^{u_{3}}(\Lambda)$
such that $l_{1}>j_{1}$,
$l_{2}=\min\{i_{2}',j_{2}'\}-1\geq\min\{i_{2},j_{2}\}$ and
$l_{4}>i_{4}$. Thus there is a maximal atom
$\nu=\uaaa{l}{\sigma}{}$ in $\Lambda$ such that $l_{3}\geq I_{3}$.
This contradicts the assumption that $\lambda$ and $\mu$ are
$(1,4)$-adjacent.

This completes the proof.
\end{proof}

\begin{prop}\label{mlc44}
Let $\Lambda$ be a pairwise molecular subcomplex of $\ubbb$. Let
$\lambda=\uaaa{i}{\alpha}{}$ and $\lambda'=\uaaa{i}{\alpha}{'}$ be
a pair of distinct maximal atoms in $\Lambda$. If $1\leq s<r<t\leq
4$ and $\lambda$ and $\lambda'$ are $(s,t)$-adjacent, and if
$i_{s}>i_{s}'$, $\min\{i_{r},i_{r}'\}>0$ and $i_{t}'>i_{t}$,  then
$\Lambda$ has a maximal atom $u_{1}[l_{1},\sigma_{1}]\times
u_{2}[l_{2},\sigma_{2}]\times u_{3}[l_{3},\sigma_{3}]\times
u_{4}[l_{4},\sigma_{4}]$ such that $l_{s}>i_{1}'$,
$l_{r}=\min\{i_{r},i_{r}'\}-1$, $l_{t}>i_{t}$ and $l_{\bar{r}}\geq
\min\{i_{\bar{r}},i_{\bar{r}}'\}$, where
$\bar{r}\in\{1,2,3,4\}\setminus\{r,s,t\}$.
\end{prop}
\begin{proof}
The arguments for various cases are similar. We only give the
proof for the case $s=1$, $r=3$ and $t=4$.

Let $I_{3}=\min\{i_{3},i_{3}'\}$. According to Lemma
\ref{adlowest42}, we may assume that $\lambda$ and $\lambda'$ are
$(u_{3},I_{3})$-projection maximal so that
$F_{I_{3}}^{u_{3}}(\lambda)$ and $F_{I_{3}}^{u_{3}}(\lambda')$ are
maximal atoms in $F_{I_{3}}^{u_{3}}(\Lambda)$.

Now $F_{I_{3}}^{u_{3}}(\lambda)=\uaa{i}{\alpha}{}{1}\times
\uaa{i}{\alpha}{}{2}\times u_{4}^{I_{3}}[i_{4},\alpha_{4}]$ and
$F_{I_{3}}^{u_{3}}(\lambda')=\uaa{i}{\alpha}{'}{1}\times
\uaa{i}{\alpha}{'}{2}\times u_{4}^{I_{3}}[i_{4}',\alpha_{4}']$. It
is easy to see that $F_{I_{3}}^{u_{3}}(\lambda)$ and
$F_{I_{3}}^{u_{3}}(\lambda')$ are adjacent in
$F_{I_{3}}^{u_{3}}(\Lambda)$. According to condition \ref{m4} in
Theorem \ref{pairwise_eq}, there is a maximal atom
$\uaa{l}{\sigma}{}{1}\times \uaa{l}{\sigma}{}{2}\times
u_{4}^{I_{3}}[l_{4},\sigma_{4}]$ such that $l_{1}>i_{1}'$,
$l_{2}=\min\{i_{2},i_{2}'\}-1$ and $l_{4}>i_{4}$. This implies
that there is a maximal atom $\uaaa{l}{\sigma}{}$ in $\Lambda$ as
required.

This completes the proof.

\end{proof}

We also need to extend the concept of projection maximal to
maximal atoms in $\ubbb$ with respect to two factors, as follows.

Let $1\leq s<t\leq 4$. Let $I_{s}$ and $I_{t}$ be  fixed
non-negative integers. A maximal atom $\uaaa{i}{\alpha}{}$ in
$\Lambda$ is {\em $(u_{s},u_{t};I_{s},I_{t})$-projection maximal}
if $i_{s}\geq I_{s}$ and $i_{t}\geq I_{t}$, and if there is no
maximal atom $\uaaa{i}{\alpha}{'}$ such that $i_{s}'\geq I_{s}$
and $i_{t}'\geq I_{t}$ and such that $i_{r}'\geq i_{r}$ for all
$r\in\{1,2,3,4\}\setminus\{s,t\}$ and $i_{r}'>i_{r}$ for some
$r\in\{1,2,3,4\}\setminus\{s,t\}$.

Evidently, if a maximal atom $\lambda$ in $\Lambda$ is
$(u_{s},u_{t};I_{s},I_{t})$-projection maximal, then
$F_{I_{s},I_{t}}^{u_{s},u_{t}}(\lambda)$ is maximal in
$F_{I_{s},I_{t}}^{u_{s},u_{t}}(\Lambda)$. Conversely, for every
maximal atom $\mu$ in $F_{I_{s},I_{t}}^{u_{s},u_{t}}(\Lambda)$,
there is a maximal atom $\mu'$ in $\Lambda$ such that
$F_{I_{s},I_{t}}^{u_{s},u_{t}}(\mu')=\mu$. The following
Proposition implies that $\mu'$ is actually
$(u_{s},u_{t};I_{s},I_{t})$-projection maximal.

\begin{prop}
Let $\Lambda$ be a pairwise molecular subcomplex of $\ubbb$ and
$\lambda$ be a maximal atom in $\Lambda$. Let $1\leq s<t\leq 4$.
Then $\lambda$ is $(u_{s},u_{t};I_{s},I_{t})$-projection maximal
if and only if $F_{I_{s},I_{t}}^{u_{s},u_{t}}(\lambda)$ is maximal
in $F_{I_{s},I_{t}}^{u_{s},u_{t}}(\Lambda)$.
\end{prop}
\begin{proof}
The argument is similar to that in Proposition
\ref{projection_reason4}.
\end{proof}

\begin{lemma}\label{adlowest4}
Let $\Lambda$ be a pairwise molecular subcomplex in $\ubbb$. Let
$\lambda=\uaaa{i}{\alpha}{}$ and $\mu=\uaaa{j}{\beta}{}$ be a pair
of $(s,t)$-adjacent maximal atoms in $\Lambda$ for some $1\leq
s<t\leq 4$. If $i_{s}>j_{s}$ and $i_{t}<j_{t}$, then there is a
pair of $(s,t)$-adjacent and
$(u_{\bar{s}},u_{\bar{t}};I_{\bar{s}},I_{\bar{t}})$-projection
maximal atoms $\uaaa{i}{\alpha}{'}$ and $\uaaa{j}{\beta}{'}$  such
that $\uaa{j}{\beta}{'}{s}=\uaa{j}{\beta}{}{s}$,
$\uaa{i}{\alpha}{'}{t}=\uaa{i}{\alpha}{}{t}$ and
$\min\{i_{\bar{s}}',j_{\bar{s}}'\}\geq I_{\bar{s}}$ and
$\min\{i_{\bar{t}}',j_{\bar{t}}'\}\geq I_{\bar{t}}$, where
$\bar{s}$ and $\bar{t}$ are distinct elements in
$\{1,2,3,4\}\setminus\{s,t\}$ and $I_{\bar{s}}=\min\{i_{\bar{s}},
j_{\bar{s}}\}$ and $I_{\bar{t}}=\min\{i_{\bar{t}}, j_{\bar{t}}\}$.
Moreover, if $s<\bar{s}<t$, then
$\min\{i_{\bar{s}}',j_{\bar{s}}'\}=I_{\bar{s}}$; if $s<\bar{t}<t$,
then $\min\{i_{\bar{t}}',j_{\bar{t}}'\}=I_{\bar{t}}$.
\end{lemma}
\begin{proof}
Let $\lambda'=\uaaa{i}{\alpha}{'}$ be the
$(u_{\bar{s}},u_{\bar{t}};I_{\bar{s}},I_{\bar{t}})$-projection
maximal atom such that $i_{s}'\geq i_{s}$ and $i_{t}'\geq i_{t}$.
Let $\mu'=\uaaa{j}{\beta}{'}$ be the
$(u_{\bar{s}},u_{\bar{t}};I_{\bar{s}},I_{\bar{t}})$-projection
maximal atom such that $j_{s}'\geq j_{s}$ and $j_{t}'\geq j_{t}$.
It follows easily from Lemma \ref{mlc34} and the adjacency of
$\lambda$ and $\mu$ that
$\uaa{j}{\beta}{'}{s}=\uaa{j}{\beta}{}{s}$,
$\uaa{i}{\alpha}{'}{t}=\uaa{i}{\alpha}{}{t}$,
$\min\{i_{\bar{s}}',j_{\bar{s}}'\}\geq I_{\bar{s}}$,
$\min\{i_{\bar{t}}',j_{\bar{t}}'\}\geq I_{\bar{t}}$ and $\lambda'$
is adjacent to $\mu'$, as required.

The second part follows easily from the adjacency of $\lambda$ and
$\mu$ and Lemma \ref{mlc44}.

This completes the proof.
\end{proof}

\begin{prop} \label{mlc24}
Let $\Lambda$ be a pairwise molecular subcomplex in $\ubbb$. Then
the following sign conditions hold.

Sign conditions: for a pair of adjacent maximal atoms
$\uaaa{i}{\alpha}{}$ and $\uaaa{i}{\alpha}{'}$ in $\Lambda$, let
$l_{r}=\min\{i_{r},i_{r}'\}$ for $1\leq r\leq 4$.
\begin{enumerate}
\item
If $\lambda$ and $\mu$ are $(1,2)$-adjacent, and if
$l_{1}=i_{1}<i_{1}'$ and $l_{2}=i_{2}'<i_{2}$, then
$\alpha_{2}'=-(-)^{l_{1}}\alpha_{1}$;
\item
If $\lambda$ and $\mu$ are $(1,3)$-adjacent, and if
$l_{1}=i_{1}<i_{1}'$ and $l_{3}=i_{3}'<i_{3}$,  then
$\alpha_{3}'=-(-)^{l_{1}+l_{2}}\alpha_{1}$;
\item\label{4glbcase}
If $\lambda$ and $\mu$ are $(1,4)$-adjacent, and if
$l_{1}=i_{1}<i_{1}'$ and $l_{4}=i_{4}'<i_{4}$,  then
$\alpha_{4}'=-(-)^{l_{1}+l_{2}+l_{3}}\alpha_{1}$;
\item
If $\lambda$ and $\mu$ are $(2,3)$-adjacent, and if
$l_{2}=i_{2}<i_{2}'$ and $l_{3}=i_{3}'<i_{3}$,  then
$\alpha_{3}'=-(-)^{l_{2}}\alpha_{2}$;
\item
If $\lambda$ and $\mu$ are $(2,4)$-adjacent, and if
$l_{2}=i_{2}<i_{2}'$ and $l_{4}=i_{4}'<i_{4}$,  then
$\alpha_{4}'=-(-)^{l_{2}+l_{3}}\alpha_{2}$;
\item
If $\lambda$ and $\mu$ are $(3,4)$-adjacent, and if
$l_{3}=i_{3}<i_{3}'$ and $l_{4}=i_{4}'<i_{4}$,  then
$\alpha_{4}'=-(-)^{l_{3}}\alpha_{3}$.
\end{enumerate}
\end{prop}
\begin{proof}
The arguments for the above cases are similar. We only give the
proof for case \ref{4glbcase}.

Let $\lambda=\uaaa{i}{\alpha}{}$ and
$\lambda'=\uaaa{i}{\alpha}{'}$ be a pair of $(1,4)$-adjacent
maximal atoms in $\Lambda$. Suppose that $i_{1}'>i_{1}$ and
$i_{4}'<i_{4}$. Let $\min\{i_{2},i_{2}'\}=l_{2}$
$\min\{i_{3},i_{3}'\}=l_{3}$. We must prove
$\alpha_{4}'=-(-)^{l_{2}+l_{3}}\alpha_{1}$. According to Lemma
\ref{adlowest4}, we may assume that $\lambda$ and $\lambda'$ are
$(u_{2},u_{3};I_{2},I_{3})$-projection maximal atoms.

It is evident that
$F_{l_{2},l_{3}}^{u_{2},u_{3}}(\lambda)=\uaa{i}{\alpha}{}{1}\times
u_{4}^{l_{2}+l_{3}}[i_{4},\alpha_{4}]$ and
$F_{l_{2},l_{3}}^{u_{2},u_{3}}(\lambda')=\uaa{i}{\alpha}{'}{1}\times
u_{4}^{l_{2}+l_{3}}[i_{4}',\alpha_{4}']$, and they are maximal
atoms in the molecule $F_{l_{2},l_{3}}^{u_{2},u_{3}}(\Lambda)$.
Moreover, by the adjacency of $\lambda$ and $\lambda'$, we can see
that $F_{l_{2},l_{3}}^{u_{2},u_{3}}(\lambda)$ and
$F_{l_{2},l_{3}}^{u_{2},u_{3}}(\lambda')$ are adjacent maximal
atoms in $F_{l_{2},l_{3}}^{u_{2},u_{3}}(\Lambda)$. According to
the formation of molecules in $u\times w^{l_{2}+l_{3}}$, we get
$\alpha_{4}'=-(-)^{l_{2}+l_{3}}\alpha_{1}$, as required.

This completes the proof.
\end{proof}

\section{An Alternative Description for Pairwise Molecular
Subcomplexes}

In this section, we give an alternative description for pairwise
molecular subcomplexes of $\ubbb$, as follows.

\begin{theorem}\label{main4}
Let $\Lambda$ be a subcomplex of $\ubbb$. Then $\Lambda$ is
pairwise molecular if and only if
\begin{enumerate}
\item
There are no distinct maximal atoms $\uaaa{i}{\alpha}{}$ and
$\uaaa{i}{\alpha}{'}$ in $\Lambda$ such that $i_{1}\leq i_{1}'$,
$i_{2}\leq i_{2}'$,  $i_{3}\leq i_{3}'$ and $i_{4}\leq i_{4}'$.
\item
Sign conditions: for a pair of maximal atoms $\uaaa{i}{\alpha}{}$
and $\uaaa{i}{\alpha}{'}$ in $\Lambda$, let
$l_{r}=\min\{i_{r},i_{r}'\}$ for $1\leq r\leq 4$.
\begin{enumerate}
\item
If $\lambda$ and $\mu$ are $(1,2)$-adjacent, and if
$l_{1}=i_{1}<i_{1}'$ and $l_{2}=i_{2}'<i_{2}$, then
$\alpha_{2}'=-(-)^{l_{1}}\alpha_{1}$;
\item
If $\lambda$ and $\mu$ are $(1,3)$-adjacent, and if
$l_{1}=i_{1}<i_{1}'$ and $l_{3}=i_{3}'<i_{3}$,  then
$\alpha_{3}'=-(-)^{l_{1}+l_{2}}\alpha_{1}$;
\item
If $\lambda$ and $\mu$ are $(1,4)$-adjacent, and if
$l_{1}=i_{1}<i_{1}'$ and $l_{4}=i_{4}'<i_{4}$,  then
$\alpha_{4}'=-(-)^{l_{1}+l_{2}+l_{3}}\alpha_{1}$;
\item
If $\lambda$ and $\mu$ are $(2,3)$-adjacent, and if
$l_{2}=i_{2}<i_{2}'$ and $l_{3}=i_{3}'<i_{3}$,  then
$\alpha_{3}'=-(-)^{l_{2}}\alpha_{2}$;
\item
If $\lambda$ and $\mu$ are $(2,4)$-adjacent, and if
$l_{2}=i_{2}<i_{2}'$ and $l_{4}=i_{4}'<i_{4}$,  then
$\alpha_{4}'=-(-)^{l_{2}+l_{3}}\alpha_{2}$;
\item
If $\lambda$ and $\mu$ are $(3,4)$-adjacent, and if
$l_{3}=i_{3}<i_{3}'$ and $l_{4}=i_{4}'<i_{4}$,  then
$\alpha_{4}'=-(-)^{l_{3}}\alpha_{3}$.
\end{enumerate}

\item\label{main4_3}
Let $\lambda=\uaaa{i}{\alpha}{}$ and
$\lambda'=\uaaa{i}{\alpha}{'}$ be a pair of distinct maximal atoms
in $\Lambda$. If $i_{s}=i_{s}'$ and $\alpha_{s}=-\alpha_{s}'$ for
some $1\leq s\leq 4$, then for every
$t\in\{1,2,3,4\}\setminus\{s\}$, there is a maximal atom
$\uaaa{l}{\sigma}{}$ in $\Lambda$ such that $l_{s}>i_{s}=i_{s}'$,
$l_{t}\geq\min\{i_{t},i_{t}'\}$ and
$\uaa{l}{\sigma}{}{r}\supset\uaa{i}{\alpha}{}{r}\cap\uaa{i}{\alpha}{'}{r}$
for all $r\in\{1,2,3,4\}\setminus\{s,t\}$.

\item \label{main4_4}
Let $\lambda=\uaaa{i}{\alpha}{}$ and
$\lambda'=\uaaa{i}{\alpha}{'}$ be a pair of distinct maximal atoms
in $\Lambda$. For $1\leq s<r<t\leq 4$, if $\lambda$ and $\lambda'$
are $(s,t)$-adjacent, and if $i_{s}>i_{s}'$,
$\min\{i_{r},i_{r}'\}>0$ and $i_{t}'>i_{t}$,  then $\Lambda$ has a
maximal atom $u_{1}[l_{1},\sigma_{1}]\times
u_{2}[l_{2},\sigma_{2}]\times u_{3}[l_{3},\sigma_{3}]\times
u_{4}[l_{4},\sigma_{4}]$ such that $l_{s}>i_{1}'$,
$l_{r}=\min\{i_{r},i_{r}'\}-1$, $l_{t}>i_{t}$ and $l_{\bar{r}}\geq
\min\{i_{\bar{r}},i_{\bar{r}}'\}$ for
$\bar{r}\in\{1,2,3,4\}\setminus\{r,s,t\}$.

\end{enumerate}
\end{theorem}

\begin{note}\label{m4_note4}
In condition \ref{main4_4}, we have a similar relations for the
signs $\sigma_{r}$, $\alpha_{s}'$ and $\alpha_{t}$ as that in Note
\ref{m4_note}. For instance, if $s=1$, $r=2$ and $t=3$, then we
have $\sigma_{2}=-(-)^{i_{1}'}\alpha_{1}'$ and
$\alpha_{t}=-(-)^{l_{2}}\sigma_{2}$.
\end{note}

In the last section, we have proved that the four conditions in
this theorem are necessary for pairwise molecular subcomplexes. We
now prove the sufficiency. The proof is separated into several
lemmas.

\begin{lemma}\label{lm34}
Let $\Lambda$ be a subcomplex satisfying the four conditions in
Theorem \ref{main4}. Let $\lambda=\uaaa{i}{\alpha}{}$ and
$\mu=\uaaa{j}{\beta}{}$ be maximal atoms in $\Lambda$. If
$i_{s}=j_{s}$ and $\alpha_{s}=-\beta_{s}$ for some $1\leq s\leq
4$, then there is a maximal atom $\nu=\uaaa{l}{\sigma}{}$ in
$\Lambda$ with $l_{s}>i_{s}=j_{s}$ such that
$l_{t}\geq\min\{i_{t},j_{t}\}$ and
$\uaa{l}{\sigma}{}{t}\supset\uaa{i}{\alpha}{}{t}\cap\uaa{j}{\beta}{}{t}$
for all  $t\in\{1,2,3,4\}\setminus\{s\}$.
\end{lemma}
\begin{proof}

We first prove the lemma when $\lambda$ and $\mu$ are adjacent.
The arguments for various cases are similar. We only give the
proof for $s=1$, $i_{2}>j_{2}$ and $i_{3}<j_{3}$.

If $i_{4}=j_{4}$ and $\alpha_{4}=-\beta_{4}$, then we can get the
required $\mu$ simply by applying condition \ref{main4_3} in
Theorem \ref{main4}. We now suppose that $i_{4}\neq j_{4}$ or that
$i_{4}=j_{4}$ and $\alpha_{4}=\beta_{4}$. In this case, we have
$\uaa{i}{\alpha}{}{4}\cap\uaa{j}{\beta}{}{4}=\uaa{i}{\alpha}{}{4}$
or
$\uaa{i}{\alpha}{}{4}\cap\uaa{j}{\beta}{}{4}=\uaa{j}{\beta}{}{4}$.
We may assume that
$\uaa{i}{\alpha}{}{4}\cap\uaa{j}{\beta}{}{4}=\uaa{i}{\alpha}{}{4}$.
According to condition \ref{main4_3} in Theorem \ref{main4}, we
can get maximal atoms $\lambda'=\uaaa{i}{\alpha}{'}$ such that
$i_{1}'>i_{1}=j_{1}$, $i_{2}'\geq j_{2}$,
$\uaa{i}{\alpha}{}{3}\subset \uaa{i}{\alpha}{'}{3}$ and
$\uaa{i}{\alpha}{}{4}\subset \uaa{i}{\alpha}{'}{4}$. Similarly, we
can get  maximal atoms $\mu'=\uaaa{j}{\beta}{'}$ such that
$j_{1}'>i_{1}=j_{1}$, $j_{3}'\geq i_{3}$,
$\uaa{j}{\beta}{}{2}\subset \uaa{j}{\beta}{'}{2}$ and
$\uaa{i}{\alpha}{}{4}\subset \uaa{j}{\beta}{'}{4}$.

Now, suppose otherwise that $\Lambda$ does not have a maximal atom
$\nu=\uaaa{l}{\sigma}{}$ as required. Then $i_{2}'=j_{2}$ and
$\alpha_{2}'=-\beta_{2}$, and  $j_{3}'=i_{3}$ and
$\beta_{3}'=-\alpha_{3}$. By applying condition \ref{main4_3} in
Theorem \ref{main4} to $\mu$ and $\lambda'$, it is easy to see
that $i_{3}'=i_{3}$; thus $\alpha_{3}'=\alpha_{3}$. Similarly, we
have $j_{2}'=j_{2}$ and $\beta_{2}'=\beta_{2}$. By applying
condition \ref{main4_3} in Theorem \ref{main4} to $\lambda'$ and
$\mu'$, we get a maximal atom $\nu''=\uaaa{l}{\sigma}{''}$ such
that $l_{1}''\geq\min\{i_{1}',j_{1}'\}>i_{1}=j_{1}$,
$l_{2}''>i_{2}'=j_{2}'=j_{2}$, $l_{3}''\geq i_{3}'=j_{3}'=i_{3}$,
$\uaa{l}{\sigma}{''}{4}\supset\uaa{i}{\alpha}{}{4}$. By the
hypothesis, we have $l_{3}''=i_{3}$ and
$\sigma_{3}''=-\alpha_{3}$. By applying condition \ref{main4_3} in
Theorem \ref{main4} to $\lambda$ and $\nu''$, we can get a maximal
atom $\lambda''=\uaaa{i}{\alpha}{''}$ such that $i_{1}''\geq
i_{1}=j_{1}$, $i_{2}''\geq\min\{i_{2},l_{2}''\}>j_{2}$,
$i_{3}''>i_{3}$ and $i_{4}''\geq i_{4}$. This contradicts the
adjacency of $\lambda$ and $\mu$.

The other cases can be proved similarly.

Now we give the proof for the general case by induction.

1. Suppose that $\sum_{r=1}^{4}\min\{i_{r},j_{r}\}$ is maximal
among the non-negative integers
$\sum_{r=1}^{4}\min\{i'_{r},j'_{r}\}$ with $\uaaa{i}{\alpha}{'}$
and $\uaaa{j}{\beta}{'}$ running over all pairs of distinct
maximal atoms in $\Lambda$. It is evident that $\lambda$ is
adjacent to $\mu$, hence the lemma holds for $\lambda$ and $\mu$.

2. Suppose that $q\geq 0$ and the lemma holds for every pair of
distinct maximal atoms $\uaaa{i}{\alpha}{'}$ and
$\uaaa{j}{\beta}{'}$ with $\sum_{r=1}^{4}\min\{i'_{r},j'_{r}\}>q$.
Suppose also that $\sum_{r=1}^{4}\min\{i_{r},j_{r}\}=q$.

If $\lambda$ and $\mu$ are adjacent, then the lemma holds by the
first part of the proof.

Suppose that $\lambda$ and $\mu$ are not adjacent. According to
Theorem \ref{adjacent_eq4}, there is a maximal atom
$\nu'=\uaaa{l}{\sigma}{'}$ with $l_{r}'\geq\min\{i_{r},j_{r}\}$
for $1\leq r\leq 4$ such that
$$\sum_{r=1}^{4}\min\{i_{r},l_{r}\}>\sum_{r=1}^{4}\min\{i_{r},j_{r}\}$$
and
$$\sum_{r=1}^{4}\min\{j_{r},l_{r}\}>\sum_{r=1}^{4}\min\{i_{r},j_{r}\}.$$
By possibly multiple applying the above argument, condition
\ref{main4_3} in Theorem \ref{main4}, and the induction
hypothesis, we can get either a maximal atom
$\mu=\uaaa{l}{\sigma}{}$ as required, or a pair of maximal atoms
$\lambda'=\uaaa{i}{\alpha}{'}$ and $\mu'=\uaaa{j}{\beta}{'}$ with
$i_{s}'=i_{s}=j_{s}=j_{s}'$ and $\alpha_{s}'=-\beta_{s}'$ such
that $\min\{i_{t},j_{t}\}\leq \min\{i_{t}',j_{t}'\}$  and
$\uaa{i}{\alpha}{}{t}\cap\uaa{j}{\beta}{}{t}\subset
\uaa{i}{\alpha}{'}{t}\cap\uaa{j}{\beta}{'}{t}$ for all
$t\in\{1,2,3,4\}\setminus\{s\}$, and such that
$\sum_{r=1}^{4}\min\{i_{r},j_{r}\}<\sum_{r=1}^{4}\min\{i_{r}',j_{r}'\}$.
It follows from the induction hypothesis that the lemma holds for
$\lambda$ and $\mu$.

This completes the proof.

\end{proof}

Note that the proof in Proposition \ref{projection_reason4} uses
only the definition for projection maximal, condition 1 for
pairwise molecular subcomplexes and Lemma \ref{mlc34}. By
condition 1 and condition 3 in Theorem \ref{main4}, we have the
following two propositions.

\begin{prop}
Let $\Lambda$ be subcomplex of $\ubbb$ satisfying the four
conditions in Theorem \ref{main4} and $\lambda$ be a maximal atom
in $\Lambda$. Let $1\leq r\leq 4$. Then $\lambda$ is
$(u_{r},I_{r})$-projection maximal if and only if
$F_{I_{r}}^{u_{r}}(\lambda)$ is maximal in
$F_{I_{r}}^{u_{r}}(\Lambda)$.
\end{prop}

\begin{prop}
Let $\Lambda$ be a  subcomplex of $\ubbb$ satisfying the four
conditions in Theorem \ref{main4} and $\lambda$ be a maximal atom
in $\Lambda$. Let $1\leq s<t\leq 4$. Then $\lambda$ is
$(u_{s},u_{t};I_{s},I_{t})$-projection maximal if and only if
$F_{I_{s},I_{t}}^{u_{s},u_{t}}(\lambda)$ is maximal in
$F_{I_{s},I_{t}}^{u_{s},u_{t}}(\Lambda)$.
\end{prop}

\begin{prop}\label{lvp14}
Let $\Lambda$ be a subcomplex of $\ubbb$ satisfying the four
conditions in Theorem \ref{main4}. Let $I_{2}$ and $I_{3}$ be
fixed non-negative integers. Then
\begin{enumerate}
\item\label{lv14}
Every maximal atom $\uaaa{i}{\alpha}{}$ in $\Lambda$ with
$i_{2}=I_{2}$ and $i_{3}=I_{3}$ is
$(u_{2},u_{3};I_{2},I_{3})$-projection maximal.

\item\label{lv24}
For every maximal atom $\uaaa{i}{\alpha}{}$ with $i_{2}\geq I_{2}$
and $i_{3}\geq I_{3}$, there is a
$(u_{2},u_{3};I_{2},I_{3})$-projection maximal atom
$\uaaa{j}{\beta}{}$  such that $u_{1}[i_{1},\alpha_{1}]\subset
u_{1}[j_{1},\beta_{1}]$ and $u_{4}[i_{4},\alpha_{4}]\subset
u_{4}[j_{4},\beta_{4}]$.

\item\label{lv34}
All  $(u_{2},u_{3};I_{2},I_{3})$-projection maximal atoms , if
exist, can be listed as $\lambda^{(1)}$, $\cdots$, $\lambda^{(S)}$
with $\lambda^{(s)}=\uaaa{i}{\alpha}{^{(s)}}$ such that
$i_{1}^{(1)}\geq\cdots\geq i_{1}^{(S)}$ and
$i_{4}^{(1)}\leq\cdots\leq i_{4}^{(S)}$. Moreover, in this list,
if $1\leq s<S$, then $i_{1}^{(s)}=i_{1}^{(s+1)}$ if and only if
$i_{4}^{(s)}=i_{4}^{(s+1)}$; in this case, one must also have
$\alpha_{1}^{(s)}=\alpha_{1}^{(s+1)}$ and
$\alpha_{4}^{(s)}=\alpha_{4}^{(s+1)}$.

\item\label{lv44}
For two consecutive $(u_{2},u_{3};I_{2},I_{3})$-projection maximal
atoms $\lambda^{(s)}$ and $\lambda^{(s+1)}$ in the above list with
$i_{1}^{(s)}>i_{1}^{(s+1)}$, one has
$\min\{i_{2}^{(s)},i_{2}^{(s+1)}\}=I_{2}$ and
$\min\{i_{3}^{(s)},i_{3}^{(s+1)}\}=I_{3}$.
\end{enumerate}
\end{prop}
\begin{proof}

Condition \ref{lv14} follows from the definition of projection
maximal and condition 1 in Theorem \ref{main4}.

To prove condition \ref{lv24}, suppose that
$\lambda=\uaaa{i}{\alpha}{}$ is not
$(u_{2},u_{3};I_{2},I_{3})$-projection maximal. Then there is a
$(u_{2},u_{3};I_{2},I_{3})$-projection maximal atom
$\lambda'=\uaaa{i}{\alpha}{'}$ such that $i_{1}'\geq i_{1}$ and
$i_{4}'\geq i_{4}$. If
$\uaa{i}{\alpha}{'}{1}\supset\uaa{i}{\alpha}{}{1}$ and
$\uaa{i}{\alpha}{'}{4}\supset\uaa{i}{\alpha}{}{4}$, then
$\lambda'$ is as required. Suppose that
$\uaa{i}{\alpha}{'}{1}\not\supset\uaa{i}{\alpha}{}{1}$. Then
$i_{1}'=i_{1}$ and $\alpha_{1}'=-\alpha_{1}$. Moreover, we have
$i_{4}'>i_{4}$ by the definition of
$(u_{2},u_{3};I_{2},I_{3})$-projection maximal. Hence, by
condition 3 in Theorem \ref{main4}, there is a maximal atom
$\mu=\uaaa{l}{\sigma}{}$ in $\Lambda$ such that
$l_{1}>i_{1}=i_{1}'$, $l_{2}\geq I_{2}$, $l_{3}\geq I_{3}$ and
$\uaa{l}{\sigma}{}{4}\supset\uaa{i}{\alpha}{'}{4}$. If $\mu$ is
$(u_{2},u_{3};I_{2},I_{3})$-projection maximal, then $\mu$ can be
taken as $\uaaa{j}{\beta}{}$, as required. If $\mu$ is not
$(u_{2},u_{3};I_{2},I_{3})$-projection maximal, then one can get a
$(u_{2},u_{3};I_{2},I_{3})$-projection maximal atom $\mu$ as
required by repeating the above argument for $\lambda$.

To prove condition \ref{lv34}, by the definition of
$(u_{2},u_{3};I_{2},I_{3})$-projection maximal, it suffices to
show that $\Lambda$ does not have
$(u_{2},u_{3};I_{2},I_{3})$-projection maximal atoms
$\lambda=\uaaa{i}{\alpha}{}$ and $\mu=\uaaa{j}{\beta}{}$ such that
$i_{1}=j_{1}$ and $\alpha_{1}=-\beta_{1}$ or such that
$i_{4}=j_{4}$ and $\alpha_{4}=-\beta_{4}$. This follows easily
from condition 3 in Theorem \ref{main4}.

Finally, condition \ref{lv44} follows easily from condition 4 in
Theorem \ref{main4}.

This completes the proof.
\end{proof}

\begin{prop}\label{lvp24}
In a subcomplex $\Lambda$ of $\ubbb$ satisfying the four
conditions in Theorem \ref{main4}, let the
$(u_{2},u_{3};I_{2},I_{3})$-projection maximal atoms be listed as
$\lambda^{(1)}$, $\cdots$, $\lambda^{(S)}$ with
$\lambda^{(s)}=\uaaa{i}{\alpha}{^{(s)}}$ such that
$i_{1}^{(1)}\geq\cdots\geq i_{1}^{(S)}$ and
$i_{4}^{(1)}\leq\cdots\leq i_{4}^{(S)}$, taking $S=0$ if there are
none. Then
\begin{enumerate}
\item
$F_{I_{2},I_{3}}^{u_{2},u_{3}}(\Lambda)=\uaa{i}{\alpha}{^{(1)}}{1}\times
u_{4}^{I_{2}+I_{3}}[i_{4}^{(1)},\alpha_{4}^{(1)}]\cup\cdots\cup
\uaa{i}{\alpha}{^{(S)}}{1}\times
u_{4}^{I_{2}+I_{3}}[i_{4}^{(S)},\alpha_{4}^{(S)}]$.
\item
$F_{I_{2},I_{3}}^{u_{2},u_{3}}(\Lambda)$ is a molecule in
$u_{1}\times u_{4}^{I_{2}+I_{3}}$ or the empty set.
\end{enumerate}
\end{prop}
\begin{proof}
The first part  is a direct consequence of the definition for
$F_{I_{2},I_{3}}^{u_{2},u_{3}}$ and conditions 2 and 3 in the
above lemma. Note that a pair of consecutive
$(u_{2},u_{3};I_{2},I_{3})$-projection maximal atoms
$\lambda^{(s)}$ and $\lambda^{(s+1)}$ with $i_{1}^{s}>i_{1}^{s+1}$
in the above lemma are (1,4)-adjacent and
$\min\{i_{2}^{(r)},i_{2}^{(r+1)}\}=I_{2}$ and
$\min\{i_{3}^{(r)},i_{3}^{(r+1)}\}=I_{3}$. It follows from the
sign conditions that
$\alpha_{4}^{(r)}=-(-)^{i_{1}^{(r+1)}+I_{2}+I_{3}}\alpha_{1}^{(r+1)}$.
This implies that $F_{I_{2},I_{3}}^{u_{2},u_{3}}(\Lambda)$ is a
molecule in $u_{1}\times u_{4}^{I_{2}+I_{3}}$ or the empty set, as
required.

\end{proof}

\bigskip

We can similarly prove that
$F_{I_{s},I_{t}}^{u_{s},u_{t}}(\Lambda)$ is a molecule or the
empty set for every pair of $s$ and $t$ with $1\leq s<t\leq 4$ in
the corresponding products of three (twisted) infinite dimensional
globes.

\begin{lemma}\label{p_condition1}
Let $\Lambda$ be a subcomplex of $\ubbb$ satisfying the four
conditions in Theorem \ref{main4}. Then
$F_{I_{r}}^{u_{r}}(\Lambda)$ satisfies condition 1 for pairwise
molecular subcomplexes.
\end{lemma}
\begin{proof}
The argument for various choices of $r$ are similar. We only give
the proof for $r=4$.

Suppose otherwise that $F_{I_{4}}^{u_{4}}(\Lambda)$ does not
satisfy condition 1 for pairwise molecular subcomplexes. The there
is a pair of maximal atoms $\lambda=\uaa{i}{\alpha}{}{1}\times
\uaa{i}{\alpha}{}{2}\times\uaa{i}{\alpha}{}{3}$ and $\lambda'=
\uaa{i}{\alpha}{'}{1}\times
\uaa{i}{\alpha}{'}{2}\times\uaa{i}{\alpha}{'}{3}$ such that
$i_{s}\leq i_{s}'$ for all $1\leq s\leq 3$. Thus we have
$i_{t}=i_{t}'$ and $\alpha_{t}=-\alpha_{t}'$ for some $1\leq t\leq
3$. Hence there is a pair of $(u_{4},I_{4})$-projection maximal
atoms of the form $\mu=\uaaa{i}{\alpha}{}$ and
$\mu'=\uaaa{i}{\alpha}{'}$. By condition 1 in Theorem \ref{main4},
we evidently have $i_{4}>i_{4}'\geq I_{4}$. It follows from
condition 3 in Theorem \ref{main4} that there is a maximal atom
$\uaaa{l}{\sigma}{}$ such that $l_{s}\geq i_{s}$ for all $1\leq
s\leq 3$, $l_{4}\geq i_{4}'\geq I_{4}$ and $l_{t}>i_{t}$. This
contradicts that $\mu$ is $(u_{4},I_{4})$-projection maximal.

This completes the proof.

\end{proof}

By the comment after Proposition \ref{lvp24} and Lemma
\ref{p_condition1}, we have now completed the proof for Theorem
\ref{main4}.

\section{Sources and Targets of Pairwise Molecular Subcomplexes}

In this section, we study source and target operators
$d_{p}^{\gamma}$ on pairwise molecular subcomplexes in $\ubbb$. We
shall prove  that $d_{p}^{\gamma}\Lambda$ is pairwise molecular
for every pairwise molecular subcomplex $\Lambda$ of $\ubbb$.

Recall that $\dpg\Lambda$ is a union of interiors of atoms for
every subcomplex $\Lambda$ of $\ubbb$. Hence the maps
$F_{I_{r}}^{u_{r}}$ and $F_{I_{s},I_{t}}^{u_{s},u_{t}}$ are
defined on $\dpg\Lambda$.

We first show that $\dpg\Lambda$ is a subcomplex for a pairwise
molecular subcomplex $\Lambda$. The proof is separated into
several Lemmas.

\begin{lemma}\label{lm5.3F}
Let $\Lambda$ be a subcomplex of $\ubbb$ and
$\lambda=\uaaa{i}{\alpha}{}$ be a $p$-dimensional atom in
$\Lambda$ with $\Int\lambda\subset\dpg\Lambda$. If there is an
atom $\lambda'=\uaaa{i}{\alpha}{'}$ in $\Lambda$ with
$\lambda'\supset\lambda$ such that $i_{s}'>i_{s}$ for some $s$,
then $\alpha_{s}=(-)^{i_{1}+\cdots+i_{s-1}}\gamma$.
\end{lemma}
\begin{proof}
The proof is similar to that in \ref{lm5.3}.
\end{proof}

\begin{lemma}\label{Int_sign4}
Let $\Lambda$ be a pairwise molecular subcomplex of $\ubbb$. Let
$\lambda=\uaaa{i}{\alpha}{}$ be a $p$-dimensional atom  with
$\Int\lambda\subset\dpg\Lambda$. For $1\leq s\leq 4$, if there is
a maximal atom $\lambda'=\uaaa{i}{\alpha}{'}$ in $\Lambda$ with
$i_{s}'>i_{s}$ such that $i_{r}'\geq i_{r}$ for all
$r\in\{1,2,3,4\}$, then
$\alpha_{s}=(-)^{i_{1}+\cdots+i_{s-1}}\gamma$.
\end{lemma}

\begin{proof}

The arguments for various choices of $s$ are similar. We only give
the proof for $s=1$.

Suppose that there is a maximal atom
$\lambda'=\uaaa{i}{\alpha}{'}$ such that $i_{1}'>i_{1}$ and
$i_{r}'\geq i_{r}$ for all $r\in\{2,3,4\}$. If $\lambda'$ can be
chosen such that $\lambda\subset\lambda'$, then we have
$\alpha_{1}=\gamma$ by Lemma \ref{lm5.3F} , as required. In the
following proof, we may assume that $\lambda'$ cannot be chosen
such that $\lambda\subset\lambda'$ so that
$\uaa{i}{\alpha}{'}{t}=u_{t}[i_{t},-\alpha_{t}]$ for some
$t\in\{2,3,4\}$. Since $\Int\lambda\subset\Lambda$, there is a
maximal atom $\mu=\uaaa{j}{\beta}{}$ such that
$\lambda\subset\mu$. By the assumption on the choice of
$\lambda'$, we have $\uaa{j}{\beta}{}{1}=\uaa{i}{\alpha}{}{1}$. By
applying Lemma \ref{lm34}, we may further assume that
$j_{t}>i_{t}$.

1. Suppose that $\lambda'$ and $\mu$ can be chosen such that
$\min\{i_{r}',j_{r}\}>i_{r}$ for two value of $r\in\{2,3,4\}$.
According to condition \ref{main4_4} in Theorem \ref{main4}, we
may suppose that $\min\{i_{r}',j_{r}\}>i_{r}$ for $r=3$ and $r=4$,
and that $j_{2}>i_{2}$. Thus $i_{2}'=i_{2}$ and
$\alpha'=-\alpha_{2}$. By Lemma \ref{lm5.3F}, we have
$\alpha_{2}=(-)^{i}\gamma$, and hence
$\alpha_{2}'=-(-)^{i}\gamma$. According to the assumption, it is
easy to see that $\lambda'$ is $(1,2)$-adjacent to $\mu$. It
follows easily that $\alpha_{1}=\gamma$, as required.

2. Suppose that $\lambda'$ and $\mu$ cannot be chosen as in case
1. Suppose also that $\lambda'$ and $\mu$ can be chosen such that
$\min\{i_{r}',j_{r}\}>i_{r}$ for only one value of
$r\in\{2,3,4\}$. According to condition \ref{main4_4} in Theorem
\ref{main4}, we may suppose that $\min\{i_{4}',j_{4}\}>i_{4}$, and
that $j_{2}>i_{2}$ or $j_{3}>i_{3}$. We may further suppose that
$j_{2}>i_{2}$ and $i_{2}'=i_{2}$ and $\alpha_{2}'=-\alpha_{2}$.
Note that, in this case, $\lambda'$ is $(1,2)$-adjacent to $\mu$.
By an argument similar to that in case 1, we have
$\alpha_{1}=\gamma$, as required.

3. Suppose that $\lambda'$ and $\mu$ cannot be chosen as in case 1
and case 2. Then $\lambda'$ is adjacent to $\mu$. By an argument
similar to that in case 1, we have $\alpha_{1}=\gamma$, as
required.

This completes the proof.

\end{proof}

\begin{lemma}\label{include_p-14}
Let $\Lambda$ be a pairwise molecular subcomplex of $\ubbb$. Let
$\lambda=\uaaa{i}{\alpha}{}$ be a $p-1$ dimensional atom  such
that $\Int\lambda\subset\dpg\Lambda$. If there is a maximal atom
$\lambda'=\uaaa{i}{\alpha}{'}$ in $\Lambda$ with
$\lambda'\supset\lambda$ such that $i_{s}'>i_{s}$ and
$i_{t}'>i_{t}$ for some $1\leq s<t\leq 4$, then
$\alpha_{s}=(-)^{i_{1}+\cdots+i_{s-1}}\gamma$ or
$\alpha_{t}=-(-)^{i_{1}+\cdots+i_{t-1}}\gamma$.
\end{lemma}
\begin{proof}
The argument is similar to that in Lemma \ref{include_p-1}.
\end{proof}

\begin{prop}\label{sign_p-14}
Let $\Lambda$ be a pairwise molecular subcomplex of $\ubbb$. Let
$\lambda=\uaaa{i}{\alpha}{}$ be a $p-1$ dimensional atom  such
that $\Int\lambda\subset\dpg\Lambda$. If there is a maximal atom
$\lambda'=\uaaa{i}{\alpha}{'}$ in $\Lambda$ with $i_{s}'>i_{s}$
and $i_{t}'>i_{t}$ for some $1\leq s<t\leq 4$ such that
$\uaa{i}{\alpha}{'}{r}\supset\uaa{i}{\alpha}{}{r}$ for at least
three values of $r\in\{1,2,3,4\}$, then
$\alpha_{s}=(-)^{i_{1}+\cdots+i_{s-1}}\gamma$ or
$\alpha_{t}=-(-)^{i_{1}+\cdots+i_{t-1}}\gamma$.
\end{prop}

\begin{proof}
The arguments for  various cases are similar. We give the proof
for the following case.

Suppose that there is a maximal atom
$\lambda'=\uaaa{i}{\alpha}{'}$ in $\Lambda$ such that
$i_{1}'>i_{1}$, $i_{2}'>i_{2}$, $i_{3}'\geq i_{3}$ and
$\uaa{i}{\alpha}{'}{4}\supset\uaa{i}{\alpha}{}{4}$. If $\lambda'$
can be chosen such that $\lambda'\supset\lambda$, then we have
$\alpha_{1}=\gamma$ or $\alpha_{2}=-(-)^{i_{1}}\gamma$, as
required, by Lemma \ref{include_p-14}. In the following, we assume
that $\lambda'$ cannot be chosen such that
$\lambda'\supset\lambda$ so that
$\uaa{i}{\alpha}{'}{3}=u_{3}[i_{3},-\alpha_{3}]$. Let
$\mu=\uaaa{j}{\beta}{}$ be a maximal atom in $\Lambda$ such that
$\lambda\subset\mu$. Then
$\uaa{j}{\beta}{}{1}=\uaa{i}{\alpha}{}{1}$ or
$\uaa{j}{\beta}{}{2}=\uaa{i}{\alpha}{}{2}$. According to Lemma
\ref{lm34}, we may assume that $j_{3}>i_{3}$. Now we consider
three cases, as follows.

1. Suppose that $\mu$ cannot be chosen such that $j_{1}>i_{1}$ or
$j_{2}>i_{2}$.  We claim that
$\alpha_{3}'=-(-)^{i_{1}+i_{2}}\alpha_{1}$ and
$\alpha_{3}'=-(-)^{i_{2}}\alpha_{2}$.

Indeed, if $\lambda'$ and $\mu$ are $(1,3)$-adjacent, then we have
$\alpha_{3}'=-(-)^{i_{1}+i_{2}}\alpha_{1}$ by sign conditions.
Suppose that $\lambda'$ and $\mu$ are not $(1,3)$-adjacent. Then
there is a maximal atom $\lambda''=\uaaa{i}{\alpha}{''}$ such that
$i_{1}''>i_{1}$, $i_{2}''\geq i_{2}$, $i_{3}''>i_{3}$ and
$i_{4}''\geq\min\{i_{4}',j_{4}\}$. According to the assumptions
and Lemma \ref{lm34}, we can see that
$\uaa{i}{\alpha}{''}{4}=u_{4}[i_{4},-\alpha_{4}]$. It follows
easily from Lemma \ref{lm34} that there is a maximal atom
$\mu'=\uaaa{j}{\beta}{'}$ such that
$\uaa{j}{\beta}{'}{1}=\uaa{i}{\alpha}{}{1}$,
$\uaa{j}{\beta}{'}{2}=\uaa{i}{\alpha}{}{2}$, $j_{3}'>i_{3}$ and
$j_{4}'>i_{4}$. By applying Lemma \ref{lm34} and the assumptions,
we can also get a maximal atom $\nu=\uaaa{k}{\varepsilon}{}$ in
$\Lambda$ such that $k_{1}>i_{1}$, $k_{2}\geq i_{2}$,
$\uaa{k}{\varepsilon}{}{3}=\uaa{i}{\alpha}{'}{3}$ and
$k_{4}>i_{4}$. It is evident that $\nu$ and $\mu'$ are
$(1,3)$-adjacent. It follows from sign conditions that
$\alpha_{3}'=-(-)^{i_{1}+i_{2}}\alpha_{1}$, as required. We can
similarly show that $\alpha_{3}'=-(-)^{i_{2}}\alpha_{2}$.

Now, if $\alpha_{3}'=-(-)^{i_{1}+i_{2}}\gamma$, then
$\alpha_{1}=\gamma$, as required; if
$\alpha_{3}=(-)^{i_{1}+i_{2}}\gamma$, then
$\alpha_{2}=-(-)^{i_{1}}\gamma$, as required.

2. Suppose that $\mu$ can be chosen such that $j_{1}>i_{1}$.
Suppose also that $\alpha_{1}=-\gamma$. Then
$\uaa{j}{\beta}{}{2}=\uaa{i}{\alpha}{}{2}$ by the assumptions.
According to Lemma \ref{include_p-14}, we have
$\alpha_{3}=-(-)^{i_{1}+i_{2}}\gamma$, and hence
$\alpha_{3}'=(-)^{i_{1}+i_{2}}\gamma$. By Lemma \ref{lm34}, it is
easy to see that $\lambda'$ and $\mu$ can be chosen such that they
are $(2,3)$-adjacent. It follows from sign conditions that
$\alpha_{2}=\beta_{2}=-(-)^{i_{1}}\gamma$, as required.

3. Suppose that $\mu$ can be chosen such that $j_{2}>i_{2}$ and
that $\alpha_{2}=(-)^{i_{1}}\gamma$. Suppose also that $\mu$
cannot be chosen such that $j_{1}>i_{1}$. According to condition
\ref{main4_4} in Theorem \ref{main4}, Lemma \ref{lm34} and the
assumptions, it is easy to see that $\lambda'$ and $\mu$ can be
chosen such that they are $(1,3)$-adjacent and
$\min\{j_{2},i_{2}'\}=i+1$. According to condition \ref{main4_4}
in Theorem \ref{main4}, there is a maximal atom
$\lambda''=\uaaa{i}{\alpha}{''}$ such that $i_{1}''>i_{1}$,
$j_{2}''=j_{2}$, $i_{3}''>i_{3}$ and
$i_{4}''\geq\min\{i_{4}',j_{4}\}\geq i_{4}$. Moreover, we have
$\alpha_{2}'=-\alpha_{2}=-(-)^{i_{1}}\gamma$ by the assumptions.
According to Note \ref{m4_note4}, we have
$\alpha_{2}''=-(-)^{j_{1}}\beta_{1}=-(-)^{i_{1}}\alpha_{1}$. This
implies that $\alpha_{1}=\gamma$, as required.

This completes the proof.

\end{proof}

\begin{lemma}\label{include_p-24}
Let $\Lambda$ be a pairwise molecular subcomplex of $\ubbb$. Let
$\lambda=\uaaa{i}{\alpha}{}$ be a $p-2$ dimensional atom  such
that $\Int\lambda\subset\dpg\Lambda$. If there is a maximal atom
$\lambda'=\uaaa{i}{\alpha}{'}$ in $\Lambda$ with
$\lambda'\supset\lambda$ such that $i_{r}'>i_{r}$, $i_{s}'>i_{s}$
and $i_{t}'>i_{t}$ for some $1\leq r<s<t\leq 4$, then
$\alpha_{r}=(-)^{i_{1}+\cdots+i_{r-1}}\gamma$ or
$\alpha_{s}=-(-)^{i_{1}+\cdots+i_{s-1}}\gamma$ or
$\alpha_{t}=(-)^{i_{1}+\cdots+i_{t-1}}\gamma$.
\end{lemma}
\begin{proof}
The argument is similar to that in Lemma \ref{include_p-1}.
\end{proof}

\begin{lemma}\label{sign_p-24}
Let $\Lambda$ be a pairwise molecular subcomplex of $\ubbb$. Let
$\lambda=\uaaa{i}{\alpha}{}$ be a $p-2$ dimensional atom  such
that $\Int\lambda\subset\dpg\Lambda$. If there is a maximal atom
$\lambda'=\uaaa{i}{\alpha}{'}$ in $\Lambda$ such that
$i_{r}'>i_{r}$, $i_{s}'>i_{s}$ and $i_{t}'>i_{t}$ for some $1\leq
r<s<t\leq 4$, and such that $i_{\bar{r}}'\geq i_{\bar{r}}$ for
$\bar{r}\in\{1,2,3,4\}\setminus\{r,s,t\}$, then
$\alpha_{r}=(-)^{i_{1}+\cdots+i_{r-1}}\gamma$ or
$\alpha_{s}=-(-)^{i_{1}+\cdots+i_{s-1}}\gamma$ or
$\alpha_{t}=(-)^{i_{1}+\cdots+i_{t-1}}\gamma$.
\end{lemma}
\begin{proof}
The arguments for various choices of $r$, $s$ and $t$ are similar.
We give the proof for $r=1$, $s=3$ and $t=4$.

Suppose that there is a maximal atom
$\lambda'=\uaaa{i}{\alpha}{'}$ such that $i_{1}'>i_{1}$,
$i_{2}'\geq i_{2}$, $i_{3}'>i_{3}$ and $i_{4}'>i_{4}$. If
$\lambda'$ can be chosen such that
$\uaa{i}{\alpha}{'}{2}\supset\uaa{i}{\alpha}{}{2}$, then we have
$\alpha_{1}=\gamma$, $\alpha_{3}=-(-)^{i_{1}+i_{2}}\gamma$ or
$\alpha_{4}=(-)^{i_{1}+i_{2}+i_{3}}\gamma$, as required, by Lemma
\ref{include_p-24}. In the following, we assume that $\lambda'$
cannot be chosen such that
$\uaa{i}{\alpha}{'}{2}\supset\uaa{i}{\alpha}{}{2}$ so that
$\uaa{i}{\alpha}{'}{2}=\uaa{i}{-\alpha}{}{2}$. Let
$\mu=\uaaa{j}{\beta}{}$ be a maximal atom in $\Lambda$ such that
$\mu\supset\lambda$. According to Lemma \ref{lm34}, we can assume
that $j_{2}>i_{2}$. We consider three cases, as follows.

1. Suppose that $\mu$ can be chosen such that there are exactly
three value of $r\in\{1,2,3,4\}$ such that $j_{r}>i_{r}$.

a. Suppose that $\mu$ can be chosen such that $j_{2}>i_{2}$,
$j_{3}>i_{3}$ and $j_{4}>i_{4}$. Then
$\uaa{j}{\beta}{}{1}=\uaa{i}{\alpha}{}{1}$ by the assumptions. We
also have $\alpha_{2}=(-)^{i_{1}}\gamma$,
$\alpha_{3}=-(-)^{i_{1}+i_{2}}\gamma$ or
$\alpha_{4}=(-)^{i_{1}+i_{2}+i_{4}}\gamma$ by Lemma
\ref{include_p-24}. If $\alpha_{3}=-(-)^{i_{1}+i_{2}}\gamma$ or
$\alpha_{4}=(-)^{i_{1}+i_{2}+i_{4}}\gamma$, then $\alpha_{3}$ or
$\alpha_{4}$ is as required. If $\alpha_{2}=(-)^{i_{1}}\gamma$,
then $\alpha_{2}'=-(-)^{i_{1}}\gamma$; moreover, it is easy to see
that $\lambda'$ is $(1,2)$-adjacent to $\mu$; it follows from sign
conditions that $\alpha_{1}=\beta_{1}=\gamma$, as required.

b. Suppose that $\mu$ can be chosen such that $j_{1}>i_{1}$,
$j_{2}>i_{2}$ and $j_{4}>i_{4}$. Then we can get
$\alpha_{1}=\gamma$, $\alpha_{3}=-(-)^{i_{1}+i_{2}}\gamma$ or
$\alpha_{4}=(-)^{i_{1}+i_{2}+i_{4}}\gamma$, as required, by
similar arguments as in case a.

c. Suppose that $\mu$ can be chosen such that $j_{1}>i_{1}$,
$j_{2}>i_{2}$ and $j_{3}>i_{3}$. Suppose also that $\mu$ cannot be
chosen such that $j_{1}>i_{1}$, $j_{2}>i_{2}$ and $j_{4}>i_{4}$.
According to the assumptions, it is easy to see that $\lambda'$
and $\mu$ are $(2,4)$-adjacent and $\min\{i_{3}',j_{3}\}=i_{3}+1$.
By condition \ref{main4_4} in Theorem \ref{main4}, there is a
maximal atom $\lambda''=\uaaa{i}{\alpha}{''}$ in $\Lambda$ such
that $i_{1}''>i$, $i_{2}''>i_{2}$, $i_{3}''=i_{3}$ and
$i_{4}''>i_{4}$. According to the assumptions, we have
$\alpha_{3}''=-\alpha_{3}$. Now, if
$\alpha_{3}=-(-)^{i_{1}+i_{2}}\gamma$, then it is as required. If
$\alpha_{3}=(-)^{i_{1}+i_{2}}\gamma$, then
$\alpha_{3}''=-(-)^{i_{1}+i_{2}}\gamma$; therefore we have
$\alpha_{4}=(-)^{i_{1}+i_{2}+i_{3}}\gamma$ by Note \ref{m4_note4},
as required.

2. Suppose that $\mu$ cannot be chosen such that there are three
value of $r\in\{1,2,3,4\}$ such that $j_{r}>i_{r}$. Suppose also
that $\mu$ can be chosen such that there are two value of
$r\in\{1,2,3,4\}$ such that $j_{r}>i_{r}$.

a. Suppose that $\mu$ can be chosen such that $j_{1}>i_{1}$ and
$j_{2}>i_{2}$. Then we have
$\uaa{j}{\beta}{}{3}=\uaa{i}{\alpha}{}{3}$ and
$\uaa{j}{\beta}{}{4}=\uaa{i}{\alpha}{}{4}$ by the assumptions.
Moreover, by Lemma \ref{lm34} and the assumptions, we can see that
$\lambda'$ is both $(2,3)$-adjacent and $(2,4)$-adjacent to $\mu$.
It follows easily that $\alpha_{3}=-(-)^{i_{1}+i_{2}}\gamma$ or
$\alpha_{4}=(-)^{i_{1}+i_{2}+i_{3}}\gamma$, as required.

b. Suppose that $\mu$ can be chosen such that $j_{2}>i_{2}$ and
$j_{4}>i_{4}$. We can get $\alpha_{1}=\gamma$ or
$\alpha_{3}=-(-)^{i_{1}+i_{2}}\gamma$ by similar arguments as that
in case a.

c. Suppose that $\mu$ can be chosen such that $j_{2}>i_{2}$ and
$j_{3}>i_{3}$. Then we have
$\uaa{j}{\beta}{}{1}=\uaa{i}{\alpha}{}{1}$ and
$\uaa{j}{\beta}{}{4}=\uaa{i}{\alpha}{}{4}$ by the assumptions.
According to condition \ref{main4_4} in Theorem \ref{main4}, Lemma
\ref{lm34} and the assumptions, it is easy to see that $\lambda'$
is both $(1,2)$-adjacent and $(2,4)$-adjacent to $\mu$, and that
$\min\{i_{3}',j_{3}\}=i_{3}+1$. It follows easily from sign
conditions that $\alpha_{1}=\gamma$ or
$\alpha_{4}=(-)^{i_{1}+i_{2}+i_{3}}\gamma$, as required.

3. There remains the case that $\mu$ cannot be chosen such that
$j_{1}>i_{1}$ or $j_{3}>i_{3}$ or $j_{4}>i_{4}$. In this case,  we
have $\uaa{j}{\beta}{}{r}=\uaa{i}{\alpha}{}{r}$ for all
$r\in\{1,3,4\}$. By Lemma \ref{lm34}, it is easy to see that
$\lambda'$ and $\mu$ are adjacent. It follows from sign conditions
that $\alpha_{1}=\gamma$ or $\alpha_{3}=-(-)^{i_{1}+i_{2}}\gamma$
or $\alpha_{4}=(-)^{i_{1}+i_{2}+i_{3}}\gamma$, as required.

This completes the proof.

\end{proof}

\begin{prop}\label{prop5.5F}
Let $\Lambda$ be a pairwise molecular subcomplex of $\ubbb$. Then
$\dpg\Lambda$ is a subcomplex.
\end{prop}

\begin{proof}
We have already seen that $\dpg\Lambda$ is a union of interiors of
atoms. By Lemma \ref{lm5.6}, it suffices to prove that for every
atom $\lambda$ with $\Int\lambda\subset\dpg\Lambda$ and every atom
$\lambda_{1}$ with $\lambda_{1}\subset\lambda$, one has
$\Int\lambda_{1}\subset\dpg\Lambda$. It is evident that there is a
sequence
$\lambda\supset\lambda_{1}^{1}\supset\lambda^{2}_{1}\supset
\cdots\supset\lambda_{1}$ such that the difference of the
dimensions of any pair of consecutive atoms is $1$. We may assume
that $\dim\lambda_{1}=\dim\lambda-1$.

Let $\lambda=\uaaa{i}{\alpha}{}$. Since
$\Int\lambda\subset\dpg\Lambda\subset\Lambda$ and $\Lambda$ is a
subcomplex, we have $\lambda_{1}\subset\lambda\subset\Lambda$ and
$\dim\lambda_{1}\leq\dim\lambda\leq p$. Suppose that
$\mu=\uaaa{l}{\sigma}{}$ is an atom with $\dim\mu=p+1$ and
$\lambda_{1}\subset\mu\subset\Lambda$. We must prove
$\lambda_{1}\subset\dpg\mu$.

If $\lambda\subset\mu$, then
$\lambda_{1}\subset\lambda\subset\dpg\mu$ since
$\lambda\subset\dpg\Lambda$. If $l_{s}>i_{s}+1$ for some $s$, then
$l_{s}>i_{s}^1+1$ for some $s$; hence we also have
$\lambda_{1}\subset\dpg\mu$ by the formation of $\dpg\mu$. In the
following, we may further assume that $\lambda\not\subset\mu$ and
that $l_{s}\leq i_{s}+1$ for every $s$. Thus
$\uaa{i}{\alpha}{}{t}\not\subset\uaa{l}{\sigma}{}{t}$ for some
$t$, and hence $\uaa{l}{\sigma}{}{t}=u_{t}[i_{t},-\alpha_{t}]$ or
$\uaa{l}{\sigma}{}{t}=u_{t}[i_{t}-1,\tilde{\alpha}_{t}]$;
moreover, we can see that
$\uaa{i}{\alpha}{}{s}\subset\uaa{l}{\sigma}{}{s}$ for every $s$
with $s\neq t$. Note that $\dim\mu=p+1$ and $\dim\lambda\leq p$,
we now have 5 cases, as follows.

1. Suppose that $\uaa{l}{\sigma}{}{t}=\uaa{i}{-\alpha}{}{t}$ for
some $t$ and $l_{s}=i_{s}+1$ for only one value of $s$. The proofs
for different choices of $t$ and the one value for $s$ are
similar. We give the proof only for the case $l_{1}=i_{1}+1$ and
$\uaa{l}{\sigma}{}{2}=\uaa{i}{-\alpha}{}{2}$. In this case, we can
see that $\lambda_{1}$ is of the form
$\lambda_{1}=\uaa{i}{\alpha}{}{1}\times
u_{2}[i_{2}-1,\tilde{\alpha}_{2}]\times \uaa{i}{\alpha}{}{3}\times
\uaa{i}{\alpha}{}{4}$; we also have $\dim\lambda=p$. According to
Proposition \ref{Int_sign4}, we have $\alpha_{1}=\gamma$. This
implies that $\lambda_{1}\subset\dpg\mu$, as required, by Lemma
\ref{dtatom}.

2. Suppose that $l_{t}=i_{t}-1$ and $l_{s}=i_{s}+1$ for only two
values of $s$. The proofs for different choices of $t$ and the two
values for $s$ are similar. We give the proof only for the case
$l_{2}=i_{2}-1$, $l_{1}=i_{1}+1$ and $l_{3}=i_{3}+1$.  In this
case, we can see that $\lambda_{1}$ is of the form
$\lambda_{1}=\uaa{i}{\alpha}{}{1}\times
u_{2}[i_{2}-1,\sigma_{2}]\times \uaa{i}{\alpha}{}{3}\times
\uaa{i}{\alpha}{}{4}$; we also have $\dim\lambda=p$ and
$\uaa{l}{\sigma}{}{4}=\uaa{i}{\alpha}{}{4}$. To get
$\lambda_{1}\subset\dpg\mu$, it suffices to prove that
$\alpha_{1}=\gamma$ or $\alpha_{3}=(-)^{i_{1}+i_{2}}\gamma$ by
Lemma \ref{dtatom}.

Let $\lambda'=\uaaa{i}{\alpha}{'}$ be a maximal atom in $\Lambda$
such that $\lambda\subset\lambda'$. Let $\mu'=\uaaa{l}{\sigma}{'}$
be a maximal atom in $\Lambda$ such that $\mu\subset\mu'$. If
$\lambda'$ can be chosen such that $i_{1}'>i_{1}$ or
$i_{3}'>i_{3}$, then we have $\alpha_{1}=\gamma$ or
$\alpha_{3}=(-)^{i_{1}+i_{2}}\gamma$ by Lemma \ref{lm5.3F} which
implies that $\lambda_{1}\subset\dpg\mu$, as required. If there is
a maximal atom $\mu''=\uaaa{l}{\sigma}{''}$ with $l_{r}''\geq
i_{r}$ for all $r\in\{1,2,3,4\}$ such that $l_{1}''>i_{1}$ or
$l_{3}''>i_{3}$, then, by Proposition \ref{Int_sign4}, we have
$\alpha_{1}=\gamma$ or $\alpha_{3}=(-)^{i_{1}+i_{2}}\gamma$ which
implies that $\lambda_{1}\subset\dpg\mu$, as required.

Now suppose that  $\lambda'$ cannot be chosen such that
$i_{1}'>i_{1}$ or $i_{3}'>i_{3}$. Suppose also that there is no
maximal atom $\mu''=\uaaa{l}{\sigma}{''}$ with $l_{r}''\geq i_{r}$
for all $r\in\{1,2,3,4\}$ such that $l_{1}''>i_{1}$ or
$l_{3}''>i_{3}$.  Then
$\uaa{i}{\alpha}{'}{1}=\uaa{i}{\alpha}{}{1}$,
$\uaa{i}{\alpha}{'}{3}=\uaa{i}{\alpha}{}{3}$ and
$\uaa{l}{\sigma}{'}{2}=\uaa{l}{\sigma}{}{2}=
u_{2}[i_{2}-1,\tilde{\alpha}_{2}]$. It is easy to see that
$\lambda'$ is both $(1,2)$-adjacent to $\mu'$ and $(2,3)$-adjacent
to $\mu'$. If $\sigma_{2}=-(-)^{i_{1}}\gamma$, then we can see
that $\alpha_{1}=\gamma$ by sign conditions; if
$\sigma_{2}=(-)^{i_{1}}\gamma$, then we can see that
$\alpha_{3}=(-)^{i_{1}+i_{2}}\gamma$ by sign conditions. These
imply that $\lambda_{1}\subset\dpg\mu$, as required.

3. Suppose that $\uaa{l}{\sigma}{}{t}=\uaa{i}{-\alpha}{}{t}$ for
some $t$ and $l_{s}=i_{s}+1$ for only 2 value of $s$. The proofs
for different choices of $t$ and the 2 values for $s$ are similar.
We give the proof only for the case
$\uaa{l}{\sigma}{}{2}=\uaa{i}{-\alpha}{}{2}$, $l_{1}=i_{1}+1$ and
$l_{3}=i_{3}+1$. In this case, we can see that $\lambda_{1}$ is of
the form $\lambda_{1}=\uaa{i}{\alpha}{}{1}\times
u_{2}[i_{2}-1,\tilde{\alpha}_{2}]\times \uaa{i}{\alpha}{}{3}\times
\uaa{i}{\alpha}{}{4}$; we also have $\dim\lambda=p-1$. According
to Proposition \ref{sign_p-14}, we have $\alpha_{1}=\gamma$ or
$\alpha_{3}=-(-)^{i_{1}+i_{2}}\gamma$. This implies that
$\lambda_{1}\subset\dpg\mu$, as required, by Lemma \ref{dtatom4}.

4. Suppose that $l_{t}=i_{t}-1$ and $l_{s}=i_{s}+1$ for the other
three values of $s$. The proofs for different choices of $t$ are
similar. We give the proof only for the case
$\uaa{l}{\sigma}{}{2}=u_{2}[i_{2}-1,\tilde{\alpha}_{2}]$,
$l_{1}=i_{1}+1$, $l_{3}=i_{3}+1$ and $l_{4}=l_{4}+1$.  In this
case, we can see that $\lambda_{1}$ is of the form
$\lambda_{1}=\uaa{i}{\alpha}{}{1}\times
u_{2}[i_{2}-1,\sigma_{2}]\times \uaa{i}{\alpha}{}{3}\times
\uaa{i}{\alpha}{}{4}$; we also have $\dim\lambda=p-1$. To get
$\lambda_{1}\subset\dpg\mu$, it suffices to prove that
$\alpha_{1}=\gamma$ or $\alpha_{3}=(-)^{i_{1}+i_{2}}\gamma$ or
$\alpha_{4}=-(-)^{i_{1}+i_{2}+i_{3}}\gamma$ by Lemma
\ref{dtatom4}.

Let $\mu'=\uaaa{l}{\sigma}{'}$ be a maximal atom in $\Lambda$ such
that $\mu\subset\mu'$. Let $\lambda'=\uaaa{i}{\alpha}{'}$ be a
maximal atom in $\Lambda$ such that $\lambda\subset\lambda'$. If
$\mu'$ can be chosen such that $l_{2}'\geq i_{2}$, then we have
$\alpha_{1}=\gamma$ or
$\alpha_{4}=-(-)^{i_{1}+i_{2}+i_{3}}\gamma$, as required, by Lemma
\ref{sign_p-14}. In the following proof, we may assume that $\mu'$
cannot be chosen such that $l_{2}'\geq i_{2}$. Now there are two
cases, as follows.

a. Suppose that $\lambda'$ can be chosen such that there are
exactly two values of $s\in\{1,3,4\}$ such that $i_{s}'>i_{s}$. If
$\lambda'$ can be chosen such that $i_{1}'>i_{1}$ and
$i_{4}'>i_{4}$, or such that $i_{3}'>i_{3}$ and $i_{4}'>i_{4}$.
Then it is easy to see that $\alpha_{1}=\gamma$ or
$\alpha_{3}=(-)^{i_{1}+i_{2}}\gamma$ or
$\alpha_{4}=-(-)^{i_{1}+i_{2}+i_{3}}\gamma$, as required, by
Proposition \ref{sign_p-14}. Suppose that  $\lambda'$ can be
chosen such that $i_{1}'>i_{1}$ and $i_{3}'>i_{3}$. Suppose also
that $\alpha_{1}=-\gamma$. By the assumptions, we have
$\uaa{i}{\alpha}{'}{4}=\uaa{i}{\alpha}{}{4}$. According to
Proposition \ref{sign_p-14}, we also have
$\alpha_{3}=-(-)^{i_{1}+i_{2}}\gamma$. By the assumptions and
condition \ref{main4_4} in Theorem \ref{main4}, it is easy to see
that $\lambda'$ is $(2,4)$-adjacent to $\mu'$ and
$\min\{i_{3}',l_{3}'\}=i_{3}+1$. According to condition
\ref{main4_4} in Theorem \ref{main4}, there exists a maximal atom
$\nu=\uaaa{j}{\beta}{}$ such that $j_{1}>i_{1}$, $j_{2}\geq
i_{2}$, $j_{3}\geq i_{3}$ and $j_{4}>i_{4}$. If
$\uaa{j}{\beta}{}{3}\supset \uaa{i}{\alpha}{}{3}$,  then it is
easy to see that $\alpha_{4}=-(-)^{i_{1}+i_{2}+i_{3}}\gamma$, as
required, by Proposition \ref{sign_p-14}. If  $j_{3}=i_{3}$ and
$\beta_{3}=-\alpha_{3}=(-)^{i_{1}+i_{2}}\gamma$, then, by Note
\ref{m4_note4}, we also have
$\alpha_{4}=\alpha_{4}'=-(-)^{j_{3}}\beta_{3}
=-(-)^{i_{1}+i_{2}+i_{3}}\gamma$, as required.

b. Suppose that $\lambda'$ cannot be chosen such that there are
two values of $s\in\{1,3,4\}$ such that $i_{s}'>i_{s}$. By our
assumptions, it is easy to see that $i_{s}'>i_{s}$ for at most one
value of $s\in\{1,3,4\}$.

Suppose that $\lambda'$ can be chosen such that $i_{1}'>i_{1}$.
Then $\uaa{i}{\alpha}{'}{3}=\uaa{i}{\alpha}{}{3}$ and
$\uaa{i}{\alpha}{'}{4}=\uaa{i}{\alpha}{}{4}$ by the assumptions.
It is also easy to see that $\mu'$ is both $(2,3)$-adjacent and
$(2,4)$-adjacent to $\lambda'$. It follows easily that
$\alpha_{3}=(-)^{i_{1}+i_{2}}\gamma$ (when
$\sigma_{2}=(-)^{i_{1}+i_{2}}\gamma$) or
$\alpha_{4}=-(-)^{i_{1}+i_{2}+i_{3}}\gamma$ (when
$\sigma_{2}=-(-)^{i_{1}+i_{2}}\gamma$), as required.

Suppose that $\lambda'$ can be chosen such that  $i_{4}'>i_{4}$ or
$\lambda'=\lambda$. By arguments similar to that in the last
paragraph, we can get $\alpha_{1}=\gamma$ or
$\alpha_{3}=(-)^{i_{1}+i_{2}}\gamma$ or
$\alpha_{4}=-(-)^{i_{1}+i_{2}+i_{3}}\gamma$, as required.

Suppose that $\lambda'$ can be chosen such that $i_{3}'>i_{3}$.
According to the assumptions and condition \ref{main4_3} in
Theorem \ref{main4}, it is easy to see that $\lambda'$ and $\mu'$
are both $(1,2)$-adjacent and $(2,4)$-adjacent. By condition
\ref{main4_3} in Theorem \ref{main4} again, we can also get
$\min\{i_{3}',l_{3}'\}=i_{3}+1$. It follows easily that
$\alpha_{1}=\gamma$ (when $\sigma_{2}=-(-)^{i_{1}+i_{2}}\gamma$)
or $\alpha_{4}=-(-)^{i_{1}+i_{2}+i_{3}}\gamma$ (when
$\sigma_{2}=(-)^{i_{1}+i_{2}}\gamma$), as required.

5. Suppose that $\uaa{l}{\sigma}{}{t}=\uaa{i}{-\alpha}{}{t}$ for
some $t$ and $l_{s}=i_{s}+1$ for 3 values of $s$. The proofs for
different choices of $t$ and the 3 values for $s$ are similar. We
give the proof only for the case
$\uaa{l}{\sigma}{}{2}=\uaa{i}{-\alpha}{}{2}$, $l_{1}=i_{1}+1$ and
$l_{3}=i_{3}+1$ and $l_{4}=i_{4}+1$. In this case, we can see that
$\lambda_{1}$ is of the form
$\lambda_{1}=\uaa{i}{\alpha}{}{1}\times
u_{2}[i_{2}-1,\tilde{\alpha}_{2}]\times \uaa{i}{\alpha}{}{3}\times
\uaa{i}{\alpha}{}{4}$; we also have $\dim\lambda=p-2$. According
to Proposition \ref{sign_p-14}, we have $\alpha_{1}=\gamma$ or
$\alpha_{3}=-(-)^{i_{1}+i_{2}}\gamma$ or
$\alpha_{4}=(-)^{i_{1}+i_{2}+i_{3}}\gamma$. This implies that
$\lambda_{1}\subset\dpg\mu$, as required, by Lemma \ref{dtatom}.

This completes the proof.
\end{proof}

\begin{prop}\label{prop5.6F}
Let $\Lambda$ be a pairwise molecular subcomplex of $\ubbb$. For
$1\leq s\leq 4$, if $p\geq I_{s}$ and
$F_{I_{s}}^{u_{s}}(\Lambda)\neq\emptyset$, then
$F_{I_{s}}^{u_{s}}(\dpg\Lambda)=
d_{p-I_{s}}^{\gamma}F_{I_{s}}^{u_{s}}(\Lambda)$.
\end{prop}

\begin{proof}

The arguments for various choices of $s$ are similar. We give the
proof for $s=2$.

We first prove that
$d_{p-I_{2}}^{\gamma}F_{I_{2}}^{u_{2}}(\Lambda)\subset
F_{I_{2}}^{u_{2}}(\dpg\Lambda)$.

Let $\uaa{i}{\alpha}{}{1}\times
u_{3}^{I_{2}}[i_{3},\alpha_{3}]\times
u_{4}^{I_{2}}[i_{4},\alpha_{4}]$ be a maximal atom in
$d_{p-I_{2}}^{\gamma}F_{I_{2}}^{u_{2}}(\Lambda)$. We must show
that $\uaa{i}{\alpha}{}{1}\times
u_{3}^{I_{2}}[i_{3},\alpha_{3}]\times
u_{4}^{I_{2}}[i_{4},\alpha_{4}]\subset
F_{I_{2}}^{u_{2}}(\dpg\Lambda)$. Since $\uaa{i}{\alpha}{}{1}\times
u_{3}^{I_{2}}[i_{3},\alpha_{3}]\times
u_{4}^{I_{2}}[i_{4},\alpha_{4}]\subset
d_{p-I_{2}}^{\gamma}F_{I_{2}}^{u_{2}}(\Lambda)\subset
F_{I_{2}}^{u_{2}}(\Lambda)$, we can see that
$\uaa{i}{\alpha}{}{1}\times u_{2}[I_{2},\alpha_{2}]\times
\uaa{i}{\alpha}{}{3}\times \uaa{i}{\alpha}{}{4}\subset\Lambda$
for some sign $\alpha_{2}$. We now consider three cases, as
follows.

1. Suppose that $\dim (\uaa{i}{\alpha}{}{1}\times
u_{3}^{I_{2}}[i_{3},\alpha_{3}]\times
u_{4}^{I_{2}}[i_{4},\alpha_{4}])<p-I_{2}$, then
$\uaa{i}{\alpha}{}{1}\times u_{3}^{I_{2}}[i_{3},\alpha_{3}]\times
u_{4}^{I_{2}}[i_{4},\alpha_{4}]$ is a maximal atom in
$F_{I_{2}}^{u_{2}}(\Lambda)$.  We claim that
$\uaa{i}{\alpha}{}{1}\times u_{2}[I_{2},\alpha_{2}]\times
\uaa{i}{\alpha}{}{3}\times
\uaa{i}{\alpha}{}{4}\subset\dpg\Lambda$.  Firstly, it is evident
that $\dim(\uaa{i}{\alpha}{}{1}\times
u_{2}[I_{2},\alpha_{2}]\times \uaa{i}{\alpha}{}{3}\times
\uaa{i}{\alpha}{}{4})<p$. Moreover, suppose that there is a
maximal atom $\uaaa{i}{\alpha}{'}$ in $\Lambda$ such that
$\uaaa{i}{\alpha}{'}\supset \uaa{i}{\alpha}{}{1}\times
u_{2}[I_{2},\alpha_{2}]\times \uaa{i}{\alpha}{}{3}\times
\uaa{i}{\alpha}{}{4}$. Then it is evident that
$\uaa{i}{\alpha}{'}{r}=\uaa{i}{\alpha}{}{r}$ for all
$r\in\{1,3,4\}$. According to Lemma \ref{dtatom4}, we have
$\uaa{i}{\alpha}{}{1}\times u_{2}[I_{2},\alpha_{2}]\times
\uaa{i}{\alpha}{}{3}\times \uaa{i}{\alpha}{}{4}\subset
\dpg(\uaaa{i}{\alpha}{'})$. It follows easily from Lemma
\ref{lm5.2} that $\uaa{i}{\alpha}{}{1}\times
u_{2}[I_{2},\alpha_{2}]\times \uaa{i}{\alpha}{}{3}\times
\uaa{i}{\alpha}{}{4}\subset\dpg\Lambda$. Hence
$\uaa{i}{\alpha}{}{1}\times u_{3}^{I_{2}}[i_{3},\alpha_{3}]\times
u_{4}^{I_{2}}[i_{4},\alpha_{4}]\subset
F_{I_{2}}^{u_{2}}(\dpg\Lambda)$.

2. Suppose that $\dim (\uaa{i}{\alpha}{}{1}\times
u_{3}^{I_{2}}[i_{3},\alpha_{3}]\times
u_{4}^{I_{2}}[i_{4},\alpha_{4}])=p-I_{2}$. Suppose also that
$\alpha_{2}$ can be chosen such that
$\alpha_{2}=(-)^{i_{1}}\gamma$. We claim that
$\uaa{i}{\alpha}{}{1}\times u_{2}[I_{2},\alpha_{2}]\times
\uaa{i}{\alpha}{}{3}\times
\uaa{i}{\alpha}{}{4}\subset\dpg\Lambda$. Indeed, it is evident
that $\dim(\uaa{i}{\alpha}{}{1}\times
u_{2}[I_{2},\alpha_{2}]\times \uaa{i}{\alpha}{}{3}\times
\uaa{i}{\alpha}{}{4})=p$. Moreover, suppose that there is a
maximal atom $\uaaa{i}{\alpha}{'}$ in $\Lambda$ such that
$\uaaa{i}{\alpha}{'}\supset \uaa{i}{\alpha}{}{1}\times
u_{2}[I_{2},\alpha_{2}]\times \uaa{i}{\alpha}{}{3}\times
\uaa{i}{\alpha}{}{4}$. Then $\uaa{i}{\alpha}{}{1}\times
u_{3}^{I_{2}}[i_{3},\alpha_{3}]\times
u_{4}^{I_{2}}[i_{4},\alpha_{4}]\subset \uaa{i}{\alpha}{'}{1}\times
u_{3}^{I_{2}}[i_{3}',\alpha_{3}']\times
u_{4}^{I_{2}}[i_{4}',\alpha_{4}']\subset
F_{I_{2}}^{u_{2}}(\Lambda)$. Thus $\uaa{i}{\alpha}{}{1}\times
u_{3}^{I_{2}}[i_{3},\alpha_{3}]\times
u_{4}^{I_{2}}[i_{4},\alpha_{4}]\subset
d_{p-I_{2}}^{\gamma}(\uaa{i}{\alpha}{'}{1}\times
u_{3}^{I_{2}}[i_{3}',\alpha_{3}']\times
u_{4}^{I_{2}}[i_{4}',\alpha_{4}'])$ by Lemma \ref{lm5.2}. It
follows easily from Lemma \ref{dtatom} and \ref{dtatom4} that
$\uaa{i}{\alpha}{}{1}\times u_{2}[I_{2},\alpha_{2}]\times
\uaa{i}{\alpha}{}{3}\times \uaa{i}{\alpha}{}{4}\subset\dpg
(\uaaa{i}{\alpha}{'})$. Therefore $\uaa{i}{\alpha}{}{1}\times
u_{2}[I_{2},\alpha_{2}]\times \uaa{i}{\alpha}{}{3}\times
\uaa{i}{\alpha}{}{4}\subset\dpg\Lambda$, and hence
$\uaa{i}{\alpha}{}{1}\times u_{3}^{I_{2}}[i_{3},\alpha_{3}]\times
u_{4}^{I_{2}}[i_{4},\alpha_{4}]\subset
F_{I_{2}}^{u_{2}}(\dpg\Lambda)$, as required.

3. Suppose that $\dim (\uaa{i}{\alpha}{}{1}\times
u_{3}^{I_{2}}[i_{3},\alpha_{3}]\times
u_{4}^{I_{2}}[i_{4},\alpha_{4}])=p-I_{2}$. Suppose also that
$\alpha_{2}$ cannot be chosen such that
$\alpha_{2}=(-)^{i_{1}}\gamma$. Then $i_{2}'=I_{2}$ and
$\alpha_{2}'=\alpha_{2}=-(-)^{i_{1}}\gamma$. By arguments similar
to those in  case 2, it is easy to see that
$\uaa{i}{\alpha}{}{1}\times u_{2}[I_{2},\alpha_{2}]\times
\uaa{i}{\alpha}{}{3}\times
\uaa{i}{\alpha}{}{4}\subset\dpg\Lambda$. Therefore
$\uaa{i}{\alpha}{}{1}\times u_{3}^{I_{2}}[i_{3},\alpha_{3}]\times
u_{4}^{I_{2}}[i_{4},\alpha_{4}]\subset
F_{I_{2}}^{u_{2}}(\dpg\Lambda)$, as required.

This completes the proof that
$d_{p-I_{2}}^{\gamma}F_{I_{2}}^{u_{2}}(\Lambda)\subset
F_{I_{2}}^{u_{2}}(\dpg\Lambda)$.

Conversely, let $\lambda=\uaa{i}{\alpha}{}{1}\times
u_{3}^{I_{2}}[i_{3},\alpha_{3}]\times
u_{4}^{I_{2}}[i_{4},\alpha_{4}]$ be an atom in $u_{1}\times
u_{3}^{I_{2}}\times u_{4}^{I_{2}}$ such that $\Int\lambda\subset
F_{I_{2}}^{u_{2}}(\dpg\Lambda)$. We must show that
$\Int\lambda\subset
d_{p-I_{2}}^{\gamma}F_{I_{2}}^{u_{2}}(\Lambda)$.

It is easy to see that there is an atom $\uaaa{i}{\alpha}{}$ in
$\Lambda$ such that $\Int(\uaaa{i}{\alpha}{})\subset\dpg\Lambda$
and $i_{2}\geq I_{2}$. Since $\dpg\Lambda$ is a subcomplex of
$\ubbb{}$, we can see that $\uaa{i}{\alpha}{}{1}\times
u_{2}[I_{2},\alpha_{2}']\times\uaa{i}{\alpha}{}{3}\times\uaa{i}{\alpha}{}{4}
\subset\dpg\Lambda$ for some sign $\alpha_{2}'$. It follows that
$\dim(\uaa{i}{\alpha}{}{1}\times
u_{3}^{I_{2}}[i_{3},\alpha_{3}]\times
u_{4}^{I_{2}}[i_{4},\alpha_{4}])\leq p-I_{2}$. Clearly, we have
$\lambda\subset F_{I_{2}}^{u_{2}}(\Lambda)$. To prove that
$\Int\lambda\subset
d_{p-I_{2}}^{\gamma}F_{I_{2}}^{u_{2}}(\Lambda)$, it suffices to
verify the second condition in the Lemma \ref{lm5.2}. Let
$\mu=\uaa{l}{\sigma}{}{1}\times
u_{3}^{I_{2}}[l_{3},\sigma_{3}]\times
u_{4}^{I_{2}}[i_{4},\sigma_{4}]$ be an atom in
$F_{I_{2}}^{u_{2}}(\Lambda)$ such that $\lambda\subset\mu$ and
$\dim\mu=p-I_{2}+1$. We must prove that $\lambda\subset
d_{p-I_{2}}^{\gamma}\mu$. It is evident that
$\uaa{l}{\sigma}{}{1}\times
u_{2}[I_{2},\sigma_{2}']\times\uaa{l}{\sigma}{}{3}\times\uaa{l}{\sigma}{}{4}
\subset\Lambda$ for some sign $\sigma_{2}'$. If $l_{1}>i_{1}+1$,
$l_{3}>i_{3}+1$ or $l_{4}>i_{4}+1$, then it is evident that
$\lambda\subset\dpg\mu$, as required. In the following prove, we
may assume that $l_{1}\leq i_{1}+1$, $l_{3}\leq i_{3}+1$ and
$l_{4}\leq i_{4}+1$  so that $\dim\lambda=p-I_{2}$ or
$\dim\lambda=p-I_{2}-1$ or $\dim\lambda=p-I_{2}-2$. Now there are
various cases, as follows.

Suppose that $\alpha_{2}'$ and $\sigma_{2}'$ can be chosen such
that $\alpha_{2}'=\sigma_{2}'$. Then $\uaa{i}{\alpha}{}{1}\times
v[I_{2},\alpha_{2}']\times\uaa{i}{\alpha}{}{3}\times\uaa{i}{\alpha}{}{4}
\subset\dpg\Lambda\cap (\uaa{l}{\sigma}{}{1}\times
u_{2}[I_{2},\sigma_{2}']\times\uaa{l}{\sigma}{}{3}\times
\uaa{l}{\sigma}{}{4}) \subset\dpg(\uaa{l}{\sigma}{}{1}\times
u_{2}[I_{2},\sigma_{2}']\times\uaa{l}{\sigma}{}{3}\times\uaa{l}{\sigma}{}{4})$.
It follows easily from Lemma \ref{dtatom4} and Lemma \ref{dtatom}
that $\lambda\subset d_{p-I_{2}}^{\gamma}\mu$, as required.

Suppose that $\alpha_{2}'$ and $\sigma_{2}'$ cannot be chosen such
that $\alpha_{2}'=\sigma_{2}'$. Suppose also that $I_{2}>0$. Since
$\dpg\Lambda$ is a subcomplex, we know that
$\uaa{i}{\alpha}{}{1}\times
u_{2}[I_{2}-1,\pm]\times\uaa{i}{\alpha}{}{3}\times
\uaa{i}{\alpha}{}{4}\subset\dpg\Lambda$. This implies that
$\uaa{i}{\alpha}{}{1}\times
u_{2}[I_{2}-1,\pm]\times\uaa{i}{\alpha}{}{3}\times
\uaa{i}{\alpha}{}{4}\subset \dpg(\uaa{l}{\sigma}{}{1}\times
u_{2}[I_{2},\sigma_{2}']\times\uaa{l}{\sigma}{}{3}\times
\uaa{l}{\sigma}{}{4})$. Hence $\uaa{i}{\alpha}{}{1}\times
u_{2}[I_{2},\sigma_{2}']\times\uaa{i}{\alpha}{}{3}\times
\uaa{i}{\alpha}{}{4}\subset \dpg(\uaa{l}{\sigma}{}{1}\times
u_{2}[I_{2},\sigma_{2}']\times\uaa{l}{\sigma}{}{3}\times
\uaa{l}{\sigma}{}{4})$. It follows easily from Lemma \ref{dtatom}
and Lemma \ref{dtatom4} that $\lambda\subset
d_{p-I_{2}}^{\gamma}\mu$, as required.

There remains the case that $J=0$ and $\alpha_{2}'$ and $\tau'$
cannot be chosen such that $\alpha_{2}'=\sigma_{2}'$. By
arguments similar to those in cases 1, 3 and 5 in the proof of
Proposition \ref{prop5.5}, we can get $\lambda\subset
d_{p-I_{2}}^{\gamma}\mu$, as required.

This completes the proof.

\end{proof}

\begin{lemma}
Let $\Lambda$ be a pairwise molecular subcomplex. Then
$\dpg\Lambda$ satisfies condition 1 for pairwise molecular
subcomplexes.
\end{lemma}
\begin{proof}

Let $\lambda=\uaaa{i}{\alpha}{}$ and
$\lambda'=\uaaa{i}{\alpha}{'}$ be a pair of maximal atom in
$\dpg\Lambda$ with $i_{s}\leq i_{s}'$ for every value of  $s$. We
must prove that $\lambda=\lambda'$. Suppose that $\dim\lambda<p$
or $\dim\lambda'<p$. Then it is easy to see that $\lambda$ is a
maximal atom in $\Lambda$ when $\dim\lambda<p$ and $\lambda'$ is a
maximal atom in $\Lambda$ when $\dim\lambda'<p$. According to
condition 1 for $\Lambda$, we can see that $\lambda=\lambda'$, as
required. In the following argument, we may assume that
$\dim\lambda=p$ and $\dim\lambda=p$ so that $i_{s}=i_{s}'$ for
every value of $s$.

Now suppose otherwise that $\lambda\neq\lambda'$. Then
$\uaa{i}{\alpha}{'}{t}=\uaa{i}{-\alpha}{}{t}$ for some $t$. Let
$r$ be such that $r\neq t$. By Proposition \ref{prop5.6F}, we have
$F_{i_{r}}^{u_{r}}(\lambda)\subset F_{i_{r}}^{u_{r}}(\dpg\Lambda)
=d_{p-i_{r}}^{\gamma}F_{i_{r}}^{u_{r}}(\Lambda)$ and similarly
$F_{i_{r}}^{u_{r}}(\lambda') \subset
d_{p-i_{r}}^{\gamma}F_{i_{r}}^{u_{r}}(\Lambda)$. Since $\dim
F_{i_{r}}^{u_{r}}(\lambda)=\dim
F_{i_{r}}^{u_{r}}(\lambda')=p-i_{r}$ and $\dim
d_{p-i_{r}}^{\gamma}F_{i_{r}}^{u_{r}}(\Lambda)\leq p-i_{r}$, we
can see that $F_{i_{r}}^{u_{r}}(\lambda)$ and
$F_{i_{r}}^{u_{r}}(\lambda')$ are maximal atoms in the molecule
$d_{p-i_{r}}^{\gamma}F_{i_{r}}^{u_{r}}(\Lambda)$. It follows from
condition 1 for $d_{p-i_{r}}^{\gamma}F_{i_{r}}^{u_{r}}(\Lambda)$
that $F_{i_{r}}^{u_{r}}(\lambda)=F_{i_{r}}^{u_{r}}(\lambda')$.
This contradicts the hypothesis that
$\uaa{i}{\alpha}{'}{t}=\uaa{i}{-\alpha}{}{t}$.

This completes the proof.

\end{proof}

\begin{prop}
Let $\Lambda$ be a pairwise molecular subcomplex. Then so is
$\dpg\Lambda$.
\end{prop}
\begin{proof}
We have shown in the last Lemma that $\dpg\Lambda$ satisfies
condition 1 for pairwise molecular subcomplexes. Moreover, by
Proposition \ref{prop5.6F}, we have
$F_{I_{s}}^{u_{s}}(\dpg\Lambda)=d_{p-I_{s}}^{\gamma}F_{I_{s}}^{u_{s}}(\Lambda)$
for all values of $s\in\{1,2,3,4\}$ and all $I_{s}$ with
$I_{s}\leq p$. Since $F_{I_{s}}^{u_{s}}(\Lambda)$ is a molecule or
empty set for every value of $s$ and every $I_{s}$, we can see
that $F_{I_{s}}^{u_{s}}(\dpg\Lambda)$ is a molecule or empty set
for every value of $s$ and every $I_{s}$. It follows from the
definition that $\dpg\Lambda$ is pairwise molecular.

This completes the proof.
\end{proof}

\begin{theorem} \label{dt4}
Let $\Lambda$ be a pairwise molecular subcomplex. Then the
dimension of every maximal atom in $d_{p}^{\gamma}\Lambda$ is not
greater than $p$. Moreover, an atom of dimension less than $p$ is
a maximal atom in $d_{p}^{\gamma}\Lambda$ if and only if it is a
maximal atom in $\Lambda$; an atom $\uaaa{i}{\alpha}{}$ of
dimension $p$ is a maximal atom in $d_{p}^{\gamma}\Lambda$ if and
only if there is a maximal atom $\uaaa{k}{\varepsilon}{}$ in
$\Lambda$ such that $k_{s}\geq i_{s}$ for $1\leq s\leq 4$, and the
signs  $\alpha_{s}$ ($1\leq s\leq 4$) satisfy the following
conditions: if $\uaaa{k}{\varepsilon}{}$ can be chosen such that
$k_{s}>i_{s}$ and $k_{r}\geq i_{r}$ for $1\leq r\leq 4$, then
$\alpha_{s}=(-)^{i_{1}+\cdots +i_{s-1}}\gamma$; otherwies
$\alpha_{s}=\varepsilon_{s}$.

\end{theorem}
\begin{note}
It follows easily from condition \ref{main4_3} in Theorem
\ref{main4} that $\alpha_{s}$ are well defined for all $1\leq
s\leq 4$.
\end{note}
\begin{proof}
By the definition of $d_{p}^{\gamma}\Lambda$, it is evident that
the dimension of every maximal atom in $d_{p}^{\gamma}\Lambda$ is
not greater than $p$.

Let $\Lambda_{1}$ be the subcomplex of $\ubbb$ as described in
this theorem. It is easy to see that $\Lambda_{1}$ satisfies
condition 1 for pairwise molecular subcomplexes. By Lemma
\ref{eq4}, it suffices to prove that
$F_{I_{s}}^{u_{s}}(\Lambda_{1})=F_{I_{s}}^{u_{s}}(\dpg\Lambda)$
for all $I_{s}$ and all  $1\leq s\leq 4$. The arguments for
different choices of $s$ are similar. We now give the proof for
$F_{I_{2}}^{u_{2}}(\Lambda_{1})=F_{I_{2}}^{u_{2}}(\dpg\Lambda)$.
If $I_{2}>p$, then it is easy to see that
$F_{I_{2}}^{u_{2}}(\Lambda_{1})=\emptyset=F_{I_{2}}^{u_{2}}(\dpg\Lambda)$,
as required. In the remaining proof, we may assume that $I_{2}\leq
p$. We have known that
$F_{I_{2}}^{u_{2}}(\dpg\Lambda)=d_{p-I_{2}}^{\gamma}F_{I_{2}}^{u_{2}}(\Lambda)$.
We need only to prove that
$F_{I_{2}}^{u_{2}}(\Lambda_{1})=d_{p-I_{2}}^{\gamma}F_{I_{2}}^{u_{2}}(\Lambda)$.

By the definition of $F_{I_{2}}^{u_{2}}$, it is easy to see that
$F_{I_{2}}^{u_{2}}(\Lambda_{1})$ and
$d_{p-I_{2}}^{\gamma}F(\Lambda)$ are subcomplexes of $u_{1}\times
u_{3}^{I_{2}}\times u_{4}^{I_{2}}$. We are going to prove that
$F_{I_{2}}^{u_{2}}(\Lambda_{1})$ and
$d_{p-I_{2}}^{\gamma}F_{I_{2}}^{u_{2}}(\Lambda)$ consists of the
same maximal atoms so that they are equal.

We first prove every maximal atom in
$F_{I_{2}}^{u_{2}}(\Lambda_{1})$ is a maximal atom in
$d_{p-I_{2}}^{\gamma}F_{I_{2}}^{u_{2}}(\Lambda)$.

Let $\mu=\uaa{i}{\alpha}{}{1}\times
u_{3}^{I_{2}}[i_{3},\alpha_{3}]\times
u_{4}^{I_{2}}[i_{4},\alpha_{4}]$ be a maximal atom in
$F_{I_{2}}^{u_{2}}(\Lambda_{1})$. Then $\Lambda_{1}$ has a
$(u_{2},I_{2})$-projection maximal atom $\lambda$  of the form
$\lambda=\uaaa{i}{\alpha}{}$ with $i_{2}\geq I_{2}$. Hence
$\Lambda$ has a maximal atom $\lambda'=\uaaa{i}{\alpha}{'}$ with
$i_{s}\leq i_{s}'$ for all $1\leq s\leq 4$.

Suppose that $i_{2}=I_{2}$ and $\dim\lambda=p$.  Since
$\uaa{i}{\alpha}{'}{1}\times
u_{3}^{I_{2}}[i_{3}',\alpha_{3}']\times
u_{4}^{I_{2}}[i_{4}',\alpha_{4}']$ is a atom in
$F_{I_{2}}^{u_{2}}(\Lambda)$ and $i_{1}+i_{3}+i_{4}=p-I_{2}$, we
know that $d_{p-I_{2}}^{\gamma}F_{I_{2}}^{u_{2}}(\Lambda)$ has a
maximal atom of the form $u_{1}[i_{1},\alpha_{1}'']\times
u_{3}^{I_{2}}[i_{3},\alpha_{3}'']\times
u_{4}^{I_{2}}[i_{4},\alpha_{4}'']$. Moreover, we can see that for
a fixed $t\in\{1,3,4\}$ there is a maximal atom
$\uaaa{j}{\beta}{}$ in $\Lambda$ with $i_{s}\leq j_{s}$ for all
$s\in\{1,3,4\}$  such that $i_{t}<j_{t}$ if and only if there is a
maximal atom $u_{1}[j_{1},\beta_{1}]\times
u_{3}^{I_{2}}[j_{2},\beta_{3}]\times
u_{4}^{I_{2}}[j_{4},\beta_{4}]$ in $F_{I_{2}}^{u_{2}}(\Lambda)$
with $i_{s}\leq j_{s}$ for all $s\in\{1,3,4\}$ such that
$i_{t}<j_{t}$. It follows that from Theorem \ref{dt} that
$\alpha_{s}=\alpha_{s}''$ for all $s\in\{1,3,4\}$, thus $\mu$ is a
maximal atom in $d_{p-I_{2}}^{\gamma}F(\Lambda)$.

Suppose that $i_{2}=I_{2}$ and $\dim\lambda<p$. Then $\lambda$ is
also a maximal atom in $\Lambda$. Therefore
$\mu=F_{I_{2}}^{u_{2}}(\lambda)$ is a maximal atom in
$F_{I_{2}}^{u_{2}}(\Lambda)$. Since $i_{1}+i_{3}+i_{4}<p-I_{2}$,
we know that $\mu$ is a maximal atom in
$d_{p-I_{2}}^{\gamma}F(\Lambda)$.

There remains the case that $i_{2}>I_{2}$. In this case, there is
no maximal atom $\uaaa{j}{\beta}{}$ in $\Lambda$ with $j_{s}\geq
i_{s}$ for all $s\in\{1,3,4\}$ such that $j_{s}>i_{s}$ for some
$s\in\{1,3,4\}$. So $i_{s}=i_{s}'$ and $\alpha_{s}=\alpha_{s}'$
for all $s\in\{1,3,4\}$. On the other hand, since
$\mu=\uaa{i}{\alpha}{}{1}\times
u_{3}^{I_{2}}[i_{3},\alpha_{3}]\times
u_{4}^{I_{2}}[i_{4},\alpha_{4}] =\uaa{i}{\alpha}{'}{1}\times
u_{3}^{I_{2}}[i_{3}',\alpha_{3}']\times
u_{4}^{I_{2}}[i_{4}',\alpha_{4}']=F_{I_{2}}^{u_{2}}(\lambda')$, we
can see that $\mu$ is a maximal atom in
$F_{I_{2}}^{u_{2}}(\Lambda)$. Because
$\dim\mu=i_{1}'+i_{3}'+i_{4}'<p-I_{2}$, it follows from Theorem
\ref{dt} that $\mu$ is a maximal atom in
$d_{p-I_{2}}^{\gamma}F_{I_{2}}^{u_{2}}(\Lambda)$.

This shows  that every maximal atom in
$F_{I_{2}}^{u_{2}}(\Lambda_{1})$ is a maximal atom in
$d_{p-I_{2}}^{\gamma}F_{I_{2}}^{u_{2}}(\Lambda)$.

We next prove that every maximal atom in
$d_{p-I_{2}}^{\gamma}F_{I_{2}}^{u_{2}}(\Lambda)$ is a maximal atom
in $F_{I_{2}}^{u_{2}}(\Lambda_{1})$.

Let $\mu=\uaa{i}{\alpha}{}{1}\times
u_{3}^{I_{2}}[i_{3},\alpha_{3}]\times
u_{4}^{I_{2}}[i_{4},\alpha_{4}]$ be a maximal atom in
$d_{p-I_{2}}^{\gamma}F_{I_{2}}^{u_{2}}(\Lambda)$. Then
$F_{I_{2}}^{u_{2}}(\Lambda)$ has a maximal atom
$\mu'=\uaa{i}{\alpha}{'}{1}\times
u_{3}^{I_{2}}[i_{3}',\alpha_{3}']\times
u_{4}^{I_{2}}[i_{4}',\alpha_{4}']$ with $i_{s}\leq i_{s}'$ for all
$s\in\{1,2,3\}$. Therefore $\Lambda$ has a
$(u_{2},I_{2})$-projection maximal atom $\lambda'$ of the form
$\lambda'=\uaaa{i}{\alpha}{'}$.

Suppose that $i_{1}+i_{3}+i_{4}=p-I_{2}$. Then $\Lambda_{1}$ has a
$(u_{2},I_{2})$-projection maximal atom of the form
$\lambda=\uaaa{i}{\alpha''}{}$. It is easy to see that for a fixed
$t\in\{1,3,4\}$ there is a maximal atom $\uaaa{j}{\beta}{}$ in
$\Lambda$ with $i_{s}\leq j_{s}$ for all $s\in\{1,3,4\}$  such
that $i_{t}<j_{t}$ if and only if there is a maximal atom
$u_{1}[j_{1},\beta_{1}]\times u_{3}^{I_{2}}[j_{2},\beta_{3}]\times
u_{4}^{I_{2}}[j_{4},\beta_{4}]$ in $F_{I_{2}}^{u_{2}}(\Lambda)$
with $i_{s}\leq j_{s}$ for all $s\in\{1,3,4\}$ such that
$i_{t}<j_{t}$. It follows  that $\alpha_{s}''=\alpha_{s}$ for
$s=1,3,4$. Therefore we have
$F_{I_{2}}^{u_{2}}(\lambda)=\uaaa{i}{\alpha''}{}=\mu$. This
implies that $\mu$ is a maximal atom in
$F_{I_{2}}^{u_{2}}(\Lambda_{1})$.

Suppose that $i_{1}+i_{3}+i_{4}<p-I_{2}$. Then
$\mu=\uaa{i}{\alpha}{}{1}\times
u_{3}^{I_{2}}[i_{3},\alpha_{3}]\times
u_{4}^{I_{2}}[i_{4},\alpha_{4}]$ is also a maximal atom in
$F_{I_{2}}^{u_{2}}(\Lambda)$. So $\Lambda$ has a
$(u_{2},I_{2})$-projection maximal atom
$\lambda'=\uaa{i}{\alpha}{}{1}\times\uaa{i}{\alpha}{'}{2}\times
\uaa{i}{\alpha}{}{3}\times \uaa{i}{\alpha}{}{4}$ with $i_{2}'\geq
I_{2}$. Now, if $i_{2}'=I_{2}$, then $\dim\lambda'<p$; hence
$\lambda'$ is also a maximal atom in $\Lambda_{1}$; it follows
that $\mu=F_{I_{2}}^{u_{2}}(\lambda')$, and hence $\mu$ is a
maximal atom in $F_{I_{2}}^{u_{2}}(\Lambda_{1})$. Suppose that
$i_{2}'>I_{2}$. Then it is easy to see that there is no maximal
atom $\uaaa{j}{\beta}{}$ in $\Lambda$ with $j_{s}\geq i_{s}$ for
all $s\in\{1,3,4\}$ such that $j_{s}>i_{s}$ for some
$s\in\{1,3,4\}$. Hence $\Lambda_{1}$ has a
$(u_{2},I_{2})$-projection maximal atom of the form
$\lambda'=\uaa{i}{\alpha}{}{1}\times\uaa{i}{\alpha}{'}{2}\times
\uaa{i}{\alpha}{}{3}\times \uaa{i}{\alpha}{}{4}$. Hence we see
that $\mu=F_{I_{2}}^{u_{2}}(\lambda')$ and $\mu$ is a maximal atom
in $F_{I_{2}}^{u_{2}}(\Lambda_{1})$.

This shows that every maximal atom in
$d_{p-I_{2}}^{\gamma}F_{I_{2}}^{u_{2}}(\Lambda)$ is a maximal atom
in $F_{I_{2}}^{u_{2}}(\Lambda_{1})$.

This completes the proof.
\end{proof}

\section{Composition of pairwise molecular subcomplexes}
In this section, we study composition of pairwise molecular
subcomplexes in $\ubbb$.

\begin{lemma}\label{lessp4}
Let $\Lambda^{-}$ and $\Lambda^{+}$ be  pairwise molecular
subcomplexes. If $d_{p}^{+}\Lambda^{-}=d_{p}^{-}\Lambda^{+}$, then
for every pair of maximal atoms
$\lambda^{-}=\uaaa{i}{\alpha}{^{-}}$ in $\Lambda^{-}$ and
$\lambda^{+}=\uaaa{i}{\alpha}{^{+}}$ in $\Lambda^{+}$ one has
$$\sum_{s=1}^{4}\min\{i_{s}^{-},i_{s}^{+}\}\leq p.$$
\end{lemma}
\begin{proof}
Let $l_{s}=\min\{i_{s}^{-},i_{s}^{+}\}$. Suppose otherwise that
$\sum_{s=1}^{4}l_{s}>p$. Then there is an ordered triple
$\{i_{1},i_{2},i_{3},i_{4}\}$ with $i_{s}\leq l_{s}$ for every
value of $s$ such that $\sum_{s}^{4}i_{s}=p$. Since
$\sum_{s=1}^{4}l_{s}>p$, we have $i_{t}<l_{t}$ for some $t$. If
$i_{1}<l_{1}$, then, by Theorem \ref{dt4}, we have
$d_{p}^{+}\Lambda^{-}$ has a maximal atom of the form
$u_{1}[i_{1},+]\times \uaa{i}{\alpha}{}{2}\times
\uaa{i}{\alpha}{}{3}\times \uaa{i}{\alpha}{}{4}$, while
$d_{p}^{-}\Lambda^{+}$ has a maximal atom of the form
$u_{1}[i_{1},-]\times \uaa{i}{\alpha}{}{2}\times
\uaa{i}{\alpha}{}{3}\times \uaa{i}{\alpha}{}{4}$. This contradicts
condition 1 for the pairwise molecular subcomplex
$d_{p}^{+}\Lambda^{-}=d_{p}^{-}\Lambda^{+}$. The argument for the
cases $i_{2}<l_{2}$, $i_{3}<l_{3}$ and $i_{4}<l_{4}$ are similar.

This completes the proof.
\end{proof}

\begin{lemma}\label{comprj}
Let $\Lambda^{-}$ and $\Lambda^{+}$ be pairwise molecular
subcomplexes in $\ubbb$. If
$d_{p}^{+}\Lambda^{-}=d_{p}^{-}\Lambda^{+}$, then
$F_{I_{s}}^{u_{s}}(\Lambda^{-})\cap
F_{I_{s}}^{u_{s}}(\Lambda^{+})=
F_{I_{s}}^{u_{s}}(\Lambda^{-}\cap\Lambda^{+})=
F_{I_{s}}^{u_{s}}(d_{p}^{+}\Lambda^{-})=
F_{I_{s}}^{u_{s}}(d_{p}^{-}\Lambda^{+})$ for every value of $s$
and every integer $I_{s}$.
\end{lemma}
\begin{proof}
The proofs for different values of $s$ are similar.  We give the
proof for $s=2$. There are two cases, as follows.

1. Suppose that $I_{2}>p$. We first claim that
$F_{I_{2}}^{u_{2}}(\Lambda^{-})\cap
F_{I_{2}}^{u_{2}}(\Lambda^{+})=\emptyset$.

Indeed, suppose otherwise that $F_{I_{2}}^{u_{2}}(\Lambda^{-})\cap
F_{I_{2}}^{u_{2}}(\Lambda^{+})\neq\emptyset$. Then it is evident
that there are atoms $\mu^{-}=\uaaa{j}{\beta}{^{-}}$ in
$\Lambda^{-}$ and $\mu^{+}=\uaaa{j}{\beta}{^{+}}$ in $\Lambda^{+}$
such that $j_{2}^{-}\geq I_{2}>p$ and $j_{2}^{+}\geq J>p$. This
implies that there are maximal atoms $u_{1}[0,\alpha_{1}']\times
u_{2}[p, +]\times u_{3}[0,\alpha_{3}']\times u_{4}[0,\alpha_{4}']$
and $u_{1}[0,\alpha_{1}'']\times u_{2}[p, -]\times
u_{3}[0,\alpha_{3}'']\times u_{4}[0,\alpha_{4}'']$ in
$d_{p}^{+}\Lambda^{-}$ and $d_{p}^{-}\Lambda^{+}$ respectively.
This contradicts the condition 1 for pairwise molecular subcomplex
$d_{p}^{+}\Lambda^{-}=d_{p}^{-}\Lambda^{+}$.

Now we have $F_{I_{2}}^{u_{2}}(d_{p}^{+}\Lambda^{-})\subset
F_{I_{2}}^{u_{2}}(\Lambda^{-}\cap\Lambda^{+})\subset
F_{I_{2}}^{u_{2}}(\Lambda^{-})\cap
F_{I_{2}}^{u_{2}}(\Lambda^{+})=\emptyset$. Therefore
$F_{I_{2}}^{u_{2}}(d_{p}^{+}\Lambda^{-})=
F_{I_{2}}^{u_{2}}(\Lambda^{-}\cap
\Lambda^{+})=F_{I_{2}}^{u_{2}}(\Lambda^{-})\cap
F_{I_{2}}^{u_{2}}(\Lambda^{+})$, as required.

2. Suppose that $I_{2}\leq p$. Since
$d_{p}^{+}\Lambda^{-}=d_{p}^{-}\Lambda^{+}$, we have
$d_{p-I_{2}}^{+}F_{I_{2}}^{u_{2}}(\Lambda^{-})
=F_{I_{2}}^{u_{2}}(d_{p}^{+}\Lambda^{-})
=F_{I_{2}}^{u_{2}}(d_{p}^{-}\Lambda^{+})
=d_{p-I_{2}}^{-}F_{I_{2}}^{u_{2}}(\Lambda^{+}).
$
Because $F_{I_{2}}^{u_{2}}(\Lambda^{-})$ and
$F_{I_{2}}^{u_{2}}(\Lambda^{+})$ are molecules, we can see that
$F_{I_{2}}^{u_{2}}(\Lambda^{-})\#_{p-J}
F_{I_{2}}^{u_{2}}(\Lambda^{+})$ is defined by Proposition
\ref{dfd}. Hence
$
 F_{I_{2}}^{u_{2}}(\Lambda^{-})\cap F_{I_{2}}^{u_{2}}(\Lambda^{+})
=d_{p-I_{2}}^{+}F_{I_{2}}^{u_{2}}(\Lambda^{-})=F_{I_{2}}^{u_{2}}(d_{p}^{+}\Lambda^{-})
\subset F_{I_{2}}^{u_{2}}(\Lambda^{-}\cap\Lambda^{+})$. Since we
automatically have
$F_{I_{2}}^{u_{2}}(\Lambda^{-}\cap\Lambda^{+})\subset
F_{I_{2}}^{u_{2}}(\Lambda^{-})\cap
F_{I_{2}}^{u_{2}}(\Lambda^{+})$, we get
$F_{I_{2}}^{u_{2}}(\Lambda^{-})\cap
F_{I_{2}}^{u_{2}}(\Lambda^{+})=F_{I_{2}}^{u_{2}}(\Lambda^{-}\cap\Lambda^{+})
=F_{I_{2}}^{u_{2}}(d_{p}^{+}\Lambda^{-})$, as required.

This completes the proof

\end{proof}

\begin{prop}\label{dfd4}
Let $\Lambda^{-}$ and $\Lambda^{+}$ be pairwise molecular
subcomplexes. If $d_{p}^{+}\Lambda^{-}=d_{p}^{-}\Lambda^{+}$, then
$\Lambda^{-}\cap\Lambda^{+}=
d_{p}^{+}\Lambda^{-}(=d_{p}^{-}\Lambda^{+})$; hence
$\Lambda^{-}\#_{p}\Lambda^{+}$ is defined.
\end{prop}

\begin{proof}
Let $M=d_{p}^{+}\Lambda^{-}=d_{p}^{-}\Lambda^{+}$. It is evident
that $M\subset\Lambda^{-}\cap\Lambda^{+}$. To prove the reverse
inclusion, it suffices to prove that every maximal atom in
$\Lambda^{-}\cap\Lambda^{+}$ is contained in $M$.

Suppose otherwise that there is a maximal atom
$\lambda=\uaaa{i}{\alpha}{}$ in $\Lambda^{-}\cap\Lambda^{+}$ such
that $\lambda\not\subset M$. Since $\uaa{i}{\alpha}{}{1}\times
\uaa{i}{\alpha}{}{2}\times \uaa{i}{\alpha}{}{3}
=F_{i_{4}}^{u_{4}}(\lambda)\subset
F_{i_{4}}^{u_{4}}(\Lambda^{-}\cap\Lambda^{+})=F_{i_{4}}^{u_{4}}(M)$,
we can see that $M$ has a maximal atom
$\lambda'=\uaaa{i}{\alpha}{'}$ such that
$\uaa{i}{\alpha}{}{s}\subset\uaa{i}{\alpha}{'}{s}$ for $s=1,2,3$.
Because $\lambda=\uaaa{i}{\alpha}{}$ is maximal in
$\Lambda^{-}\cap\Lambda^{+}$ and $M\subset
\Lambda^{-}\cap\Lambda^{+}$, we have $i_{4}'=i_{4}$ and
$\alpha_{4}'=-\alpha_{4}$. Now we know that
$\lambda\cup\lambda'\subset\Lambda^{-}$ and
$\lambda\cup\lambda'\subset\Lambda^{+}$. By Lemma \ref{lm34}, it
is easy to see that there are maximal atoms
$\lambda^{-}=\uaaa{i}{\alpha}{^{-}}$ in $\Lambda^{-}$ and
$\lambda^{+}=\uaaa{i}{\alpha}{^{+}}$ in $\Lambda^{+}$ such that
$\uaa{i}{\alpha}{^{-}}{s}\cap\uaa{i}{\alpha}{^{+}}{s}\supset
\uaa{i}{\alpha}{}{s}$ for $s=1,2,3$ and
$\min\{i_{4}^{-},i_{4}^{+}\}>i_{4}$. Since $\lambda$ is a maximal
atom in $\Lambda^{-}\cap\Lambda^{+}$, we have
$i_{4}^{-}=i_{4}^{+}=i_{4}+1$ and
$\alpha_{4}^{-}=-\alpha_{4}^{+}$.

Now, we have $\uaa{i}{\alpha}{}{1}\times
\uaa{i}{\alpha}{}{2}\times \uaa{i}{\alpha}{}{3}\subset
F_{i_{4}+1}^{u_{4}}(\Lambda^{-})\cap
F_{i_{4}+1}^{u_{4}}(\Lambda^{+})
=F_{i_{4}+1}^{u_{4}}(\Lambda^{-}\cap\Lambda^{+})$. Therefore
$\Lambda^{-}\cap\Lambda^{+}$ has a maximal atom
$\lambda''=\uaaa{i}{\alpha}{''}$ with
$\uaa{i}{\alpha}{''}{s}\supset\uaa{i}{\alpha}{'}{s}$ for $s=1,2,3$
and $i_{4}''>i_{4}$. This contradicts the assumption that
$\lambda$ is a maximal atom in $\Lambda^{-}\cap\Lambda^{+}$.

This completes the proof.

\end{proof}

\begin{prop}\label{comp4}
Let $\Lambda^{-}$ and $\Lambda^{+}$ be pairwise molecular
subcomplexes of $\ubbb$. If
$d_{p}^{+}\Lambda^{-}=d_{p}^{-}\Lambda^{+}$, then the maximal
atoms in the composite $\Lambda^{-}\#_{p}\Lambda^{+}$ are the
common $q$-dimensional maximal atoms of $\Lambda^{-}$ and
$\Lambda^{+}$ for $q\leq p$ together with the $r$-dimensional
maximal atoms in either $\Lambda^{-}$ or $\Lambda^{+}$ for $r>p$.
\end{prop}
\begin{proof}

 Let $\Lambda$ be the union of the maximal atoms described in the
proposition. We must prove that
$\Lambda=\Lambda^{-}\cup\Lambda^{+}$. Clearly, we have
$\Lambda\subset\Lambda^{-}\cup\Lambda^{+}$; it suffices to prove
that $\Lambda^{-}\cup\Lambda^{+}\subset\Lambda$. By the formation
of $\Lambda$, we must prove that  $\lambda\subset \Lambda$ for
every maximal atom $\lambda=\uaaa{i}{\alpha}{}$ in either
$\Lambda^{-}$ or $\Lambda^{+}$ with $\dim\lambda\leq p$ and such
that $\lambda$ is not a common maximal atom in $\Lambda^{-}$ and
$\Lambda^{+}$. In this case, it is easy to see that
$\dim\lambda=p$. Suppose that $\lambda$ is a maximal atom in
$\Lambda^{\gamma}$ which is not a maximal atom in
$\Lambda^{-\gamma}$. Then $\lambda$ must be a maximal atom in
$d_{p}^{+}\Lambda^{-}=d_{p}^{-}\Lambda^{+}$ which implies that
$\lambda\subset \lambda^{-\gamma}$ for some maximal atom
$\lambda^{-\gamma}=\uaaa{i}{\alpha}{^{-\gamma}}$ with
$\dim\lambda^{-\gamma}>p$. Thus $\lambda\subset \Lambda$.
Therefore, we have $\Lambda^{-}\cup\Lambda^{+}\subset\Lambda$.

This completes the proof.

\end{proof}

\begin{prop}
Let $\Lambda^{-}$ and $\Lambda^{+}$ be pairwise molecular
subcomplexes. If $d_{p}^{+}\Lambda^{-}=d_{p}^{-}\Lambda^{+}$, then
$\Lambda^{-}\#_{p}\Lambda^{+}$ is a pairwise molecular subcomplex
of $\ubbb$.
\end{prop}
\begin{proof}
Let $\Lambda=\Lambda^{-}\#_{p}\Lambda^{+}$. According to Lemma
\ref{lessp4} and Proposition \ref{comp4}, it is easy to see that
$\Lambda$ satisfies condition 1 for pairwise molecular
subcomplexes. Moreover, we have
$F_{I_{s}}^{u_{s}}(\Lambda^{-}\#_{p}\Lambda^{+})
=F_{I_{s}}^{u_{s}}(\Lambda^{-}\cup\Lambda^{+})
=F_{I_{s}}^{u_{s}}(\Lambda^{-})\cup
F_{I_{s}}^{u_{s}}(\Lambda^{+})$ for every value of $s$.

Now suppose that $p\geq I_{s}$. We  have
$d_{p-I_{s}}^{+}F_{I_{s}}^{u_{s}}(\Lambda^{-})
=F_{I_{s}}^{u_{s}}(d_{p}^{+}\Lambda^{-})
=F_{I_{s}}^{u_{s}}(d_{p}^{-}\Lambda^{+})
=d_{p-I_{s}}^{-}F_{I_{s}}^{u_{s}}(\Lambda^{+})$. Thus
$F_{I_{s}}^{u_{s}}(\Lambda^{-}\#_{p}\Lambda^{+})
=F_{I_{s}}^{u_{s}}(\Lambda^{-})\#_{p-I_{s}}
F_{I_{s}}^{u_{s}}(\Lambda^{+})$. Therefore
$F_{I_{s}}^{u_{s}}(\Lambda^{-}\#_{p}\Lambda^{+})$ is a molecule.

Suppose that $p<I_{s}$. Then it is easy to see that
$F_{I_{s}}^{u_{s}}(\Lambda^{-})=\emptyset$ or
$F_{I_{s}}^{u_{s}}(\Lambda^{+})=\emptyset$. (Otherwise, we have
$F_{I_{s}}^{u_{s}}(\Lambda^{-}\cap\Lambda^{+})\neq\emptyset$. This
would lead to a contradiction to Lemma \ref{lessp4}.) Therefore
$F_{I_{s}}^{u_{s}}(\Lambda^{-}\#_{p}\Lambda^{+})$ is a molecule or
empty set.

We have now proved that
$F_{I_{s}}^{u_{s}}(\Lambda^{-}\#_{p}\Lambda^{+})$ is a molecule or
empty set for every value of $s$ and every $I_{s}$. Evidently,
$\Lambda$ satisfies condition 1 for pairwise molecular
subcomplexes. It follows from the definition that $\Lambda$ is a
pairwise molecular subcomplex of $\ubbb$.

\end{proof}

\section{Decomposition of Pairwise Molecular  subcomplexes}

 The aim of this section is to prove the following theorem.

\begin{theorem}\label{pairwise_mole4}
If  $\Lambda$ is a pairwise molecular subcomplex, then  $\Lambda$
is a molecule.
\end{theorem}

It is trivial that the theorem holds when $\Lambda$ is an atom.
Thus we may assume that $\Lambda$ is a pairwise molecular
subcomplex in $\ubbb$ which is not an atom throughout this
section. We are going to show that $\Lambda$ is a molecule.

Let
\begin{center}
$p=\max\{\dim(\lambda\cap\mu)$:\,\,$\lambda$ and $\mu$ are
distinct maximal atoms in $\Lambda\}$.
\end{center}
It is evident that there are at least two maximal atoms $\lambda$
and $\mu$ in $\Lambda$ with $\dim\lambda>p$ and $\dim\mu>p$. By
Theorem \ref{lm34} for pairwise molecular subcomplex $\Lambda$, it
is easy to see that $p$ is the maximal number among the numbers
$\sum_{s=1}^{4}\min\{i_{s},j_{s}\}$, where
$\lambda=\uaaa{i}{\alpha}{}$ and $\mu=\uaaa{j}{\beta}{}$ run over
all pairs of  distinct maximal atoms in $\Lambda$.

\begin{lemma}\label{n_order_p4}
Let $\lambda=\uaaa{i}{\alpha}{}$ and $\mu=\uaaa{j}{\beta}{}$ be
maximal atoms in $\Lambda$ with
$\sum_{s=1}^{4}\min\{i_{s},j_{s}\}=p$.
\begin{enumerate}
\item
Let $i_{1}=j_{1}$ and $\alpha_{1}=-\beta_{1}$. If $i_{2}<j_{2}$,
then $\alpha_{2}=(-)^{i_{1}}\alpha_{1}$; if $i_{3}<j_{3}$,  then
$\alpha_{3}=(-)^{i_{1}+\min\{i_{2},j_{2}\}}\alpha_{1}$; if
$i_{4}<j_{4}$, then
$\alpha_{4}=(-)^{i_{1}+\min\{i_{2},j_{2}\}+\min\{i_{3},j_{3}\}}\alpha_{1}$.
\item
Let $i_{2}=j_{2}$ and $\alpha_{2}=-\beta_{2}$. If $i_{1}<j_{1}$,
then $\alpha_{2}=(-)^{i_{1}}\alpha_{1}$; if $i_{3}<j_{3}$,  then
$\alpha_{3}=(-)^{i_{2}}\alpha_{2}$; if $i_{4}<j_{4}$,  then
$\alpha_{4}=(-)^{i_{2}+\min\{i_{3},j_{3}\}}\alpha_{2}$.
\item
Let $i_{3}=j_{3}$ and $\alpha_{3}=-\beta_{3}$. If $i_{1}<j_{1}$,
then $\alpha_{3}=(-)^{i_{1}+\min\{i_{2},j_{2}\}}\alpha_{1}$; if
$i_{2}<j_{2}$, then $\alpha_{3}=(-)^{i_{2}}\alpha_{2}$; if
$i_{4}<j_{4}$, then $\alpha_{4}=(-)^{i_{3}}\alpha_{3}$.

\item
Let $i_{4}=j_{4}$ and $\alpha_{4}=-\beta_{4}$. If $i_{1}<j_{1}$,
then
$\alpha_{4}=(-)^{i_{1}+\min\{i_{2},j_{2}\}+\min\{i_{3},j_{3}\}}\alpha_{1}$;
if $i_{2}<j_{2}$,  then
$\alpha_{4}=(-)^{i_{2}+\min\{i_{3},j_{3}\}}\alpha_{1}$; if
$i_{3}<j_{3}$, then $\alpha_{4}=(-)^{i_{3}}\alpha_{1}$.
\end{enumerate}

\end{lemma}
\begin{proof}
The arguments for various cases are similar. We give the proof
only for the first case.

Suppose that $i_{1}=j_{1}$, $\alpha_{1}=-\beta_{1}$ and
$i_{2}<j_{2}$. According to Theorem \ref{lm34} for pairwise
molecular subcomplexes, we can get a maximal atom
$\nu=\uaaa{k}{\varepsilon}{}$ with $k_{1}>i_{1}=j_{1}$,
$\uaa{k}{\varepsilon}{}{2}\supset\uaa{i}{\alpha}{}{2}$, $k_{3}\geq
\min\{i_{3},j_{3}\}$ and $k_{4}\geq \min\{i_{4},j_{4}\}$. Since
$\sum_{s=1}^{4}\min\{i_{s},j_{s}\}=p$, we have
$k_{s}=\min\{i_{s},j_{s}\}$ for $s=2,3,4$. Hence
$\uaa{k}{\varepsilon}{}{2}=\uaa{i}{\alpha}{}{2}$. Moreover, it is
easy to see that $\lambda$, $\mu$ and $\nu$ are pairwise adjacent
by the choice of $p$. It follows easily from the sign conditions
that $\alpha_{2}=(-)^{i_{1}}\alpha_{1}$, as required. The
arguments for other cases are similar.

This completes the proof.

\end{proof}

To decompose $\Lambda$ into atoms, we need a total order $<$ on
the atoms in the product of four globes analogous to that on the
atoms in the product of three globes. For a pair of atom atoms
$\lambda=\uaaa{i}{\alpha}{}$ and $\mu=\uaaa{j}{\beta}{}$ in
$\ubbb$, we write $\lambda<\mu$ if one of the following holds:
\begin{itemize}
\item
$\alpha_{1}=\beta_{1}=-$ and $i_{1}<j_{1}$;
\item
$\alpha_{1}=\beta_{1}=+$ and $i_{1}>j_{1}$;
\item
$\alpha_{1}=-$ and $\beta_{1}=+$;
\item
$i_{1}=j_{1}$ are even, $\alpha_{1}=\beta_{1}$,
$\alpha_{2}=\beta_{2}=-$ and $i_{2}<j_{2}$;
\item
$i_{1}=j_{1}$ are even, $\alpha_{1}=\beta_{1}$,
$\alpha_{2}=\beta_{2}=+$ and $i_{2}>j_{2}$;
\item
$i_{1}=j_{1}$ are even, $\alpha_{1}=\beta_{1}$, $\alpha_{2}=-$ and
$\beta_{2}=+$.
\item
$i_{1}=j_{1}$ are odd, $\alpha_{1}=\beta_{1}$,
$\alpha_{2}=\beta_{2}=+$ and $i_{2}<j_{2}$;
\item
$i_{1}=j_{1}$ are odd, $\alpha_{1}=\beta_{1}$,
$\alpha_{2}=\beta_{2}=-$ and $i_{2}>j_{2}$;
\item
$i_{1}=j_{1}$ are odd, $\alpha_{1}=\beta_{1}$, $\alpha_{2}=+$ and
$\beta_{2}=-$.
\item
$i_{1}=j_{1}$, $\alpha_{1}=\beta_{1}$, $j_{1}=j_{2}$,
$\beta_{1}=\beta_{2}$, $i_{1}+i_{2}$ is even,
$\alpha_{3}=\beta_{3}=-$ and $i_{3}<j_{3}$;
\item
$i_{1}=j_{1}$, $\alpha_{1}=\beta_{1}$, $j_{1}=j_{2}$,
$\beta_{1}=\beta_{2}$, $i_{1}+i_{2}$ is even,
$\alpha_{3}=\beta_{3}=+$ and $i_{3}>j_{3}$;
\item
$i_{1}=j_{1}$, $\alpha_{1}=\beta_{1}$, $j_{1}=j_{2}$,
$\beta_{1}=\beta_{2}$, $i_{1}+i_{2}$ is even, $\alpha_{3}=-$ and
$\beta_{3}=+$;
\item
$i_{1}=j_{1}$, $\alpha_{1}=\beta_{1}$, $j_{1}=j_{2}$,
$\beta_{1}=\beta_{2}$, $i_{1}+i_{2}$ is odd,
$\alpha_{3}=\beta_{3}=+$ and $i_{3}<j_{3}$;
\item
$i_{1}=j_{1}$, $\alpha_{1}=\beta_{1}$, $j_{1}=j_{2}$,
$\beta_{1}=\beta_{2}$, $i_{1}+i_{2}$ is odd,
$\alpha_{3}=\beta_{3}=-$ and $i_{3}>j_{3}$;
\item
$i_{1}=j_{1}$, $\alpha_{1}=\beta_{1}$, $j_{1}=j_{2}$,
$\beta_{1}=\beta_{2}$, $i_{1}+i_{2}$ is even, $\alpha_{3}=+$ and
$\beta_{3}=-$.
\end{itemize}

It is evident that the relation $<$ is a total order on the set of
atoms in $\ubbb$.

\begin{lemma}
For any pair of maximal atoms $\lambda$ and $\mu$ in $\Lambda$
with $\dim\lambda>p$ and $\dim\mu>p$, if $\lambda<\mu$, then
$\lambda\cap\mu\subset d_{p}^{+}\lambda\cap d_{p}^{-}\mu$.
\end{lemma}
\begin{proof}
In the proof of this lemma, we use Lemma \ref{dtatom4} without
comments.

Let $\lambda=\uaaa{i}{\alpha}{}$ and $\mu=\uaaa{j}{\beta}{}$. We
consider several cases, as follows.

1. Suppose that $\sum_{s=1}^{4}\min\{i_{s},j_{s}\}=p$. Then
$\lambda$ and $\mu$ are adjacent by the choice of $p$.  According
to  Lemma \ref{n_order_p4} and sign conditions for pairwise
molecular subcomplexes, it is easy to see that
$\lambda\cap\mu\subset d_{p}^{+}\lambda\cap d_{p}^{-}\mu$, as
required.

2. Suppose that $\sum_{s=1}^{4}\min\{i_{s},j_{s}\}<p-2$. According
to condition 1 for pairwise molecular subcomplexes, it is evident
that $\lambda\cap\mu\subset d_{p}^{+}\lambda\cap d_{p}^{-}\mu$, as
required.

3. Suppose that $\sum_{s=1}^{4}\min\{i_{s},j_{s}\}=p-1$ and that
$\lambda$ and $\mu$ are adjacent. There are several case, as
follows: (1) $i_{1}=j_{1}$ and $\alpha_{1}=\beta_{1}$; (2)
$i_{1}=j_{1}$, $\alpha_{1}=-\beta_{1}$, $i_{2}<j_{2}$,
$\alpha_{2}=(-)^{i_{1}}\alpha_{1}$; (3)  $i_{1}=j_{1}$,
$\alpha_{1}=-\beta_{1}$, $i_{2}>j_{2}$ and
$\beta_{2}=(-)^{i_{1}}\beta_{1}$; (4) $i_{1}\neq j_{1}$; (5)
$i_{1}=j_{1}$, $\alpha_{1}=-\beta_{1}$, $i_{2}<j_{2}$,
$\alpha_{2}=-(-)^{i_{1}}\alpha_{1}$, $i_{3}>j_{3}$ and
$i_{4}<j_{4}$; (6) $i_{1}=j_{1}$, $\alpha_{1}=-\beta_{1}$,
$i_{2}<j_{2}$ $\alpha_{2}=-(-)^{i_{1}}\alpha_{1}$, $i_{3}>j_{3}$
and $i_{4}<j_{4}$; (7) $i_{1}=j_{1}$, $\alpha_{1}=-\beta_{1}$,
$i_{2}>j_{2}$ $\beta_{2}=-(-)^{i_{1}}\beta_{1}$, $i_{3}>i_{3}$ and
$i_{4}<j_{4}$; (8) $i_{1}=j_{1}$, $\alpha_{1}=-\beta_{1}$,
$i_{2}>j_{2}$ $\beta_{2}=-(-)^{i_{1}}\beta_{1}$, $i_{3}<i_{3}$ and
$i_{4}>j_{4}$. In the first 4 cases, it follows from the sign
conditions that $\lambda\cap\mu\subset d_{p-1}^{+}\lambda\cap
d_{p-1}^{-}\mu$, thus $\lambda\cap\mu\subset d_{p}^{+}\lambda\cap
d_{p}^{-}\mu$, as required. The arguments for cases (5) to (8) are
similar, we give the proof  for only case (5). In this case, we
have $\alpha_{1}=-$ and hence $\beta_{1}=+$,
$\alpha_{2}=(-)^{i_{1}}$, $\beta_{3}=-(-)^{i_{1}+i_{2}}$ and
$\alpha_{4}=(-)^{i_{1}+i_{2}+j_{3}}$, thus $\lambda\cap\mu\subset
(\uaa{i}{\alpha}{}{1}\times\uaa{i}{\alpha}{}{2}\times
u_{3}[j_{3}+1, \tilde{\alpha_{3}}]\times \uaa{i}{\alpha}{}{4})\cap
(\uaa{j}{\beta}{}{1}\times u_{2}[i_{2}+1,\tilde{\beta_{2}}]\times
\uaa{j}{\beta}{}{3}\times u_{4}[i_{4},(-)^{i_{1}+i_{2}+j_{3}}])
\subset d_{p}^{+}\lambda\cap d_{p}^{-}\mu$, as required.

4. Suppose that $\sum_{s=1}^{4}\min\{i_{s},j_{s}\}=p-1$ and that
$\lambda$ and $\mu$ are not adjacent. Suppose also that
$i_{s}=j_{s}$ for two values of $s$. Then it is easy to see that
$\lambda\cap\mu\subset d_{p}^{+}\lambda\cap d_{p}^{-}\mu$, as
required.

5. Suppose that $\sum_{s=1}^{4}\min\{i_{s},j_{s}\}=p-1$ and that
$\lambda$ and $\mu$ are not adjacent. Suppose also that
$i_{1}=j_{1}$, $i_{2}\neq j_{2}$, $i_{3}\neq j_{3}$ and $i_{4}\neq
j_{4}$. There are various cases: (1) $i_{2}<j_{2}$, $i_{3}<j_{3}$
and $i_{4}>j_{4}$; (2) $i_{2}<j_{2}$, $i_{3}>j_{3}$ and
$i_{4}<j_{4}$; (3) $i_{2}<j_{2}$, $i_{3}>j_{3}$ and $i_{4}>j_{4}$;
(4) $i_{2}>j_{2}$, $i_{3}>j_{3}$ and $i_{4}<j_{4}$; (5)
$i_{2}>j_{2}$, $i_{3}<j_{3}$ and $i_{4}>j_{4}$; (3) $i_{2}>j_{2}$,
$i_{3}<j_{3}$ and $i_{4}<j_{4}$. The arguments for cases (1), (2),
(4) and (5) are similar, and the arguments for cases (3) and (6)
are similar. We give the proof for cases (1) and (3).

In case (1), we have $i_{1}=j_{1}$, $i_{2}<j_{2}$, $i_{3}<j_{3}$
and $i_{4}>j_{4}$. We claim that $\alpha_{2}=-(-)^{i_{1}}$ or
$\alpha_{3}=(-)^{i_{1}+i_{2}}$ which implies that
$\lambda\cap\mu\subset d_{p}^{+}\lambda\cap d_{p}^{-}\mu$, as
required. Indeed, suppose otherwise that $\alpha_{2}=(-)^{i_{1}}$
and $\alpha_{3}=-(-)^{i_{1}+i_{2}}$. By the definition of $<$, we
must have $\alpha_{1}=-$ and $\beta_{1}=+$. If
$\beta_{4}=-(-)^{i_{1}+i_{2}+i_{3}}$, then $\Lambda$ has a maximal
atom $\nu=\uaaa{k}{\varepsilon}{}$ with $k_{1}\geq i_{1}=i_{2}$,
$k_{2}\geq i_{2}$, $k_{3}>i_{3}$ and $k_{4}>j_{4}$; moreover, we
can see that $k_{2}=i_{2}$, $\min\{k_{3},j_{3}\}=i_{3}+1$ and
$\min\{k_{4},i_{4}\}=j_{4}+1$ by the definition of $p$;
furthermore, we can see that $\nu$ is adjacent to both $\lambda$
and $\mu$; it follows that $\varepsilon_{2}=-(-)^{i_{1}}$; this
contradicts the sign condition for $\alpha_{3}$ and
$\varepsilon_{2}$. If $\beta_{4}=(-)^{i_{1}+i_{2}+i_{3}}$, then
$\Lambda$ has a maximal atom $\nu'=\uaaa{k}{\varepsilon}{'}$ with
$k_{1}'\geq i_{1}=i_{2}$, $k_{2}'> i_{2}$, $k_{3}'\geq i_{3}$ and
$k_{4}'>j_{4}$; moreover, we can see that
$\min\{k_{2},j_{2}\}=i_{2}+1$, $k_{3}=i_{3}$ and
$\min\{k_{4},i_{4}\}=j_{4}+1$ by the definition of $p$;
furthermore, we can see that $\nu$ is adjacent to both $\lambda$
and $\mu$; it follows that $\varepsilon_{1}=+=-\alpha_{1}$ when
$k_{1}=i_{1}$; this contradicts the sign condition for
$\alpha_{1}$ and $\alpha_{2}$.

In case (3), we have $i_{1}=j_{1}$,  $i_{2}<j_{2}$, $i_{3}>j_{3}$
and $i_{4}>j_{4}$. We claim that $\beta_{3}=(-)^{i_{1}+i_{2}}$ or
$\beta_{4}=-(-)^{i_{1}+i_{2}+j_{3}}$ which implies that
$\lambda\cap\mu\subset d_{p}^{+}\lambda\cap d_{p}^{-}\mu$, as
required. Indeed, suppose otherwise that
$\beta_{3}=-(-)^{i_{1}+i_{2}}$ and
$\beta_{4}=(-)^{i_{1}+i_{2}+i_{3}}$. If $\alpha_{2}=-(-)^{i_{1}}$,
then $\Lambda$ has a maximal atom $\nu=\uaaa{k}{\varepsilon}{}$
adjacent to both $\lambda$ and $\mu$ with $k_{1}=i_{1}$,
$\min\{k_{2},j_{2}\}=i_{2}+1$, $\min\{k_{3},i_{3}\}=j_{3}+1$ and
$k_{4}=j_{4}$. By the sign condition for $\mu$ and $\nu$, we have
$\varepsilon_{4}=(-)^{i_{1}+i_{2}+i_{3}}$. This contradicts the
sign condition for $\alpha_{2}$ and $\varepsilon_{4}$. If
$\alpha_{2}=(-)^{i_{1}}$, then $\alpha_{1}=-$ and $\beta_{1}=+$ by
the definition of $<$. According to the sign conditions for $\mu$
and $\nu$, we get $\varepsilon_{1}=+$. This contradicts the sign
condition for $\beta_{1}$ and $\alpha_{2}$.

6. Suppose that $\sum_{s=1}^{4}\min\{i_{s},j_{s}\}=p-1$ and that
$\lambda$ and $\mu$ are not adjacent. Suppose also that $i_{1}\neq
j_{1}$, $i_{2}=j_{2}$, $i_{3}\neq j_{3}$ and $i_{4}\neq j_{4}$.
There are various cases: (1) $i_{1}<j_{1}$, $i_{3}<j_{3}$ and
$i_{4}>j_{4}$; (2) $i_{1}<j_{1}$, $i_{3}>j_{3}$ and $i_{4}<j_{4}$;
(3) $i_{1}<j_{1}$, $i_{3}>j_{3}$ and $i_{4}>j_{4}$; (4)
$i_{1}>j_{1}$, $i_{3}>j_{3}$ and $i_{4}<j_{4}$; (5) $i_{1}>j_{1}$,
$i_{3}<j_{3}$ and $i_{4}>j_{4}$; (6) $i_{1}>j_{1}$, $i_{3}<j_{3}$
and $i_{4}<j_{4}$. The arguments for cases (1), (2), (4) and (5)
are similar, and the arguments for cases (3) and (6) are similar.
We give the proof for cases (1) and (3).

In case (1), we have  $i_{1}<j_{1}$, $i_{2}=j_{2}$, $i_{3}<j_{3}$
and $i_{4}>j_{4}$. According to the definition of $<$, we get
$\alpha_{1}=-$. It follows easily that $\lambda\cap\mu\subset
d_{p}^{+}\lambda\cap d_{p}^{-}\mu$, as required.

In case (3), we have  $i_{1}<j_{1}$, $i_{2}=j_{2}$, $i_{3}>j_{3}$
and $i_{4}>j_{4}$. According to the definition of $<$, we get
$\alpha_{1}=-$. We claim that $\beta_{3}=(-)^{i_{1}+i_{2}}$ or
$\beta_{4}=-(-)^{i_{1}+i_{2}+j_{3}}$ which implies that
$\lambda\cap\mu\subset d_{p}^{+}\lambda\cap d_{p}^{-}\mu$, as
required. Indeed, suppose otherwise that
$\beta_{3}=-(-)^{i_{1}+i_{2}}$ and
$\beta_{4}=(-)^{i_{1}+i_{2}+j_{3}}$. Then $\Lambda$ has a maximal
atom $\nu=\uaaa{k}{\varepsilon}{}$ adjacent to both $\lambda$ and
$\mu$ with $\min\{k_{1},j_{1}\}=i_{1}+1$, $k_{2}\geq i_{2}$,
$\min\{k_{3},i_{3}\}=j_{3}+1$ and $k_{4}=j_{4}$. According to the
sign conditions for $\lambda$ and $\nu$, we have
$\varepsilon_{4}=-(-)^{i_{1}+i_{2}+j_{3}}$ which contradicts the
sign condition for $\beta_{3}$ and $\varepsilon_{4}$.

7. Suppose that $\sum_{s=1}^{4}\min\{i_{s},j_{s}\}=p-1$ and that
$\lambda$ and $\mu$ are not adjacent. Suppose also that $i_{1}\neq
j_{1}$, $i_{2}\neq j_{2}$, $i_{3}=j_{3}$ and $i_{4}\neq j_{4}$. By
similar arguments as in case 6, we can get $\lambda\cap\mu\subset
d_{p}^{+}\lambda\cap d_{p}^{-}\mu$, as required.

8. Suppose that $\sum_{s=1}^{4}\min\{i_{s},j_{s}\}=p-1$ and that
$\lambda$ and $\mu$ are not adjacent. Suppose also that $i_{1}\neq
j_{1}$, $i_{2}\neq j_{2}$, $i_{3}\neq j_{3}$ and $i_{4}=j_{4}$. By
similar arguments as in case 6, we can get $\lambda\cap\mu\subset
d_{p}^{+}\lambda\cap d_{p}^{-}\mu$, as required.

9. Suppose that $\lambda$ and $\mu$ are not adjacent. Suppose also
that $\sum_{s=1}^{4}\min\{i_{s},j_{s}\}=p-1$ and $i_{s}\neq j_{s}$
for all values of $s$. There are various cases. The arguments for
these cases are similar. We give the proof for two cases.

Suppose that $i_{1}<j_{1}$, $i_{2}<j_{2}$, $i_{3}>j_{3}$ and
$i_{4}>j_{4}$. Then $\alpha_{1}=-$ by the definition of $<$. We
claim that $\beta_{3}=(-)^{i_{1}+i_{2}}$ or
$\beta_{4}=-(-)^{i_{1}+i_{2}+i_{3}}$ which implies that
$\lambda\cap\mu\subset d_{p}^{+}\lambda\cap d_{p}^{-}\mu$, as
required. Indeed, suppose otherwise that
$\beta_{3}=-(-)^{i_{1}+i_{2}}$ and
$\beta_{4}=(-)^{i_{1}+i_{2}+i_{3}}$. Then $\Lambda$ has a maximal
atom $\nu=\uaaa{k}{\varepsilon}{}$ adjacent to both $\lambda$ and
$\mu$ with $k_{1}=i_{1}+1$, $k_{2}=i_{2}$, $k_{3}=j_{3}+1$ and
$k_{4}=j_{4}$; it follow from the sign conditions for $\nu$ and
$\lambda$ that $\varepsilon_{4}=-(-)^{i_{1}+i_{2}+j_{3}}$ which
contradicts the sign condition for $\beta_{3}$ and
$\varepsilon_{4}$.

Suppose that $i_{1}<j_{1}$, $i_{2}>j_{2}$, $i_{3}>j_{3}$ and
$i_{4}>j_{4}$. By similar arguments as that in the above case, we
can prove that $\beta_{2}=(-)^{i_{1}}$ and
$\beta_{3}=(-)^{i_{1}+i_{2}}$, or $\beta_{2}=(-)^{i_{1}}$ and
$\beta_{4}=-(-)^{i_{1}+i_{2}+i_{3}}$, or
$\beta_{3}=-(-)^{i_{1}+i_{2}}$ and
$\beta_{4}=-(-)^{i_{1}+i_{2}+i_{3}}$, which implies that
$\lambda\cap\mu\subset d_{p}^{+}\lambda\cap d_{p}^{-}\mu$, as
required.

10. Suppose that $\sum_{s=1}^{4}\min\{i_{s},j_{s}\}=p-2$. If
$i_{s}=j_{s}$ for some value of $s$, then it is evident that
$\lambda\cap\mu\subset d_{p}^{+}\lambda\cap d_{p}^{-}\mu$, as
required.

Now suppose that $i_{s}\neq j_{s}$ for every value of $s$. If
$\lambda$ and $\mu$ are adjacent, then we have
$\lambda\cap\mu\subset d_{p-2}^{+}\lambda\cap
d_{p-2}^{-}\mu\subset d_{p}^{+}\lambda\cap d_{p}^{-}\mu$, as
required.

If $i_{1}<j_{1}$ and if $i_{s}<j_{s}$ for some value of $s$ with
$s=2,3,4$, then we have $\alpha_{1}=-$; it follows easily that
$\lambda\cap\mu\subset d_{p}^{+}\lambda\cap d_{p}^{-}\mu$, as
required. If $i_{1}>j_{1}$ and if $i_{s}>j_{s}$ for some value of
$s$ with $s=2,3,4$, then we have $\beta_{1}=+$; it follows easily
that $\lambda\cap\mu\subset d_{p}^{+}\lambda\cap d_{p}^{-}\mu$, as
required.

There remain two cases: (1) $\lambda$ and $\mu$ are not adjacent
and $i_{1}<j_{1}$ and $i_{s}>j_{s}$ for $s=2, 3,4$; (2) $\lambda$
and $\mu$ are not adjacent and $i_{1}>j_{1}$ and $i_{s}<j_{s}$ for
$s=2, 3,4$.  The arguments for the two cases are similar. We give
the proof for the first case.

In the first case, we have $\alpha_{1}=-$ by the definition of
$<$. We claim that $\beta_{2}=(-)^{i_{1}}$ or
$\beta_{3}=-(-)^{i_{1}+j_{2}}$ or
$\beta_{4}=(-)^{i_{1}+j_{2}+j_{3}}$ which implies that
$\lambda\cap\mu\subset d_{p}^{+}\lambda\cap d_{p}^{-}\mu$, as
required. Indeed, suppose otherwise that $\beta_{2}=-(-)^{i_{1}}$
and $\beta_{3}=(-)^{i_{1}+j_{2}}$ and
$\beta_{4}=-(-)^{i_{1}+j_{2}+j_{3}}$. Then there is a maximal atom
$\nu=\uaaa{k}{\varepsilon}{}$ such that $k_{1}>i_{1}$, $k_{2}\geq
j_{2}$, $k_{3}\geq j_{3}$ and $k_{4}\geq j_{4}$. According to sign
conditions and condition \ref{main4_4} in Theorem \ref{main4}, we
can see that, for each fixed value of $s$ with $s=2,3,4$, $\nu$
can be chosen such that $k_{s}>j_{s}$. Moreover, by the choice of
$p$, there are at most two values of $s$ with $s=2,3,4$ such that
$k_{s}>j_{s}$. Now there are several cases, as follows.

(a). Suppose that $\nu$ can be chosen such that there are two
values of $s$ with $s=2,3,4$ such that $k_{s}>j_{s}$. The
arguments for various choices of the two values are similar. We
give the proof for $k_{2}>j_{2}$ and $k_{3}>j_{3}$. In this case,
we have $k_{2}=j_{2}+1$ and $k_{3}=j_{3}+1$ and $k_{4}=j_{4}$ by
the choice of $p$, and $\lambda$ and $\nu$ are adjacent. It
follows from sign conditions that
$\varepsilon_{4}=(-)^{i_{1}+j_{2}+j_{3}}=-\beta_{4}$. According to
Lemma \ref{lm34}, we can get a maximal atom
$\mu'=\uaaa{j}{\beta}{'}$ such that
$j_{1}'\geq\min\{j_{1},k_{1}\}>i_{1}$,
$\uaa{j}{\beta}{'}{2}\supset\uaa{j}{\beta}{}{2}$,
$\uaa{j}{\beta}{'}{3}\supset\uaa{j}{\beta}{}{3}$ and
$j_{4}'>j_{4}$. If $j_{2}'>j_{2}$, then it is evident that
$j_{2}'=j_{2}+1$ and $\uaa{j}{\beta}{'}{3}=\uaa{j}{\beta}{}{3}$ by
the choice of $p$, and $\mu'$ is adjacent to $\lambda$; this gives
a contradiction to the sign condition for $\beta_{3}'$ and
$\alpha_{1}$. Suppose that $j_{2}'=j_{2}$. Then
$\beta_{2}'=\beta_{2}=-(-)^{i_{1}}$. Thus $\lambda$ and $\mu'$ are
not $(1,2)$-adjacent. It follows that there is a maximal atom
$\lambda'=\uaaa{i}{\alpha}{'}$ in $\Lambda$ such that
$i_{1}'>i_{1}$, $i_{2}'>j_{2}'=j_{2}$, $i_{3}'\geq j_{3}$ and
$i_{4}'\geq \min\{j_{4}',i_{4}\}>j_{4}$. By the choice of $p$, it
is easy to see that $i_{2}'=j_{2}+1$, $i_{3}'=j_{3}$ and
$i_{4}'=j_{4}+1$, and that $\lambda'$ is adjacent to both
$\lambda$ and $\mu'$. This leads to a contradiction to the sign
conditions for $\alpha_{1}$, $\alpha_{3}'$ and $\beta_{2}'$.

(b). Suppose that $\nu$ cannot be chosen such that there are two
values of $s$ with $s=2,3,4$ such that $k_{s}>j_{s}$. Then, for
each value of $s$ with $s=2,3,4$, $\nu$ can be chosen such that
$k_{s}>j_{s}$. In particular, $\nu$ can be chosen such that
$k_{1}>i_{1}$ and $k_{2}>j_{2}$. Moreover, we have $k_{2}=j_{2}+1$
or $k_{2}=j_{2}+2$ by the choice of $p$. By the assumption, we can
see that $\lambda$ is both $(1,3)$-adjacent and $(1,4)$-adjacent
to $\nu$.

Suppose that $k_{2}=j_{2}+1$.  It follows from sign conditions
that $\varepsilon_{3}=-(-)^{i_{1}+j_{2}}=-\beta_{3}$ and
$\varepsilon_{4}=-(-)^{i_{1}+j_{2}+j_{3}}=\beta_{4}$. According to
Lemma \ref{lm34} and the assumptions, there is a maximal atom
$\nu'=\uaaa{k}{\varepsilon}{'}$ such that
$k_{1}'\geq\min\{j_{1},k_{1}\}>i_{1}$,
$\uaa{k}{\varepsilon}{'}{2}=\uaa{j}{\beta}{'}{2}$, $k_{3}'>j_{3}$
and $\uaa{k}{\varepsilon}{'}{4}=\uaa{j}{\beta}{}{4}$. It follows
that $\lambda$ and $\nu'$ are not $(1,2)$-adjacent. Thus $\Lambda$
has a maximal atom $\nu''=\uaaa{k}{\varepsilon}{''}$ such that
$k_{1}''>i_{1}$, $k_{2}''>k_{2}'=j_{2}$,
$k_{3}''\geq\min\{k_{3}',i_{3}\}>j_{3}$ and $k_{4}''\geq
k_{4}'=j_{4}$. This contradicts to the assumption on the choice of
$\nu$.

Suppose that $k_{2}=j_{2}+2$. Then one can get a contradiction by
a similar argument.

This completes the proof.

\end{proof}

By this lemma, we can arrange all the maximal atoms in $\Lambda$
with dimension greater than $p$ as $$\lambda_{1},
\lambda_{2},\cdots,\lambda_{n}$$ such $\lambda_{i}\cap
\lambda_{j}\subset d_{p}^{+}\lambda_{i}\cap d_{p}^{-}\lambda_{j}$
for $i<j$. We denote $\lambda_{k}$ in the list by
$$\lambda_{k}=\uaa{i}{\alpha}{^{(k)}}{1}\times
\uaa{i}{\alpha}{^{(k)}}{2}\times \uaa{i}{\alpha}{^{(k)}}{3}\times
\uaa{i}{\alpha}{^{(k)}}{4}.$$

Let $\Lambda^{-}=d_{p}^{-}\Lambda\cup\lambda_{1}$ and
$\Lambda^{+}=d_{p}^{+}\Lambda\cup\lambda_{2}\cdots\cup\lambda_{n}$.
We are going to prove that $\Lambda^{-}$ and $\Lambda^{+}$ are
pairwise molecular subcomplexes and $\Lambda$ can be decomposed
into $\Lambda^{-}$ and $\Lambda^{+}$.

\begin{lemma} \label{L4-c1}
$\Lambda^{-}$ satisfies condition 1 for pairwise molecular
subcomplexes.
\end{lemma}
\begin{proof}
We first prove that $d_{p}^{-}\lambda_{1}\subset
d_{p}^{-}\Lambda$. Suppose that $\xi\in d_{p}^{-}\lambda_{1}$.
Then, for every maximal atom $\lambda'$ in $\Lambda$ with
$\xi\in\lambda'$, if $\lambda'=\lambda_{t}$ for some $t>1$, then
$\xi\in \lambda_{1}\cap\lambda_{t}\subset
d_{p}^{-}\lambda_{t}=d_{p}^{-}\lambda'$; if $\dim\lambda'\leq p$,
then we automatically have $\xi\in d_{p}^{-}\lambda'$.  It follows
from Lemma \ref{xi_dpg} that $d_{p}^{-}\lambda_{1}\subset
d_{p}^{-}\Lambda$, as required.

We now verify that $\Lambda^{-}$ satisfies condition 1 for
pairwise molecular subcomplexes. It suffices to  prove that any
maximal atom $\lambda=\uaaa{i}{\alpha}{}$ in $d_{p}^{-}\Lambda$
with $i_{s}\leq i_{s}^{(1)}$ for $s=1,2,3,4$ is contained in
$\lambda_{1}$. By the formation of $d_{p}^{-}\lambda_{1}$ and
$d_{p}^{-}\Lambda$, it is easy to see that $\lambda$ is a maximal
atom in $d_{p}^{-}\lambda_{1}$, and hence
$\lambda\subset\lambda_{1}$, as required.
\end{proof}

\begin{lemma}\label{L4+c1}
$\Lambda^{+}$ satisfies condition 1 for pairwise molecular
subcomplexes.
\end{lemma}
\begin{proof} It suffices to prove that any maximal atom
$\lambda=\uaaa{i}{\alpha}{}$ in $d_{p}^{+}\Lambda$ with $i_{s}\leq
i_{s}^{(t)}$ for $s=1,2,3,4$ and some $2\leq t\leq n$ is contained
in some $\lambda_{r}$ with $2\leq r\leq n$. It is evident that
$\dim\lambda=p$.

Let $r$ be the maximal integer $t$ between $2$ and $n$ such that
with $i_{s}\leq i_{s}^{(t)}$ for $s=1,2,3,4$. Then
$d_{p}^{+}\lambda_{r}$ has a maximal atom of the form
$\lambda'=\uaaa{i}{\alpha'}{}$. By the choice of $r$, it is
evident that $\Int\lambda'\cap\lambda_{t}=\emptyset$ for any
$t>r$. Moreover, for any $1\leq s<r$, we have
$\lambda'\cap\lambda_{s}\subset\lambda_{r}\cap\lambda_{s}\subset
d_{p}^{+}\lambda_{s}$. By Lemma \ref{xi_dpg}, we can see that
$\Int\lambda'\subset d_{p}^{+}\Lambda$ and hence $\lambda'\subset
d_{p}^{+}\Lambda$. So, by condition 1 for the pairwise molecular
subcomplex $d_{p}^{+}\Lambda$, we can see that
$\lambda=\lambda'\subset\lambda_{r}$, as required.

This completes the proof.
\end{proof}

\begin{lemma}\label{Fi1u}
Let $1\leq r\leq 4$. If $p\geq I_{r}$ and $\lambda_{1}$ is  a
$(u_{r},I_{r})$-projection maximal atom, then
\begin{enumerate}
\item
$F_{I_{r}}^{u_{r}}(\Lambda^{-})$ and
$F_{I_{r}}^{u_{r}}(\Lambda^{+})$ are molecules;
\item
$d_{p-I_{r}}^{+}F_{I_{r}}^{u_{r}}(\Lambda^{-})=d_{p-I_{r}}^{-}
F_{I_{r}}^{u_{r}}(\Lambda^{+})$, hence
$F_{I_{r}}^{u_{r}}(\Lambda^{-})\#_{p-I_{r}}F_{I_{r}}^{u_{r}}(\Lambda^{+})$
is defined;
\item
$F_{I_{r}}^{u_{r}}(\Lambda)=F_{I_{r}}^{u_{r}}(\Lambda^{-})\#_{p-I_{r}}F_{I_{r}}^{u_{r}}(\Lambda^{+})$.
\end{enumerate}
\end{lemma}
\begin{proof}
The arguments for various choices of $r$ are similar. We prove
only for $r=1$.

Since $F_{I_{1}}^{u_{1}}$ preserves unions, we have
$F_{I_{1}}^{u_{1}}(\Lambda^{-})
=F_{I_{1}}^{u_{1}}(d_{p}^{-}\Lambda\cup\lambda_{1})
=F_{I_{1}}^{u_{1}}(d_{p}^{-}\Lambda)\cup
F_{I_{1}}^{u_{1}}(\lambda_{1})$ and
$F_{I_{1}}^{u_{1}}(\Lambda^{+})
=F_{I_{1}}^{u_{1}}(d_{p}^{+}\Lambda\cup\lambda_{2}\cup\cdots\cup\lambda_{n})
=F_{I_{1}}^{u_{1}}(d_{p}^{+}\Lambda)\cup
F_{I_{1}}^{u_{1}}(\lambda_{2})\cup\cdots\cup
F_{I_{1}}^{u_{1}}(\lambda_{n})$. If $\dim
F_{I_{1}}^{u_{1}}(\lambda_{1})\leq p-I_{1}$, then it is evident
that
$F_{I_{1}}^{u_{1}}(\Lambda^{-})=F_{I_{1}}^{u_{1}}(d_{p}^{-}\Lambda)
=d_{p-I_{1}}^{-}F_{I_{1}}^{u_{1}}(\Lambda)$ and
$F_{I_{1}}^{u_{1}}(\Lambda^{+})=F_{I_{1}}^{u_{1}}(\Lambda)$; it
follows easily that $F_{I_{1}}^{u_{1}}(\Lambda^{-})$ and
$F_{I_{1}}^{u_{1}}(\Lambda^{+})$ are molecules and
$d_{p-I_{1}}^{+}F_{I_{1}}^{u_{1}}(\Lambda^{-})=d_{p-I_{1}}^{-}F_{I_{1}}^{u_{1}}(\Lambda^{+})$,
as required. If $F_{I_{1}}^{u_{1}}(\Lambda)$ consists of only one
maximal atom, then
$F_{I_{1}}^{u_{1}}(\Lambda)=F_{I_{1}}^{u_{1}}(\lambda_{1})$; it
follows that $F_{I_{1}}^{u_{1}}(\Lambda^{-})=
F_{I_{1}}^{u_{1}}(\lambda_{1})=F_{I_{1}}^{u_{1}}(\Lambda)$ and
$F_{I_{1}}^{u_{1}}(\Lambda^{+})=F_{I_{1}}^{u_{1}}(d_{p}^{+}\Lambda)
=d_{p-I_{1}}^{+}F_{I_{1}}^{u_{1}}(\Lambda)$; hence
$F_{I_{1}}^{u_{1}}(\Lambda^{-})$ and
$F_{I_{1}}^{u_{1}}(\Lambda^{+})$ are molecules and
$d_{p-I_{1}}^{+}F_{I_{1}}^{u_{1}}(\Lambda^{-})=d_{p-I_{1}}^{-}F_{I_{1}}^{u_{1}}(\Lambda^{+})$,
as required. In the following proof, we may assume that $\dim
F_{I_{1}}^{u_{1}}(\lambda_{1})>p-I_{1}$ and
$F_{I_{1}}^{u_{1}}(\Lambda)$ consists of at least two distinct
maximal atoms.

Let $$q=\max\{\dim(\mu\cap\mu'):\text{$\mu$ and $\mu'$ are
distinct maximal atoms in $F_{I_{1}}^{u_{1}}(\Lambda)$}\}.$$ It is
clear that $q\leq p-I_{1}$ by the choice of $p$. Let
$\mu=u_{2}^{I_{1}}[j_{2},\beta_{2}]\times
u_{3}^{I_{1}}[j_{3},\beta_{3}]\times
u_{4}^{I_{1}}[j_{4},\beta_{4}]$ be a maximal atom in
$F_{I_{1}}^{u_{1}}(\Lambda)$ distinct from
$F_{I_{1}}^{u_{1}}(\lambda_{1})$ such that $\dim\mu> p-I_{1}$.  We
first prove that $F_{I_{1}}^{u_{1}}(\lambda_{1})\cap\mu\subset
d_{p-I_{1}}^{+}F_{I_{1}}^{u_{1}}(\lambda_{1})\cap
d_{p-I_{1}}^{-}\mu$.

Since $\mu$ is a maximal atom in $F_{I_{1}}^{u_{1}}(\Lambda)$,
there is a $(u_{1},I_{1})$-projection maximal atom $\tilde{\mu}$
of the form $\tilde{\mu}=\uaaa{j}{\beta}{}$. We consider several
cases, as follows.

1. Suppose that $\min\{i_{1}, j_{1}\}=I_{1}$. Since
$\lambda_{1}\cap\tilde{\mu}\subset d_{p}^{+}\lambda_{1}\cap
d_{p}^{-}\tilde{\mu}$, it is easy to see that
$F_{I_{1}}^{u_{1}}(\lambda_{1})\cap\mu\subset
d_{p-I_{1}}^{+}F_{I_{1}}^{u_{1}}(\lambda_{1})\cap
d_{p-I_{1}}^{-}\mu$, as required.

2. Suppose that $\min\{i_{1},j_{1}\}>I_{1}+1$. Then
$\min\{i_{2},j_{2}\}+\min\{i_{3},j_{3}\}+\min\{i_{4},j_{4}\}\leq
p-I_{1}-2$. It follows easily that
$F_{I_{1}}^{u_{1}}(\lambda_{1})\cap\mu\subset
d_{p-I_{1}}^{+}F_{I_{1}}^{u_{1}}(\lambda_{1})\cap
d_{p-I_{1}}^{-}\mu$, as required.

3. Suppose that $\min\{i_{1},j_{1}\}=I_{1}+1$. Then
$\min\{i_{2},j_{2}\}+\min\{i_{3},j_{3}\}+\min\{i_{4},j_{4}\}\leq
p-I_{1}-1$. If
$\min\{i_{2},j_{2}\}+\min\{i_{3},j_{3}\}+\min\{i_{4},j_{4}\}<p-I_{1}-1$,
then it is evident that
$F_{I_{1}}^{u_{1}}(\lambda_{1})\cap\mu\subset
d_{p-I_{1}}^{+}F_{I_{1}}^{u_{1}}(\lambda_{1})\cap
d_{p-I_{1}}^{-}\mu$, as required. If
$\min\{i_{2},j_{2}\}+\min\{i_{3},j_{3}\}+\min\{i_{4},j_{4}\}
=p-I_{1}-1$, and if $i_{s}=j_{s}$ for some value of $s$ with
$s=2,3,4$, then it is evident that
$F_{I_{1}}^{u_{1}}(\lambda_{1})\cap\mu\subset
d_{p-I_{1}}^{+}F_{I_{1}}^{u_{1}}(\lambda_{1})\cap
d_{p-I_{1}}^{-}\mu$, as required. If
$\min\{i_{2},j_{2}\}+\min\{i_{3},j_{3}\}+\min\{i_{4},j_{4}\}
=p-I_{1}-1$, and if $i_{s}\neq j_{s}$ for $s=2,3,4$, then
$\min\{i_{1},j_{1}\}+
\min\{i_{2},j_{2}\}+\min\{i_{3},j_{3}\}+\min\{i_{4},j_{4}\}=p$;
thus $\lambda_{1}$ and $\tilde{\mu}$ are adjacent; it follows
easily from the sign condition for $\lambda_{1}$ and $\tilde{\mu}$
that $F_{I_{1}}^{u_{1}}(\lambda_{1})\cap\mu\subset
d_{p-I_{1}}^{+}F_{I_{1}}^{u_{1}}(\lambda_{1})\cap
d_{p-I_{1}}^{-}\mu$, as required.

Now, we have
$F_{I_{1}}^{u_{1}}(\Lambda^{-})=d_{p-I_{1}}^{-}F_{I_{1}}^{u_{1}}(\Lambda)\cup
F_{I_{1}}^{u_{1}}(\lambda_{1})$ and
$$F_{I_{1}}^{u_{1}}(\Lambda^{+})=d_{p-I_{1}}^{+}F_{I_{1}}^{u_{1}}(\Lambda)\cup\bigcup\{\mu:
\mu\text{ is a maximal atom in }F_{I_{1}}^{u_{1}}(\Lambda)\text{
with } \mu\neq F_{I_{1}}^{u_{1}}(\lambda_{1})\}$$ (Note that it is
possible that
$F_{I_{1}}^{u_{1}}(\Lambda^{+})=d_{p-I_{1}}^{+}F_{I_{1}}^{u_{1}}(\Lambda)$).
It follows from Theorem \ref{decom_mole} that
$F_{I_{1}}^{u_{1}}(\Lambda^{-})$ and
$F_{I_{1}}^{u_{1}}(\Lambda^{+})$ are molecules in
$u_{2}^{I_{1}}\times u_{3}^{I_{1}}\times u_{4}^{I_{1}}$,
$d_{p-I_{1}}^{+}F_{I_{1}}^{u_{1}}(\Lambda^{-})=
d_{p-I_{1}}^{-}F_{I_{1}}^{u_{1}}(\Lambda^{+})$ and
$F_{I_{1}}^{u_{1}}(\Lambda)=F_{I_{1}}^{u_{1}}(\Lambda^{-})\#_{p-I_{1}}F_{I_{1}}^{u_{1}}(\Lambda^{+})$,
as required.

This completes the proof.

\end{proof}

\begin{prop}
Let $\Lambda$ be a pairwise molecular subcomplex. Then
\begin{enumerate}
\item
$\Lambda^{-}$ and $\Lambda^{+}$ are pairwise molecular
subcomplexes.
\item
$d_{p}^{+}\Lambda^{-}=d_{p}^{-}\Lambda^{+}$, hence  the composite
$\Lambda^{-}\#_{p}\Lambda^{+}$ is defined.
\item
$\Lambda=\Lambda^{-}\#_{p}\Lambda_{+}$.
\end{enumerate}
\end{prop}
\begin{proof}
We first prove that $\Lambda^{-}$ and $\Lambda^{+}$ are pairwise
molecular subcomplexes. If $\lambda_{1}$ is not a
$(u_{1},I_{1})$-projection maximal atom in $\Lambda$, then it is
easy to see that
$F_{I_{1}}^{u_{1}}(\Lambda^{-})=F_{I_{1}}^{u_{1}}(d_{p}^{-}\Lambda)$
and $F_{I_{1}}^{u_{1}}(\Lambda^{+})=F_{I_{1}}^{u_{1}}(\Lambda)$ by
the choice of $p$ and Lemmas \ref{L4-c1} and \ref{L4+c1}; hence
$F_{I_{1}}^{u_{1}}(\Lambda^{-})$ and
$F_{I_{1}}^{u_{1}}(\Lambda^{+})$ are the empty set or molecules in
$u_{2}^{I_{1}}\times u_{3}^{I_{1}}\times u_{4}^{I_{1}}$. If
$\lambda_{1}$ is a $(u_{1},I_{1})$-projection maximal atom in
$\Lambda$, then we have already seen that
$F_{I_{1}}^{u_{1}}(\Lambda^{-})$ and
$F_{I_{1}}^{u_{1}}(\Lambda^{+})$ are molecules in
$u_{2}^{I_{1}}\times u_{3}^{I_{1}}\times u_{4}^{I_{1}}$ from Lemma
\ref{Fi1u}. Consequently, $F_{I_{1}}^{u_{1}}(\Lambda^{-})$ and
$F_{I_{1}}^{u_{1}}(\Lambda^{+})$ are the empty set or molecules in
$u_{2}^{I_{1}}\times u_{3}^{I_{1}}\times u_{4}^{I_{1}}$ for every
integer $I_{1}$. Similarly, $F_{I_{s}}^{u_{s}}(\Lambda^{-})$ and
$F_{I_{s}}^{u_{s}}(\Lambda^{+})$ are the empty set or molecules in
the corresponding $\omega$-complex for every value of $s$ and
every integer $I_{s}$. It follows that $\Lambda^{-}$ and
$\Lambda^{+}$ are pairwise molecular subcomplex of $\ubbb$.

Now, if $p\geq I_{1}$ and $\lambda_{1}$ is not
$(u_{1},I_{1})$-projection maximal, then we can see that $$
\begin{array}{rl}
&F_{I_{1}}^{u_{1}}(d_{p}^{+}\Lambda^{-})\\
=&d_{p-I_{1}}^{+}F_{I_{1}}^{u_{1}}(\Lambda^{-})\\
=&d_{p-I_{1}}^{+}F_{I_{1}}^{u_{1}}(d_{p}^{-}\Lambda))\\
=&d_{p-I_{1}}^{-}F_{I_{1}}^{u_{1}}(\Lambda)\\
=&d_{p-I_{1}}^{-}F_{I_{1}}^{u_{1}}(\Lambda^{+})\\
=&F_{I_{1}}^{u_{1}}(d_{p}^{-}\Lambda^{+});
\end{array}
$$ if $p<I_{1}$, then
$F_{I_{1}}^{u_{1}}(d_{p}^{+}\Lambda^{-})=\emptyset=F_{I_{1}}^{u_{1}}(d_{p}^{+}\Lambda^{-})$;
if $p\geq I_{1}$ and $\lambda$ is $(u_{1},I_{1})$-projection
maximal, then $F_{I_{1}}^{u_{1}}(d_{p}^{+}\Lambda^{-})
=F_{I_{1}}^{u_{1}}(d_{p}^{-}\Lambda^{+})$ by Proposition
\ref{prop5.6F}. Consequently, we have
$F_{I_{1}}^{u_{1}}(d_{p}^{+}\Lambda^{-})
=F_{I_{1}}^{u_{1}}(d_{p}^{+}\Lambda^{-})$ for every value of
$I_{1}$. Similarly, we can see that we have
$F_{I_{s}}^{u_{s}}(d_{p}^{+}\Lambda^{-})
=F_{I_{s}}^{u_{s}}(d_{p}^{-}\Lambda^{+})$ for every value of $s$
and every value of $I_{s}$. If follows from Proposition \ref{eq4}
that $d_{p}^{+}\Lambda^{-}=d_{p}^{-}\Lambda^{+}$. Clearly, we have
$\Lambda=\Lambda^{-}\cup\Lambda^{+}$. Therefore
$\Lambda=\Lambda^{-}\#_{p}\Lambda^{+}$, as required. This
completes the proof.
\end{proof}

We have now proved that a pairwise molecular subcomplex $\Lambda$
in $\ubbb$ can be decomposed into pairwise molecular subcomplexes
$\Lambda=\Lambda^{-}\#_{p}\Lambda^{+}$. It is evident that this is
a proper decomposition. By induction, we can see that $\Lambda$
can be eventually decomposed into atoms. Thus $\Lambda$ is a
molecule. So we get the proof for Theorem \ref{pairwise_mole4}.


\printindex
\end{document}